\documentclass[12pt, a4paper]{article}
\usepackage{amssymb,amsmath,amsthm}
\usepackage[all]{xy}
\usepackage{amssymb,amsmath,amsthm}
\usepackage[mathscr]{eucal}\usepackage[all]{xy}
\usepackage{lmodern}
\usepackage[T1]{fontenc}
\usepackage[utf8]{inputenc}
\usepackage{tikz-cd}
\usepackage[shortlabels]{enumitem}
\usepackage{relsize}
\usepackage{geometry}
\usepackage{mathtools}
\usepackage{adjustbox}
\usepackage{lscape}
\usepackage{diagbox}
\usepackage{slashbox}
\usepackage{bm}
\usepackage{caption}
\usepackage{float}

\usepackage{stmaryrd}
\usepackage{mathrsfs}  
\usepackage{multirow}
\usepackage[hyperfootnotes=false]{hyperref}
\usetikzlibrary{positioning}
\usepackage{tikz-3dplot}
\usepackage{xifthen}
\usetikzlibrary{quotes,arrows.meta}
\tdplotsetmaincoords{60}{125}
\tdplotsetrotatedcoords{0}{0}{0}
\usepackage{xcolor}
\definecolor{mygreen}{rgb}{0.16,.55,0.0}
\usepackage[nottoc]{tocbibind}

\theoremstyle{plain}
\newtheorem{prop}{Proposition}[section]
\newtheorem{conj}[prop]{Conjecture}

\newtheorem{thm}[prop]{Theorem}

\theoremstyle{definition}

\newtheorem{ques}[prop]{Question}

\theoremstyle{plain} 
\newtheorem{prop1}{Proposition}[subsection]
\newtheorem{conj1}[prop1]{Conjecture}
\newtheorem{lem1}[prop1]{Lemma}
\newtheorem{sublem1}[prop1]{Sublemma}
\newtheorem{thm1}[prop1]{Theorem}
\newtheorem{cor1}[prop1]{Corollary}
\theoremstyle{definition}

\newtheorem{rem1}[prop1]{Remark}

\theoremstyle{plain}

\theoremstyle{definition}

\theoremstyle{plain}

\theoremstyle{definition}

\def\plim#1{\displaystyle \lim_{\stackrel{\longleftarrow}{#1}}}
\def\varddots{\mathinner{\raise7pt\vbox{\kern3pt\hbox{.}}\mkern1mu\smash{\raise4pt\hbox{.}}\mkern1mu\smash{\raise1pt\hbox{.}}}}

\newcommand{\vplim}{\varprojlim}
\newcommand{\vilim}{\varinjlim}

\DeclareMathAlphabet{\mathpzc}{OT1}{pzc}{m}{it}
\DeclarePairedDelimiter{\Iouv}{]}{[}

\DeclarePairedDelimiter{\vabs}{\lvert}{\rvert}

\DeclarePairedDelimiter{\scalar}{\langle}{\rangle}
\DeclarePairedDelimiter{\set}{\{}{\}}

\DeclareMathOperator{\Hom}{Hom}
\DeclareMathOperator{\cont}{cont}
\DeclareMathOperator{\Ind}{Ind}

\DeclareMathOperator{\GL}{GL}

\DeclareMathOperator{\Gal}{Gal}

\DeclareMathOperator{\tr}{tr}

\DeclareMathOperator{\Id}{Id}

\DeclareMathOperator{\Rep}{Rep}

\DeclareMathOperator{\gr}{gr}
\DeclareMathOperator{\ad}{ad}

\newcommand{\p}{\mathfrak{p}}

\newcommand{\m}{\mathfrak{m}}

\newcommand{\Z}{{\mathbb Z}}

\newcommand{\cC}{\mathcal{C}}

\newcommand{\QK}[1]{I_{1}}

\newcommand{\Qp}{\mathbb{Q}_{p}}

\newcommand{\Zp}{\mathbb{Z}_{p}}
\newcommand{\F}{\mathbb{F}}
\newcommand{\Fp}{\mathbb{F}_{p}}

\newcommand{\Fq}{\mathbb{F}_{q}}

\newcommand{\Qpbar}{\overline{\mathbb{Q}}_p}
\newcommand{\Cp}{\mathbb{C}_p}
\newcommand{\oCp}{\cO_{\mathbb{C}_p}}

\newcommand{\Fpbar}{\overline{\mathbb{F}}_p}

\newcommand{\val}{\mathrm{val}}

\newcommand{\Q}{\mathbb{Q}}

\newcommand{\R}{\mathbb{R}}

\newcommand{\ang}[1]{\langle #1 \rangle}

\newcommand{\onto}{\twoheadrightarrow}
\newcommand{\into}{\hookrightarrow}
\newcommand{\congto}{\xrightarrow{\,\sim\,}}
\newcommand{\congfrom}{\xleftarrow{\,\sim\,}}

\newcommand{\cO}{\mathcal{O}}
\newcommand{\et}{\acute{\mathrm{e}}\mathrm{t}}

\DeclareMathOperator{\ind}{ind}

\newcommand{\rbar}{\overline{r}}
\newcommand{\rhobar}{\overline{\rho}}

\newcommand{\gp}{{\Gal}(\Qpbar/\Qp)}

\newcommand{\gK}{{\Gal}(\Qpbar/K)}

\newcommand{\gFw}{{\Gal}(\overline F_w/F_w)}
\newcommand{\gF}{{\Gal}(\overline F/F)}
\DeclareMathOperator{\Spf}{Spf}

\newcommand{\oK}{{\mathcal O}_K}

\newcommand{\Ker}{\mathrm{Ker}}\newcommand{\brho}{\overline{\rho}}\newcommand{\smatr}[4]{\bigl(\begin{smallmatrix} {#1}& {#2}\\ {#3}&{#4}\end{smallmatrix}\bigr)}\newcommand{\smat}[1]{\left( \begin{smallmatrix} #1 \end{smallmatrix} \right)}

\newcommand{\defeq}{\stackrel{\textrm{\tiny{\upshape{def}}}}{=}}

\newcommand{\un}[1]{\underline{#1}}

\newcommand{\rig}{\mathrm{rig}}

\newcommand{\ra}{\rightarrow}

\newcommand{\ppar}[1]{(\mkern-3mu(#1)\mkern-3mu)}
\newcommand{\bbra}[1]{\llbracket #1\rrbracket}

\newcommand{\fm}{\mathfrak{m}}

\DeclareMathOperator{\id}{id}

\DeclareMathOperator{\alg}{alg}
\DeclareMathOperator{\Spec}{Spec}

\DeclareMathOperator{\diag}{diag}
\DeclareMathOperator{\rank}{rk}

\renewcommand{\subset}{\subseteq}
\renewcommand{\simeq}{\cong}

\newcommand{\sumfi}{\sum_{i=0}^{f-1}}

\newcommand{\vp}{\varphi}
\newcommand{\okt}{\cO_K^\times}

\newcommand{\ve}{\varepsilon}
\newcommand{\zp}{\Z_p}
\newcommand{\dual}{^\vee}
\newcommand{\rhob}{\overline{\rho}} 
\renewcommand{\o}[1]{\overline{#1}}
\newcommand{\s}{^\times}

\newcommand{\abs}{\vabs{\cdot}}
\newcommand{\B}{\mathbf{B}}
\DeclareMathOperator{\Div}{\mathrm{Div}}
\newcommand{\gen}{\mathrm{gen}}
\newcommand{\LT}{\mathrm{LT}}
\newcommand{\Perf}{\mathrm{Perf}}
\newcommand{\pr}{\mathrm{pr}}
\newcommand{\ronron}{\circ\circ}
\DeclareMathOperator{\Spa}{\mathrm{Spa}}
\DeclareMathOperator{\Tr}{\mathrm{Tr}} 
\theoremstyle{plain}

\title{Multivariable $(\varphi,\cO_K^\times)$-modules and local-global compatibility}

\usepackage{bigfoot}

\DeclareNewFootnote{AAffil}[alph]

\usepackage{etoolbox}
\makeatletter
\patchcmd\maketitle{\def\@makefnmark{\rlap{\@textsuperscript{\normalfont\@thefnmark}}}}{}{}{}
\makeatother

\makeatletter
\def\thanksAAffil#1{\footnotemarkAAffil\protected@xdef\@thanks{\@thanks \protect\footnotetextAAffil[\the \c@footnoteAAffil]{#1}}}
\def\thanksANote#1{\footnotemarkANote \protected@xdef\@thanks{\@thanks \protect\footnotetextANote[\the \c@footnoteANote]{#1}}}
\makeatother

\author{Christophe Breuil\thanksAAffil{CNRS, B\^atiment 307, Facult\'e d'Orsay, Universit\'e Paris-Saclay, 91405 Orsay Cedex, France}\\
\and
Florian Herzig\thanksAAffil{Dept.\ of Math., Univ.\ of Toronto, 40 St.\ George St., BA6290, Toronto, ON M5S 2E4, Canada}\\
\and
Yongquan Hu\thanksAAffil{Morningside Center of Math., No.\ 55, Zhongguancun East Road, Beijing, 100190, China}\\
\and
Stefano Morra\thanksAAffil{Lab.\ d'Analyse, G\'eom\'etrie, Alg\`ebre, 99
  Av.\ Jean Baptiste Cl\'ement, 93430 Villetaneuse, France }$^{,}$\footnotemarkAAffil[6]\\
\and
Benjamin Schraen\thanksAAffil{Institut Camille Jordan, Universit\'e Claude Bernard Lyon I, 69622 Villeurbanne, France}$^{,}$\thanksAAffil{Institut
  Universitaire de France (IUF)}}

\date{ }

\begin{document} 

\maketitle

\setcounter{tocdepth}{3}

\begin{abstract}
Let $p$ be a prime number, $K$ a finite unramified extension of $\Qp$ and $\F$ a finite extension of $\Fp$. Using perfectoid spaces we associate to any finite-dimensional continuous representation $\rhobar$ of ${\rm Gal}(\overline K/K)$ over $\F$ an \'etale $(\varphi,\cO_K^\times)$-module $D_A^\otimes(\rhobar)$ over a completed localization $A$ of $\F\bbra{\cO_K}$. We conjecture that one can also associate an \'etale $(\varphi,\cO_K^\times)$-module $D_A(\pi)$ to any smooth representation $\pi$ of $\GL_2(K)$ occurring in some Hecke eigenspace of the mod $p$ cohomology of a Shimura curve, and that moreover $D_A(\pi)$ is isomorphic (up to twist) to $D_A^\otimes(\rhobar)$, where $\rhobar$ is the underlying $2$-dimensional representation of ${\rm Gal}(\overline K/K)$. Using previous work of the same authors, we prove this conjecture when $\rhobar$ is semi-simple and sufficiently~generic. 
\end{abstract}

\tableofcontents

\newpage

\clearpage{}\section{Introduction}

  Let $p$ be a prime number. The main motivation of this work
  is the investigation of the (hoped for) mod $p$ Langlands
  correspondence for $\GL_2(K)$, where $K$ is a finite unramified
  extension of $\Q_p$. The case $K=\Q_p$ is now well known
  (\cite{breuilI}, \cite{Colmez}, \cite{emerton-local-global}), whereas the case $K\neq\Q_p$ is
  still resisting after more than 10 years (\cite{BP}). An important
  aspect of the $\GL_2(\Q_p)$-case is the construction by Colmez in 
  {\it loc.~cit.}~of an exact functor from the category of admissible finite length mod
  $p$ representations of $\GL_2(\Q_p)$ to the category of finite-dimensional continuous mod $p$ representations of
  $\Gal(\overline{\Q_p}/\Q_p)$. The construction of this functor uses, as an intermediate step, Fontaine's category of
  $(\varphi,\Gamma)$-modules. In a previous article (\cite{BHHMS2}), we
  constructed an exact functor $D_A^{\et}$ from a ``good''
  subcategory of admissible mod $p$ representations of $\GL_2(K)$ to a
  category of \'etale multivariable $(\varphi,\cO_K^\times)$-modules. These multivariable
  $(\varphi,\cO_K^\times)$-modules are $A$-modules with additional
  structures, where $A$ is a ring obtained as a completed localization
  of the Iwasawa algebra of $\cO_K$. In this work we propose a
  construction of a functor $D_A^{\otimes}$ from the category of
  continuous mod $p$ representations of $\Gal(\overline{K}/K)$ to the
  category of \'etale multivariable
  $(\varphi,\cO_K^\times)$-modules. This construction is based on the
  equivalence, also due to Fontaine (\cite{Fo}), between mod $p$
  representations of $\Gal(\overline{K}/K)$ and Lubin--Tate \'etale
  $(\varphi,\cO_K^\times)$-modules. One of the main obstructions to pass from
  Lubin--Tate $(\varphi,\cO_K^\times)$-modules to multivariable
  $(\varphi,\cO_K^\times)$-modules over $A$ lies in the comparison
  between the $\cO_K^\times$-action on $A$ and the
  $\cO_K^\times$-action on (some tensor power of) the structural ring
  of the Lubin--Tate group. To solve this problem, we need to work at a perfectoid level and 
  use the ``Abel--Jacobi map'' considered by Fargues in \cite{FarguesAJ}. We then
  prove, under some conditions, that the two functors $D_A^{\et}$ and
  $D_A^{\otimes}$ satisfy a local-global compatibility property in
  the completed cohomology of a tower of Shimura curves. 
  
  We now describe in more detail the content of this article.

Let $F$ be a totally real number field and let $X_U$ be the smooth projective Shimura curve over $F$ associated to a quaternion algebra $D$ of center $F$ (which splits at one infinite place) and to a compact open subgroup $U$ of $(D\otimes_F{\mathbb A}_F^\infty)^\times$. For $v$ a place of $F$ above $p$ which splits $D$ and $\F$ a finite extension of $\Fp$ (``sufficiently large'', as usual), consider the admissible smooth representation of $\GL_2(F_v)$ over $\F$
\begin{equation}\label{goal}
\pi\defeq \varinjlim_{U_v}\Hom_{{\rm Gal}(\overline F/F)}\!\big(\rbar,H^1_{{\rm \acute et}}(X_{U^vU_v} \times_F \overline F, \F)\big),
\end{equation}
where $U^v$ is a fixed compact open subgroup of $(D\otimes_F{\mathbb A}_F^{\infty,v})^\times$, $\rbar : {\rm Gal}(\overline F/F)\rightarrow \GL_2(\F)$ is an absolutely irreducible continuous Galois representation such that $\pi\ne 0$, and where the inductive limit runs over compact open subgroups $U_v$ of $(D\otimes_FF_v)^\times\cong \GL_2(F_v)$. In this introduction, we moreover assume for simplicity that $v$ is the only $p$-adic place of $F$ and that we are in a ``multiplicity $1$'' situation, which then roughly means that $U^v$ is ``as big as possible'' (in general, one needs to take into account the action of certain operators, which requires mild assumptions on $F$, $D$ and $\rbar$, see (\ref{locpiindef})).

We know that the isomorphism class of $\pi$ always determines the one of $\rbar_v\defeq \rbar\vert_{{\rm Gal}(\overline F_v/F_v)}$, see \cite{BD}, \cite{ScholzeLT}. We also expect that $\pi$ is always of finite length, which is known in several cases, see \cite{HuWang2}, \cite{BHHMS2}. However, the representation $\pi$ is still not understood when $F_v\ne \Qp$, in particular we have the key question:

\begin{ques}\label{hard}
Assume $F_v\ne \Qp$, does $\pi$ only depend on $\rbar_v$?
\end{ques}

Question \ref{hard}, as routine as it may seem at first, has
unfortunately proven to be surprisingly difficult, and there is not
one single instance of a $\pi$ as in (\ref{goal}) for which we know
the answer. For instance the mod $p$ \'etale cohomology of the
Drinfeld tower in dimension $1$, which provides a smooth
representation of $\GL_2(F_v)$ only depending on $\rbar_v$, cannot
give rise to representations like $\pi$ as soon as $F_v\ne \Qp$, see
\cite{CDN} (together with \cite{Benj}, \cite{Wu}). On the other hand,
we know that, for $F_v$ unramified and most $\rbar_v$, the diagram
$(\pi^{I_1}\hookrightarrow \pi^{K_1})$ (where $K_1\defeq
1+pM_2(\cO_{F_v})\subset I_1\defeq$ pro-$p$-Iwahori) only depends on
$\rbar_v$, and this is a really non-trivial fact, see \cite{DoLe}. We
do not answer Question \ref{hard} in this work, but we provide one
further step towards the understanding of the representation $\pi$,
and certainly Question \ref{hard} was a motivation. More precisely, we
completely describe the multivariable \'etale
$(\varphi,\cO_{F_v}^\times)$-module $D_A(\pi)$ associated to $\pi$ in
\cite[\S3]{BHHMS2} when $F_v$ is unramified and $\rbar_v$ is
semi-simple sufficiently generic, in particular we prove that it only
depends on $\rbar_v$, and we provide a precise conjecture on what
$D_A(\pi)$ should be for all $\rbar_v$ (and $F_v$ unramified),
crucially using perfectoid spaces. As an intermediate result, we
construct a new fully faithful functor from continuous representations
of ${\rm Gal}(\overline F_v/F_v)$ over $\F$ to a certain category of
multivariable \'etale $(\varphi_q,\cO_{F_v}^\times)$-modules: this is
the functor $D_A^{(0)}$ constructed in Corollary \ref{descenti}.

Let us first recall the definition of these modules. Let $K$ be a finite unramified extension of $\Qp$ of degree $f\geq 1$, then we can write the Iwasawa algebra $\F\bbra{\cO_K}$ as $\F\bbra{Y_{\sigma},\ \sigma\!:\!\Fq\hookrightarrow \F}$ for $Y_{\sigma}\defeq \sum_{\lambda\in \Fq^\times}\sigma(\lambda)^{-1}[\lambda]\in \F\bbra{\cO_K}$, where $q\defeq p^f$ and $[\lambda]\in \cO_K$ is the multiplicative representative of $\lambda$ (seen in $\F\bbra{\cO_K}$). We then define $A$ to be the completion of $\F\bbra{\cO_K}[1/Y_{\sigma},\ \sigma\!:\!\Fq\hookrightarrow \F]$ for the $(Y_{\sigma})_{\sigma}$-adic topology (in a suitable sense), see \eqref{eq:def:fil:Aq} for the precise definition. In fact $A$ is isomorphic to the Tate algebra $\F\ppar{Y_{\sigma}}\langle (Y_{\sigma'}/Y_{\sigma})^{\pm 1},\sigma'\ne \sigma\rangle$ for any choice of $\sigma$, see Lemma \ref{tatefora}. It is endowed with an $\F$-linear Frobenius $\varphi$ coming from the multiplication by $p$ on $\cO_K$ and with a commuting continuous action of $\cO_K^\times$ coming from its action on $\F\bbra{\cO_K}$ (by multiplication on $\cO_K$). Then an \'etale $(\varphi,\cO_K^\times)$-module over $A$ is by definition a finite free $A$-module endowed with a semi-linear Frobenius $\varphi$ whose image generates everything and a commuting continuous semi-linear action of $\cO_K^\times$. Replacing $\varphi$ on $A$ by $\varphi_q\defeq \varphi^f$, we define in the same way \'etale $(\varphi_q,\cO_K^\times)$-modules over $A$. When $f=1$, the two definitions recover Fontaine's classical $(\varphi,\Zp^\times)$-modules (or $(\varphi,\Gamma)$-modules) in characteristic $p$.

Now let $\pi$ be an admissible smooth representation of ${\rm GL}_2(\cO_K)$ over $\F$. We endow $\pi^\vee\defeq \Hom_\F(\pi,\F)$ with the $\m_{I_1}$-adic topology, where $\m_{I_1}$ is the maximal ideal of the Iwasawa algebra $\F\bbra{I_1}$. In particular we can see $\pi^\vee$ as an $\F\bbra{\cO_K}$-module via $\F\bbra{\cO_K}\simeq \F\bbra{\smatr {1}{\mathcal{O}_K}{0}{1}}\subset \F\bbra{I_1}$. We define
\[D_A(\pi)\defeq \big(\F\bbra{\cO_K}[1/Y_{\sigma},\sigma\!:\!\Fq\hookrightarrow \F]\otimes_{\F\bbra{\cO_K}}\pi^\vee\big)^\wedge,\]
where the completion is for the tensor product topology, see \cite[\S 3.1.1]{BHHMS2} or \S\ref{conjecture}. Even though $D_A(\pi)$ is an $A$-module endowed with a semi-linear action of $\cO_K^\times$ (coming from the action of $\smatr {\cO_K^\times}{0}{0}{1}$ on $\pi^\vee$), it is not clear if it has good properties in general (it might not have a Frobenius $\varphi$, it might not be of finite type, etc.). But we know that $D_A(\pi)$ is an \'etale $(\varphi,\cO_K^\times)$-module of rank $2^f$ for some of the $\pi$ in (\ref{goal}) when $K\defeq F_v$ is unramified, see \cite[\S 1.3]{BHHMS2}\footnote{Note that, with the notation of \cite[\S 3.1.2]{BHHMS2}, $D_A(\pi)$ is equal to its \'etale quotient $D_A(\pi)^{\et}$ in our case, see \cite[Rem.~3.3.5.4(ii)]{BHHMS2}.} together with Remark \ref{freeness}. In fact we conjecture in this paper that $D_A(\pi)$ is always an \'etale $(\varphi,\cO_K^\times)$-module over $A$ (hence equal to $D_A(\pi)^{\et}$) of rank $2^f$ for all representations $\pi$ in (\ref{goal}) (when $F_v$ is unramified).

On the Galois side, for $\rhobar: {\rm Gal}(\overline
K/K)\!\rightarrow \GL_n(\F)$ ($n\geq 1$) a continuous representation
and $\sigma\!:\!\Fq\hookrightarrow \F$ we can associate to $\rhobar$ a
Lubin--Tate $(\varphi_q,\cO_K^\times)$-module. Recall that it is an
$n$-dimensional $\F\ppar{T_{K,\sigma}}$-vector space
$D_{K,\sigma}(\rhobar)$ equipped with a semi-linear endomorphism
$\varphi_q$ whose image generates $D_{K,\sigma}(\rhobar)$ and a
commuting continuous action of $\cO_K^\times$. Here $\varphi_q$ is
$\F$-linear and satisfies $\varphi_q(T_{K,\sigma})=T_{K,\sigma}^q$,
and the action of $\cO_K^\times$ on $\F\ppar{T_{K,\sigma}}$ is given
by the Lubin--Tate power series associated to the choice of logarithm
$\sum_{n\geq 0}p^{-n}T_{K,\sigma}^{q^n}$ composed with
$\sigma\!:\!\Fq\hookrightarrow \F$ on the coefficients. Recall we have
\[ \F\ppar{T_{K,\sigma}}\otimes_{\F\ppar{T_{K,\sigma}^{q-1}}}~D_{K,\sigma}(\rhobar)^{[\Fq^\times]}\buildrel\sim\over\rightarrow ~D_{K,\sigma}(\rhobar).\]

Assume now that $\rhobar$ is a direct sum of absolutely irreducible representations and define
\begin{equation}\label{recipeintro}
D_{A,\sigma}(\rhobar)\defeq A\otimes_{\F\ppar{T_{K,\sigma}^{q-1}}} D_{K,\sigma}(\rhobar)^{[\Fq^\times]},
\end{equation}
where the embedding $\F\ppar{T_{K,\sigma}^{q-1}}\hookrightarrow A$ sends $T_{K,\sigma}^{q-1}$ to $\varphi(Y_{\sigma})/Y_{\sigma}\in A$. We endow $D_{A,\sigma}(\rhobar)$ with $\varphi_q\defeq \varphi^f\otimes \varphi_q$. The embedding $\F\ppar{T_{K,\sigma}^{q-1}}\hookrightarrow A$ does not commute with $\cO_K^\times$, but one easily checks that, when $\rhobar$ is a direct sum of absolutely irreducible representations, there exists a unique (in a certain sense) continuous semi-linear action of $\cO_K^\times$ on $D_{A,\sigma}(\rhobar)$ which commutes with $\varphi_q$ and makes $D_{A,\sigma}(\rhobar)$ an \'etale $(\varphi_q,\cO_K^\times)$-module over $A$ of rank $\dim_{\F}\rhobar$, see Lemma \ref{PSnatural}. Moreover there is a canonical isomorphism $\id\otimes \varphi:A\otimes_{\varphi,A}D_{A,\sigma\circ \varphi}(\rhobar)\buildrel\sim\over\longrightarrow D_{A,\sigma}(\rhobar)$ of \'etale $(\varphi_q,\cO_K^\times)$-modules over $A$, where $\sigma\circ \varphi\defeq \sigma((-)^p)$. We then define:
\begin{equation}\label{daotimesintro}
D_A^\otimes(\rhobar)\defeq \bigotimes_{A,\sigma:\Fq\hookrightarrow \F} D_{A,\sigma}(\rhobar)
\end{equation}
endowed with the ``diagonal'' action of $\cO_K^\times$. Using the isomorphism $\id\otimes \varphi$, we can define a canonical endomorphism $\varphi:D_A^\otimes(\rhobar)\rightarrow D_A^\otimes(\rhobar)$ which cyclically permutes the factors $D_{A,\sigma}(\rhobar)$, is semi-linear with respect to $\varphi$ on $A$ and is such that $\varphi^f=\varphi_q\otimes\cdots\otimes\varphi_q$. It is then clear that $D_A^\otimes(\rhobar)$ is an \'etale $(\varphi,\cO_K^\times)$-module over $A$ of rank $(\dim_{\F}\rhobar)^f$. The following theorem is our main result:

\begin{thm}[Corollary \ref{maincor}]\label{mainintro}
Assume that $\rbar_v$ is semi-simple and sufficiently gene\-ric {\upshape(}see (\ref{stronggen}){\upshape)}, and assume standard technical assumptions on the global setting {\upshape(}see \S\ref{conjecture} for precise statements{\upshape)}. Then there is an isomorphism of \'etale $(\varphi,\cO_K^\times)$-modules $D_A(\pi)\simeq D_A^\otimes(\rbar_v(1))$ over $A$, where $\rbar_v(1)$ is the usual Tate twist of $\rbar_v$.
\end{thm}

The proof of Theorem \ref{mainintro} is a long explicit computation of the dual \'etale $(\varphi,\cO_K^\times)$-module $\Hom_A(D_A(\pi),A)$. Let us briefly indicate the various steps. We first describe $\Hom_\F^{\cont}(D_A(\pi),\F)$, which is not so hard, see Proposition \ref{prop:F-dual}. We then prove that there is a canonical injection
\[\Hom_A(D_A(\pi),A)\hookrightarrow \Hom_\F^{\cont}(D_A(\pi),\F)\]
induced by a nonzero continuous morphism $\mu:A\rightarrow \F$ uniquely determined (up to scalar in $\F^\times$) by the condition $\mu\circ \psi\in \F^\times\mu$, where $\psi:A\rightarrow A$ is a certain canonical left inverse of $\varphi$, see Lemma \ref{lm:psi}, Proposition \ref{prop:mu} and (\ref{eq:mu-star}). To each Serre weight $\sigma$ of $\rbar_v^\vee$ we then associate in (\ref{eq:13f}) a certain projective system $x_\sigma=(x_{\sigma,k})_{k\geq 0}$, where $x_{\sigma,k}\in \pi[\m_{I_1}^{kf+1}]$, and we prove via Proposition \ref{prop:F-dual} that $x_\sigma$ lies in $\Hom_\F^{\cont}(D_A(\pi),\F)$, see Lemma \ref{lem:a(n)} and Proposition \ref{prop:degree-of-phi-n-v}. Then the key calculation is to prove that $x_\sigma$ actually also lies in the submodule $\Hom_A(D_A(\pi),A)$, and that the $2^f$-tuple $(x_\sigma)_{\sigma\in W(\rbar_v^\vee)}$ even forms an $A$-basis of the free $A$-module $\Hom_A(D_A(\pi),A)$, see Theorem \ref{thm:descent}. For that we prove a crucial finiteness result (Proposition \ref{prop:cond-ii}) using the technical -- but important -- computations in \cite[\S 3.2]{BHHMS2} that we need to strengthen, see \S\ref{ss:degree}. Once all this is done, it is easy to derive the explicit actions of $\varphi$ and $\cO_K^\times$ on $\Hom_A(D_A(\pi),A)$, see Propositions~\ref{prop:phi-x} and \ref{prop:gamma-action}. We can then at last compare the two $(\varphi,\cO_K^\times)$-modules $D_A(\pi)$ and $D_A^\otimes(\rbar_v(1))$ and prove that they are isomorphic, see Theorem \ref{main}. The same proof works verbatim for quaternion algebras $D$ which are definite at all infinite places (and split at $v$) and the representations $\pi$ of $\GL_2(K)=\GL_2(F_v)$ defined analogously to (\ref{goal}). 

There is no doubt to us that there should exist a more conceptual proof of Theorem \ref{mainintro} which will hopefully avoid both the genericity assumptions on $\rbar_v$ and the technical computations. At present however, we do not know how to do this. But the first issue is to find a more conceptual definition of $D_{A,\sigma}(\rhobar)$ and of $D_A^\otimes(\rhobar)$. Indeed, when $\rhobar$ is not semi-simple, the recipe (\ref{recipeintro}) does not work in general because there might not always exist a continuous semi-linear action of $\cO_K^\times$ on $A\otimes_{\F\ppar{T_{K,\sigma}^{q-1}}} D_{K,\sigma}(\rhobar)^{[\Fq^\times]}$ which commutes with $\varphi^f\otimes \varphi_q$ (or such an action might not be unique) {, see for instance \cite[\S 4]{YitongWang3}}. Using perfectoid spaces we give below a functorial construction of an \'etale $(\varphi_q,\cO_K^\times)$-module $D_{A,\sigma}(\rhobar)$, and subsequently of an \'etale $(\varphi,\cO_K^\times)$-module $D_A^\otimes(\rhobar)$, which works for all $\rhobar$.

The first step is to replace the ring $A$ by its perfectoid version
\begin{equation}\label{tateintro}
A_\infty\defeq \F\ppar{Y_{\sigma}^{1/p^{\infty}}}\left\langle (Y_{\sigma'}/Y_{\sigma})^{\pm {1/p^{\infty}}},\sigma'\ne \sigma\right\rangle
\end{equation}
which is a perfectoid Tate algebra over the perfectoid field
$\F\ppar{Y_{\sigma}^{1/p^{\infty}}}$ (for any $\sigma$). Using the
equivalence between finite \'etale $A$-algebras and finite \'etale
$A_\infty$-algebras together with the equivalence between  {locally constant \'etale sheaves of finite-dimensional $\Fq$-vector spaces on $\Spec(R)$ and finite projective $R$-modules with an action of Frobenius for perfect rings $R$ over $\Fq$,}
it is not hard to check that the extension of scalars $(-)\mapsto (-)\otimes_AA_\infty$ induces an equivalence of categories between \'etale $(\varphi_q,\cO_K^\times)$-modules over $A$ and \'etale $(\varphi_q,\cO_K^\times)$-modules over $A_\infty$, and similarly with $(\varphi,\cO_K^\times)$ instead of $(\varphi_q,\cO_K^\times)$, see Corollary \ref{descentoK}. Hence we may as well look for a definition of $D_{A_\infty,\sigma}(\rhobar)$ and $D_{A_\infty}^\otimes(\rhobar)$.

It is now convenient to fix an embedding $\sigma_0:\Fq\hookrightarrow \F$ and set $\sigma_i\defeq \sigma_0\circ\varphi^i$ for $i\in \Z$. The second step is to consider the two perfectoid spaces 
\begin{eqnarray*}
Z_{\LT}&\defeq &\underbrace{\Spa\big(\F\ppar{T_{K,\sigma_0}^{1/p^\infty}},\F\bbra{T_{K,\sigma_0}^{1/p^\infty}}\big)\times_{\Spa(\F)}\cdots \times_{\Spa(\F)}\Spa\big(\F\ppar{T_{K,\sigma_0}^{1/p^\infty}},\F\bbra{T_{K,\sigma_0}^{1/p^\infty}}\big)}_{f\mathrm{\ times}}\\
Z_{\cO_K}&\defeq &\Spa\Big(\F\bbra{Y_{\sigma_0}^{1/p^\infty},\dots,Y_{\sigma_{f-1}}^{1/p^\infty}},\F\bbra{Y_{\sigma_0}^{1/p^\infty},\dots,Y_{\sigma_{f-1}}^{1/p^\infty}}\Big)\setminus V(Y_{\sigma_0},\dots, Y_{\sigma_{f-1}}),
\end{eqnarray*}
where $Z_{\LT}$ is endowed with an obvious action of $(K^\times)^f\rtimes\mathfrak{S}_f$ ($p\in K^\times$ acting via $\varphi_q$ which is now bijective) and $Z_{\cO_K}$ is endowed with an action of $K^\times$ ($p$ acting via $\varphi$). It turns out that there is a morphism of perfectoid spaces (see the beginning of \S\ref{mapm})
\[m:Z_{\LT}\longrightarrow Z_{\cO_K}\]
such that $m \circ ((a_0,\dots,a_{f-1}),w) = (\prod_i a_i) \circ m$ for $a_i\in K^\times$ and $w\in \mathfrak{S}_f$, a crucial fact that we learnt from \cite{FarguesAJ}. Indeed, the sheaf on the perfectoid $v$-site over $\F$ represented by $Z_{\LT}$ sends a perfectoid $\F$-algebra $R$ to a subset of $(\B^+(R)^{\varphi_q=p})^f$ stable under multiplication, where $\B^+(R)$ is the (relative version of the) ring defined in \cite[\S 1.10]{FF} (a certain completion of $W(R^\circ)[1/p]$, where $R^\circ\subseteq R$ is the subring of power-bounded elements). Likewise, the sheaf represented by $Z_{\cO_K}$ sends $R$ to a subset of $\B^+(R)^{\varphi_q=p^f}$ stable under multiplication, see \S\ref{reminderp}. The map $m$ then is induced by the product map $(\B^+(R)^{\varphi_q=p})^f\rightarrow \B^+(R)^{\varphi_q=p^f}$ in the ring $\B^+(R)$, which satisfies the above relation with respect to the various group actions.

Note that $\Spa(A_\infty, A_\infty^\circ)$ is an affinoid open subspace of $Z_{\cO_K}$ by (\ref{tateintro}). Let $\Delta\defeq \{(a_0,\dots,a_{f-1})\in (K^\times)^f,\ \prod a_i =1\}$ and $\Delta_1\defeq \Delta \cap (\cO_K^\times)^f$. The third step is to prove that the morphism $m$ induces a commutative diagram of perfectoid spaces over $\F$:
\begin{equation*}
\begin{gathered}
\xymatrix{Z_{\LT}\ar@{<-^{)}}[r]\ar^{m}[d] & \ m^{-1}(\Spa(A_\infty, A_\infty^\circ))\ar^{m}[d] \simeq & \!\!\!\!\!\!\!\!\!\!\!\!\!\underline{(\Delta/\Delta_1)\rtimes \mathfrak{S}_f}\times \Spa(A'_\infty,(A'_\infty)^\circ)\ar[ld]\\
Z_{\cO_K}\ar@{<-^{)}}[r] & \ \Spa(A_\infty, A_\infty^\circ)}, 
\end{gathered}
\end{equation*}
where \ the \ middle \ vertical \ morphism \ is \ a \ pro-\'etale \ $\Delta\rtimes \mathfrak{S}_f$-torsor \ and \ where $\Spa(A'_\infty,(A'_\infty)^\circ)$ is an explicit affinoid open subspace of $Z_{\LT}$ preserved by the action of $\Delta_1$ which is itself a pro-\'etale $\Delta_1$-torsor over $\Spa(A_\infty, A_\infty^\circ)$, see Proposition \ref{prop:Delta1_torsor}, Corollary \ref{prop:DeltaS_torsor} and Lemma \ref{lemm:norm}.

Now let $\rhobar$ be any finite-dimensional continuous representation of ${\rm Gal}(\overline K/K)$ over $\F$, then $\F\ppar{T_{K,\sigma_0}^{1/p^\infty}}\otimes_{\F\ppar{T_{K,\sigma_0}}}D_{K,\sigma_0}(\rhobar)$ is the space of global sections of a $K^\times$-equivariant vector bundle $\mathcal{V}_{\rhobar}$ on
$\Spa(\F\ppar{T_{K,\sigma_0}^{1/p^\infty}},\F\bbra{T_{K,\sigma_0}^{1/p^\infty}})$. For $i\in \{0,\dots, f-1\}$ we define $\mathcal{V}_{\rhobar}^{(i)}\defeq \pr_i^*\mathcal{V}_{\rhobar}$, where $\pr_i:Z_{\LT}\rightarrow \Spa(\F\ppar{T_{K,\sigma_0}^{1/p^\infty}},\F\bbra{T_{K,\sigma_0}^{1/p^\infty}})$ is the $i$-th projection. Then $\mathcal{V}_{\rhobar}^{(i)}$ is a $(K^\times)^f$-equivariant vector bundle on $Z_{\LT}$, and thus $\mathcal{V}_{\rhobar}^{(i)}\vert_{\Spa(A'_\infty,(A'_\infty)^\circ)}$ is a $\Delta_1$-equivariant vector bundle on $\Spa(A'_\infty,(A'_\infty)^\circ)$. By the third step and using \cite[Lemma 17.1.8]{SWBerkeley}, we deduce that $\Gamma(\Spa(A'_\infty,(A'_\infty)^\circ),\mathcal{V}_{\rhobar}^{(i)})^{\Delta_1}$ is an \'etale $(\varphi_q,\cO_K^\times)$-module over $A_\infty$ of rank $\dim_{\F}\rhobar$, see Theorem \ref{descentdelta1} and \S\ref{arbitrary}. Hence by the first step $\Gamma(\Spa(A'_\infty,(A'_\infty)^\circ),\mathcal{V}_{\rhobar}^{(i)})^{\Delta_1}$ is the extension of scalars of a unique \'etale $(\varphi_q,\cO_K^\times)$-module $D_A^{(i)}(\rhobar)$ over $A$ of rank $\dim_{\F}\rhobar$.

The following theorem sums up the main properties of the functor $\rhobar\mapsto D_{A}^{(i)}(\rhobar)$.

\begin{thm}\label{sumsupintro}
Let $i\in \{0,\dots,f-1\}$.
\begin{enumerate}\item There is a functorial $A$-linear isomorphism $\phi_i:A\otimes_{\varphi,A}D_{A}^{(i)}(\rhobar)\buildrel\sim\over\longrightarrow D_{A}^{(i+1)}(\rhobar)$ which commutes with $(\varphi_q,\cO_K^\times)$ and is such that $\phi_{f-1}\circ \phi_{f-2}\circ \cdots \circ \phi_0:A\otimes_{\varphi^f,A}D_{A}^{(0)}(\rhobar)\buildrel\sim\over\longrightarrow D_{A}^{(0)}(\rhobar)$ is $\id\otimes \varphi_q$, see Corollary \ref{descenti}.
\item The functor $\rhobar\mapsto D_{A}^{(i)}(\rhobar)$ from finite-dimensional continuous representations of ${\rm Gal}(\overline K/K)$ over $\F$ to \'etale $(\varphi_q,\cO_K^\times)$-modules over $A$ is exact and fully faithful, see Corollary \ref{fully}.
\item There is $d\in \{0,\cdots, f-1\}$ such that the surjection $A\twoheadrightarrow\F\ppar{T}$ induced by the trace $\F\bbra{\cO_K}\twoheadrightarrow\F\bbra{\Zp}\simeq\F\bbra{T}$ gives a functorial isomorphism of $(\varphi_q,\Zp^\times)$-modules
\[\F\ppar{T}\otimes_{A}D_{A}^{(i)}(\rhobar)\simeq D_{\sigma_{d-i}}(\rhobar),\]
where $D_{\sigma_{d-i}}(\rhobar)$ is the usual {\upshape(}cyclotomic{\upshape)} $(\varphi_q,\Zp^\times)$-module over $\F\ppar{T}$ associated to $\rhobar$ using $\sigma_{d-i}$ to embed $\Fq$ into $\F$, see Proposition \ref{lift} and Remark \ref{shift}.
\item If $\rhobar$ is a direct sum of absolutely irreducible representations then there is an isomorphism of $(\varphi_q,\cO_K^\times)$-modules over $A$
\[D_A^{(i)}(\rhobar)\simeq D_{A,\sigma_{f-i}}(\rhobar),\]
where $D_{A,\sigma_{f-i}}(\rhobar)$ is as in (\ref{recipeintro}), see Theorem \ref{explicitK}.
\end{enumerate}
\end{thm}

Because of Theorem \ref{sumsupintro}(iv) it is natural to rename $D_A^{(i)}(\rhobar)$ as $D_{A,\sigma_{f-i}}(\rhobar)$ for any $\rhobar$. Using Theorem \ref{sumsupintro}(i) we can then associate to any $\rhobar$ an \'etale $(\varphi,\cO_K^\times)$-module $D_A^\otimes(\rhobar)$ over $A$ of rank $(\dim_{\F}\rhobar)^f$ by exactly the same formula as in (\ref{daotimesintro}). Note that by Theorem \ref{sumsupintro}(iii) $\F\ppar{T}\otimes_{A}D_{A}^\otimes(\rhobar)$ can be identified with the $(\varphi,\Gamma)$-module of the tensor induction from $K$ to $\Qp$ of $\rhobar$. 

We can now state our conjecture:

\begin{conj}[Conjecture \ref{dAPi}]\label{conjintro}
For any $\pi$ as in (\ref{goal}) {\upshape(}with $F_v=K$ unramified{\upshape)} there is an isomorphism of \'etale $(\varphi,\cO_K^\times)$-modules $D_A(\pi)\simeq D_A^\otimes(\rbar_v(1))$ over $A$.
\end{conj}

By Theorem \ref{sumsupintro}(iv) we see that Theorem \ref{mainintro} proves special cases of Conjecture \ref{conjintro} (but recall that our somewhat technical proof of Theorem \ref{mainintro} does not use perfectoids). Note that Conjecture \ref{conjintro} implies (the analogue of) \cite[Conjecture 1.2.5]{BHHMS2} for the representations $\pi$ in (\ref{goal}). It is also reminiscent of the plectic structure of the local Galois action at $p$ on the $\ell$-adic cohomology ($\ell\ne p$) of certain Shimura varieties recently proven in \cite{LiHu}, where the above map $m$ also plays a key role.

We finish this introduction by going back to Question \ref{hard} assuming Conjecture \ref{conjintro}. The image of the natural map $\pi^\vee\rightarrow D_A(\pi)\simeq D_A^\otimes(\rbar_v(1))$ is a compact $\F\bbra{\cO_K}$-submodule $D_A(\pi)^\natural$ which generates $D_A^\otimes(\rbar_v(1))$ over $A$ and is preserved by $\cO_{K}^\times$ and the operator $\psi$, with moreover $\psi:D_A(\pi)^\natural\twoheadrightarrow D_A(\pi)^\natural$ surjective. Assuming there is an admissible smooth representation of $\GL_2(K)$ naturally associated to $\rbar_v$, and that this representation is $\pi$ (as is the case when $K=\Qp$), one could hope to ``guess'' what $D_A(\pi)^\natural$ is inside $D_A^\otimes(\rbar_v(1))$, as the latter is pretty explicit, at least when $\rbar_v$ is semi-simple and sufficiently generic. However, even in the simplest case where $K$ is quadratic (unramified) and $\rbar_v$ is the direct sum of two characters, where we know that $\pi$ is semi-simple (\cite{BHHMS2}), it seems impossible to find $D_A(\pi)^\natural$ ``by hand'' (there exists a natural explicit generating compact $\F\bbra{\cO_{K}}$-submodule in $D_A^\otimes(\rbar_v(1))$ which is preserved by $\cO_K^\times$ and $\psi$ with $\psi$ surjective, but we can prove that it cannot be $D_A(\pi)^\natural$). Going back to perfectoids, one could hope to find instead a natural $\F\bbra{Y_{\sigma_0}^{1/p^\infty},\dots,Y_{\sigma_{f-1}}^{1/p^\infty}}$-submodule $D_{A_\infty}(\pi)^\natural$ inside $D_{A_\infty}^\otimes(\rbar_v(1))=A_\infty\otimes_AD_{A}^\otimes(\rbar_v(1))$ and from there go to $D_A(\pi)^\natural$ in a similar way as what was done by Colmez when $K=\Qp$ in \cite[\S IV.2]{Colmez2}. However, even though there is a natural candidate, namely the $\F\bbra{Y_{\sigma_0}^{1/p^\infty},\dots,Y_{\sigma_{f-1}}^{1/p^\infty}}$-submodule
\begin{multline*}
\Gamma\big(Z_{\LT}, \mathcal{V}_{\rbar_v(1)}^{(0)}\otimes_{\cO_{Z_{\LT}}}\cdots\otimes_{\cO_{Z_{\LT}}}\mathcal{V}_{\rbar_v(1)}^{(f-1)}\big)^{\Delta\rtimes \mathfrak{S}_f}\\
\subseteq \Gamma\big(m^{-1}(\Spa(A_\infty, A_\infty^\circ)), \mathcal{V}_{\rbar_v(1)}^{(0)}\otimes_{\cO_{Z_{\LT}}}\cdots\otimes_{\cO_{Z_{\LT}}}\mathcal{V}_{\rbar_v(1)}^{(f-1)}\big)^{\Delta\rtimes \mathfrak{S}_f}\simeq D_{A_\infty}^\otimes(\rbar_v(1)),
\end{multline*}
computations for $f=2$ show no evidence for this submodule to be large enough (or even nonzero when $\rbar_v$ is irreducible).

We fix some general notation (most of which has already been introduced above, but we remind the reader). We fix $K$ a finite unramified extension of $\Qp$ of residue field $\Fq= \F_{p^f}$, so $\cO_K=W(\F_q)$ and $K=\cO_K[1/p]$. We normalize the local reciprocity map so that it sends $p\in K^\times$ to (the image of) the geometric Frobenius $x\mapsto x^{-q}$. We fix an algebraic closure $\overline K$ of $K$ with ring of integers ${\mathcal O}_{\overline K}$ and maximal ideal $\m_{\overline K}$. We denote by $\F$ the coefficients, which is a finite extension of $\Fq$ that we always tacitly assume to be ``large enough''. We fix an embedding $\sigma_0 : \Fq\hookrightarrow\F$ (which is sometimes omitted from the notation when the context is clear) and we let $\sigma_i\defeq \sigma_0\circ\varphi^i$ for $\varphi$ the Frobenius on $\Fq$ (i.e.\ $\varphi(x)=x^p$) and $i\in \Z$.

\textbf{Acknowledgements}: We thank Laurent Berger, Xavier Caruso, K{\k e}stutis {\v C}esnavi{\v c}ius and Ofer Gabber for their answers to our questions,  {Mohamed Boubakeur for his help in computer-assisted computations} and Laurent Berger,  {Changjiang Du}, Dat Pham, Yitong Wang for their comments. We heartily thank Laurent Fargues for many discussions and patient explanations on the perfectoid world and \cite{FarguesAJ}.

Y.\;H.\ is \ partially \ supported by \ National \ Key \ R$\&$D \ Program \ of \ China \ 2020YFA0712600, National Natural Science Foundation of China Grants 12288201 and 12425103, National Center for Mathematics and Interdisciplinary Sciences and Hua Loo-Keng Key Laboratory of Mathematics, Chinese Academy of Sciences. F.\;H.\ is partially supported by an NSERC grant. S.\;M.\ and B.\;S.\ are partially supported by Institut Universitaire de France.  {During the genesis and writing of this paper, C.\;B., S.\;M.\ and B.\;S.\ were} members of the A.N.R.\ project CLap-CLap ANR-18-CE40-0026.

\clearpage{}

\clearpage{}\section{\'Etale \texorpdfstring{$(\varphi,\!\oK^\times)$}{(phi,O\_K\^{}x)}-modules and Galois representations}\label{chapterGalois}

In this section we functorially associate to any finite-dimensional continuous representation of ${\rm Gal}(\overline K/K)$ over $\F$ an \'etale $(\varphi_q,\!\oK^\times)$-module $D_{A,\sigma}(\rhobar)$ of rank $\dim_{\F}\rhobar$ over the ring $A$ of \cite[\S 3.1.1]{BHHMS2} (depending on an embedding $\sigma:\Fq\hookrightarrow \F$) and an \'etale $(\varphi,\!\oK^\times)$-module $D_{A}^\otimes(\rhobar)$ of rank $(\dim_{\F}\rhobar)^f$ over $A$. We prove various properties of these modules and we make them explicit when $\rhobar$ is a direct sum of absolutely irreducible representations.

If $X$ is an adic space over $\F$, we denote by $h_X$ the functor $\Hom_{\Spa(\F,\F)}(-,X)$ from the category of adic spaces over $\F$ to the category of sets. If $R$ is an adic Huber ring, i.e.\ a topological ring whose topology is $I$-adic for a finitely generated ideal $I$ (see for instance \cite[\S 2.2]{SWBerkeley}), we use the shorthand $\Spa(R)$ for the adic spectrum $\Spa(R,R)$. We denote by $\Perf_{\F}$ the category of perfectoid spaces over $\F$. For background on adic spaces or perfectoid spaces we refer (mostly without comment) to \cite{Huberbook}, \cite{Scholzeperfectoid} or \cite{SWBerkeley}.

Let $A$ be a (commutative) ring and let $\varphi$ be a ring
endomorphism of $A$. We
define a \emph{$\varphi$-module over $A$} as a finite free
$A$-module $D$ endowed with a $\varphi$-semi-linear map
$\varphi : D\rightarrow D$. We say that a $\varphi$-module over $A$ is
\emph{\'etale} if the $A$-linear map
$\id_A\otimes\varphi : A\otimes_{\varphi,A}D\rightarrow D$ is an
isomorphism. Assume moreover that $A$ is a topological ring and that
there exists a continuous action of an abelian topological group
$\Gamma$ on $A$ via endomorphisms commuting with $\varphi$. We
define a \emph{$(\varphi,\Gamma)$-module over $A$} as a $\varphi$-module
$D$ over $A$ endowed with a continuous semi-linear action of
$\Gamma$ such that, for $a\in A$, $v\in D$ and $\gamma\in\Gamma$:
\[
\varphi(\gamma(v))=\gamma(\varphi(v)).
\] 
Moreover we say that a $(\varphi,\Gamma)$-module is
\emph{\'etale} if its underlying $\varphi$-module over $A$ is so.

Let $\Rep_{\F}\gK$ denote the category of $\F$-linear continuous
representations of the topological group $\gK$ on finite-dimensional
$\F$-vector spaces.

\subsection{Review of Lubin--Tate and classical \texorpdfstring{$(\varphi,\Gamma)$}{(phi,Gamma)}-modules}\label{LT}

We \ review Lubin--Tate \ and \ classical \ $(\varphi,\Gamma)$-modules
\ associated \ to \ an \ object \ of $\Rep_{\F}\gK$.

Let $\oCp$ be the $p$-adic completion of
${\mathcal O}_{\overline K}$ and let $\Cp\defeq\oCp[1/p]$. Set
\[\oCp^\flat\defeq \plim{x\mapsto x^p}\oCp/(p)\simeq \plim{x\mapsto
    x^p}{\mathcal O}_{\overline K}/(p)\simeq \plim{x\mapsto
    x^q}{\mathcal O}_{\overline K}/(p)\]
(which Fontaine used to denote by $R$ in \cite[\S2.1]{Fontaine_potBT}). Note that there is an isomorphism of (multiplicative) monoids $\oCp^\flat\buildrel\sim\over\rightarrow\plim{x\mapsto
  x^p}\oCp$. This allows us to define a map $v :
\oCp^\flat\rightarrow\Z\cup\set{+\infty}$ by $v((x_m)_{m\geq
  1})\defeq \val(x_1)$, where $x_m\in \Cp$, $x_m^p=x_{m-1}$ (for
$m>1$)
and val is the usual $p$-adic valuation on $\Cp$ normalized by
$\val(p)=1$. Then $v$ is a valuation on $\oCp^\flat$ and extends
therefore to a valuation on $\Cp^\flat\defeq \plim{x\mapsto x^p}\Cp$. Then $\Cp^\flat$ is an
algebraically closed field of characteristic $p$ which is complete
with respect to the valuation $v$. Moreover its ring of integers
$\{x\in \Cp^\flat,\ v(x)\geq 0\}$ is $\oCp^\flat$  and $\Cp^\flat
\simeq {\rm Frac}(\oCp^\flat)$. There is an action of $\gp$, hence of $\gK$, on $\Cp^\flat$ which preserves $\oCp^\flat$.

We denote by $G_{\LT}$ the unique (up to isomorphism) Lubin--Tate
formal $\cO_K$-module over $\cO_K$ associated to the uniformizer
$p$. Let $T_K$ be a formal variable of $G_{\LT}$. The structure of
$\cO_K$-module on $G_{\LT}$ is given by power series:
\[a_{\rm LT}(T_K)\in aT_K+T_K^2\oK\bbra{T_K}\ {\rm for}\ a\in \oK,\]
and recall that $p_{\LT}(T_K)\in T_K^q+p\cO_K\bbra{T_K}$. Let
$T_pG_{\LT}\defeq\varprojlim_{m\geq1}G_{\LT}[p^m](\cO_{\overline{K}})$
be the Tate module of $G_{\LT}$, which is a free $\cO_K$-module
of rank $1$. Let $u$ be a generator of the $\cO_K$-module $T_p G_{\LT}$.
We can write $u=(u_m)_{m\geq1}$ with
$u_m\in G_{\LT}[p^m](\cO_{\overline{K}})\subset\oCp$ for $m\geq1$,
where we embed $G_{\LT}(\oCp)$ into $\oCp$ using $T_K$. For $m\geq1$
let $\overline{u}_m$ be the image of $u_m$ in $\oCp/(p)$ and
$\overline{u}\defeq(\overline{u}_m)_{m\geq1}\in\oCp^\flat$. The map
$\F_q\bbra{T_K}\rightarrow\oCp^\flat$ sending $T_K$ to $\overline{u}$
is injective, and we use it to identify $\F_q\bbra{T_K}$ with a subring of
$\oCp^\flat$.

We denote by $K_\infty$ the abelian extension of $K$ generated by all
the elements $u_m$ and recall that we have the commutative diagram:
\begin{equation}\label{LTp}
\begin{gathered}
\xymatrix{\Gal(\overline K/K)\ar@{^{}->>}[r] &\Gal(\overline K/K)^{\rm ab}\ar@{^{}->>}[r] & \Gal(K_\infty/K)\ar@{^{}->>}[r] &\Gal(K(\!\sqrt[p^\infty]{1})/K)\\
&K^\times\simeq p^{\Z}\times \oK^\times\ar@{^{(}->}[u] \ar@{^{}->>}[r]&\oK^\times\ar^{\wr}[u]\ar@{^{}->>}[r]&\Zp^\times\ar^{\wr}[u]}
\end{gathered}
\end{equation}
where the left vertical injection is the local reciprocity map, the
bottom left horizontal surjection is the projection sending $p$ to
$1$, and the bottom right horizontal surjection is the norm
map. 

{We endow the topological ring $\F\otimes_{\Fp}\Fq\ppar{T_K}$ with a continuous $\F$-linear endomorphism $\varphi$ and a continuous
$\F\otimes_{\Fp}\Fq$-linear action of $\oK^\times$ commuting with $\varphi$ and satisfying the following conditions for $\lambda\in \F$, $f\in \Fq\bbra{T_K}$, and $a\in \oK^\times$:}
\begin{equation}\label{phiok}
\left\{\begin{array}{l}
\varphi(\lambda\otimes f)=\lambda\otimes f^p,\\
a(\lambda\otimes f)=\lambda\otimes (f\circ a_{\rm LT}),
\end{array}\right.
\end{equation}
where we still denote by
$a_{\rm LT}(T_K)\in \Fq\bbra{T_K}$ the reduction mod $p$ of
$a_{\rm LT}(T_K)\in \oK\bbra{T_K}$. Lubin--Tate Theory implies that $\F_q\bbra{T_K}\subset\oCp^\flat$ is
stable under the action of $\gK$, and moreover that the action of $\gK$ on
$\F_q\bbra{T_K}$ factors through $\Gal(K_\infty/K)$ and coincides with
action of $\cO_K^\times$ in (\ref{phiok}) via the local reciprocity map.

Denote \ by \ $\Fq\ppar{T_{K}}^{\rm sep}$ \ the \ separable \ closure \ of
\ $\Fq\ppar{T_{K}}$ \ in \ $\Cp^\flat$. \ If \ $\overline{\rho}\in\Rep_{\F}\gK$, define
\[ D_K(\rhobar)\defeq\big(\Fq\ppar{T_K}^{\rm
    sep}\otimes_{\Fp}\rhobar\big)^{\Gal(\overline K/K_\infty)}. \]

Then $D_K(\rhobar)$ is an \'etale $(\varphi,\oK^\times)$-modules over
$\F\otimes_{\Fp}\Fq\ppar{T_K}$ and it follows from
Fontaine's theory of $(\varphi,\Gamma)$-modules (\cite{Fo}) that $D_K$
is a (covariant) rank-preserving {$\otimes$-}equivalence \ of \ categories \ between \ $\Rep_{\F}\gK$ \ and \ the \ category \ of
\ \'etale \ $(\varphi,\oK^\times)$-modules over
$\F\otimes_{\Fp}\Fq\ppar{T_K}$. Note that the injectivity of
$\id_{\F_{q}\ppar{T_K}}\otimes\varphi$
implies that the endomorphism
$\varphi$ of an \'etale $(\varphi,\oK^\times)$-module over
$\F\otimes_{\Fp}\Fq\ppar{T_K}$ is automatically injective.

The isomorphism
\begin{equation}\label{decomp}
\begin{array}{rcl}
\F\otimes_{\Fp}\Fq\ppar{T_K}\!\!&\!\!\buildrel\sim\over\longrightarrow \!\!&\!\!\F\ppar{T_{K,\sigma_0}}\times \F\ppar{T_{K,\sigma_1}}\times \cdots \times \F\ppar{T_{K,\sigma_{f-1}}}\\
\lambda\otimes (\sum_{n\gg -\infty}c_nT_K^n)\!\!&\!\!\longmapsto \!\!&\!\! \big(\sum_{n\gg -\infty}\lambda\sigma_0(c_n)T_{K,\sigma_0}^n,\dots,\sum_{n\gg -\infty}\lambda\sigma_{f-1}(c_n)T_{K,\sigma_{f-1}}^n\big)
\end{array}
\end{equation}
induces an analogous decomposition for any $\F\otimes_{\Fp}\Fq\ppar{T_K}$-module $D_K$:
\[D_K\buildrel\sim\over\longrightarrow D_{K,\sigma_0}\times \cdots\times D_{K,\sigma_{f-1}}.\]
If $D_K$ is an \'etale $\varphi$-module over
$\F\otimes_{\Fp}\Fq\ppar{T_K}$, then $\varphi$ induces a morphism
(still denoted by) $\varphi:D_{K,\sigma_i}\rightarrow D_{K,\sigma_{i-1}}$
such that
$\varphi(\sum_{n\gg -\infty}c_nT_{K,\sigma_{i}}^nv)=\sum_{n\gg -\infty}c_nT_{K,\sigma_{i-1}}^{pn}\varphi(v)$
for $c_n\in \F$ and $v\in D_{K,\sigma_i}$. By a standard argument, the
functor $D_K\longmapsto D_{K,\sigma_0}$ induces an equivalence of
categories (compatible with tensor products) between the category of
\'etale $(\varphi,\oK^\times)$-modules over
$\F\otimes_{\Fp}\Fq\ppar{T_K}$ and the category of {\'etale}
$(\varphi_q,\oK^\times)$-modules over $\F\ppar{T_{K,\sigma_0}}$, where
$\F\ppar{T_{K,\sigma_0}}$ is endowed with a continuous $\F$-linear endomorphism
$\varphi_q\ (=\varphi^f)$ and a continuous $\F$-linear action of $\oK^\times$ that commutes with $\varphi_q$ and satisfies the following conditions for $f\in \F\bbra{T_{K,\sigma_0}}$ and $a\in \oK^\times$:
\begin{equation*}\left\{\begin{array}{l}
\varphi_q(f(T_{K,\sigma_0}))=f(T_{K,\sigma_0}^q),\\
a(f(T_{K,\sigma_0}))=f(a_{\rm LT}(T_{K,\sigma_0})).
\end{array}\right.
\end{equation*}
{To be precise, in the last formula,}
\begin{equation*}
a_{\rm LT}(T_{K,\sigma_0})\defeq\sigma_0(a_{\rm LT}(T_K))\in \sigma_0(\overline a)T_{K,\sigma_0}+T_{K,\sigma_0}^2\F\bbra{T_{K,\sigma_0}},
\end{equation*}
where $\sigma_0(a_{\rm LT}(T_K))$ is the image of $a_{\rm LT}(T_K)\in \Fq\bbra{T_K}$ via
\[\Fq\bbra{T_K}\hookrightarrow \F\bbra{T_{K,\sigma_0}}, \sum_{n\gg -\infty}c_nT_K^{n}\mapsto \sum_{n\gg -\infty}\sigma_0(c_n)T_{K,\sigma_0}^{n}.\]
If one chooses the embedding $\sigma_i$ for some $i\in \{1,\dots,f-1\}$ instead of $\sigma_0$, one goes from $D_{K,\sigma_0}$ to $D_{K,\sigma_i}$ by the isomorphism
\[\Id \otimes \varphi^{f-i}:\F\bbra{T_{K,\sigma_i}}\otimes_{\varphi^{f-i},\F\bbra{T_{K,\sigma_0}}}D_{K,\sigma_0}\buildrel{\sim}\over\longrightarrow D_{K,\sigma_i}.\]

We can also work with the infinite Galois extension
$K(\!\sqrt[p^\infty]{1})$ instead of $K_\infty$ (see (\ref{LTp})). Let
$T$ be a coordinate of the formal group ${\mathbb G}_{\rm m}$. We
endow the topological ring $\F\otimes_{\Fp}\Fq\ppar{T}$ with a continuous $\F$-linear endomorphism $\varphi$ and a continuous
$\F\otimes_{\Fp}\Fq$-linear action of $\Zp^\times$ commuting with $\varphi$ and satisfying the following conditions for $\lambda\in \F$, $f\in \Fq\bbra{T}$, and $a\in \Zp^\times$:
\begin{equation}\label{phizp}
\left\{\begin{array}{l}
\varphi(\lambda\otimes f)=\lambda\otimes f^p,\\
a(\lambda\otimes f)=\lambda\otimes (f\circ a).
\end{array}\right.
\end{equation}
The choice of a generator of the Tate module of $\mathbb{G}_m$ and the
choice of $T$ induce an embedding $\F_q\bbra{T}\hookrightarrow\oCp^\flat$
whose image is stable under $\gK$ and on which the action of $\gK$
factors through {$\Gal(K(\!\sqrt[p^{\infty}]{1})/K)$} with action given by (\ref{phizp}) {(via local class field theory)}.

If $\overline{\rho}\in\Rep_{\F}\gK$, define
\[ D(\rhobar)\defeq\big(\Fq\ppar{T}^{\rm
    sep}\otimes_{\Fp}\rhobar\big)^{\Gal(\overline
    K/K(\!\sqrt[p^\infty]{1}))}. \] The functor $D$ is, as before, a
(covariant) rank-preserving {$\otimes$-}equivalence of categories between the category $\Rep_{\F}(\gK)$ and the category
of \'etale $(\varphi,\Z_p^\times)$-modules over
$\F\otimes_{\Fp}\Fq\ppar{T}$.

Here a standard choice is to take $T$ such that $a(T)\in \overline
aT+T^2\Fp\bbra{T}\subseteq \overline aT+T^2\Fq\bbra{T}$
is the reduction mod $p$ of $(1+T)^a-1\in \Zp\bbra{T}$. Using a
decomposition analogous to (\ref{decomp}) and choosing the embedding $\sigma_0$, we again have an
equivalence (compatible with tensor products) $D\longmapsto D_{\sigma_0}$ between the category of
\'etale $(\varphi,\Z_p^\times)$-modules over $\F\otimes_{\Fp}\Fq\ppar{T}$ and the category of \'etale
$(\varphi_q,\Z_p^\times)$-modules over $\F\ppar{T}$, where
$\F\ppar{T}$ is endowed with an $\F$-linear endomorphism $\varphi_q\
(=\varphi^f)$ and a continuous $\F$-linear action of $\Zp^\times$ that commutes with
$\varphi_q$ and satisfies the following conditions for $f\in\F\bbra{T}$ and $a\in \Zp^\times$:
\begin{equation}\label{phizp0}
\left\{\begin{array}{l}
\varphi_q(f(T))=f(T^q),\\
a(f(T))=f(a(T)).
\end{array}\right.
\end{equation}

We will mostly use $D_{K,\sigma_0}(\rhobar) \!{{}\defeq \!(D_{K}(\rhobar))_{\sigma_0}}${, an \'etale $(\varphi_q,\oK^\times)$-module over $\F\ppar{T_{K,\sigma_0}}$,} and $D_{\sigma_0}(\rhobar) {{}\defeq (D(\rhobar))_{\sigma_0}}${, an \'etale $(\varphi_q,\Zp^\times)$-module over $\F\ppar{T}$,} in the sequel.

We now relate $D_{K}(\rhobar)$ and $D(\rhobar)$, $D_{K,\sigma_0}(\rhobar)$ and $D_{\sigma_0}(\rhobar)$. In order to do so, we have to use the perfectoid versions of $\Fq\ppar{T_{K}}$, $\Fq\ppar{T_{K,\sigma_0}}$, etc.

We let $\Fq\bbra{T_{K}^{p^{-\infty}}}$ be the completion of the perfection $\bigcup_{n\geq 0}\Fq\bbra{T_{K}^{p^{-n}}}$ of $\Fq\bbra{T_{K}}$ with respect to the $T_{K}$-adic topology and $\Fq\ppar{T_{K}^{p^{-\infty}}}$ the fraction field of $\Fq\bbra{T_{K}^{p^{-\infty}}}$. Concretely:
\[\Fq\bbra{T_{K}^{p^{-\infty}}}\simeq \left\{\sum_{n\geq 0}c_nT_{K}^{\frac{d_n}{p^n}},\ c_n\in \Fq,\ d_n\in\Z_{\geq 0},\ \frac{d_n}{p^n}\rightarrow +\infty\textrm{ in $\Q$ when }n\rightarrow +\infty\right\}\]
and $\Fq\ppar{T_K^{p^{-\infty}}}=\Fq\bbra{T_{K}^{p^{-\infty}}}[\frac{1}{T_K}]$. We define in a similar way  $\Fq\bbra{T^{p^{-\infty}}}$ and $\Fq\ppar{T^{p^{-\infty}}}$. 

As $\oCp^\flat$ is perfect and complete for the $T_K$-adic
(resp.~$T$-adic) topology, we
have morphisms of $\F_q$-algebras
\begin{equation}\label{eq:embeddings_OCPflat}
\F_q\bbra{T_K^{p^{-\infty}}}\rightarrow\oCp^\flat, \qquad
\F_q\bbra{T^{p^{-\infty}}}\rightarrow\oCp^\flat.
\end{equation}
The following well-known theorem follows from the work of Wintenberger (\cite{wintenberger}) and the Ax--Sen--Tate Theorem, see for instance \cite[Cor.~3.4]{CE}:

\begin{thm1}\label{axsentate}
The morphisms (\ref{eq:embeddings_OCPflat}) induce isomorphisms of topological rings compatible with the action of $\oK^\times$ {\upshape(}via (\ref{LTp}){\upshape)}:
\[\Fq\bbra{T_{K}^{p^{-\infty}}}\simeq {\oCp^\flat}^{\!\!\!\Gal(\overline K/K_\infty)}\textrm{  and  \ }\Fq\ppar{T_{K}^{p^{-\infty}}}\simeq {\Cp^\flat}^{\Gal(\overline K/K_\infty)}\]
and isomorphisms of topological rings compatible with the action of $\Zp^\times$ {\upshape(}via (\ref{LTp}){\upshape)}:
\[\Fq\bbra{T^{p^{-\infty}}}\simeq {\oCp^\flat}^{\!\!\!\Gal(\overline K/K(\!\sqrt[p^\infty]{1}))}\textrm{  and  \ }\Fq\ppar{T^{p^{-\infty}}}\simeq {\Cp^\flat}^{\Gal(\overline K/K(\!\sqrt[p^\infty]{1}))}.\]
In particular, $\Fq\bbra{T^{p^{-\infty}}}\simeq \Fq\bbra{T_{K}^{p^{-\infty}}}^{\Gal(K_\infty/K(\!\sqrt[p^\infty]{1}))}\hookrightarrow \Fq\bbra{T_{K}^{p^{-\infty}}}$
and $\Fq\ppar{T^{p^{-\infty}}}\simeq \Fq\ppar{T_{K}^{p^{-\infty}}}^{\Gal(K_\infty/K(\!\sqrt[p^\infty]{1}))}\hookrightarrow \Fq\ppar{T_{K}^{p^{-\infty}}}$.
\end{thm1}

By Theorem \ref{axsentate} we have in particular embeddings
$\Fq\bbra{T}\hookrightarrow\Fq\bbra{T^{p^{-\infty}}}\hookrightarrow
\Fq\bbra{T_{K}^{p^{-\infty}}}$. Applying $\F\otimes_{\sigma_0,\Fq}(-)$ to Theorem \ref{axsentate}, we deduce embeddings $\F\ppar{T}\hookrightarrow\F\ppar{T^{p^{-\infty}}}\hookrightarrow \F\ppar{T_{K,\sigma_0}^{p^{-\infty}}}$ and $\F\bbra{T}\hookrightarrow\F\bbra{T^{p^{-\infty}}}\hookrightarrow \F\bbra{T_{K,\sigma_0}^{p^{-\infty}}}$.

\begin{prop1}\label{compare}
Let $\rhobar\in\Rep_{\F}\gK$. There is a canonical $\Fq\ppar{T_K^{p^{-\infty}}}$-linear isomorphism which commutes with the actions of $\oK^\times$ and $\varphi$:
\begin{equation}\label{eq:compare}\Fq\ppar{T_K^{p^{-\infty}}}\otimes_{\Fq\ppar{T}}D(\rhobar)\buildrel\sim\over\longrightarrow \Fq\ppar{T_K^{p^{-\infty}}}\otimes_{\Fq\ppar{T_K}}D_K(\rhobar)\end{equation}
where $\oK^\times$, $\varphi$ act
diagonally on each side,
$\oK^\times$ acting on $D(\rhobar)$ via the norm map $\oK^\times
\twoheadrightarrow \Zp^\times$. Moreover there is a canonical $\F\ppar{T_{K,\sigma_0}^{p^{-\infty}}}$-linear isomorphism which commutes with the actions of $\oK^\times$ and $\varphi_q$:
\[\F\ppar{T_{K,\sigma_0}^{p^{-\infty}}}\otimes_{\F\ppar{T}}D_{\sigma_0}(\rhobar)\buildrel\sim\over\longrightarrow \F\ppar{T_{K,\sigma_0}^{p^{-\infty}}}\otimes_{\F\ppar{T_{K,\sigma_0}}}D_{K,\sigma_0}(\rhobar)\]
where $\oK^\times$, $\varphi_q$ act diagonally on each side, $\oK^\times$ acting on $D_{\sigma_0}(\rhobar)$ via the norm map $\oK^\times \twoheadrightarrow \Zp^\times$.

\end{prop1}
\begin{proof}
  From Wintenberger's theory of the field of norms
  (\cite{wintenberger}), recall that we have topological isomorphisms
\[\Gal(\Fq\ppar{T_K}^{\rm sep}/\Fq\ppar{T_K}\!)\simeq \Gal(\overline K/K_\infty),\ \Gal(\Fq\ppar{T}^{\rm sep}/\Fq\ppar{T}\!)\simeq \Gal(\overline K/K(\!\sqrt[p^\infty]{1})).\]
Since we have for any integer $n\geq 1$:
\[H^1\big(\Gal(\Fq\ppar{T_K}^{\rm sep}/\Fq\ppar{T_K}\!),\GL_n(\Fq\ppar{T_{K}}^{\rm sep})\big)=1\]
as follows by taking inductive limit from \cite[Prop.~X.1.3]{Serre68}, we have a canonical isomorphism
\begin{equation}\label{isosep}
\Fq\ppar{T_K}^{\rm sep}\otimes_{\Fq\ppar{T_K}}D_K(\rhobar)\buildrel\sim\over\longrightarrow \Fq\ppar{T_K}^{\rm sep}\otimes_{\Fp}\rhobar
\end{equation}
that is compatible with the actions of $\vp$ and  {$\Gal(\overline K/K)$ ($\Gal(\overline K/K)$ acting on $D_K(\rhobar)$ via $\Gal(K_\infty/K)$)}, and likewise with $\Fq\ppar{T}^{\rm sep}$, $\Fq\ppar{T}$ and $D(\rhobar)$. Tensoring (\ref{isosep}) by $\Cp^\flat$ over $\Fq\ppar{T_K}^{\rm sep}$, resp.~its analogue over $\Fq\ppar{T}^{\rm sep}$, we obtain a canonical isomorphism
\begin{equation*}
\Cp^\flat\otimes_{\Fq\ppar{T_K}}D_K(\rhobar)\congto \Cp^\flat\otimes_{\Fp}\rhobar \congfrom \Cp^\flat\otimes_{\Fq\ppar{T}}D(\rhobar)
\end{equation*}
compatible \ with \ the \ actions \ of \ $\vp$ \ and \  {$\Gal(\overline K/K)$}. \ Taking \ invariants \ under $\Gal(\overline K /K_\infty)$, which acts trivially on $D_K(\rhobar)$,
$D(\rhobar)$, and remembering ${\Cp^\flat}^{\Gal(\overline K/K_\infty)} = \Fq\ppar{T_{K}^{p^{-\infty}}}$ from Theorem~\ref{axsentate} we obtain the desired isomorphism \eqref{eq:compare}.
The last assertion follows from an analogous discussion~(the details of which are left to the reader).
\end{proof}

\begin{rem1}
Arguing as in the proofs of Theorem \ref{descent} and Corollary \ref{descentoK} below, the functor $D_{K,\sigma_0}\mapsto \F\ppar{T_{K,\sigma_0}^{p^{-\infty}}}\otimes_{\F\ppar{T_{K,\sigma_0}}}D_{K,\sigma_0}$ in fact still induces an equivalence of categories from the category of \'etale $(\varphi_q,\oK^\times)$-modules over $\F\ppar{T_{K,\sigma_0}}$ to the category of \'etale $(\varphi_q,\oK^\times)$-modules over $\F\ppar{T_{K,\sigma_0}^{p^{-\infty}}}$. Likewise with the functor $D_{\sigma_0}\mapsto \F\ppar{T^{p^{-\infty}}}\otimes_{\F\ppar{T}}D_{\sigma_0}$ and \'etale $(\varphi_q,\Zp^\times)$-modules.
\end{rem1}

We finally recall a convenient explicit presentation of $D_{K,\sigma_0}(\rhobar)$ for $\rhobar$ absolutely irreducible.

For simplicity, we now choose the formal variable $T_K$ such that $a_{\rm LT}(T_K)=\overline aT_K$ when $a\in [\Fq]$ (so $a(T_{K,\sigma_0})=\sigma_0(\overline a)T_{K,\sigma_0}$ for $a\in [\Fq]$); for instance, this holds if $T_K$ is such that the logarithm of the Lubin--Tate group $G_{\LT}$ (\cite[ch.8~\S 6]{Lang}) is the series $\sum_{n\geq0}p^{-n}T_K^{q^n}$. Note that in that case $\F\ppar{T_{K,\sigma_0}}^{[\Fq^\times]}=\F\ppar{T_{K,\sigma_0}^{q-1}}$ and that the commutativity of the action of $a\in \oK$ with $[\Fq]$ implies:
\begin{equation}\label{inq-1}
a(T_{K,\sigma_0})\in \sigma_0(\overline a)T_{K,\sigma_0}+T_{K,\sigma_0}^q\oK\bbra{T_{K,\sigma_0}^{q-1}}.
\end{equation}
We recall the following straightforward lemma.

\begin{lem1}\label{triv}
Let $\rhobar$ be a finite-dimensional continuous representation of $\Gal(\overline K/K)$ over $\F$. Denote by $D_{K,\sigma_0}(\rhobar)^{[\Fq^\times]}$ the $\F\ppar{T_{K,\sigma_0}^{q-1}}$-vector subspace of $D_{K,\sigma_0}(\rhobar)$ fixed by $[\Fq^\times]\subseteq \oK^\times$. Then $D_{K,\sigma_0}(\rhobar)^{[\Fq^\times]}$ is preserved by $\varphi_q$ and the action of $\oK^\times$, and we have an $\F\ppar{T_{K,\sigma_0}}$-linear isomorphism compatible with $\varphi_q$ and $\oK^\times$:
\[\F\ppar{T_{K,\sigma_0}}\otimes_{\F\ppar{T_{K,\sigma_0}^{q-1}}}D_{K,\sigma_0}(\rhobar)^{[\Fq^\times]}\buildrel\sim\over\longrightarrow D_{K,\sigma_0}(\rhobar)\]
where the actions of $\varphi_q$ and $\oK^\times$ on the left-hand side are the diagonal ones.
\end{lem1}
\begin{proof}
It is enough to prove that the morphism in the statement is an isomorphism, everything else being trivial. It is enough to prove
\[H^1\big([\Fq^\times],\GL_n(\F\ppar{T_{K,\sigma_0}})\big)=1\]
for any integer $n\geq 1$. But this is again the generalization of Hilbert $90$ applied to the Galois extension $\F\ppar{T_{K,\sigma_0}}/\F\ppar{T_{K,\sigma_0}^{q-1}}$ (which has Galois group $[\Fq^\times]$), see for instance \cite[Prop.~X.1.3]{Serre68}.
\end{proof}

We now give explicitly $D_{K,\sigma_0}(\rhobar)^{[\Fq^\times]}$ for an absolutely irreducible $\rhobar$.

For $\lambda\in \F^\times$ denote by ${\rm unr}(\lambda)$ the unramified character of $\Gal(\overline K/K)$ sending the Frobenius $x\mapsto x^q$ to $\lambda^{-1}$. For $f'\geq 1$ denote by $\omega_{f'}:I_K\rightarrow \F_{p^{f'}}^\times$ Serre's fundamental character of level $f'$, where $I_K\subseteq \Gal(\overline K/K)$ is the inertia subgroup. We also denote by $\omega_f$ (instead of $\sigma_0\circ \omega_f$) the composition
\begin{equation}\label{omegaf}
I_K\buildrel\omega_f\over \longrightarrow \Fq^\times\buildrel\sigma_0\over\hookrightarrow \F^\times
\end{equation}
and again $\omega_f$ its unique extension to $\Gal(\overline K/K)$ such that $\omega_f(p)=1$ (via local class field theory). Recall that $\omega_f:\Gal(\overline K/K)\rightarrow \F^\times$ is the composition by $\sigma_0:\Fq\hookrightarrow \F$ of the mod $p$ Lubin--Tate character of $\Gal(\overline K/K)$. For $d\in \Z_{\geq 1}$, it goes back to Serre that any absolutely irreducible $d$-dimensional representation of $\Gal(\overline K/K)$ over $\F$ is isomorphic to $(\ind\omega_{df}^h)\otimes {\rm unr}(\lambda)$ for some $\lambda\in \F^\times$ and some positive integer $h$ which is not of the form $m\frac{q^d-1}{q^{d'}-1}$ for some $m\in \Z_{\geq 1}$ and some $d'\in \{1,\dots,d-1\}$, where $\ind\omega_{df}^h$ is the induction from $\Gal(\overline K/K_d)$ to $\Gal(\overline K/K)$ of the mod $p$ Lubin--Tate character of $\Gal(\overline K/K_d)$ (seen with values in $\F$ via any embedding $\F_{q^d}\hookrightarrow \F$ lifting $\sigma_0$), where $K_d$ is the unramified extension of $K$ of degree $d$. Equivalently $\ind\omega_{df}^h$ is the unique representation of $\Gal(\overline K/K)$ over $\F$ with determinant $\omega_{f}^h\cdot {\rm unr}(-1)^{d-1}$ such that $(\ind\omega_{df}^h)\vert_{I_K}\cong \omega_{df}^h\oplus \omega_{df}^{qh}\oplus \cdots \oplus \omega_{df}^{q^{d-1}h}$ (for any choice of embedding $\F_{q^{d}}\hookrightarrow \F$). Note that
\[(\ind\omega_{df}^h)\otimes {\rm unr}(\lambda)\cong (\ind\omega_{df}^{h'})\otimes {\rm unr}(\lambda')\]
if and only if $h'\equiv q^ih\!\!\mod q^d-1$ for some $i\in \{0,\dots,d-1\}$ and $\lambda^d=\lambda'^d$.

For $a\in \oK^\times$, we set:
\[f_a^{\rm LT}\defeq f_a^{\rm LT}(T_{K,\sigma_0})\defeq\frac{\sigma_0(\overline a)T_{K,\sigma_0}}{a(T_{K,\sigma_0})}\in 1+T_{K,\sigma_0}\F\bbra{T_{K,\sigma_0}}.\]
Note that $f_a^{\rm LT}=1$ if $a\in [\Fq^\times]$ and that (\ref{inq-1}) implies
\[f_a^{\rm LT}\in 1+T_{K,\sigma_0}^{q-1}\F\bbra{T_{K,\sigma_0}^{q-1}}.\]

\begin{lem1}\label{PS}
Let $\rhobar\in\Rep_{\F}\gK$ and write  $\rhobar=(\ind\omega_{df}^h)\otimes {\rm unr}(\lambda)$ for some $d,h,\lambda$ as above. Then $D_{K,\sigma_0}(\rhobar)\simeq \F\ppar{T_{K,\sigma_0}}\otimes_{\F\ppar{T_{K,\sigma_0}^{q-1}}}D_{K,\sigma_0}(\rhobar)^{[\Fq^\times]}$ {\upshape(}Lemma \ref{triv}{\upshape)}, where $D_{K,\sigma_0}(\rhobar)^{[\Fq^\times]}$ is explicitly described as follows:
\begin{equation}\label{D0rhobar}
\left\{\begin{array}{cll}
D_{K,\sigma_0}(\rhobar)^{[\Fq^\times]}&=&\bigoplus_{i=0}^{d-1}\F\ppar{T_{K,\sigma_0}^{q-1}}e_i\\
\varphi_q(e_i)&=& e_{i+1},\ i<d-1\\
\varphi_q(e_{d-1})&=&\displaystyle{\frac{\lambda^d}{T_{K,\sigma_0}^{h(q-1)}} e_0}\\
a(e_i)&=&{{\left(f_a^{\rm LT}\right)}^{\frac{hq^i(q-1)}{q^d-1}}e_i,\ a\in \oK^\times.}\\
\end{array}\right.
\end{equation}
Moreover a basis $(e_0,\dots,e_{d-1})=(e_0,\varphi_q(e_0),\dots,\varphi_q^{d-1}(e_0))$ as in (\ref{D0rhobar}) is uniquely determined up to a scalar in $\F^\times$. Finally, if $h'=q^jh+m(q^d-1)$ for some $j\in \{0,\dots,d-1\}$ and some $m\in \Z$, then the unique basis $(e'_i)_i=(e'_0,\varphi_q(e'_0),\dots,\varphi_q^{d-1}(e'_0))$ in (\ref{D0rhobar}) corresponding to $h'$ is given by $e'_0= \frac{1}{T_{K,\sigma_0}^{m(q-1)}}e_j$ {\upshape(}up to a scalar in $\F^\times${\upshape)}.
\end{lem1}
\begin{proof}
The first statement is \cite[Cor.~10.10]{PS2}. We prove the uniqueness of the basis $(e_i)$ in (\ref{D0rhobar}) (up to scalar). Let $(f_0,\dots,f_{d-1})$ be another basis of $D_{K,\sigma_0}(\rhobar)^{[\Fq^\times]}$ satisfying (\ref{D0rhobar}), it is enough to prove that $f_0\in \F e_0$. Write $f_0=\sum_{i=0}^{d-1}x_ie_i$ for some $x_i\in \F\ppar{T_{K,\sigma_0}^{q-1}}$. Since $\varphi_q^d(f_0)=\frac{\lambda^d}{T_{K,\sigma_0}^{h(q-1)}} f_0$ and $\varphi_q^d(e_i)=\frac{\lambda^d}{T_{K,\sigma_0}^{q^ih(q-1)}}e_i$, we deduce that $\frac{1}{T_{K,\sigma_0}^{h(q-1)}}x_i=\frac{1}{T_{K,\sigma_0}^{q^ih(q-1)}}\varphi_q^d(x_i)$ for $i\in \{0,\dots,d-1\}$, i.e.\ $\varphi_q^d(x_i)=T_{K,\sigma_0}^{h(q-1)(q^i-1)}x_i$. This easily implies $x_i\in \F T_{K,\sigma_0}^{m_i}$, where $m_i\defeq\frac{h(q-1)(q^i-1)}{q^d-1}\in \Z_{\geq 0}$. If $x_i\ne 0$, since $(q-1)\vert m_i$ in $\Z$, we obtain $h=\frac{m_i}{q-1}\frac{q^d-1}{q^i-1}$ for some $i\in \{1,\dots,d-1\}$ which contradicts the assumption on $h$. Hence $x_i=0$ for all $i\ne 0$ and thus $f_0\in \F e_0$. The last statement is an easy check that is left to the reader.
\end{proof}

\begin{rem1}
One can prove that the action of $a\in \oK^\times$ in (\ref{D0rhobar}) is the unique semi-linear action on $D_{K,\sigma_0}(\rhobar)^{[\Fq^\times]}$ which commutes with $\varphi_q$ and is such that $a(e_i)\in e_i+T_{K,\sigma_0}^{q-1}\sum_{j=0}^{d-1}\F\bbra{T_{K,\sigma_0}^{q-1}}e_j$ for all $i$. The argument is the same as in the proof of Lemma \ref{PSnatural} below.
\end{rem1}

As a special case of Lemma \ref{PS} we have:

\begin{lem1}\label{twist0}
Let $\chi:\Gal(\overline K/K)\rightarrow \F^\times$ be a continuous character and write $\chi=\omega_f^{h_{\chi}}{\rm unr}(\lambda_{\chi})$ for $h_\chi\in \Z_{\geq 0}$ and $\lambda_\chi\in \F^\times$, then {\upshape(}for $a\in \oK^\times${\upshape)}:
\begin{equation*}
\left\{\begin{array}{cll}
D_{K,\sigma_0}(\chi)^{[\Fq^\times]}&=&\F\ppar{T_{K,\sigma_0}^{q-1}}e_\chi\\
\varphi_q(e_{\chi})&=&\displaystyle{\frac{\lambda_{\chi}}{T_{K,\sigma_0}^{h_\chi(q-1)}}e_\chi}\\
a(e_\chi)&=&({f_a^{\rm LT}})^{h_\chi}e_\chi
\end{array}\right.
\end{equation*}
\end{lem1}

If $\rhobar$ is any finite-dimensional continuous representation of $\Gal(\overline K/K)$ over $\F$, write $D_{K,\sigma_0}(\rhobar)(\chi)$ for $D_{K,\sigma_0}(\rhobar)\otimes_{\F\ppar{T_{K,\sigma_0}}}D_{K,\sigma_0}(\chi)$ with tensor product structures. Then we have $D_{K,\sigma_0}(\rhobar\otimes \chi)\simeq D_{K,\sigma_0}(\rhobar)(\chi)$ as follows from the compatibility of $D_{K,\sigma_0}(-)$ with tensor products.

\subsection{The \texorpdfstring{$(\varphi_q,\cO_K^\times)$}{(phi\_q,O\_K\^{}x)}-module over \texorpdfstring{$A$}{A} of a semi-simple Galois representation}\label{oK}

To an arbitrary semi-simple $\rhobar$ we associate by an elementary recipe an \'etale $(\varphi_q,\oK^\times)$-module $D_{A,\sigma_0}(\rhobar)$ over $A$ (depending on the fixed choice of the embedding $\sigma_0$).

Let $N_0\defeq\smat{1&\oK\\0&1}\subseteq \GL_2(\oK)$ and $\m_{N_0}$ the maximal ideal of $\Fq\bbra{N_0}$. Recall that $\Fq\bbra{N_0}=\Fq\bbra{Y_{0},\dots,Y_{{f-1}}}$ and $\m_{N_0}=(Y_{0},\dots,Y_{{f-1}})$, where
\[Y_{i}\defeq\sum_{\lambda\in
    \Fq^\times}\lambda^{-p^i}\begin{pmatrix}1&[\lambda]\\0&1\end{pmatrix}\in
  \Fq\bbra{N_0}.\]
(Namely, as $N_0$ is a uniform pro-$p$-group isomorphic to $\Zp^f$, it
follows from \cite[Thm.~7.23(i)]{DDMS} that $\Fq\bbra{N_0}$ is
isomorphic to a power series ring in $f$ variables over $\F_q$. This is a local ring. We easily check that the images of $Y_0,\dots,Y_{f-1}$
in $\m_{N_0}/\m_{N_0}^2$ form a basis of this $\F_q$-vector
space, so that $\F_q\bbra{N_0}=\F_q\bbra{Y_0,\dots,Y_{f-1}}$.)
  
As in \cite[\S 3.1.1]{BHHMS2} consider the multiplicative system
\[S\defeq\{(Y_{0}\cdots Y_{f-1})^k,\ k\geq 0\}\subseteq \Fq\bbra{N_0}\]
and $A_q\defeq\widehat{\Fq\bbra{N_0}_S}$ the completion of the localization $\Fq\bbra{N_0}_S$ with respect to the ascending filtration ($n\in \Z$):
\begin{equation}
\label{eq:def:fil:Aq}
F_n(\Fq\bbra{N_0}_S)\defeq\sum_{k\geq0}\frac{1}{(Y_{0}\cdots Y_{f-1})^{k}}\m_{N_0}^{kf-n}=\bigcup_{k\geq0}\frac{1}{(Y_{0}\cdots Y_{f-1})^{k}}\m_{N_0}^{kf-n}
\end{equation}
where $\m_{N_0}^m\defeq\Fq\bbra{N_0}$ if $m\leq 0$ (see \cite[\S 3.1.1]{BHHMS2}). We denote by $F_nA_q$ ($n\in \Z$) the induced ascending filtration on $A_q$ and endow $A_q$ with the associated topology (\cite[\S I.3]{LiOy}). The ring $A_q$ contains $\Fq\bbra{N_0}$ and the $\Fq$-linear action of $\oK^\times$ on $\Fq\bbra{N_0}$ (induced by the multiplication on $\oK\simeq N_0$) canonically extends by continuity to $A_q$ (but not to $\Fq\bbra{N_0}_S$ as it does not preserve $S$). We will write this action of $\oK^\times$ on $\Fq\bbra{N_0}$ and $A_q$ as $a(x)$ for $(a,x)\in \oK^\times\times A_q$. In fact using $a-[\overline a]\in p\oK$ one has for $a\in \oK^\times$ and $i\in \Z$:
\[a(Y_{i})\in \overline a^{p^i}Y_{i}+\m_{N_0}^p\subseteq \m_{N_0}\]
which implies
\begin{equation}\label{natural}
a(Y_{i})\in \overline a^{p^i}Y_{i}\big(1+\frac{1}{Y_{i}}\m_{N_0}^p\big)\subseteq \overline a^{p^i}Y_{i}(1+F_{1-p}A_q)\subseteq A_q^\times.
\end{equation}
We define $\varphi$ as the Frobenius endomorphism of $\Fq\bbra{N_0}$,
i.e.\ by $\varphi(f)=f^p$ for $f\in\F_q\bbra{N_0}$. It canonically extends by continuity to $A_q$ and obviously commutes with the action of $\oK^\times$ on $\Fq\bbra{N_0}$, hence on $A_q$. 

Let $A$ be the complete filtered ring in \cite[\S 3.1.1]{BHHMS2}. Recall that $A$ is defined similarly to $A_q$ replacing $\Fq\bbra{N_0}$ by $\F\bbra{N_0}$ {\it except} that the Frobenius $\varphi$ on $\F\bbra{N_0}$ is now $\F$-linear. As in (\ref{decomp}), we have an isomorphism $\F\otimes_{\Fp}A_q\buildrel\sim\over\longrightarrow \underbrace{A\times A \times \cdots \times A}_{f\textrm{ times}}$ which sends $\lambda\otimes \sum_{\underline n}c_{\underline n}Y_0^{n_0}\cdots Y_{f-1}^{n_{f-1}}\in \F\otimes_{\Fp}A_q$ to:
\[\big(\lambda\sum_{\underline n}\sigma_0(c_{\underline n})Y_{\sigma_0}^{n_0}\cdots Y_{\sigma_{f-1}}^{n_{f-1}},\lambda\sum_{\underline n}\sigma_1(c_{\underline n})Y_{\sigma_1}^{n_0}\cdots Y_{\sigma_{0}}^{n_{f-1}},\dots,\lambda\sum_{\underline n}\sigma_{f-1}(c_{\underline n})Y_{\sigma_{f-1}}^{n_0}\cdots Y_{\sigma_{f-2}}^{n_{f-1}}\big)\]
where we set for $\sigma:\Fq\hookrightarrow \F$:
\begin{equation}\label{variableysigma}
Y_\sigma\defeq\sum_{\lambda\in \Fq^\times}\sigma(\lambda)^{-1}\begin{pmatrix}1&[\lambda]\\0&1\end{pmatrix}\in \F\bbra{N_0}\subseteq A.
\end{equation}
It induces an analogous decomposition for any $\F\otimes_{\Fp}A_q$-module $D_{A_q}$:
\begin{equation*}D_{A_q}\buildrel\sim\over\longrightarrow D_{A,\sigma_0}\times \cdots\times D_{A,\sigma_{f-1}}.
\end{equation*}
We extend $\F$-linearly the Frobenius $\varphi$ and the action of
$\oK^\times$ from $A_q$ to $\F\otimes_{\Fp}A_q$. Note that we have
$\varphi(Y_{\sigma_i})=Y_{\sigma_{i-1}}^p$ for $i\in \Z$ (see \cite[\S
3.1.1]{BHHMS2}, where $\varphi$ on $A$ is denoted by $\phi$). We let
$\varphi_q\defeq\varphi^f$ on $A$. As in \S\ref{LT}, the functor
$D_{A_q}\mapsto D_{A,\sigma_0}$ induces an equivalence of categories
between the category of \'etale $(\varphi,\oK^\times)$-modules over
$\F\otimes_{\Fp}A_q$ and the category of \'etale
$(\varphi_q,\oK^\times)$-modules over $A$.

The embedding of $\F$-algebras
\begin{equation}\label{recette}
\F\ppar{T_{K,\sigma_0}^{q-1}}\hookrightarrow A,\ \ \ \sum_{n\gg -\infty} c_nT_{K,\sigma_0}^{n(q-1)}\longmapsto \sum_{n\gg -\infty} c_n\Big(\frac{\varphi(Y_{\sigma_0})}{Y_{\sigma_0}}\Big)^n
\end{equation}
trivially commutes with $\varphi_q$ and $[\Fq^\times]$ (the latter acting trivially on both sides). When $\rhobar$ is a direct sum of absolutely irreducible finite-dimensional continuous representations of $\Gal(\overline K/K)$ over $\F$, we define:
\begin{equation}\label{DAq0}
D_{A,\sigma_0}(\rhobar) \defeq A\otimes_{\F\ppar{T_{K,\sigma_0}^{q-1}}}D_{K,\sigma_0}(\rhobar)^{[\Fq^\times]}
\end{equation}
where $D_{K,\sigma_0}(\rhobar)^{[\Fq^\times]}$ is as in Lemma
\ref{triv}. It follows from its definition that
$D_{A,\sigma_0}(\rhobar)$ is an \'etale $\varphi_q$-module over
$A$ if it is endowed with the endomorphism
$\varphi_q\defeq\varphi_q\otimes \varphi_q$.

\begin{rem1}
  Definition (\ref{DAq0}) does not need the semi-simplicity of
  $\rhobar$, but we will only use it in that case, see also Remark
  \ref{nonsemisimplebad} below.
\end{rem1}

For $a\in \oK^\times$, we set (see (\ref{natural})):
\begin{equation}\label{1-p}
f_{a,\sigma_0}\defeq f_{a,\sigma_0}(Y_{\sigma_0},\dots,Y_{\sigma_{f-1}})\defeq\frac{\sigma_0(\overline a)Y_{\sigma_0}}{a(Y_{\sigma_0})}\in 1+F_{1-p}A.
\end{equation}

\begin{lem1}\label{PSnatural}
Let \ $\rhobar$ \ be \ an \ absolutely \ irreducible \ continuous \ representation \ of $\Gal(\overline K/K)$ over $\F$ and $(e_0,\dots,e_{d-1})$ a basis of $D_{K,\sigma_0}(\rhobar)^{[\Fq^\times]}$ as in Lemma \ref{PS}. Then we have:
\begin{equation*}
\left\{\begin{array}{cll}
D_{A,\sigma_0}(\rhobar)&=&\bigoplus_{i=0}^{d-1}A(1\otimes e_i)\\
\varphi_q(1\otimes e_i)&=& 1\otimes e_{i+1},\ i<d-1\\
\varphi_q(1\otimes e_{d-1})&=&\displaystyle{\lambda^d\Big(\frac{Y_{\sigma_0}}{\varphi(Y_{\sigma_0})}\Big)^h(1\otimes e_0).}
\end{array}\right.
\end{equation*}
Moreover there is a unique structure of
$(\varphi_q,\oK^\times)$-module over $A$ on $D_{A,\sigma_0}(\rhobar)$ such that
\[a(1\otimes e_i)\in 1\otimes e_i+\sum_{j=0}^{d-1}(F_{1-p}A)(1\otimes e_j)\textrm{ for all $i$ and $a\in \oK^\times$.}\]
This action of $\oK^\times$ is explicitly given by {\upshape(}$i\in \{0,\dots,d-1\}$, $a\in \oK^\times${\upshape)}:
\begin{equation}\label{naturalaction}
a(1\otimes e_i)=\displaystyle{\Big(\frac{f_{a,\sigma_0}}{\varphi(f_{a,\sigma_0})}\Big)^{\frac{hq^i}{1-q^d}}(1\otimes e_i)\in (1+F_{q^i(1-p)}A)(1\otimes e_i)}
\end{equation}
and \ does \ not \ depend \ {\upshape(}up \ to \ isomorphism{\upshape)} \ on \ the \ choice \ of \ the \ basis $(e_i)_i$ \ of \ $D_{K,\sigma_0}(\rhobar)^{[\Fq^\times]}$.
\end{lem1}
\begin{proof}
The first part of the statement follows from the definition of $D_{A,\sigma_0}(\rhobar)$ in (\ref{DAq0}). Fix $a\in \oK^\times$ and write $a(1\otimes e_0)=\sum_{i=0}^{d-1}C_i(1\otimes e_i)$ for some $C_0\in 1+F_{1-p}A$ and $C_i\in F_{1-p}A$ if $i\ne 0$. Assume $C_i\ne 0$ for some $i\ne 0$ and let $m_i\geq p-1$ be the maximal integer such that $C_i\in F_{-m_i}A\setminus F_{-(m_i+1)}A$. Since $a(1\otimes e_0)$ and the $1\otimes e_j$ are fixed by $[\Fq^\times]$, the constants $C_j$ are also fixed by $[\Fq^\times]$ in $A$ for all $j$, and thus in particular by $[\Fp^\times]$, from which it is an exercise to deduce that we must have $(p-1)\vert m_i$. Since $\varphi_q^d(1\otimes e_j)=\lambda^d\Big(\frac{Y_{\sigma_0}}{\varphi(Y_{\sigma_0})}\Big)^{q^jh}\!(1\otimes e_j)$ for all $j$, the equality $a(\varphi_q^d(e_0))=\varphi_q^d(a(e_0))$ yields for $j\in \{0,\dots,d-1\}$ (using $\sigma_0(\overline a)^{q-1}=1$):
\begin{equation}\label{cj}
C_j=\Big(\frac{f_{a,\sigma_0}}{\varphi(f_{a,\sigma_0})}\Big)^{h}\Big(\frac{Y_{\sigma_0}}{\varphi(Y_{\sigma_0})}\Big)^{(q^j-1)h}\varphi_q^d(C_j)
\end{equation}
which implies in particular $-m_i=(q^i-1)h(p-1)-q^dm_i$, i.e.\
$(q^i-1)h(p-1)=(q^d-1)m_i$, i.e.\
$h=\frac{q^d-1}{q^i-1}\frac{m_i}{p-1}$, which contradicts the
assumption on $h$ since $\frac{m_i}{p-1}\in \Z$. Hence we must have
$C_i=0$ if $i\ne 0$. When $i=0$, (\ref{cj}) is just
$C_0=\big(\frac{f_{a,\sigma_0}}{\varphi(f_{a,\sigma_0})}\big)^{h}\varphi_q^d(C_0)$. The
equation $C_0=\big(\frac{f_{a,\sigma_0}}{\varphi(f_{a,\sigma_0})}\big)^{h}\varphi_q^d(C_0)$ has a solution in $1+F_{1-p}A$ given by
\begin{multline*}
C_0=\prod_{n=0}^{+\infty}\varphi_q^{nd}\Big(\Big(\frac{f_{a,\sigma_0}}{\varphi(f_{a,\sigma_0})}\Big)^{h}\Big)=\prod_{n=0}^{+\infty}\Big(\frac{f_{a,\sigma_0}}{\varphi(f_{a,\sigma_0})}\Big)^{q^{nd}h}=\Big(\frac{f_{a,\sigma_0}}{\varphi(f_{a,\sigma_0})}\Big)^{h(1+q^d+q^{2d}+\cdots)}\\
=\Big(\frac{f_{a,\sigma_0}}{\varphi(f_{a,\sigma_0})}\Big)^{\frac{h}{1-q^d}},
\end{multline*}
where the second equality uses $x^{q^{nd}}=x$ if $x\in \Fq$. This
solution is unique in $1+F_{1-p}A (\subset A^\times)$: the quotient of two
solutions is an element of $1+F_{1-p}A$ fixed by $\varphi_q^d$ and the
map $1-\varphi_q^d$ induces an automorphism of $F_{1-p}A$ (with inverse
$\sum_{n\geq0}\varphi_q^{dn}$) so that only $1$ is fixed by
$\varphi_q^d$ in $1+F_{1-p}A$. Then (\ref{naturalaction}) immediately
follows, from which the continuity of the action of $\oK^\times$ is
clear (as it is continuous on $A$). If one changes the basis
$(e_i)_i$, or equivalently by (the last statement in) Lemma \ref{PS}
changes the integer $h$, the last statement easily follows from the
last statement of Lemma \ref{PS}.
\end{proof}

\begin{rem1} The uniqueness of the action of $\oK^\times$ in the proof of Lemma \ref{PSnatural} works just assuming $a(1\otimes e_i)\in 1\otimes e_i+\sum_{j=0}^{d-1}(F_{-1}A)(1\otimes e_j)$ (and lands automatically in $1\otimes e_i+\sum_{j=0}^{d-1}(F_{1-p}A)(1\otimes e_j)$).
\end{rem1}

Let us finally make twists explicit. Let $\chi:\Gal(\overline K/K)\rightarrow \F^\times$ be a continuous character and write $\chi=\omega_f^{h_{\chi}}{\rm unr}(\lambda_{\chi})$ for $h_\chi\in \Z_{\geq 0}$ and $\lambda_\chi\in \F^\times$, then (using Lemma \ref{twist0}) the \'etale $(\varphi_q,\cO_K^\times)$-module $D_{A,\sigma_0}(\chi)$ is explicitly given by ($a\in \oK^\times$):
\begin{equation}\label{dachi}
\left\{\begin{array}{cll}
D_{A,\sigma_0}(\chi)&=&A(1\otimes e_\chi)\\
\varphi_q(1\otimes e_\chi)&=& \displaystyle{\lambda_\chi\Big(\frac{Y_{\sigma_0}}{\varphi(Y_{\sigma_0})}\Big)^{h_\chi}(1\otimes e_\chi)}\\
a(1\otimes e_\chi)&=&\displaystyle{\Big(\frac{f_{a,\sigma_0}}{\varphi(f_{a,\sigma_0})}\Big)^{\frac{h_\chi}{1-q}}(1\otimes e_\chi).}
\end{array}\right.
\end{equation}
One has an action of $\oK^\times$ on $D_{A,\sigma_0}(\rhobar\otimes \chi)\defeq D_{A,\sigma_0}(\rhobar)\otimes_A D_{A,\sigma_0}(\chi)$ by taking the tensor product action. We leave to the reader the exercise to check that, when $\rhobar\otimes \chi\simeq \rhobar'\otimes \chi'$, then $D_{A,\sigma_0}(\rhobar\otimes \chi)\simeq D_{A,\sigma_0}(\rhobar'\otimes \chi')$ as $(\varphi_q,\oK^\times)$-modules over $A$.

\subsection{A reminder on \texorpdfstring{$p$}{p}-divisible groups and \texorpdfstring{$K$}{K}-vector spa\-ces}\label{reminderp}

We review some results on constructions of Fargues and Fontaine (\cite{FF}) related to $p$-divisible groups (in a relative context, see for example \cite[\S 5.1]{LeBras}) and we define the important perfectoid spaces $Z_{\LT}$ and $Z_{\cO_K}$ over $\F$.

Let $R$ be a perfectoid
$\F$-algebra and $\varpi$ a
pseudo-uniformizer of $R$. As usual we denote by $R^\circ$ the subring of power-bounded elements in $R$ and by $R^{\ronron}\subset R^\circ$ the subset of topologically nilpotent elements (i.e.\ those $a\in R$ such that ${a}^n$ converges to $0$ in $R$). We fix a power-multiplicative norm
$\abs$ on $R$ defining the topology of $R$. Such a norm exists and can
be explicitly given by
\begin{equation}\label{normul} 
\vabs{a}=\inf\set{2^{\frac{m}{n}}, (m,n)\in\Z\times\Z_{>0}, \,
    \varpi^m a^n\in R^\circ}\in \R_{\geq 0}
\end{equation}
(so $\vabs{a}\leq 1\Leftrightarrow a\in R^\circ$). We endow the Witt vectors $W(R^\circ)$ with the
$(p,[\varpi])$-adic topology (where $[\cdot]$ is the multiplicative representative). Let $\B^+(R)$ be the Fr\'echet $K$-algebra defined as the completion
of $W(R^\circ)[1/p]$ for the family of norms $\vabs{\cdot}_\rho$,
$0<\rho<1$ given by
\begin{equation}\label{rho}
  \vabs*{\sum_{n>\!>-\infty}[x_n]p^n}_\rho\defeq\sup_{n\in\Z}(\vabs{x_n}\rho^n).
  \end{equation}
It is endowed with a continuous $K$-semi-linear endomorphism defined by
\[ \varphi\left(\sum_n [x_n]p^n\right)\defeq\sum_n[x_n^p]p^n\]
and we define $\varphi_q\defeq\varphi^f$ which is $K$-linear.

Let $\F\bbra{x_0^{1/p^\infty},\dots,x_{d-1}^{1/p^\infty}}$ be the
completion of the perfection of $\F[x_0,\dots,x_{d-1}]$ for the
$(x_0,\dots,x_{d-1})$-topology. If $R$ is a perfectoid $\F$-algebra
and $(r_0,\dots,r_{d-1})\in(R^{\ronron})^d$, let
\begin{equation*}F(r_0,\dots,r_{d-1})\defeq\sum_{n\in\Z}\sum_{i=0}^{d-1}[r_i^{p^{-i-nf}}]p^{i+nd}\in\B^+(R)^{\varphi_q=p^d}
\end{equation*}
then we have:

\begin{lem1}\label{lemm:BW}
  Let $1\leq d \leq f$. For each perfectoid $\F$-algebra $R$, the
  following functorial map is a bijection:
  \[
    \begin{array}{rcc}
      \Hom_{\F-\alg}^{\cont}(\F\bbra{x_0^{1/p^\infty},\dots,x_{d-1}^{1/p^\infty}},R)\simeq(R^{\ronron})^d&\longrightarrow&\B^+(R)^{\varphi_q=p^d}\\
      (r_0,\dots,r_{d-1})&\longmapsto & F(r_0,\dots,r_{d-1}).
    \end{array}
  \]
\end{lem1}

\begin{proof}
  This follows from \cite[Prop.~II.2.5(iv)]{FS}. See also
  \cite[Prop.~4.2.1]{FF} for the case where $R$ is an algebraically closed
  perfectoid field.
\end{proof}

\begin{rem1}
  If $R$ is a Huber ring over $\F$ and $R^+\subset R^\circ$ is an open and integrally closed subring ($(R,R^+)$ is then called a \emph{Huber pair}), we have
  $R^{\ronron}\subset R^+$ so that, by \cite[Prop.~2.1(i)]{HuberFS}
  \begin{multline*}
    \Hom_{\Spa(\F)}\big(\Spa(R,R^+),\Spa(\F\bbra{x_0^{1/p^\infty},\dots,x_{d-1}^{1/p^\infty}})\big)\\
    \simeq\Hom_{\F-\alg}^{\cont}(\F\bbra{x_0^{1/p^\infty},\dots,x_{d-1}^{1/p^\infty}},R).
    \end{multline*}
  Thus Lemma \ref{lemm:BW} and \cite[Lemma 18.1.1]{SWBerkeley} imply
  that the functor $(R,R^+)\!\mapsto\!\B^+(R)^{\varphi_q=p^d}$ can be
  extended to a sheaf on the site $\Perf_{\F}$ of perfectoid spaces
  over $\F$ endowed with either the pro-\'etale topology or the
  $v$-topology.
\end{rem1}

\begin{rem1}\label{rem:BW_epointe}
  Let $(R,R^+)$ be a perfectoid Huber pair over $\F$. If
  $x\in \Spa(R,R^+)$, then its residue field $k(x)$ is a perfectoid
  field containing $\F$ (see for example
  \cite[Cor.~6.7(ii)]{Scholzeperfectoid}). If
  $z\in \B^+(R)^{\varphi_q=p^d}$, we let $z_x$ be its image in
  $\B^+(k(x))^{\varphi_q=p^d}$. Then the functorial bijection of Lemma
  \ref{lemm:BW} induces a functorial bijection:
  \[
    \big(\Spa(\F\bbra{x_0^{1/p^\infty},\dots,x_{d-1}^{1/p^\infty}})\setminus
    V(x_0,\dots,x_{d-1})\big)(R,R^+) \!\simeq \!\set{z\in B^+(R)^{\varphi_q=p^{d}}\!,\ z_x \neq0 \ \forall \ x}. \]
\end{rem1}

The following remark will be used in \S\ref{explicitss}.

\begin{rem1}\label{rema:norm1}
  There exists a norm $\vabs{\cdot}_1$ on $\B^+(R)$ which induces
  on $W(R^\circ)[1/p]$ the norm
  \[ \vabs*{\sum_{n>\!>-\infty}[x_n]p^n}_1=\sup\set{\vabs{x_n},\ n\in\Z}\in [0,1]\subset \R_{\geq 0}\] and is such that
  $\vabs{x}_1=\lim_{\substack{\rho <1 \\
        \rho\rightarrow1}}\vabs{x}_\rho$ (see
  \cite[Prop.~1.10.5 \& Prop.~1.6.16]{FF}). Equivalently, there exists a valuation
  $v_0 : \B^+(R)\rightarrow [0,+\infty]$ such that
  \[ \forall \ x=\sum_{n>\!>-\infty}[x_n]p^n\in W(R^\circ)[1/p], \quad
    v_0(x)=\inf\set{v(x_n),\ n\in\Z}, \] where $v$ is the valuation
  $-\log\abs $ on $R$. This description implies that if $(x_{-n})_{n\geq0}$ is a sequence
  of elements of $R^{\ronron}$ such that $\sum_{n\leq 0}[x_{n}]p^n\in \B^+(R)$ and such that there exists $0\leq c<1$ with
  $\vabs{x_{-n}} \leq c$ for all $n\geq0$, then
  \[ \vabs*{\sum_{n\leq 0}[x_{n}]p^n}_1 \leq c.\]
  Note that $\abs_1 : \B^+(R)\rightarrow [0,1]$ is not continuous since, for instance, $\vabs{p^n}_1=1$
  for $n\in \Z$ although $p^n\rightarrow 0$ in $\B^+(R)$ when $n\rightarrow +\infty$
  (in fact $\abs_1$ induces the discrete topology on $K\subset \B^+(R)$). 
\end{rem1}

Now we review the interpretation of $\B^+(-)^{\varphi_q=p^d}$ in
terms of $p$-divisible groups in the two extreme cases $d=1$ and $d=f$.

\paragraph{The case $d=1$}

Let $G_{\LT}$ be the Lubin--Tate formal group of \S\ref{LT}. As at the end of {\it loc.~cit.}~we choose an isomorphism
$G_{\LT}\simeq\Spf(\cO_K\bbra{T_K})$ such that the logarithm map
$\log_{G_{\LT}} : G_{\LT}^{\rig}\rightarrow\mathbb{G}_{a,K}^{\rig}$ (where $\mathbb{G}_{a,K}$ is the additive formal group over $\cO_K$ and ``rig'' the rigid analytic generic fiber) is given by the series $\sum_{n\geq0}p^{-n}T_K^{q^n}$. Let $\widetilde{G}_{\LT}\defeq \varprojlim_p
(G_{\LT}\times_{\Spf(\cO_K)}\Spf(\Fq))\simeq\Spf(\Fq\bbra{T_K^{1/p^\infty}})$ be the universal
cover of $G_{\LT}\times_{\Spf(\cO_K)}\Spf(\Fq)$ (see for instance \cite[\S 3.1]{SWdivisible}). The action of $\cO_K$ on $G_{\LT}$ extends to an
action of $K$ on $\widetilde{G}_{\LT}$. Note that if $R$ is a perfectoid $\F$-algebra, we have $\widetilde{G}_{\LT}(R)\buildrel\sim\over\rightarrow (G_{\LT}\times_{\Spf(\cO_K)}\Spf(\Fq))(R)\simeq G_{\LT}(R)$ (see for example \cite[Prop.~3.1.3(iii)]{SWdivisible}) so that $G_{\LT}(R)$ already has a
  structure of a $K$-vector space. We also have $\widetilde{G}_{\LT}(R^\circ)\buildrel\sim\over\rightarrow \widetilde{G}_{\LT}(R)$. By \cite[Prop.~4.4.5]{FF} or \cite[Prop.~II.2.2]{FS},
for each perfectoid $\F$-algebra $R$, we have an isomorphism of
$K$-vector spaces
$\widetilde{G}_{\LT}(R^\circ)\xrightarrow{\sim}
\B^+(R)^{\varphi_q=p}$ given by
\begin{equation}\label{eq:F(r)}
 r\in R^{\ronron}\simeq\widetilde{G}_{\LT}(R^\circ)\longmapsto
  F(r)\defeq \sum_{n\in\Z}[r^{q^{-n}}]p^n\in \B^+(R)^{\varphi_q=p}
\end{equation}
(this is the isomorphism of Lemma \ref{lemm:BW} when $d=1$, where the variable $x_0$ in {\it loc.~cit.} is denoted by $T_K$). Note that on the left-hand side, the $K$-linear structure is given by (for $r\in R^{\ronron}$)
\begin{equation}\label{actionLT}
  \begin{cases}
    \forall \ n\in\Z, & p^n(r)=r^{q^n}, \\
    \forall \ a\in\cO_K, & a(r)=a_{\LT}(r),
  \end{cases}
\end{equation}
where we view the coefficients of the power series $a_{\LT}$ in $\F$
via $\cO_K\twoheadrightarrow \Fq\buildrel\sigma_0\over\hookrightarrow
\F$. We let
\begin{equation*}Z_{\LT}\defeq((\widetilde{G}_{\LT}\times_{\Spf(\Fq)}\Spf(\F))^{\ad}\setminus\set{0})^f,\end{equation*}
where $(\widetilde{G}_{\LT}\times_{\Spf(\Fq)}\Spf(\F))^{\ad}$ is the adic space associated to the
formal scheme
$\widetilde{G}_{\LT}\times_{\Spf(\Fq)}\Spf(\F)$ and $\set{0}$ is
the closed analytic subspace image of the $0$-section, i.e.\ $f$-times the fiber product of $(\widetilde{G}_{\LT}\times_{\Spf(\Fq)}\Spf(\F))^{\ad}\setminus\set{0}$ over $\Spa(\F)$ (still using $\sigma_0$). Using obvious notation, we have an
isomorphism of adic spaces
\begin{equation*}Z_{\LT}\simeq\Spa(\F\bbra{T_{K,0}^{1/p^\infty},\dots,
  T_{K,f-1}^{1/p^\infty}})\setminus V(T_{K,0}\cdots T_{K,f-1})\simeq\prod_{i=0}^{f-1}\Spa(\F\ppar{T_{K,i}^{1/p^\infty}},\F\bbra{T_{K,i}^{1/p^\infty}})
\end{equation*}
and there is an action of $(K^\times)^f$ on $Z_{\LT}$ given by
\[ \forall \ \underline{a}=(a_0,\dots,a_{f-1})\in (K^\times)^f, \quad
  \underline{a}(T_{K,i})=a_{i,\LT}(T_{K,i}).\]

\paragraph{The case $d=f$}

Let $\mathcal{G}_{f,f}$ be the $p$-divisible group over $\F_p$ defined
in \cite[\S 4.3.2]{FF} (with $\cO=\Zp$) as the kernel of $V^f-1$ on the group scheme
of Witt covectors $CW$ (we use without comment the notation of {\it loc.~cit.}, for instance $V$ is the Verschiebung, $F$ is the Frobenius, see \cite[\S 1.10.2]{FF} for $CW$, etc.). The base change of $\mathcal{G}_{f,f}$ to $\F$ is endowed with an additional structure of functor in
$\cO_K$-modules. Namely if $R$ is an $\F$-algebra, then $CW(R)$ is an
$\cO_K=W(\F_q)$-module via $\sigma_0:\F_q\hookrightarrow \F$ and the action of $\cO_K$ on $CW(R)$ commutes
with $V^f$ and $F^f$ (but not with $V$ and $F$).

As $\ker(V-1)\subset\ker(V^f-1)$, there is a natural injection of
$p$-divisible groups
$\mathcal{G}_{1,1}\hookrightarrow\mathcal{G}_{f,f}$ which induces a
morphism of $p$-divisible groups over $\F$ with $\cO_K$-action
\begin{equation}\label{eq:map_pdvi_OK}
  \cO_K\otimes_{\Z_p}\mathcal{G}_{1,1,\F}\longrightarrow\mathcal{G}_{f,f,\F}.
\end{equation}

\begin{lem1}\label{lemm:isoG_dd}
  The map \eqref{eq:map_pdvi_OK} is an isomorphism of $p$-divisible
  groups over $\F$ with $\cO_K$-action.
\end{lem1}

\begin{proof}
 In this proof we will use (contravariant) Dieudonn\'e Theory $\mathbb{D}(-)$ for $p$-divisible groups over $\Fp$. Recall that it yields free $\Zp$-modules, and that when the $p$-divisible group is over $\F$ it yields free $W(\F)$-modules. The map (\ref{eq:map_pdvi_OK}) corresponds to a nonzero map of Dieudonn\'e modules which is both $\cO_K$-linear and $W(\F)$-linear:
  \begin{multline}\label{mapdieud}
    \mathbb{D}(\mathcal{G}_{f,f,\F})\longrightarrow\mathbb{D}(\cO_K\otimes_{\Z_p}\mathcal{G}_{1,1,\F})\simeq\Hom_{\Z_p-{\rm mod}}(\cO_K,\mathbb{D}(\mathcal{G}_{1,1,\F}))\\
    \simeq \Hom_{\Z_p-{\rm mod}}(\cO_K,W(\F))\otimes_{W(\F)}\mathbb{D}(\mathcal{G}_{1,1,\F})
    \end{multline}
    where $\cO_K$ acts on the right-hand side via its natural action on $\Hom_{\Z_p-{\rm mod}}(\cO_K,W(\F))$. Note that $\mathbb{D}(\mathcal{G}_{f,f,\F})=W(\F)\otimes_{\Z_p}\mathbb{D}(\mathcal{G}_{f,f})$, where the Dieudonn\'e module $\mathbb{D}(\mathcal{G}_{f,f})$ has a
  $\Z_p$-basis $(e_0,e_1\defeq V(e_0),\dots,e_{f-1}\defeq V^{f-1}(e_0))$ such that $F(e_i)=pe_{i-1}$
  for all $0\leq i\leq f-1$ (see \cite[\S 4.3.2]{FF}, we write $e_0:\mathcal{G}_{f,f}\hookrightarrow CW$ for the canonical embedding $e$ of {\it loc.~cit.}~and we use the convention that
  $i=i+f$). Moreover the action of $\cO_K$ on $\mathcal{G}_{f,f,\F}$
  induces an action of $\cO_K$ on $\mathbb{D}(\mathcal{G}_{f,f,\F})=W(\F)\otimes_{\Z_p}\mathbb{D}(\mathcal{G}_{f,f})$
  such that $a(1\otimes e_i)=\varphi^{-i}(a)\otimes e_i$ for $a\in\cO_K$, where $\varphi$ is the absolute Frobenius on $W(\F)$ and $\cO_K$ is seen in $W(\F)$ via $\sigma_0:\F_q\hookrightarrow \F$. Using the $\cO_K$- and $W(\F)$-linearities, and the commutativity with $F$, one checks that there is an isomorphism $W(\F)\simeq\mathbb{D}(\mathcal{G}_{1,1,\F})$ such that the map (\ref{mapdieud}) is given by
  \[\sum_{i=0}^{f-1}\lambda_i\otimes e_i\longmapsto\left(a\longmapsto
      \sum_{i=0}^{f-1}\lambda_i\varphi^{-i}(a)\right) \in \Hom_{\Z_p-{\rm mod}}(\cO_K,W(\F))\]
      (in particular, $e_0$ maps to the inclusion $\cO_K\hookrightarrow W(\F)$). To conclude the proof we need to show that the elements $a\mapsto \varphi^{-i}(a)$, $i\in \{0,\dots, f-1\}$, generate the $W(\F)$-module $\Hom_{\Z_p-{\rm mod}}(\cO_K,W(\F))$. This can be
  checked after reduction mod $p$ and we have to prove that the
  elements $a\mapsto a^{p^i}$, $i\in \{0,\dots, f-1\}$, generate the $\F$-vector space
  $\Hom_{\F_p-{\rm vs}}(\F_q,\F)$, which is a consequence of the linear
  independence of characters.
\end{proof}

By \cite[Prop~4.4.5]{FF} (replacing $\Fpbar$ by $\F$, the field $F$ by a perfectoid $\F$-algebra $R$ and where the variable $x_i$ of {\it loc.~cit.}~is reindexed $x_{f-i}$ here for $i\in \{1,\dots, f-1\}$, $x_0$ being unchanged) there exists a coordinate $z$ (resp.~coordinates $x_0,\dots,x_{f-1}$) on
the formal group $\mathcal{G}_{1,1,\F}$ (resp.~$\mathcal{G}_{f,f,\F}$) such that the following map is an isomorphism
of $\Qp$-vector spaces (resp.~$K$-vector spaces) for any perfectoid $\F$-algebra $R$:
\begin{equation}\label{formulaf}
\begin{gathered}
  \gamma_1 :
  \begin{array}{ccc}
  {\mathcal{G}}_{1,1,\F}(R) & \xrightarrow{\sim} & \B^+(R)^{\varphi=p}\\
  z & \xrightarrow{\sim} & \sum_{n\in\Z}[z^{p^{-n}}]p^n
  \end{array}\\
 \left(\textrm{resp. }\gamma_f : 
  \begin{array}{ccc}
    {\mathcal{G}}_{f,f,\F}(R) &\xrightarrow{\sim}&\B^+(R)^{\varphi_q=p^f}\\
    (x_0,\dots,x_{f-1}) & \xrightarrow{\sim} & \sum_{i=0}^{f-1}\sum_{n\in\Z}[x_{i}^{p^{-i-nf}}]p^{i+nf}
  \end{array}\right)
\end{gathered}
\end{equation}
(we use $\widetilde{\mathcal{G}}_{1,1,\F}(R)\buildrel\sim\over\rightarrow {\mathcal{G}}_{1,1,\F}(R)$ by \cite[Prop.~3.1.3(iii)]{SWdivisible} for the structure of $\Qp$-vector space on ${\mathcal{G}}_{1,1,\F}(R)$, likewise with ${\mathcal{G}}_{f,f,\F}(R)$).
Moreover these isomorphisms are given by the composition of the isomorphisms in
the following diagram (we only give $\gamma_f$ and refer to {\it loc.~cit.} for the notation):
\begin{multline}\label{diagr}
      \mathcal{G}_{f,f,\F}(R)\xrightarrow{\sim}
      \Hom_{W(\F)[F]}(\mathbb{D}(\mathcal{G}_{f,f,\F}),CW(R))\xleftarrow{\sim}
      \Hom_{W(\F)[F]}(\mathbb{D}(\mathcal{G}_{f,f,\F}),BW(R)) \\
      \xrightarrow{\sim}\Hom_{W(\F)[F]}(\mathbb{D}(\mathcal{G}_{f,f,\F}),\B^+(R))
    =  \B^+(R)^{\varphi_q=p^f},
\end{multline}
where the second isomorphism is a consequence of
\cite[Prop.~4.4.2]{FF} and the third a
consequence of \cite[Prop.~4.2.1]{FF}. We deduce from (\ref{diagr}) the commutativity of the following
diagram of $\Qp$-vector spaces:
\[\begin{tikzcd}
    {\mathcal{G}}_{1,1,\F}(R) \ar[r,"\gamma_1"] \ar[d] &
    \B^+(R)^{\varphi=p} \ar[d] \\
    {\mathcal{G}}_{f,f,\F}(R) \ar[r,"\gamma_f"] & \B^+(R)^{\varphi_q=p^f}
  \end{tikzcd}\]
and thus the commutativity of the following diagram of $K$-vector spaces:
\begin{equation}\label{commuteR}
  \begin{tikzcd}
    \cO_K\otimes_{\Z_p}{\mathcal{G}}_{1,1,\F}(R)
    \ar[r,"\Id_{\cO_K}\!\otimes\gamma_1"] \ar[d] &
    \cO_K\otimes_{\Z_p}\B^+(R)^{\varphi=p} \ar[d,"\simeq"] \\
    {\mathcal{G}}_{f,f,\F}(R) \ar[r,"\gamma_f"] & \B^+(R)^{\varphi_q=p^f}.
  \end{tikzcd}
\end{equation}
Let $\widehat{\mathbb{G}}_{m,\F_p}$ be the multiplicative formal group over $\F_p$ and $\widehat{\mathbb{G}}_{m,\F}$ its base change over $\F$, we have $\mathcal{G}_{1,1}\simeq\widehat{\mathbb{G}}_{m,\F_p}$ (see
\cite[Ex.~4.4.7]{FF}) and isomorphisms of $p$-divisible groups over $\F$ with $\cO_K$-action
\begin{equation}\label{11ff}
\Hom_{\Z_p-{\rm mod}}(\cO_K,\Z_p)\otimes_{\Z_p}\widehat{\mathbb{G}}_{m,\F}\simeq \cO_K\otimes_{\Z_p}{\mathcal{G}}_{1,1,\F}\buildrel\sim\over \rightarrow {\mathcal{G}}_{f,f,\F}
\end{equation}
using the isomorphism of $\cO_K$-modules
\begin{equation}\label{traceK}
\cO_K\buildrel\sim\over\longrightarrow \Hom_{\Z_p-{\rm mod}}(\cO_K,\Z_p),\ a\mapsto\Tr_{K/\Q_p}(a\cdot)
\end{equation}
and Lemma \ref{lemm:isoG_dd}. Here, the $\cO_K$-action on the left-hand side of (\ref{11ff}) is via the action of $\cO_K$ on $\Hom_{\Z_p-{\rm mod}}(\cO_K,\Z_p)$ given by $a(\lambda)\!=\!\lambda(a-)$ ($a\!\in\!\cO_K$, $\lambda\!\in\!\Hom_{\Z_p-{\rm mod}}(\cO_K,\Z_p)$). Using
\[\Hom_{\F-\alg}^{\cont}(\F\bbra{\cO_K},A)\simeq \Hom_{\Z_p-{\rm mod}}(\cO_K,A^{\ronron})\simeq \Hom_{\Z_p-{\rm mod}}(\cO_K,\Z_p)\otimes_{\Z_p}A^{\ronron}\]
for any complete topological $\F$-algebra $A$, we deduce from (\ref{11ff}) an isomorphism of formal modules over $\F$ with $\cO_K$-action
\begin{equation}\label{gff}
{\mathcal{G}}_{f,f,\F}\simeq\Spf(\F\bbra{\cO_K}),
\end{equation}
where $\cO_K$ acts (continuously) on $\F\bbra{\cO_K}$ by multiplication on itself. It follows that $\widetilde{\mathcal{G}}_{f,f,\F}\defeq \varprojlim_p {\mathcal{G}}_{f,f,\F}$ is represented by the formal scheme $\Spf(\F\bbra{K})$, where $\F\bbra{K}$ is the $\mathfrak{m}_{\cO_K}$-adic completion of $\F[K]\otimes_{\F[\cO_K]}\F\bbra{\cO_K}$
($\mathfrak{m}_{\cO_K}$ being the maximal ideal of $\F\bbra{\cO_K}$). It also follows from the formula for $\gamma_f$ in (\ref{formulaf}) that there exist elements $X_0,\dots,X_{f-1}\in\F\bbra{\cO_K}$ satisfying
$\F\bbra{\cO_K}=\F\bbra{X_0,\dots,X_{f-1}}$ such that we have isomorphisms ${\mathcal{G}}_{f,f,\F}(R)\simeq \Hom_{\F-\alg}^{\cont}(\F\bbra{\cO_K},R)\simeq \B^+(R)^{\varphi_q=p^f}$ for any perfectoid $\F$-algebra $R$, where the second isomorphism is given by (where $X_i\mapsto r_i\in R^{\ronron}$)
\begin{equation}\label{eq:F(r,..)}
  (r_0,\dots,r_{f-1})\in (R^{\ronron})^f\longmapsto
F(r_0,\dots,r_{f-1}) \defeq \sum_{i=0}^{f-1}\sum_{n\in\Z}[r_{i}^{p^{-i-nf}}]p^{i+nf}\in  \B^+(R)^{\varphi_q=p^f}.
\end{equation}
We then easily check that, in the coordinates $X_i$, the action of $K^\times$ on $\F\bbra{K}$ has
the following properties
\begin{equation}\label{frobdown}
  \begin{cases}
    \forall\ 0\leq i\leq f-1,\, \forall\ n\in\Z, & p^n(X_i)=X_{i-n}^{p^n} \\
    \forall\ 0\leq i\leq f-1,\, \forall\ a\in\F_q^\times,&
    [a](X_i)=\sigma_0(a)^{p^i}X_i
  \end{cases}
  \end{equation}
(with the usual convention that $X_{i+f}=X_i$). Finally, we let
\begin{equation*}Z_{\cO_K} \defeq
\widetilde{\mathcal{G}}_{f,f,\F}^{\ad}\setminus\set{0},\end{equation*} 
where
$\widetilde{\mathcal{G}}_{f,f,\F}^{\ad}$ is the adic space over $\F$ associated to the formal scheme
$\widetilde{\mathcal{G}}_{f,f,\F}$. We have an isomorphism
\begin{equation*}Z_{\cO_K}\simeq
\Spa(\F\bbra{X_0^{1/p^\infty},\dots,X_{f-1}^{1/p^\infty}})\setminus V(X_0,\dots,X_{f-1}).
\end{equation*}
Note that the adic spaces $Z_{\LT}$ and $Z_{\cO_K}$ are both in
$\Perf_{\F}$.

We fix now $C$ a perfectoid field containig $\F$ and $v$ a continuous
rank $1$ valuation on $F$. If $x\in\B(C)$,
the \emph{Newton polygon} of $x$, defined in \cite[D\'ef.~1.5.2, D\'ef.~1.6.18,
D\'ef.~1.6.21]{FF} is a decreasing convex function from
$\R\rightarrow\R\cup\set{+\infty}$. From \cite[Ex.~1.6.22]{FF}, it can
be computed as the inverse Legendre transform of the function
$\lambda\mapsto v_\lambda(x)$, $\lambda\in\left[0,+\infty\right[$
(from which we remove the zero slope if it appears), where $v_\lambda$ is the continuous
extension to $\B(C)$ of the valuation $v_\lambda$ defined in
\cite[D\'ef.~1.4.1]{FF}.

\begin{lem1}\label{lemm:polygone_Newton}
  Let $(x_n)_{n\in\Z}$ be a family of nonzero elements of $\cO_C$ such
  that, for any $\lambda\in\Iouv{0,+\infty}$,
  $v(x_n)+\lambda n\rightarrow+\infty$ when $n\rightarrow-\infty$ and
  ${v}(x_n)\rightarrow0$ when $n\rightarrow+\infty$. Let
  $x {\defeq}\sum_{n\in\Z}[x_n]p^n\in\B^+(C)$.
  \begin{enumerate}\item\label{lemm:polygone_Newton1} Assume that
    $x\in\B^+(C)^{\varphi_q=p^f}$. Then the set of slopes of the
    Newton polygon of $x$ is of the form $ap^{\Z}$
for some $a>0$ if
    and only if {$v(x_n)=v(x_0)p^{-n}$ for all $n \in \Z$}. In this case we can choose $a=(p-1)v(x_0)$.
  \item\label{lemm:polygone_Newton2} Assume that
    $x\in\B^+(C)^{\varphi_q=p}$. Then the set of slopes of the Newton
    polygon of $x$ is $(q-1)v(x_0)q^{\Z}$.
  \end{enumerate}
\end{lem1}

\begin{proof}
  We prove \ref{lemm:polygone_Newton1}, \ref{lemm:polygone_Newton2}
  being similar and simpler. For $\lambda\in\Iouv{0,+\infty}$ we
  define $f_x(\lambda)\defeq v_\lambda(x)$. By \cite[\S1.5.1]{FF}, the
  set of slopes of the Newton polygon of $x$ is the set of breakpoints
  of $f_x$. For any compact interval $[a,b]\subset\Iouv{0,+\infty}$,
  there exists an integer $N>0$ such that
  \[ \forall\lambda\in[a,b],\qquad
    f_x(\lambda)=v_\lambda(x)=\inf_{n\in [-N,N] \cap \Z}(v(x_n)+n\lambda). \] This shows that the function $f_x$
  has only finitely  {many} slopes in the interval $[a,b]$ and that these
  slopes are in the set $[-N,N]\cap\Z$. Therefore the function $f_x$
  has integral slopes. Moreover, as the breakpoints of $f_x$ in
  $[a,b]$ are coordinates of the intersection points of finitely many
  lines, $f_x$ has finitely many breakpoints in $[a,b]$. The
  relation $\varphi_q(x)=p^fx$ implies
  $qv_{\lambda/q}(x)=v_\lambda(x)+f\lambda$ for all $\lambda >0$. As a
  consequence, the set of breakpoints of the function $f_x$ is
  stable under multiplication by $q$ and $q^{-1}$. Moreover, if
  $\lambda$ is a regular point of $f_x$, so is $q\lambda$ and
  $f'_x(q\lambda)=f'_x(\lambda)-f$. Let us fix $\lambda_0 >0$ some
  regular point of $f_x$ and let $s {\defeq}f'_x(\lambda_0)\in\Z$. As
  $f'_x(q\lambda_0)=s-f$, between $\lambda_0$ and $q\lambda_0$, $f_x$
  can have at most $f+1$ different slopes and so at most $f$ breakpoints. 
  These breakpoints are representatives of the quotient of
  the set of breakpoints of $f_x$ by $q^{\Z}$. Therefore if the
  set of breakpoints of $f_x$ is of the form $ap^{\Z}$, for some
  $a>0$, then $f_x$ has exactly $f$ breakpoints in
  $[\lambda_0,q\lambda_0]$ which are of the form
  $\lambda_1,p\lambda_1,\dots,p^{f-1}\lambda_1$ and the successive
  slopes in this interval are $s,s-1,\dots,s-(f-1),s-f$. This implies that
  $\lambda_1$ is the coordinate of the intersection point of the
  graphs of the functions $\lambda\mapsto v(x_s)+s\lambda$ and
  $\lambda\mapsto v(x_{s-1})+(s-1)\lambda$, that is
  $\lambda_1=v(x_{s-1})-v(x_s)$. Similarly, we have
  $p^j\lambda_1=v(x_{s-j-1})-v(x_{s-j})$ for all $0\leq j\leq
  f-1$. Thus we have the relation
  $v(x_{n-1})-v(x_n)=p(v(x_n)-v(x_{n+1}))$ for all $n\in\Z$. As the
  sequence $(v(x_n))_{n\in\Z}$ is not constant and has limit $0$ as $n \to +\infty$, we easily deduce that $v(x_n)=v(x_0)p^{-n}$ for all
  $n\in\Z$. Moreover, as the computation shows, $v(x_0)-v(x_1)$ is the
  coordinate of a breakpoint up to some power of $p$, {hence} we can
  choose $a=p(v(x_0)-v(x_1))=(p-1)v(x_0)$. A direct computation
  following the same lines shows, conversely, that if
  $v(x_n)=p^{-n}v(x_0)$ for all $n\in\Z$, then the set of breakpoints 
  of $f_x$ is $(p-1)v(x_0)p^{\Z}$.
\end{proof}

{Recall the series $F(r) \in \B^+(R)^{\vp_q=p}$ from~\eqref{eq:F(r)} and $F(x_0,\dots,x_{f-1}) \in \B^+(R)^{\vp_q=p^f}$ from~\eqref{eq:F(r,..)}}

\begin{cor1}\label{coro:polygone_Newton}\
  \begin{enumerate}\item\label{coro:polygone_Newton1} Let $t_0,\dots,t_{f-1}\in C$ such
    that $v(t_0),\dots,v(t_{f-1})>0$. The set of slopes of the Newton
    polygon of $F(t_0,\dots,t_{f-1})$ is of the form $ap^{\Z}$ for
    some $a>0$ if and only if $v(t_0)=v(t_1)=\cdots=v(t_{f-1})$. 
    In this case we can choose $a=(p-1)v(t_0)$.
  \item\label{coro:polygone_Newton2} Let $t\in C$. The set of slopes
    of the Newton polygon of $F(t)$ is $(q-1)v(t_0)q^{\Z}$.
  \end{enumerate}
\end{cor1}

\begin{proof}
  This is a direct consequence of Lemma \ref{lemm:polygone_Newton}.
\end{proof}

\subsection{A ``sum of divisors'' map}\label{mapm}

We define and study certain open subspaces of the perfectoid spaces $Z_{\LT}$ and $Z_{\cO_K}$ of \S\ref{reminderp}, as well as a canonical map $m:Z_{\LT}\longrightarrow Z_{\cO_K}$ preserving these subspaces.

For any perfectoid $\F$-algebra $R$, the product in the ring $\B^+(R)$
induces a functorial map:
\begin{equation}\label{mR} m_R:
  \begin{array}{ccc}
    (\B^+(R)^{\varphi_q=p})^f&\longrightarrow& \B^+(R)^{\varphi_q=p^f}
    \\
    (z_1,\dots,z_f)&\longmapsto & z_1\cdots z_f.
  \end{array}
\end{equation}

Using Remark \ref{rem:BW_epointe}, the fact that each $\B^+(k)$ is a
domain for $k$ a perfectoid field (see \cite[Thm.~6.2.1 \& Thm.~3.6.1]{FF}) and \cite[Prop.~8.2.8(2)]{SWBerkeley}, the family of maps
$(m_R)$ induces a morphism of perfectoid spaces over $\F$
\begin{equation}\label{mapmzlt}
m : Z_{\LT}\longrightarrow Z_{\cO_K}.
\end{equation}
The map $m_R$ being compatible with the
actions of $(K^\times)^f$ (on the source) and $K^\times$ (on the target), we deduce that
$m$ is compatible with the
actions of $(K^\times)^f$ and $K^\times$ on $Z_{\LT}$ and $Z_{\cO_K}$, i.e.~$m\circ (a_0,\dots,a_{f-1})=(\prod_i a_i)\circ m$. For $0\leq i\leq f-1$ let $j_i$ be the morphism $K^\times\rightarrow(K^\times)^f$ sending $a$ to the $f$-uple with $1$ at all entries
except at the $i$-th entry where it is $a$, then for all
$a\in K^\times$ and $0\leq i\leq f-1$, we have in particular
\begin{equation}\label{eq:ji}
m\circ j_i(a)=a\circ m:Z_{\LT}\rightarrow Z_{\cO_K}.
\end{equation}

\begin{rem1}
  The map $m$ can be reinterpreted using the Abel--Jacobi map (cf.\ \cite{FarguesAJ}). Namely the sheaf on the
  pro-\'etale site of $\Perf_{\F}$ associated to the quotient presheaf
  $(\B^+(-)^{\varphi_q=p}\setminus \{0\})/K^\times$ is isomorphic to the pro-\'etale sheaf
  $\Div_{\F}^1$ of degree $1$ divisors on the relative
  Fargues--Fontaine curve over $\F$ and likewise
  $(\B^+(-)^{\varphi_q=p^f}\setminus \{0\})/K^\times$ is isomorphic to the pro-\'etale sheaf
  $\Div^f_{\F}$ of degree $f$ divisors. The map $m$ induces a morphism
  of pro-\'etale sheaves $(\Div_{\F}^1)^f\rightarrow\Div^f_{\F}$ which
  is given by the sum of divisors, cf.~\cite[\S 2.4]{FarguesAJ}.
\end{rem1}

The group $\mathfrak{S}_f$ acts on the left on $(K^\times)^f$ by permutation of
coordinates:
\[\forall\  \sigma\in\mathfrak{S}_f, \, \forall\ (a_i)_{0\leq i\leq
    f-1}\in (K^\times)^f, \quad
  \sigma(a_i)\defeq (a_{\sigma^{-1}(i)}).\]
The group $\mathfrak{S}_f$ acts likewise on $Z_{\LT}$ by permuting the factors
$(\widetilde{G}_{\LT}\times_{\Spf(\cO_K)}\Spf(\F)\setminus\set{0})^{\ad}$ so that the action
of $(K^\times)^f$ on $Z_{\LT}$ extends to an action of the semi-direct product
$(K^\times)^f\rtimes\mathfrak{S}_f$. Let $\Delta$ be the kernel of the
multiplication $(K^\times)^f\rightarrow K^\times$ and
$\Delta_1\defeq\Delta\cap(\cO_K^\times)^f$. Then
$\Delta\rtimes\mathfrak{S}_f$ is a subgroup of
$(K^\times)^f\rtimes\mathfrak{S}_f$ and the map $m$ is invariant under
the action of $\Delta\rtimes\mathfrak{S}_f$. By
\cite[Lemme 7.6]{FarguesAJ}\footnote{Note that \cite[Lemme 7.6]{FarguesAJ} extends scalars to $\overline{\mathbb F}_q$, however the proof works the same without extending scalars as it is based on the proof of \cite[Prop.~2.18]{FarguesAJ} where one does not extend scalars.}, the map $m$ induces an isomorphism of
pro-\'etale sheaves on $\Perf_{\F}$
\begin{equation}\label{fargues}
\Delta\rtimes\mathfrak{S}_f\backslash Z_{\LT} \xrightarrow{\sim} Z_{\cO_K}.
\end{equation}
  
We let $Z_{\cO_K}^{\gen}$ be the open subspace of $\Spa(\F\bbra{K})$ defined by the
relations
\[ \vabs{X_0}=\dots=\vabs{X_{f-1}}\neq0\]
(it is open as it is the intersection over $i\in \{0,\dots,f-1\}$ of the rational open subsets $U(\frac{X_0,\dots,X_{f-1}}{X_i})$ of $\Spa(\F\bbra{K})$).
Note that we have $Z_{\cO_K}^{\gen}\subset Z_{\cO_K}$. We also define
$Z_{\LT}^{\gen}\defeq m^{-1}(Z_{\cO_K}^{\gen})$, an open subspace of $Z_{\LT}$, so that $Z_{\cO_K}^{\gen}$ and $Z_{\LT}^{\gen}$ are both in $\Perf_{\F}$. We now give explicit descriptions of $Z_{\cO_K}^{\gen}$ and $Z_{\LT}^{\gen}$.
  
We start with $Z_{\cO_K}^{\gen}$. We denote by $A_\infty\defeq\cO_{Z_{\cO_K}}(Z_{\cO_K}^{\gen})$ the ring of global
sections on $Z_{\cO_K}^{\gen}$.

\begin{lem1}\label{lemm:A_infty_expl}
The following statements hold.
\begin{enumerate}\item\label{lemm:A_infty_expl(i)}
    The ring $A_\infty$ is the perfectoid $\F$-algebra \[\F\ppar{X_0^{1/p^\infty}}\scalar*{\left(\frac{X_i}{X_0}\right)^{\pm
        1/p^\infty}\!\!,\ 1\leq i\leq f-1}.\]
  \item\label{lemm:A_infty_expl(ii)}
    We have $Z_{\cO_K}^{\gen}=\Spa(A_\infty,A_\infty^\circ)$, in particular $Z_{\cO_K}^{\gen}$ is affinoid perfectoid.
  \item\label{lemm:A_infty_expl(iibis)} There exists a multiplicative
    norm $\abs$ on $A_\infty$ such that $\vabs{X_0}=p^{-1}$ inducing
    the topology of $A_\infty$.
  \item\label{lemm:A_infty_expl(iii)} Any quasi-compact open subset of $Z_{\cO_K}$ whose points of rank $1$ are exactly the points of $Z_{\cO_K}^{\gen}$ of rank $1$ is necessarily $Z_{\cO_K}^{\gen}$ itself. 
  \end{enumerate}
\end{lem1}

\begin{proof}
Define the adic spaces
    \[T^{\gen}\defeq\set{\vabs{X_0}=\cdots=\vabs{X_{f-1}}\neq0}\subset
    T\defeq\Spa(\F\bbra{X_0,\dots,X_{f-1}}).\]
It is enough to prove (i), (ii) and (iii) replacing everywhere
$Z_{\cO_K}^{\gen}\subset \Spa(\F\bbra{K})$ by $T^{\gen}\subset T$
(i.e.\ completed perfection will not change the arguments in the proof
below). Moreover, as the map
$T=\Spa(\F\bbra{X_0,\dots,X_{f-1}})\rightarrow
\Spa(\F\bbra{X_0^{1/p^\infty},\dots,X_{f-1}^{1/p^\infty}})$ is a
homeomorphism, it is also enough to prove (iv) with $T^{\gen}$ and
$T\setminus V(X_0,\dots,X_{f-1})$.

We first show the analogue (iii). Let
  $S\defeq
  \F\ppar{X_0}\scalar*{\left(\frac{X_i}{X_0}\right)^{\pm1},\ 1\leq
    i\leq f-1}$ that we endow with the $X_0$-adic topology (it is a Tate algebra), then the norm in (iii) is the unique multiplicative extension to $S$ of the Gauss norm on the restricted power series $\F\ppar{X_0}\scalar*{\left(\frac{X_i}{X_0}\right),\ 1\leq i\leq f-1}$ (which is well-known to be multiplicative). Note that $S^\circ=\F\bbra{X_0}\scalar*{\left(\frac{X_i}{X_0}\right)^{\pm1},\ 1\leq i\leq f-1}$ is the unit ball for this norm.
    
Let us prove (the analogues of) (i) and (ii). Looking at continuous valuations, it is clear that the morphism of adic spaces $\Spa(S,S^\circ)\rightarrow T$ factors as $\Spa(S,S^\circ)\rightarrow T^{\gen}\subset T$. In order to prove that the morphism of adic spaces $\Spa(S,S^\circ)\rightarrow T^{\gen}$ is an isomorphism, it is enough to prove that it induces an isomorphism $\Spa(S,S^\circ)(W)\buildrel\sim\over\rightarrow T^{\gen}(W)$ for any analytic adic space $W$ over $\F$, and it is enough to take $W=\Spa(R,R^+)$ for an arbitrary complete analytic Huber pair $(R,R^+)$ over $\F$ (the case $R$ Tate would be enough). Then this easily follows from the definitions of $T$ and $S$.

Let us finally prove (the analogue of) (iv). First note that $T\setminus V(X_0,\dots,X_{f-1})$
  is the analytic locus of the adic space $T$, the only non-analytic
  point of $T$ being the unique (rank $0$) valuation with kernel the
  maximal ideal of the local ring $\F\bbra{X_0,\dots,X_{f-1}}$. Let
  $U$ be a quasi-compact open subset of $T\setminus V(X_0,\dots,X_{f-1})$ whose points of rank $1$ are the points of $T^{\gen}$ of rank $1$. For $i\in \{0,\dots,f-1\}$ consider the open subset $U_i$ of $T$ defined by $|X_j|\leq |X_i|\ne 0$ for all $j$, or equivalently (by the same argument as for the proof of (i))
\[U_i=\mathrm{Spa}\Big(\F\ppar{X_i}\left\langle \frac{X_j}{X_i}, j\ne i\right\rangle, \F\bbra{X_i}\left\langle \frac{X_j}{X_i}, j\ne i\right\rangle\Big)\subset T\setminus V(X_0,\dots,X_{f-1}).\]
Then $U\cap U_i$ and $T^{\gen}=\bigcap_j U_j$ are two open subsets of $U_i$ with the same points of rank $1$, and thus {\it a fortiori} with the same points with residue field being a finite extension of $\F\ppar{X_i}$. Let $U_i^{\rig}\subset U_i$ (resp.\ $(T^{\gen})^{\rig}\subset T^{\gen}$) be the subset of points of $U_i$ (resp.\ $T^{\gen}$) with residue field being a finite extension of $\F\ppar{X_i}$, then $U_i^{\rig}$ (resp.\ $(T^{\gen})^{\rig}$) can be identified with the affinoid rigid analytic space over $\F\ppar{X_i}$ corresponding to $U_i$ (resp. $T^{\gen}$) by \cite[(1.1.11)(a)]{Huberbook}, and we have $U \cap U_i^{\rig}=(T^{\gen})^{\rig}$. Note that $U$, $U_i$ and $T^{\gen}$ are quasi-compact ($U$ by assumption, $U_i$, $T^{\gen}$ as they are affinoid). As $T$ is a quasi-separated adic space (being spectral as the adic space associated to a Huber pair, see for instance \cite[Cor.~III.2.4]{Morel}), the open subset $U\cap U_i$ is still quasi-compact. As $U_i^{\rig}$ is quasi-separated, we deduce $U\cap U_i=T^{\gen}$ from $U \cap U_i^{\rig}=(T^{\gen})^{\rig}$ by \cite[(1.1.11)]{Huberbook} (see also \cite[Thm.~2.21]{Scholzeperfectoid}). Since $U=\bigcup_i (U\cap U_i)$ (as $U\subset T\setminus V(X_0,\dots,X_{f-1})$), we finally obtain $U= T^{\gen}$ in $T$.
\end{proof}

\begin{lem1}\label{lemm:stability}
The following statements hold.
  \begin{enumerate}\item\label{lemm:stab1} The open subset $Z_{\cO_K}^{\gen}$ of
    $Z_{\cO_K}$ is stable under the action of $K^\times$.
  \item\label{lemm:stab2} The open subset $Z_{\LT}^{\gen}$ of $Z_{\LT}$ is stable
    under the action of $(K^\times)^f\rtimes\mathfrak{S}_f$.
  \end{enumerate}
\end{lem1}

\begin{proof}
  (ii) can be easily deduced from (i) and $Z_{\LT}^{\gen}\defeq m^{-1}(Z_{\cO_K}^{\gen})$, so we only prove (i).
  
  The fact that
  $Z_{\cO_K}^{\gen}$ is stable under the actions of $p$ and $p^{-1}$ on $Z_{\cO_K}$
  is a direct computation on $\F\bbra{X_0^{1/p^\infty},\dots,X_{f-1}^{1/p^\infty}}\hookrightarrow A_\infty^\circ$ using (\ref{frobdown}). Let us show that $Z_{\cO_K}^{\gen}$ is stable under the action of
  $\cO_K^\times$. It follows from Lemma \ref{lemm:A_infty_expl}(iv)
  that it is sufficient to check that $Z_{\cO_K}^{\gen}(C,\cO_C)$ is
  stable under the action of $\cO_K^\times$ on $\Spa(\F\bbra{K})(C,\cO_C)=\B^+(C)^{\varphi_q=p^f}$ for $C$ a perfectoid field containing $\F$ (using $Z_{\cO_K}(C,C^+)\buildrel\sim\over\rightarrow Z_{\cO_K}(C,\cO_C)$ for any open bounded valuation  subring $C^+\subset C$). Recall that $\B^+(C)^{\varphi_q=p^f}$ is the set of converging power series in $\B^+(C)$:
\[F(x_0,\dots,x_{f-1})=\sum_{n\in\Z}\sum_{i=0}^{f-1}[x_i^{p^{-i-nf}}]p^{i+nf}\]
where $\vabs{x_{i}}<1$ for all $i$ with $\vabs{\cdot}$ a fixed power-multiplicative norm on $C$ (e.g.~as in (\ref{normul})). A point $x\in \B^+(C)^{\varphi_q=p^f}$ is in $Z_{\cO_K}^{\gen}(C,\cO_C)$ if and
  only if $0\ne \vabs{x_{0}}=\cdots = \vabs{x_{f-1}}<1$, equivalently
  if and only if its Newton polygon has slopes $\{cp^n,\ n\in \Z\}$
  for some $c>0$ by
  Corollary \ref{coro:polygone_Newton}. As the Newton polygon of $x$ only depends on the norms $\vabs{x}_\rho$ for $0<\rho <1$ (see \cite[Ex.~1.6.22]{FF} and (\ref{rho}) for $\vabs{\cdot}_\rho$), it is enough to show that
  $\vabs{x}_{\rho}$ does not change if we multiply $x$ by an element of
  $\cO_K^\times$. This follows from the multiplicativity of
  $\abs_\rho$ (see \cite[Prop.~1.4.9]{FF}) and the fact that
  $\abs_\rho$ induces the $p$-adic norm on $K$. 
\end{proof}

From Lemma \ref{lemm:stability} we deduce a continuous action of
$K^\times$ on the topological $\F$-algebra $A_\infty$. We denote by
$\varphi$ the endomorphism of $A_\infty$ induced by the action of
$p\in K^\times$. It is $\F$-linear and satisfies (see (\ref{frobdown}))
\begin{equation}\label{phionainfini}
\varphi(X_i)=X_{i-1}^p\ {\rm for}\ 0\leq i\leq f-1
\end{equation}
(with $X_{-1}=X_{f-1}$ as usual). We also note $\varphi_q\defeq\varphi^f$ (which coincides with $x\mapsto x^q$ on $A_\infty$ when $\Fq=\F$).

We now give an explicit description of $Z_{\LT}^{\gen}$.

Recall first that if a locally profinite group $H$ acts continuously on a perfectoid space $X'$ over $\F$, a morphism $X'\rightarrow X$ in $\Perf_{\F}$ ($H$ acting trivially on $X$) is a pro-\'etale $H$-torsor if there exists a pro-\'etale cover $Y\rightarrow X$ in $\Perf_{\F}$ such that there is an isomorphism $X'\times_XY\simeq \underline H \times Y$ in $\Perf_{\F}$, where $\underline H$ is the sheaf on $\Perf_{\F}$ defined by $\underline H(T)\defeq {\rm Cont}(\vert T\vert, H)$, $\vert T\vert$ being the underlying topological space of the perfectoid space $T$ (note that $\underline H \times Y$ is perfectoid by \cite[Lemma 10.13]{ScholzeEtcohDia}).

Let $\Z^f/\Z$ be the additive group quotient of $\Z^f$ by the diagonal embedding of $\Z$
into $\Z^f$. If $\underline{n}=(n_0,\dots,n_{f-1})\in\Z^f/\Z$ we let
$U_{\underline{n}}$ be the open affinoid perfectoid subspace of $Z_{\LT}\subset\Spa(\F\bbra{T_{K,0}{^{1/p^\infty}},\dots,T_{K,f-1}{^{1/p^\infty}}})$ defined
by the relations
\[ \vabs{T_{K,i}}^{p^{n_j}}=\vabs{T_{K,j}}^{p^{n_i}}\neq0, \qquad \forall\ 
  0\leq i,j\leq f-1\]
or equivalently $\vabs{T_{K,i}}=\vabs{T_{K,0}}^{p^{n_i-n_0}}$ for $0\leq i\leq f-1$. Note that $U_{\underline{n}}$ is well-defined as it only depends on the class of $\underline{n}$ in $\Z^f/\Z$, and that $U_{\underline{n}}$ is disjoint from $U_{\underline{n}'}$ if $\underline{n}\ne \underline{n}'$ in $\Z^f/\Z$. The group $\mathfrak{S}_f$ acts on $\Z^f/\Z$ by
permutation, for $\sigma\in\mathfrak{S}_f$ and
$\underline{n}\in\Z^f/\Z$ we have
\[ \sigma(\underline{n})\defeq (n_{\sigma^{-1}(i)})_{0\leq i\leq f-1}\] and
we check that
$\sigma(U_{\underline{n}})=U_{\sigma(\underline{n})}$. Moreover, if
$\underline{a}=(a_0,\dots,a_{f-1})\in (K^\times)^f$, we also easily check that (where $v_p$ is the unique valuation on $K$ with $v_p(p)=1$):
\[
  \underline{a}(U_{\underline{n}})=U_{\underline{n}+fv_p(\underline{a})}.\]

\begin{prop1}\label{prop:Delta1_torsor}
  Let $\underline{n}_0\defeq (0,1,\dots,f-1)$ and let $(\Z^f/\Z)_0$ be the image in $\Z^f/\Z$ of the subgroup of $\Z^f$ of $\underline{m}=(m_0,\dots,m_{f-1})$ such that $\sum_{i=0}^{f-1}m_i=0$. We have in $Z_{\LT}$
  \begin{equation}\label{eq:explicite}
    Z_{\LT}^{\gen}=\bigcup_{\sigma\in\mathfrak{S}_f}\bigcup_{\underline{m}\in\Z^f/\Z}U_{\sigma(\underline{n}_0)+f\underline{m}}=\coprod_{\sigma\in\mathfrak{S}_f}\coprod_{\underline{m}\in(\Z^f/\Z)_0}U_{\sigma(\underline{n}_0)+f\underline{m}}= \!\!\coprod_{\gamma \in (\Delta\rtimes\mathfrak{S}_f)/\Delta_1}\gamma(U_{\underline{n}_0}).
  \end{equation}
  Moreover for each $U_{\underline{n}}$ in (\ref{eq:explicite}) the map $m : Z_{\LT}^{\gen}\rightarrow Z_{\cO_K}^{\gen}$ restricts to a pro-\'etale $\Delta_1$-torsor $m|_{U_{\underline{n}}} : U_{\underline{n}}\rightarrow Z_{\cO_K}^{\gen}$.
\end{prop1}

\begin{proof}
One first easily checks that any element in $\Z^f/\Z$ of the form $\sigma(\underline{n}_0)+f\underline{m}$ can uniquely be written (in $\Z^f/\Z$) as $\sigma'(\underline{n}_0)+f\underline{m'}$ for a unique $\sigma'\in \mathfrak{S}_f$ and a unique $\underline{m'}\in (\Z^f/\Z)_0$. Assuming $\gamma(U_{\underline{n}_0})=U_{\underline{n}_0}$ when $\gamma\in \Delta_1$ this gives the last two equalities in (\ref{eq:explicite}) (recall that $\Delta_1$ is normal in $\Delta \rtimes\mathfrak{S}_f$).

We check that $Z_{\LT}^{\gen}$ and
$\bigcup_{\sigma\in\mathfrak{S}_f}\bigcup_{\underline{m}\in\Z^f/\Z}U_{\sigma(\underline{n}_0)+f\underline{m}}$
have the same rank $1$ points, i.e.~the same \ $(C,\cO_C)$-points, \
where \ $C$ \ is \ a \ perfectoid \ field \ containing $\F$ \ (recall
that $Z_{\LT}(C,C^+)\buildrel\sim\over\rightarrow Z_{\LT}(C,\cO_C)$
for any open bounded valuation subring $C^+\subset C$). We use Newton
polygons and notation as in the proof of Lemma
\ref{lemm:stability}. If
$(F(t_0),\dots,F(t_{f-1}))\in(\B^+(C)^{\varphi_q=p})^f$, the element
$F(t_i)$ has slopes $\{(q-1)v(t_i)q^n,\ \!n\!\in\Z\}$ by Corollary
\ref{coro:polygone_Newton}, where $v$ is the valuation of $C$ such
that $\vabs{\cdot}=q^{-v(\cdot)}$, and recall that
$(F(t_0),\dots,F(t_{f-1}))\in Z_{\LT}^{\gen}(C,\cO_C)$ if and only if
$F(t_0)\cdots F(t_{f-1})\in \B^+(C)^{\varphi_q=p^f}$ lies in
$Z_{\cO_K}^{\gen}(C,\cO_C)$. As the slopes of the Newton polygon of a
product $ab$ in $\B^+(C)$ is the union of the slopes of the Newton
polygons of $a$ and $b$ (see \cite[Prop.~1.6.20]{FF} for instance), we
see that $F(t_0)\cdots F(t_{f-1})\in Z_{\cO_K}^{\gen}(C,\cO_C)$ if and
only if there exists $c >0$ such that
$\bigcup_i \{(q-1)v(t_i)q^n,\ n\in\Z\} = \{cp^n,\ n\in \Z\}$ (see the proof of
Lemma \ref{lemm:stability}). Equivalently
$F(t_0)\cdots F(t_{f-1})\in Z_{\cO_K}^{\gen}(C,\cO_C)$ if and only if
there exist $c>0$, $\sigma\in\mathfrak{S}_f$ and
$m_0,\dots,m_{f-1}\in\Z$ such that $v(t_i)=cp^{\sigma^{-1}(i)+fm_i}$
for $0\leq i\leq f-1$ if and only if there exist
$\sigma\in\mathfrak{S}_f$ and $m_0,\dots,m_{f-1}\in\Z$ such that
$v(t_i)=p^{(\sigma^{-1}(i)+fm_i)-(\sigma^{-1}(0)+fm_0)}v(t_0)$ for
$0\leq i\leq f-1$ if and only if
$F(t_0)\cdots F(t_{f-1})\in U_{\sigma(\underline{n}_0)+f\underline
  m}$. This proves our statement on rank $1$ points.

  For a point $x$ of the analytic adic space $Z_{\LT}$ define $\tilde{x}\in Z_{\LT}$ as its maximal generization, then the corresponding valuation $\abs_{\tilde{x}}$ is of rank $1$, i.e.\ real valued (see for instance \cite[Lemma 1.1.10]{Huberbook} or \cite[Cor.~II.2.4.8]{Morel}). Thus one can define continuous maps as in \cite[proof of Prop.~4.2.6]{SWBerkeley}:
 \[\kappa_{i,j} : Z_{\LT}\rightarrow ]0,+\infty[,\ x\mapsto \kappa_{i,j}(x)\defeq \tfrac{\log(\vabs{T_{K,i}}_{\tilde{x}})}{\log(\vabs{T_{K,j}}_{\tilde{x}})}.\]
 For $\underline{n}\in\Z^f/\Z$, define the closed subset of $Z_{\LT}$
  \[
    V_{\underline{n}}\defeq\kappa^{-1}(p^{n_0-n_1},\dots,p^{n_0-n_{f-1}}),\]
  where $\kappa=(\kappa_{0,1},\dots,\kappa_{0,f-1})$. For $x\in U_{\underline{n}}$, we still have $\tilde{x}\in U_{\underline{n}}$ by \cite[Lemma 1.1.10(v)]{Huberbook} applied to $X\defeq U_{\underline{n}}\hookrightarrow Y\defeq   Z_{\LT}$, hence we have an inclusion of topological spaces $U_{\underline{n}}\subset V_{\underline{n}}$. Let us prove that the open subspace $Z_{\LT}^{\gen}$ of $Z_{\LT}$ is contained in $V\defeq \bigcup_{\sigma\in\mathfrak{S}_f}\bigcup_{\underline{m}\in\Z^f/\Z}V_{\sigma(\underline{n}_0)+f\underline{m}}$. Let $x\in Z_{\LT}^{\gen}$ of rank $1$, then $x\in U_{\underline{n}}\subset V_{\underline{n}}$ for some $\underline{n}$ of the form $\underline{n}_0)+f\underline m$ by the second paragraph. As $V_{\underline{n}}$ is closed, we have $\overline{\{x\}}\subset V_{\underline{n}}$. Now let $x \in Z_{\LT}^{\gen}$ be any point and $\tilde{x}$ its maximal generization (which is in $Z_{\LT}^{\gen}$ by \cite[Lemma 1.1.10(v)]{Huberbook} applied to $Z_{\LT}^{\gen}\hookrightarrow Z_{\LT}$), then $\tilde{x}$ is of rank $1$ and $x\in \overline{\{\tilde x\}}$, which implies $x\in V_{\underline{n}}$ for some $\underline{n}$, i.e.\ $Z_{\LT}^{\gen}\subset V$. As $Z_{\LT}^{\gen}$ is open in $Z_{\LT}$, we have $Z_{\LT}^{\gen}\subset \mathring{V}\subset V$, where $\mathring{V}$ is the interior of the topological space $V$ in $Z_{\LT}$ ($\mathring{V}$ is then open in the perfectoid space $Z_{\LT}$, hence itself a perfectoid space). Let $x\in \mathring{V}$, then $x\in V_{\underline{n}}$ for some ${\underline{n}}$. But $V_{\underline{n}}$ is open in $V$ as $V$ is the inverse image by $\kappa$ of a discrete set and $V_{\underline{n}}$ is the inverse image of a single, hence open, element in this discrete set. Hence there exists an open subset $U$ of $Z_{\LT}$ such that $V_{\underline{n}}=U\cap V$. As $x\in U\cap \mathring{V}$ which is open in $Z_{\LT}$, we deduce $x\in \mathring{V_{\underline{n}}}$ which proves that $\mathring{V} = \bigcup_{\sigma\in\mathfrak{S}_f}\bigcup_{\underline{m}\in\Z^f/\Z}\mathring{V}_{\sigma(\underline{n}_0)+f\underline{m}}$. Thus we finally have $Z_{\LT}^{\gen}\subset \bigcup_{\sigma\in\mathfrak{S}_f}\bigcup_{\underline{m}\in\Z^f/\Z}\mathring{V}_{\sigma(\underline{n}_0)+f\underline{m}}$ which implies (using the first sentence of the proof)
\begin{equation}\label{Zgendecomp}
Z_{\LT}^{\gen} = \coprod_{\sigma\in\mathfrak{S}_f}\coprod_{\underline{m}\in(\Z^f/\Z)_0}(Z_{\LT}^{\gen}\cap \mathring{V}_{\sigma(\underline{n}_0)+f\underline{m}})
\end{equation}
 as open (perfectoid) subspaces of $Z_{\LT}$.
 
Now we go into group actions. It is not hard to check that $\Delta\rtimes\mathfrak{S}_f$ stabilizes $V$ (inside $Z_{\LT}$), more precisely $\sigma\in \mathfrak{S}_f$ sends $V_{\underline{n}}$ to $V_{\sigma(\underline n)}$, $(p^{d_0},\dots,p^{d_{f-1}})\in \Delta\cap (p^{\Z})^f$ sends $V_{\underline{n}}$ to $V_{\underline{n} + f(d_0,\dots,d_{f-1})}$ and $\Delta_1$ preserves each $V_{\underline{n}}$ (indeed, using that $f(\widetilde x)=\widetilde{f(x)}$ for any $x\in Z_{\LT}$ and any endomorphism $f$ of $Z_{\LT}$ by \cite[Lemma 1.1.10(iv)\&(v)]{Huberbook}, it is enough to check this for rank $1$ points, i.e.\ $(C,\cO_C)$-points for perfectoid fields $C$ containing $\F$, which is an easy exercise left to the reader). Then by continuity of the action of $\Delta\rtimes\mathfrak{S}_f$ the same holds for the interiors $\mathring V_{\underline{n}}$, and thus also for $Z_{\LT}^{\gen}\cap\mathring{V}_{\underline{n}}$ by Lemma \ref{lemm:stability}(ii). In particular, the group $(\Delta\cap (p^{\Z})^f)\rtimes\mathfrak{S}_f$ permutes transitively the perfectoid spaces $Z_{\LT}^{\gen}\cap\mathring{V}_{\underline{n}}$ for $\underline{n}\in\Z^f/\Z$ of the form $\sigma(\underline{n}_0)+f\underline{m}$ as in (\ref{Zgendecomp}), and the group $\Delta_1$ preserves each $Z_{\LT}^{\gen}\cap\mathring{V}_{\underline{n}}$. Thus the associated sheaf $\underline{\Delta_1}$ acts on (the sheaf corresponding to) $Z_{\LT}^{\gen}\cap\mathring{V}_{\underline{n}}$, and one easily checks that the group $\Delta_1$ moreover acts freely on the $(C,\cO_C)$-points of $Z_{\LT}^{\gen}\cap\mathring{V}_{\underline{n}}$. By the proof of \cite[Prop.~4.3.2]{WeinsteinGalois}, $Z_{\LT}^{\gen}\cap\mathring{V}_{\underline{n}}$ is a pro-\'etale $\Delta_1$-torsor over $\Delta_1 \backslash (Z_{\LT}^{\gen}\cap\mathring{V}_{\underline{n}})$, seen as a pro-\'etale sheaf on $\Perf_{\F}$. Since $\Delta_1$ is a normal subgroup in $\Delta \rtimes\mathfrak{S}_f$, we deduce with (\ref{Zgendecomp}) that $Z_{\LT}^{\gen}$ is a pro-\'etale $\Delta \rtimes\mathfrak{S}_f$-torsor over
\begin{eqnarray*}
\Delta \rtimes\mathfrak{S}_f \backslash Z_{\LT}^{\gen}&\simeq &((\Delta\cap (p^{\Z})^f)\rtimes\mathfrak{S}_f) \backslash (\Delta_1 \backslash Z_{\LT}^{\gen})\\
&\simeq &((\Delta\cap (p^{\Z})^f)\rtimes\mathfrak{S}_f)\backslash \big(\coprod_{\sigma, \underline{m}} \Delta_1 \backslash (Z_{\LT}^{\gen}\cap \mathring{V}_{\sigma(\underline{n}_0)+f\underline{m}})\big)\\
&\simeq & \Delta_1 \backslash (Z_{\LT}^{\gen}\cap\mathring{V}_{\underline{n}})
\end{eqnarray*}
for each $\underline n=\sigma(\underline{n}_0)+f\underline{m}$ as in (\ref{Zgendecomp}). Now, it follows from (\ref{fargues}) (and Lemma \ref{lemm:stability}) that we have an isomorphism $\Delta\rtimes\mathfrak{S}_f\backslash Z_{\LT}^{\gen} \xrightarrow{\sim} Z_{\cO_K}^{\gen}$ of pro-\'etale sheaves, hence $\Delta_1 \backslash (Z_{\LT}^{\gen}\cap\mathring{V}_{\underline{n}})\simeq Z_{\cO_K}^{\gen}$ for each $\underline n$ as above.

We now finish the proof. As $Z_{\cO_K}^{\gen}$ is affinoid perfectoid by Lemma \ref{lemm:A_infty_expl}(ii), each $Z_{\LT}^{\gen}\cap\mathring{V}_{\underline{n}}$ is affinoid perfectoid by \cite[Prop.~9.3.1]{SWBerkeley}, in particular is a quasi-compact open subset of $Z_{\LT}$. The quasi-compact open subspaces $Z_{\LT}^{\gen}\cap\mathring{V}_{\underline{n}}$ and $U_{\underline{n}}$ of $Z_{\LT}\subset \Spa(\F\bbra{K})\setminus V(T_{K,0})$ have the same points of rank $1$ by the second paragraph of this proof, and we can then argue in a similar way as for the proof of Lemma \ref{lemm:A_infty_expl}(iv), applying the results in \cite[(1.1.11)]{Huberbook} (or \cite[Thm.~2.21]{Scholzeperfectoid}) to the affinoid rigid analytic space over $\F\ppar{T_{K,0}}$ associated to $\Spa(\F\bbra{T_{K,0},\dots,T_{K,{f-1}}})\setminus V(T_{K,0})$ (recalling that $\Spa(\F\bbra{T_{K,0},\dots,T_{K,{f-1}}})\rightarrow \Spa(\F\bbra{T_{K,0}^{1/p^\infty},\dots,T_{K,{f-1}}^{1/p^\infty}})=Z_{\LT}$ is a homeomorphism). In particular, we obtain $U_{\underline{n}}=Z_{\LT}^{\gen}\cap\mathring{V}_{\underline{n}}$ for all $\underline{n}\in\Z^f/\Z$ of the form $\sigma(\underline{n}_0)+f\underline{m}$, which finishes the proof.
\end{proof}

As a consequence of the above proof and of \cite[Lemma 10.13]{ScholzeEtcohDia}, we also obtain:

\begin{cor1}\label{prop:DeltaS_torsor}
  The map $m : Z_{\LT}^{\gen}\rightarrow Z_{\cO_K}^{\gen}$ is a pro-\'etale $\Delta\rtimes\mathfrak{S}_f$-torsor, in particular is a pro-\'etale cover.
\end{cor1}

\begin{rem1}
  Note that $Z_{\LT}^{\gen}$ is not affinoid (contrary to $Z_{\cO_K}^{\gen}$) as it is not quasi-compact.
\end{rem1}

Let us denote by $A'_\infty\defeq\cO(U_{\underline{n}_0})$ the ring of global
sections on $U_{\underline{n}_0}$. The following result on $A'_\infty$ can be proved exactly as Lemma \ref{lemm:A_infty_expl}, and we leave the details to the reader.

\begin{lem1}\label{lemm:norm}
The following statements hold.
  \begin{enumerate}\item\label{lemm:norm(i)}
    The ring $A'_\infty$ is the perfectoid $\F$-algebra \[\F\ppar{T_{K,0}^{1/p^\infty}}\scalar*{\left(\frac{T_{K,i}}{T_{K,0}^{p^i}}\right)^{\pm
        1/p^\infty}\!\!,\ 1\leq i\leq f-1}.\]
  \item\label{lemm:norm(ii)}
    We have $U_{\underline{n}_0}=\Spa(A'_\infty,(A'_\infty)^\circ)$.
  \item\label{lemm:norm(iibis)} There exists a multiplicative norm
    $\abs$ on $A'_\infty$ such that $\vabs{T_{K,0}}=p^{-1}$ inducing
    the topology of $A'_\infty$.
  \item\label{lemm:norm(iii)} Any quasi-compact open subset of $Z_{\LT}$ whose points of rank $1$ are exactly the points of $U_{\underline{n}_0}$ of rank $1$ is necessarily $U_{\underline{n}_0}$ itself.
    \end{enumerate}
\end{lem1}

\subsection{Equivariant vector bundles on \texorpdfstring{$Z_{\cO_K}^{\gen}$}{Z\_{\{O\_K\}}\^{}gen} and \texorpdfstring{$Z_{\LT}^{\gen}$}{Z\_{LT}\^{}gen}}\label{equibundle}

We show that continuous $(\!K^\times\!)^f\!\rtimes\mathfrak{S}_f$-equivariant vector bundles on $Z_{\LT}^{\gen}$ and \'etale $(\varphi,\cO_K^\times)$-modules over $A_\infty$ are the same thing.

Recall first that if $X$ is an adic space with a left action of a group $H$,
an $H$-equivariant vector bundle on $X$ is a locally finite free ${\mathcal O}_X$-module $\mathcal{V}$ with a collection of $\cO_X$-linear isomorphisms $(c_h :h^*\mathcal{V}\xrightarrow{\sim}
\mathcal{V})_{h\in H}$ satisfying
the relation $c_{h_2h_1}=c_{h_1}\circ h_1^*(c_{h_2})$ for all
$h_1,h_2\in H$. This induces a {\it right} action of $H$ on
$\Gamma(X,\mathcal{V})$ given by
\[ c_h^* :
  \Gamma(X,\mathcal{V})=\Gamma(X,h^*\mathcal{V})\buildrel\sim\over
  \longrightarrow \Gamma(X,\mathcal{V}).\]
Now assume that $X$ is perfectoid space (the only case we will use) and that $H$ is a locally profinite topological group acting continuously on $X$. Let $\mathcal{V}$ be a vector bundle on $X$, for an open affinoid perfectoid subspace $U=\Spa(A,A^+)\subseteq X$, the finite projective $A$-module $\mathcal{V}(U)$ is endowed with the Banach topology given by the quotient topology of any surjection of $A$-modules $A^{\oplus d}\twoheadrightarrow \mathcal{V}(U)$. If $U\subseteq X$ is any open subspace, we endow $\mathcal{V}(U)\simeq \displaystyle{\lim_{\substack{\longleftarrow \\ U'\subseteq U}}\mathcal{V}(U')}$ with the projective limit topology, where $U'$ ranges over open affinoid subspaces of $U$, and we define $H_U\defeq \{h\in H,\ h(U)=U\}$, which is a closed subgroup of $H$ by continuity of the action of $H$ on $X$. We then define a {\it continuous} $H$-equivariant vector bundle on $X$ as an $H$-equivariant vector bundle $\mathcal{V}$ on $X$ such that for any open subspace $U\subseteq X$ the natural map $H_U\times \mathcal{V}(U)\rightarrow \mathcal{V}(U)$, $(h,s)\mapsto c_h^*(s)$ is continuous (for the product topology on the left).
  
By Lemma \ref{lemm:A_infty_expl}(i),(ii) and \cite[Thm.~2.7.7]{KL1} or \cite[Thm.~3.5.8]{KL2}, the functor of global sections induces an
equivalence of categories from the category of vector bundles on
$Z_{\cO_K}^{\gen}$ to the category of finite projective
$A_\infty$-modules. This equivalence is exact, rank preserving and compatible
with tensor products. As a finite projective $A_\infty$-module is in fact always free
(see \cite[Thm.~2.19]{DH}) and as the action of $K^\times$ on $Z_{\cO_K}^{\gen}$ is continuous, we see that the functor of global sections induces a {rank-preserving $\otimes$-}equivalence of categories from the category of continuous $K^\times$-equivariant vector bundles on
$Z_{\cO_K}^{\gen}$ to the category of \'etale $(\varphi,\cO_K^\times)$-modules
over $A_\infty$, where $\varphi$ on $A_\infty$ is given by (\ref{phionainfini}).

As $Z_{\LT}^{\gen}$ is perfectoid and as the fibered category of vector bundles on $\Perf_{\F}$ is a $v$-stack by \cite[Lemma 17.1.8]{SWBerkeley}, we easily deduce from Corollary \ref{prop:DeltaS_torsor} an equivalence of categories between the category of (continuous) $\Delta\rtimes\mathfrak{S}_f$-equivariant vector bundles on $Z_{\LT}^{\gen}$ and the category of vector bundles on $Z_{\cO_K}^{\gen}$ (the continuity condition is then automatic in that case, as $\Delta\rtimes\mathfrak{S}_f$ acts continuously on $Z_{\LT}^{\gen}$), hence also between the category of continuous $(K^\times)^f\rtimes\mathfrak{S}_f$-equivariant vector bundles on
$Z_{\LT}^{\gen}$ and the category of continuous $K^\times$-equivariant vector
bundles on $Z_{\cO_K}^{\gen}$. In both cases this equivalence is given by the two
functors
$\mathcal{V}\mapsto (m_*\mathcal{V})^{\Delta\rtimes\mathfrak{S}_f}$
and $\mathcal{W}\mapsto m^*\mathcal{W}$, where $m : Z_{\LT}^{\gen}\rightarrow Z_{\cO_K}^{\gen}$. If $\mathcal{V}$ is
$(K^\times)^f\rtimes \mathfrak{S}_f$-equivariant, the
$K^\times$-equivariant structure on
$(m_*\mathcal{V})^{\Delta\rtimes\mathfrak{S}_f}$ can be made explicit as follows. For $a\in K^\times$ and any $i\in \{0,\dots,f-1\}$ we have an isomorphism using the notation in (\ref{eq:ji})
\[
  a^*m_*\mathcal{V}\simeq(a^{-1})_*m_*\mathcal{V}\simeq(a^{-1}m)_*\mathcal{V}\simeq
  (m j_i(a)^{-1})_*\mathcal{V}\simeq
  m_*(j_i(a)^{-1})_*\mathcal{V}\simeq m_*j_i(a)^*\mathcal{V}\]
(where the first isomorphism is $\id\in \Hom(m_*\mathcal{V},m_*\mathcal{V})=\Hom((a^{-1})^*a^*m_*\mathcal{V},m_*\mathcal{V}) \simeq \Hom(a^*m_*\mathcal{V},(a^{-1})_*m_*\mathcal{V})$, the third comes from (\ref{eq:ji}) and the last is analogous to the first). We then obtain an isomorphism of sheaves for $a\in K^\times$ and any $i\in \{0,\dots,f-1\}$:
\[m_*(c_{j_i(a)}) : a^*m_*\mathcal{V}\simeq
  m_*j_i(a)^*\mathcal{V}\buildrel\sim\over\rightarrow m_*\mathcal{V}  \]
which preserves the subsheaf $(m_*\mathcal{V})^{\Delta\rtimes\mathfrak{S}_f}$ (as
$\Delta\rtimes\mathfrak{S}_f$ is a normal subgroup of
$(K^\times)^f\rtimes\mathfrak{S}_f$) and induces an isomorphism $m_*(c_{j_i(a)}):a^*(m_*\mathcal{V})^{\Delta\rtimes\mathfrak{S}_f}\buildrel\sim\over\rightarrow (m_*\mathcal{V})^{\Delta\rtimes\mathfrak{S}_f}$ which does not depend on $i$.

We deduce from Proposition \ref{prop:Delta1_torsor} that we have an
isomorphism of perfectoid $\F$-algebras
$A_\infty\xrightarrow{\sim}(A'_\infty)^{\Delta_1}$, and as above using
\cite[Lemma 17.1.8]{SWBerkeley} that there is also an equivalence of
categories between the category of $\Delta_1$-equivariant vector
bundles on $U_{\underline{n}_0}$ and the category of vector bundles
on $Z_{\cO_K}^{\gen}$. Using again \cite[Thm.~2.7.7]{KL1} (or \cite[Thm.~3.5.8]{KL2}) and
\cite[Thm.~2.19]{DH}, we deduce:

\begin{thm1}\label{descentdelta1}
The functor
$D_{A_\infty}\mapsto A'_\infty\otimes_{A_\infty}D_{A_\infty}$ induces
an exact {rank-preserving $\otimes$-}equivalence of categories from the category of finite free
$A_\infty$-modules to the category of finite free $A'_\infty$-modules
with a semi-linear action of $\Delta_1$. A quasi-inverse is
given by $D_{A'_\infty}\mapsto D_{A'_\infty}^{\Delta_1}$.
\end{thm1}

Let $\delta \in\mathfrak{S}_f$ be the cyclic permutation $i\mapsto i+1$ (with $f-1\mapsto 0$). If
$\sigma\in\mathfrak{S}_f$, let
$p_\sigma\defeq (1,\dots,p,\dots,1)\in (K^\times)^f$ with $p$ at the
$\sigma(0)$-th entry. From the discussion before Proposition \ref{prop:Delta1_torsor} we get
\[ (p_\sigma\circ \sigma) (U_{\underline{n}_0})=p_\sigma(U_{\sigma(\underline{n}_0)})=U_{\sigma \delta
    (\underline{n}_0)}.\]
In particular, $p_{\delta^{-1}}\circ \delta^{-1} : U_{\underline{n}_0}\xrightarrow{\sim} U_{\underline{n}_0}$ and we define an $\F$-linear continuous automorphism $\varphi$ of
$A'_\infty=\cO(U_{\underline{n}_0})$ by
\begin{equation}\label{frob'} 
\varphi\defeq (p_{\delta^{-1}}\circ \delta^{-1})^*=(\delta^{-1})^*\circ p_{\delta^{-1}}^*.
\end{equation}
Using (\ref{actionLT}) and since $\delta^{-1}(0)=f-1$ this automorphism
is easily checked to satisfy
\begin{equation}\label{phionT}
\varphi(T_{K,i})=T_{K,i+1}\ {\rm for\ }i\neq f-1\ {\rm and}\ \varphi(T_{K,f-1})=T_{K,0}^q.
\end{equation}
In particular, $\varphi^f$ on $A'_\infty$ is $\F$-linear and such that $\varphi^f(T_{K,i})=T_{K,i}^q$ for all $i$. Moreover if
$\underline{a}\in(\cO_K^\times)^f$, we have $\varphi\circ\underline{a}=\delta(\underline{a})\circ\varphi$, where
$\delta(\underline{a})= (a_{i-1})_{0\leq i\leq f-1}$ (with $a_{-1}=a_{f-1}$), in particular $\varphi^f$ commutes with $(\cO_K^\times)^f$. As $m : Z_{\LT}^{\gen}\rightarrow Z_{\cO_K}^{\gen}$ is $K^\times$-equivariant and $\mathfrak{S}_f$-equivariant, the action of $K^\times$ on $Z_{\LT}^{\gen}$ being through $j_i$ for any $0\leq i\leq f-1$ and the action of $\mathfrak{S}_f$ on $Z_{\cO_K}^{\gen}$ being trivial, the isomorphism $A_\infty\xrightarrow{\sim}(A'_\infty)^{\Delta_1}$
commutes with the actions of $\varphi$ and $\cO_K^\times$ on both sides (see (\ref{phionainfini}) for $\varphi$ on $A_\infty$).

The following result sums up the previous discussion and gives a more explicit way to compute
the $(\varphi,\cO_K^\times)$-module over $A_\infty$ associated to a
continuous $(K^\times)^f\rtimes\mathfrak{S}_f$-equivariant vector bundle on
$Z_{\LT}^{\gen}$.

\begin{cor1}\label{coro:concrete}
  There is an equivalence of categories between the category of continuous $(\!K^\times\!)^f\!\rtimes\mathfrak{S}_f$-equivariant vector bundles on $Z_{\LT}^{\gen}$ and the category of \'etale $(\varphi,\cO_K^\times)$-modules over $A_\infty$. If $\mathcal{V}$ is a
  continuous $(K^\times)^f\rtimes\mathfrak{S}_f$-equivariant vector bundle on
  $Z_{\LT}^{\gen}$, its associated $A_\infty$-module is
  $\Gamma(Z_{\cO_K}^{\gen},(m_*\mathcal{V})^{\Delta \rtimes\mathfrak{S}_f})$ which is
  isomorphic to
  $\Gamma(U_{\underline{n}_0},\mathcal{V}\vert_{U_{\underline{n}_0}})^{\Delta_1}$. The action of
  $a\in\cO_K^\times$ on
  $\Gamma(Z_{\cO_K}^{\gen},(m_*\mathcal{V})^{\Delta \rtimes\mathfrak{S}_f})$ is
  induced by the action of $(a,1\dots,1)=j_0(a)$ on
  $\Gamma(U_{\underline{n}_0},\mathcal{V}\vert_{U_{\underline{n}_0}})$ and the action of
  $\varphi$ on $\Gamma(\!Z_{\cO_K}^{\gen},(m_*\mathcal{V})^{\Delta \rtimes\mathfrak{S}_f}\!)$ is induced by
  \begin{multline*}
  (\delta^{-1})^*\circ p_{\delta^{-1}}^*:\Gamma(U_{\underline{n}_0},\mathcal{V}\vert_{U_{\underline{n}_0}}) = \Gamma(U_{\delta^{-1}(\underline{n}_0)},(p_{\delta^{-1}}^*\mathcal{V})\vert_{U_{\delta^{-1}(\underline{n}_0)}})\\
  \simeq \Gamma(U_{\underline{n}_0},((p_{\delta^{-1}}\circ \delta^{-1})^*\mathcal{V})\vert_{U_{\underline{n}_0}})
  \buildrel\sim\over\longrightarrow \Gamma(U_{\underline{n}_0},\mathcal{V}\vert_{U_{\underline{n}_0}}).
  \end{multline*}
\end{cor1}

\subsection{The \texorpdfstring{$(\varphi_q,\cO_K^\times)$}{(phi\_q,O\_K\^{}x)}-module over \texorpdfstring{$A$}{A} of an arbitrary Galois representation}\label{arbitrary}

To an arbitrary $\rhobar$ we functorially associate an \'etale $(\varphi_q,\cO_K^\times)$-module $D_{A}^{(i)}(\rhobar)$ over $A$ for $i\in \{0,\dots,f-1\}$.

Let $\rhobar$ be a continuous representation of $\Gal(\overline{K}/K)$ on a finite-dimensional $\F$-vector space and $D_{K,\sigma_0}(\rhobar)$ its
Lubin--Tate $(\varphi_q,\cO_K^\times)$-module (see \S\ref{LT}). The \'etale $(\varphi_q,\cO_K^\times)$-module
$\F\ppar{T_{K,\sigma_0}^{1/p^\infty}}\otimes_{\F\ppar{T_{K,\sigma_0}}}D_{K,\sigma_0}(\rhobar)$
is the space of global sections of a continuous $K^\times$-equivariant vector
bundle $\mathcal{V}_{\rhobar}$ on
$\Spa(\F\ppar{T_{K,\sigma_0}^{1/p^\infty}},\F\bbra{T_{K,\sigma_0}^{1/p^\infty}})$, where $p$ acts by $\varphi^f\otimes\varphi_q$. For $i\in \{0,\dots, f-1\}$ we define
\[\mathcal{V}_{\rhobar}^{(i)}\defeq \cO_{Z_{\LT}}\otimes_{\F\ppar{T_{K,\sigma_0}^{1/p^\infty}},\iota_i}\mathcal{V}_{\rhobar}\simeq \cO_{Z_{\LT}}\otimes_{\F\ppar{T_{K,\sigma_0}},\iota_i}D_{K,\sigma_0}(\rhobar),\]
where $\iota_i$ denotes the $\F$-linear embedding $\F\ppar{T_{K,\sigma_0}^{1/p^\infty}}\hookrightarrow\cO_{Z_{\LT}}$ corresponding to $T_{K,\sigma_0}\mapsto T_{K,i}$. Each $\mathcal{V}_{\rhobar}^{(i)}$ is a $\Delta$-equivariant vector bundle on $Z_{\LT}$ with $(a_0,\dots,a_{f-1})\in \Delta\subset (K^\times)^f$ acting on $\F\ppar{T_{K,\sigma_0}^{1/p^\infty}}\otimes_{\F\ppar{T_{K,\sigma_0}}}D_{K,\sigma_0}(\rhobar)$ via $a_i$. In particular, $\mathcal{V}_{\rhobar}^{(i)}\vert_{U_{\underline{n}_0}}$ is a $\Delta_1$-equivariant vector
bundle on $U_{\underline{n}_0}$ and $\Gamma(U_{\underline{n}_0},\mathcal{V}_{\rhobar}^{(i)}\vert_{U_{\underline{n}_0}})=A'_\infty\otimes_{\F\ppar{T_{K,\sigma_0}},\iota_i}D_{K,\sigma_0}(\rhobar)$. We define for $i\in \{0,\dots,f-1\}$
\begin{equation}\label{dai}
D_{A_\infty}^{(i)}(\rhobar)\defeq \Gamma(U_{\underline{n}_0},\mathcal{V}_{\rhobar}^{(i)}\vert_{U_{\underline{n}_0}})^{\Delta_1}=
  (A'_\infty\otimes_{\F\ppar{T_{K,\sigma_0}},\iota_i}D_{K,\sigma_0}(\rhobar))^{\Delta_1}
\end{equation}
which is a finite free $A_\infty$-module of rank
$\dim_{\F}\rhobar$ by Theorem \ref{descentdelta1}.

The endomorphism
$\varphi^f\otimes\varphi_q$ on
$A'_\infty\otimes_{\F\ppar{T_{K,\sigma_0}},\iota_i}D_{K,\sigma_0}(\rhobar)$ (see below (\ref{phionT}) for $\varphi^f$ on $A'_\infty$) commutes with the action of $\Delta_1$ and
induces a $\varphi_q$-semi-linear automorphism of
$D_{A_\infty}^{(i)}(\rhobar)$, which is thus naturally a
$\varphi_q$-module (see below (\ref{phionainfini}) for $\varphi_q$ on $A_\infty$). The action of $\cO_K^\times$ on
$A'_\infty\otimes_{\F\ppar{T_{K,\sigma_0}},\iota_i}D_{K,\sigma_0}(\rhobar)$ defined by
$a(x\otimes v)\defeq j_i(a)(x)\otimes a(v)$ induces a continuous semi-linear
action of $\cO_K^\times$ on $D_{A_\infty}^{(i)}(\rhobar)$ (with respect to the action of $\cO_K^\times$ on $A_\infty$) which commutes with $\varphi_q$. In particular, $D_{A_\infty}^{(i)}(\rhobar)$ is naturally an \'etale $(\varphi_q,\cO_K^\times)$-module over $A_\infty$. Note that the functor
$\rhobar\mapsto D_{A_\infty}^{(i)}(\rhobar)$ from continuous representations of $\Gal(\overline{K}/K)$ on finite-dimensional $\F$-vector spaces to \'etale $(\varphi_q,\cO_K^\times)$-modules over $A_\infty$ is exact and $\F$-linear. We also have isomorphisms of functors for $0\leq i\leq f-1$:
\begin{equation}\label{phii}
\phi_i :
    D_{A_\infty}^{(i)}(-)  \buildrel\sim\over\longrightarrow
    D_{A_\infty}^{(i+1)}(-), \qquad \phi_i=
                            \begin{cases}
                            \varphi\otimes \Id & \text{if } i<
                            f-1 \\
                            \varphi\otimes\varphi_q & \text{if }i=f-1.
                          \end{cases}
\end{equation}

We now show that \'etale $\varphi_q$-modules over $A_\infty$, and hence \'etale $(\varphi_q,\cO_K^\times)$-modules over $A_\infty$, canonically descend to the ring $A$ of \S\ref{oK}. First we need an easy lemma.

\begin{lem1}\label{tatefora}
The ring $A$ of \S\ref{oK} can be identified with the ring of global sections of the structure sheaf $\cO$ on the
rational open subset of the adic space $\Spa(\F\bbra{\cO_K})$ defined by the
relations
\[ \vabs{Y_{\sigma_0}}=\cdots=\vabs{Y_{\sigma_{f-1}}}\neq0,\]
where the variables $Y_{\sigma_i}\in \F\bbra{\cO_K}$ are defined in (\ref{variableysigma}).
\end{lem1}
\begin{proof}
Recall that $A$ is by definition the completed localization $(\F\bbra{\cO_K}_{(Y_{\sigma_0}\cdots Y_{\sigma_{f-1}})})^\wedge=(\F\bbra{Y_{\sigma_0},\dots,Y_{\sigma_{f-1}}}_{(Y_{\sigma_0}\cdots Y_{\sigma_{f-1}})})^{\wedge}$, where the completion is for the $(Y_{\sigma_0},\dots ,Y_{\sigma_{f-1}})$-adic to\-po\-logy. Then using (for instance) \cite[Rk.~3.1.1.3(iii)]{BHHMS2} one easily checks that
\begin{eqnarray*}
A&\simeq &\F[(Y_{\sigma_1}/Y_{\sigma_0})^{\pm 1},\dots,(Y_{\sigma_{f-1}}/Y_{\sigma_0})^{\pm 1}]\bbra{Y_{\sigma_0}}[1/Y_{\sigma_0}]\\
&\simeq & \F\ppar{Y_{\sigma_0}}\langle (Y_{\sigma_1}/Y_{\sigma_0})^{\pm 1},\dots,(Y_{\sigma_{f-1}}/Y_{\sigma_0})^{\pm 1}\rangle,
\end{eqnarray*}
where $\langle \ \rangle$ means, as usual, the corresponding Tate algebra with respect to the non-archimedean local field $\F\ppar{Y_{\sigma_0}}$. This is exactly the Tate algebra of the statement.
\end{proof}

Note that the open subset of Lemma \ref{tatefora} is stable under the endomorphisms deduced from the actions of $p$ and
$\cO_K^\times$ on $\cO_K$ by multiplication, in particular the $\F$-linear endomorphism $\varphi$ on $A$ sending $Y_{\sigma_i}$ to $Y_{\sigma_{i-1}}^p$ (see \S\ref{oK}) is the one deduced from the action of $p$.

\begin{rem1}\label{freeness}
It follows from Lemma \ref{tatefora} and \cite[Satz~3, p.~131]{Lut} (we thank Ofer Gabber for pointing out this reference to us) that any projective $A$-module of finite type is actually free.
\end{rem1}

Let $X_0,\dots,X_{f-1}$ be as at the end of \S\ref{reminderp}, we have $\F\bbra{\cO_K}=\F\bbra{X_0,\dots,X_{f-1}}=\F\bbra{Y_{\sigma_0},\dots,Y_{\sigma_{f-1}}}$, and from the  {equalities in (\ref{frobdown}) we deduce that there is $\lambda\in \Fq^\times$ such that} for $i\in \{0,\dots,f-1\}$
\begin{equation}\label{xiyi}
X_i= {\sigma_0(\lambda)^{p^i}}Y_{\sigma_i}+({\rm degree}\geq 2\ {\rm in\ the\ variables}\ Y_{\sigma_j}).
\end{equation}
This easily implies an isomorphism of completed localized rings
\[\big(\F\bbra{X_0,\dots,X_{f-1}}_{(X_0\cdots X_{f-1})}\big)^{\wedge}\simeq \big(\F\bbra{Y_{\sigma_0},\dots,Y_{\sigma_{f-1}}}_{(Y_{\sigma_0}\cdots Y_{\sigma_{f-1}})}\big)^{\wedge}=A,\]
where the completion on the left-hand side is for the $(X_0,\dots, X_{f-1})$-adic topology. In other words we can use the variables $X_i$ defined in \S\ref{reminderp} instead of the variables $Y_{\sigma_{i}}$ to define the ring $A$. In particular, the perfectoid Tate algebra $A_\infty$ in Lemma \ref{lemm:A_infty_expl} is the completion of the perfection of $A$ and the action of $\varphi$ and $\cO_K^\times$ on $A_\infty$ are compatible with the corresponding actions on $A$.

We will use the following result:

\begin{prop1}\label{prop:Katz}
{For $R$ a perfect $\Fq$-algebra there is an equivalence of categories between the category of locally constant sheaves $L$ of finite-dimensional $\Fq$-vector spaces on $\Spec(R)_{\acute et}$ and the
  category of pairs $(M,\phi)$ where $M$ is a finite projective $R$-module and $\phi$ is an isomorphism
  $\varphi_q^*M\xrightarrow{\sim}M$ {\upshape(}where $\varphi_q(-)=(-)^q${\upshape)}. This equivalence is given by the two inverse functors $L\mapsto (\Gamma(L\otimes_{\F_q}\mathcal{O}_{\Spec(R)}),\mathrm{Id}\otimes\varphi_q)$ and $(M,\phi)\mapsto (S\mapsto (M\otimes_RS)^{\phi\otimes \varphi_q=1}$, $S$ \'etale $R$-algebra).}
\end{prop1}

\begin{proof}
   {This is (a trivial variant of) \cite[Prop.~3.2.7]{KL1}.}
\end{proof}

We let $\displaystyle{A^{1/p^\infty}=\varinjlim_{x\mapsto x^p}A}=\bigcup_{n\geq 0}\big(\F\bbra{X_0^{1/p^n},\dots,X_{f-1}^{1/p^n}}_{(X_0\cdots X_{f-1})}\big)$ be the perfection of the ring $A$.  {It is used in the next proof}.

\begin{thm1}\label{descent}
  The functor $D_A\mapsto A_\infty \otimes_A D_A$ induces an exact equivalence of categories
  from the category of \'etale $\varphi_q$-modules over $A$ to the
  category of \'etale $\varphi_q$-modules over $A_\infty$.
\end{thm1}

\begin{proof}
 {Note first that $D_A\mapsto A^{1/p^\infty} \otimes_A D_A$ induces an exact equivalence of categories
  from the category of \'etale $\varphi_q$-modules over $A$ to the
  category of \'etale $\varphi_q$-modules over $A^{1/p^\infty}$ (use that any \'etale $\varphi_q$-module over $A^{1/p^\infty}$ comes by extension of scalars from an \'etale $\varphi_q$-module over $A^{1/q^N}$ for some $N\gg 0$ and apply $\varphi_q^N$). Hence we can replace $A$ by $A^{1/p^\infty}$. It follows from
\cite[Thm.~7.4.8]{SWBerkeley} (more precisely the discussion following {\it loc.~cit.})~that there is an equivalence of categories between the
category of finite \'etale $A^{1/p^\infty}$-algebras and the category
of finite \'etale $A_\infty$-algebras. Hence, when $\F=\Fq$, the result follows from Proposition \ref{prop:Katz} applied to both $R=A^{1/p^\infty}$ and $R=A_\infty$.}
In general, let $A_q$ be the ring of \S\ref{oK}, i.e.~$A_q$ is $A$ but with $\Fq$ instead of $\F$,  {$A_q^{1/p^\infty}$ its perfection and $A_{q,\infty}$ the completion of $A_q^{1/p^\infty}$. Then one can see an \'etale $\varphi_q$-module over $A^{1/p^\infty}$ (resp.\ $A_\infty$) as an \'etale $\varphi_q$-module over $A_{q}^{1/p^\infty}$ (resp.\ $A_{q,\infty}$) together with the structure of an $\F$-vector space compatible with the action of $\Fq$ (seen in $\F$ via $\sigma_0$). We only prove essential surjectivity (full faithfulness being easy). Let $D_{A_\infty}$ be an \'etale $\varphi_q$-module over $A_\infty$. By the equivalence of categories for $\F=\Fq$, there is an \'etale $\varphi_q$-module $D_{A^{1/p^\infty}}$ over $A_q^{1/p^\infty}$, which is also an $\F\otimes_{\Fq}A_q^{1/p^\infty}=A^{1/p^\infty}$-module, such that
\[A_{q,\infty}\otimes_{A_q^{1/p^\infty}}D_{A^{1/p^\infty}}\simeq (\F\otimes_{\Fq}A_{q,\infty})\otimes_{\F\otimes_{\Fq}A_q^{1/p^\infty}}D_{A^{1/p^\infty}}=A_\infty\otimes_{A^{1/p^\infty}}D_{A^{1/p^\infty}}\buildrel\sim\over\rightarrow D_{A_\infty}.\]
 {We need to prove that $D_{A^{1/p^\infty}}$ is finite projective over $A^{1/p^\infty}$ (or equivalently free by faithfully flat descent with Remark \ref{freeness}). The following argument is due to Ch.~Du. Let $0\ra N'\ra N\ra N''\ra 0$ be a short exact sequence of $A^{1/p^\infty}$-modules. Since $D_{A^{1/p^\infty}}$ is free over $A_q^{1/p^\infty}$ there is a short exact sequence $0\ra \Hom_{A_q^{1/p^\infty}}(D_{A^{1/p^\infty}},N'')\ra \Hom_{A_q^{1/p^\infty}}(D_{A^{1/p^\infty}},N)\ra \Hom_{A_q^{1/p^\infty}}(D_{A^{1/p^\infty}},N')\ra 0$. Making $\F^\times$ act on $\Hom_{A_q^{1/p^\infty}}(D_{A^{1/p^\infty}},(-))$ by $(\lambda f)(x)\defeq \lambda f(\lambda^{-1}x)$ and taking $\F^\times$-invariants we deduce a short exact sequence
\[0\ra \Hom_{A^{1/p^\infty}}(D_{\!A^{1/p^\infty}},N'')\ra \Hom_{A^{1/p^\infty}}(D_{\!A^{1/p^\infty}},N)\ra \Hom_{A^{1/p^\infty}}(D_{\!A^{1/p^\infty}},N')\ra 0.\]
Hence $D_{A^{1/p^\infty}}$ is finite projective over $A^{1/p^\infty}$}.}
\end{proof}

\begin{rem1}
We thank Laurent Berger for a discussion around Theorem \ref{descent}, and Laurent Fargues for suggesting to use Proposition \ref{prop:Katz} in its proof. Note that one can characterize the subspace $\displaystyle{A^{1/p^\infty}}\otimes_A D_A$ of an \'etale $\varphi_q$-module $D_{A_\infty}$ over $A_\infty$ as the $A$-submodule of $D_{A_\infty}$ of elements $d\in D_{A_\infty}$ such that $\sum_{n\geq 0}A\varphi_q^n(d)$ is a finite type $A$-module.
\end{rem1}

\begin{cor1}\label{descentoK}
  The functor $D_A\mapsto A_\infty \otimes_A D_A$ induces a rank-preserving $\otimes$-equivalence \ between \ the \ category \ of \ \'etale
  \ $(\varphi_q,\cO_K^\times)$-modules\ {\upshape(}resp. \'etale
  \ $(\varphi,\cO_K^\times)$-modules{\upshape)} over $A$ and the
  category of \'etale $(\varphi_q,\cO_K^\times)$-modules
  {\upshape(}resp.~\'etale $(\varphi,\cO_K^\times)$-modules{\upshape)}
  over $A_\infty$.
\end{cor1}

\begin{proof}
Let $D_{A_\infty}$ be an \'etale $(\varphi_q,\cO_K^\times)$-module over $A_\infty$. Any $a\in \cO_K^\times$ gives an isomorphism of \'etale $\varphi_q$-modules $\id \otimes a:a^*D_{A_\infty}\buildrel\sim\over\rightarrow D_{A_\infty}$ which canonically descends to an isomorphism of \'etale $\varphi_q$-modules $a^*D_{A}\buildrel\sim\over\rightarrow D_{A}$ by Theorem \ref{descent}. Now let $D_{A_\infty}$ be an \'etale $(\varphi,\cO_K^\times)$-module over $A_\infty$, then replacing $\varphi$ by $\varphi_q\defeq \varphi^f$, it is also an \'etale $(\varphi_q,\cO_K^\times)$-module over $A_\infty$. Let $\varphi^*D_{A_\infty}\defeq A_\infty\otimes_{\varphi,A_\infty}D_{A_\infty}$, then $\id\otimes \varphi$ induces an isomorphism of \'etale $\varphi_q$-modules $\varphi^*D_{A_\infty}\buildrel\sim\over\rightarrow D_{A_\infty}$ which canonically descends to an isomorphism $\varphi^*D_{A}\buildrel\sim\over\rightarrow D_{A}$ by Theorem \ref{descent}, giving the endomorphism $\varphi$ on $D_A$. The action of $\cO_K^\times$ canonically descends too by the first case of the proof and commutes with $\varphi$ (using Theorem \ref{descent} again). The rest of the statement is easy and left to the reader.
\end{proof}

From (\ref{dai}), (\ref{phii}) and Corollary \ref{descentoK}, we deduce:

\begin{cor1}\label{descenti}
For $i\in \{0,\dots ,f-1\}$ there is a covariant exact $\F$-linear
functor $\rhobar\mapsto D_{A}^{(i)}(\rhobar)$ compatible with tensor
products from $\Rep_{\F}\gK$ to \'etale
$(\varphi_q,\cO_K^\times)$-modules over $A$ and an isomorphism
$A_\infty\otimes_AD_A^{(i)}(-)\buildrel\sim\over\longrightarrow
D_{A_\infty}^{(i)}(-)$ between
functors from $\Rep_{\F}\gK$ to the category of \'etale
$(\varphi_q,\cO_K^\times)$-modules over $A_\infty$. These functors are related by functorial $A$-linear isomorphisms $\phi_i:A\otimes_{\varphi,A}D_{A}^{(i)}(\rhobar)\buildrel\sim\over\longrightarrow D_{A}^{(i+1)}(\rhobar)$ which commute with $(\varphi_q,\cO_K^\times)$ and are such that $\phi_{f-1}\circ \phi_{f-2}\circ \cdots \circ \phi_0:A\otimes_{\varphi^f,A}D_{A}^{(0)}(\rhobar)\buildrel\sim\over\longrightarrow D_{A}^{(0)}(\rhobar)$ is $\id\otimes \varphi_q$.
 \end{cor1}
 
 \begin{rem1}
One can check that $D_{A}^{(0)}(\rhobar)\times D_{A}^{(f-1)}(\rhobar)\times D_{A}^{(f-2)}(\rhobar)\times \cdots \times D_{A}^{(1)}(\rhobar)$ can be given the structure of an \'etale $(\varphi,\cO_K^\times)$-module over $\F\otimes_{\Fp}A_q$ in the sense of \S\ref{oK}.
\end{rem1}

\subsection{The \texorpdfstring{$(\varphi,\cO_K^\times)$}{(phi,O\_K\^{}x)}-module over \texorpdfstring{$A$}{A} associated to a Galois representation}\label{otimesarbitrary}

To an arbitrary $\rhobar$ we associate an \'etale $(\varphi,\cO_K^\times)$-module $D_{A}^{\otimes}(\rhobar)$ (which will be particularly important when $\dim_{\F}\rhobar=2$).

Keep the notation of \S\ref{arbitrary} and let
$\mathcal{V}_{\rhobar}^{\boxtimes
  f}\defeq\bigotimes_{i=0}^{f-1}\pr_i^*\mathcal{V}_{\rhobar}$ be
the $f$-th ``exterior tensor product'' of $\mathcal{V}_{\rhobar}$ on
$Z_{\LT}\!=\!(\Spa(\F\ppar{T_{K,\sigma_0}^{1/p^\infty}},\F\bbra{T_{K,\sigma_0}^{1/p^\infty}}))^f,$
where
\[\pr_i:(\Spa(\F\ppar{T_{K,\sigma_0}^{1/p^\infty}},\F\bbra{T_{K,\sigma_0}^{1/p^\infty}}))^{\!f}\onto 
\Spa(\F\ppar{T_{K,\sigma_0}^{1/p^\infty}},\F\bbra{T_{K,\sigma_0}^{1/p^\infty}})\]
is the $i$-th projection (so $\pr_i^*\mathcal{V}_{\rhobar}$ is the sheaf $\mathcal{V}_{\rhobar}^{(i)}$ of \S\ref{arbitrary}). As
$\mathcal{V}_{\rhobar}$ is a continuous $K^\times$-equivariant vector bundle,
$\mathcal{V}_{\rhobar}^{\boxtimes f}$ is naturally a
continuous $(K^\times)^f$-equivariant vector bundle. We promote it to a
(continuous) $(K^\times)^f\rtimes\mathfrak{S}_f$-equivariant vector bundle using
the commutativity of the tensor product (where $\sigma\in\mathfrak{S}_f$):
\[ c_{\sigma}: \sigma^*\mathcal{V}_{\rhobar}^{\boxtimes
    f}=\sigma^*\Big(\bigotimes_{i=0}^{f-1}\pr_i^*\mathcal{V}_{\rhobar}\Big)
\simeq\bigotimes_{i=0}^{f-1}\sigma^*\pr_i^*\mathcal{V}_{\rhobar}\simeq\bigotimes_{i=0}^{f-1}\pr_{\sigma^{-1}(i)}^*\mathcal{V}_{\rhobar}\buildrel {\sigma\atop \sim}\over \longrightarrow \bigotimes_{i=0}^{f-1}\pr_{i}^*\mathcal{V}_{\rhobar}=\mathcal{V}_{\rhobar}^{\boxtimes f}.\]
We define
$D_{A_\infty}^{\otimes}(\rhobar)$ as the $A_\infty$-module with a
continuous semi-linear action of $K^\times$ obtained as the global
sections of the continuous $K^\times$-equivariant vector bundle on $Z_{\cO_K}^{\gen}$
corresponding to $\mathcal{V}^{\boxtimes f}|_{Z_{\LT}^{\gen}}$, more concretely (see \S\ref{equibundle}):
\[
  D_{A_\infty}^{\otimes}(\rhobar)\defeq\Gamma\big(Z_{\cO_K}^{\gen},(m_*(\mathcal{V}_{\rhobar}^{\boxtimes
    f}|_{Z_{\LT}^{\gen}}))^{\Delta\rtimes\mathfrak{S}_f}\big)=\Gamma(Z_{\LT}^{\gen},\mathcal{V}_{\rhobar}^{\boxtimes
    f})^{\Delta\rtimes\mathfrak{S}_f}.\]
This is an \'etale $(\varphi,\cO_K^\times)$-module over $A_\infty$ (recall $\varphi$ is bijective).

Using Corollary \ref{coro:concrete} and \S\ref{arbitrary}, we can give a more explicit description of $D_{A_\infty}^{\otimes}(\rhobar)$. Note that we have: 
\[
  D_{A_\infty}^{\otimes}(\rhobar)=\Gamma(U_{\underline{n}_0},\mathcal{V}_{\rhobar}^{\boxtimes
    f}\vert_{U_{\underline{n}_0}})^{\Delta_1}\]
and that the vector bundle $\mathcal{V}_{\rhobar}^{\boxtimes f}$ is isomorphic to the tensor
product
\[\mathcal{V}_{\rhobar}^{(0)}\otimes_{\cO_{Z_{\LT}}} \cdots\otimes_{\cO_{Z_{\LT}}}\mathcal{V}_{\rhobar}^{(f-1)}.\]
As the equivalence with vector bundles on $Z_{\cO_K}^{\gen}$, i.e.~finite free
$A_\infty$-modules, is compatible with tensor products (see \S\ref{equibundle}), we deduce an isomorphism of
$A_\infty$-modules
\[
  D_{A_\infty}^{\otimes}(\rhobar)\simeq(A'_\infty\otimes_{\F\ppar{T_{K,\sigma_0}},\iota_0}D_{K,\sigma_0}(\rhobar))^{\Delta_1}\otimes_{A_\infty}\cdots\otimes_{A_\infty}(A'_\infty\otimes_{\F\ppar{T_{K,\sigma_0}},\iota_{f-1}}D_{K,\sigma_0}(\rhobar))^{\Delta_1}.\]

\begin{lem1}\label{dotimes}
  There is a functorial isomorphism of \'etale $(\varphi,\cO_K^\times)$-modules over
  $A_\infty$
  \[
    D_{A_\infty}^{\otimes}(\rhobar)\buildrel\sim\over\longrightarrow\bigotimes_{i=0}^{f-1}D_{A_\infty}^{(i)}(\rhobar),\]
  where the automorphism $\varphi$ on the right-hand side is given by {\upshape(}see (\ref{phii}) for $\phi_i${\upshape)}
  \[ \varphi(v_0\otimes \cdots \otimes
    v_{f-1})=\phi_{f-1}(v_{f-1})\otimes\phi_0(v_0)\otimes\cdots\otimes\phi_{f-2}(v_{f-2})\]
    {\upshape(}and the action of $\cO_K^\times$ is as defined in \S\ref{arbitrary} on each factor $D_{A_\infty}^{(i)}(\rhobar)${\upshape)}.
\end{lem1}

\begin{proof}
  Recall that $\delta \in\mathfrak{S}_f$ sends $i$ to $i+1$. Let
  $\alpha_i :
  (\delta^{-1})^*\mathcal{V}_{\rhobar}^{(i-1)}\xrightarrow{\sim}\mathcal{V}_{\rhobar}^{(i)}$
  be the tautological isomorphism deduced from the identifications
  \[
    (\delta^{-1})^*\mathcal{V}_{\rhobar}^{(i-1)}=(\delta^{-1})^*\pr_{i-1}^*\mathcal{V}_{\rhobar}\simeq(\pr_{i-1}\circ
    \delta^{-1})^*\mathcal{V}_{\rhobar}=\pr_{i}^*\mathcal{V}_{\rhobar}=\mathcal{V}_{\rhobar}^{(i)}.\]
    Recall that $p_{\delta^{-1}}\in (K^\times)^f$ is defined in \S\ref{equibundle} and let
  $\beta_i :
  p_{\delta^{-1}}^*\mathcal{V}_{\rhobar}^{(i)}\xrightarrow{\sim}\mathcal{V}_{\rhobar}^{(i)}$
  be the isomorphism of sheaves on $Z_{\LT}$ defined by ($f\in \cO_{Z_{\LT}}$, $v\in \mathcal{V}_{\rhobar}$ and compare with (\ref{phionT})):
  \[ f\otimes v \longmapsto
    \begin{cases}
      f(p_{\delta^{-1}}(-))\otimes v & \text{if } i\neq f-1 \\
      f(p_{\delta^{-1}}(-))\otimes \varphi_q(v) & \text{if } i=f-1.
    \end{cases}\]
  We obtain isomorphisms of sheaves on $Z_{\LT}$ for $i\in \{0,\dots,f-1\}$
  \begin{equation*}\alpha_i\circ(\delta^{-1})^*(\beta_{i-1}) :
    \varphi^*\mathcal{V}_{\rhobar}^{(i-1)}\buildrel (\ref{frob'})\over
    \simeq((\delta^{-1})^*\circ p_{\delta^{-1}}^*)\mathcal{V}_{\rhobar}^{(i-1)}\xrightarrow{\sim}(\delta^{-1})^*\mathcal{V}_{\rhobar}^{(i-1)}\xrightarrow{\sim}\mathcal{V}_{\rhobar}^{(i)}.
    \end{equation*}
  The isomorphism $c_{p_{\delta^{-1}}\circ \delta^{-1}} : \varphi^*\mathcal{V}_{\rhobar}^{\boxtimes f}\xrightarrow{\sim}\mathcal{V}_{\rhobar}^{\boxtimes
    f}$ (with the notation as at the beginning of \S\ref{equibundle}) is easily checked to decompose as a tensor product
  \[ \bigotimes_{i=0}^{f-1}\big(\alpha_i\circ(\delta^{-1})^*(\beta_i)\big) :
    \varphi^*\mathcal{V}_{\rhobar}^{\boxtimes
      f}\simeq\bigotimes_{i=0}^{f-1}\varphi^*\mathcal{V}_{\rhobar}^{(i-1)}\buildrel\sim\over\longrightarrow\bigotimes_{i=0}^{f-1}\mathcal{V}_{\rhobar}^{(i)}.\]
Taking global sections on $U_{\underline{n}_0}$ and $\Delta_1$-invariants, we obtain the desired formula.
\end{proof}

From Lemma \ref{dotimes} and Corollary \ref{descentoK} we deduce $D_{A_\infty}^{\otimes}(\rhobar)\simeq A_\infty\otimes_AD_{A}^{\otimes}(\rhobar)$ for a unique \'etale $(\varphi,\cO_K^\times)$-module $D_{A}^{\otimes}(\rhobar)$ over $A$ such that
 \begin{equation}\label{daotimesbis}
 D_{A}^{\otimes}(\rhobar)\simeq \bigotimes_{i=0}^{f-1}D_{A}^{(i)}(\rhobar)
 \end{equation}
 with the same $\varphi$ and action of $\cO_K^\times$ on the right-hand side as in Lemma \ref{dotimes}.

\begin{rem1}
Note that, for $0\leq i < f-1$, the isomorphism $\phi_i$ in (\ref{phii}) is induced by the natural $A_\infty$-linear isomorphism $\varphi^*D_{A_\infty}^{(i)}(-)\simeq D_{A_\infty}^{(i+1)}(-)$, whereas $\phi_{f-1}$ coincides with the
  $A_\infty$-linear isomorphism
  \[
    \varphi^*D_{A_\infty}^{(f-1)}(-)\simeq \varphi^*((\varphi^{f-1})^*D_{A_\infty}^{(0)}(-))=\varphi_q^*D_{A_\infty}^{(0)}(-)\longrightarrow
    D_{A_\infty}^{(0)}(-)\] induced by the $\varphi_q$-semi-linear
  automorphism $\varphi_q$ of $D_{A_\infty}^{(0)}(-)$. Therefore the
  isomorphism class of the $(\varphi,\cO_K^\times)$-module
  $D_{A_\infty}^{\otimes}(\rhobar)$ (equivalently of $D_{A}^{\otimes}(\rhobar)$) is completely characterized by the
  isomorphism class of the $(\varphi_q,\cO_K^\times)$-module
  $D_{A_\infty}^{(0)}(\rhobar)$ (equivalently of $D_{A}^{(0)}(\rhobar)$).
\end{rem1}

\subsection{Relation to classical \texorpdfstring{$(\varphi,\Gamma)$}{(phi,Gamma)}-modules}
 
We show that the \'etale $(\varphi_q,\cO_K^\times)$-module $D_{A}^{(0)}(\rhobar)$ is related in a simple way to the (usual) \'etale $(\varphi_q,\Zp^\times)$-module $D_{\sigma_0}(\rhobar)$ of \S\ref{LT} and derive some consequences.
 
 As in \cite[\S 3.1.3]{BHHMS2}, let $\tr : A\twoheadrightarrow\F\ppar{T}$ be the ring surjection induced by the trace $\tr : \F\bbra{\cO_K}\rightarrow\F\bbra{\Zp}\simeq\F\bbra{T}$. Since the map $\tr$ commutes with $\varphi$ (hence $\varphi_q$) and the action of $\Zp^\times$, we deduce that $\F\ppar{T}\otimes_AD_{A}^{(0)}(\rhobar)$ is an \'etale $(\varphi_q,\Zp^\times)$-module.
 
\begin{prop1}\label{lift}
There is $d\in \{0,\dots,f-1\}$ such that we have a functorial isomorphism of $(\varphi_q,\Zp^\times)$-modules
\[\F\ppar{T}\otimes_{A}D_{A}^{(0)}(\rhobar)\simeq D_{\sigma_d}(\rhobar),\]
where $D_{\sigma_d}(\rhobar)$ is as above (\ref{phizp0}) choosing the embedding $\sigma_d:\Fq\hookrightarrow \F$ instead of $\sigma_0$.
\end{prop1}
\begin{proof}
The trace $\tr : \F\bbra{K}\rightarrow\F\bbra{\Qp}\simeq\F\bbra{T^{p^{-\infty}}}$ induces a ring surjection $\tr : A_\infty \!\twoheadrightarrow \F\ppar{T^{p^{-\infty}}}$ commuting (in an obvious way) with $\tr : A\twoheadrightarrow\F\ppar{T}$. Using Corollary \ref{descentoK} it is enough to prove $\F\ppar{T^{p^{-\infty}}}\otimes_{A_\infty}D_{A_\infty}^{(0)}(\rhobar)\simeq \F\ppar{T^{p^{-\infty}}}\otimes_{\F\ppar{T}}D_{\sigma_d}(\rhobar)$.

For any perfectoid $\F$-algebra $R$ we have a commutative diagram
\begin{equation}\label{allB+}
\begin{gathered}
\xymatrix{\B^+(R)^{\varphi_q=p}\ar@{^{(}->}[r]\ar[d] & (\B^+(R)^{\varphi_q=p})^f\ar^{m_R}@{^{}->>}[d] \\
\B^+(R)^{\varphi=p}\ar@{^{(}->}[r] & \B^+(R)^{\varphi_q=p^f}}
\end{gathered}
\end{equation}
where the top horizontal injection sends $x\in \B^+(R)^{\varphi_q=p}$ to $(x,\varphi(x),\dots,\varphi^{f-1}(x))\in (\B^+(R)^{\varphi_q=p})^f$, the left vertical map sends $x\in \B^+(R)^{\varphi_q=p}$ to $x\varphi(x)\cdots\varphi^{f-1}(x)\in \B^+(R)^{\varphi=p}$ and where the bottom horizontal injection is the canonical injection. Note that the left vertical map commutes with the action of $K$, where $K$ acts on $\B^+(R)^{\varphi=p}$ via ${\rm Norm}_{K/\Q_p}:K\rightarrow \Qp$. As at the beginning of \S\ref{mapm}, using Remark \ref{rem:BW_epointe} and \cite[Prop.~8.2.8(2)]{SWBerkeley}, we deduce from (\ref{allB+}) a corresponding commutative diagram of perfectoid spaces over $\F$:
 \begin{equation}\label{allrings0}
\begin{gathered}
\xymatrix{(\widetilde{G}_{\LT}\times_{\Spf(\F_q)}\Spf(\F)\setminus\set{0})^{\ad} \ar@{^{(}->}[r] \ar[d] & Z_{\LT} \ar[d]\\
Z_{\Z_p} \ar[r] & Z_{\cO_K} }
\end{gathered}
\end{equation}
where the top horizontal map is $r\mapsto (r,r^p,\dots,r^{p^{f-1}})$ on the coordinates and the right vertical map is the map $m$ in (\ref{mapmzlt}). From the discussion above, the map $Z_{\Z_p} \rightarrow Z_{\cO_K}$ commutes with the action of $K^\times$, where $K^\times$ acts on $Z_{\Z_p}$ via ${\rm Norm}_{K/\Q_p}$. Also, it follows from the end of \S\ref{reminderp} (see in particular (\ref{commuteR}), (\ref{11ff}), (\ref{traceK}) and (\ref{gff})) that the bottom horizontal map is induced by the morphism $\F\bbra{K}\rightarrow\F\bbra{\Q_p}$ deduced from the trace $\Tr_{K/\Q_p}:K\rightarrow \Qp$. Hence we deduce from (\ref{allrings0}) a commutative diagram of perfectoid rings over $\F$:
\begin{equation}\label{allrings}
\begin{gathered}
\xymatrix{A'_\infty\ar@{^{}->>}[r] & \F\ppar{T_{K,\sigma_0}^{p^{-\infty}}}\\
A_\infty \ar@{^{(}->}[u] \ar^{\!\!\!\!\!\!\!\!\!\!\tr}@{^{}->>}[r] & \F\ppar{T^{p^{-\infty}}} \ar@{^{(}->}[u]}
\end{gathered}
\end{equation}
where the top horizontal surjection sends $T_{K,i}^{p^{-n}}$ to $T_{K,\sigma_0}^{p^{i-n}}$ for $i\in \{0,\dots, f-1\}$. The right vertical injection commutes with $\cO_K^\times$ (acting on $\F\ppar{T^{p^{-\infty}}}$ via the norm $\oK^\times \twoheadrightarrow \Zp^\times$), hence we deduce from Theorem \ref{axsentate} (and (\ref{LTp})) that it induces an injection of perfectoid fields $\iota:\F\ppar{T^{p^{-\infty}}} \hookrightarrow \F\ppar{T_{K,\sigma_0}^{p^{-\infty}}}^{\Gal(K_\infty/K(\!\sqrt[p^\infty]{1}))}\simeq \F\ppar{T^{p^{-\infty}}}$.
Since $\iota$ commutes with the action of $\Zp^\times$, one easily checks that it must be an isomorphism, as any continuous $\F$-algebra homomorphism $\F\ppar{\Qp}\rightarrow \F\ppar{\Qp}$ commuting with the action of $\Zp^\times$ sends $[1]\in \F\ppar{\Qp}$ to $[\lambda]\in \F\ppar{\Qp}$ for some $\lambda\in \Qp^\times$ by \cite[Thm.\ 3.1]{MR4526253} (and continuity).

Now let $\rhobar$ be a continuous representation of $\Gal(\overline{K}/K)$ on a finite-dimensional $\F$-vector space, using the isomorphism $A'_\infty\otimes_{\F\ppar{T_{K,\sigma_0}},\iota_0}D_{K,\sigma_0}(\rhobar)\simeq A'_\infty \otimes_{A_\infty}D_{A_\infty}^{(0)}(\rhobar)$ from Theorem \ref{descentdelta1} we deduce from (\ref{allrings}):
\begin{eqnarray*}
\F\ppar{T_{K,\sigma_0}^{p^{-\infty}}}\otimes_{\F\ppar{T_{K,\sigma_0}}}D_{K,\sigma_0}(\rhobar)&\simeq &\F\ppar{T_{K,\sigma_0}^{p^{-\infty}}}\otimes_{A'_\infty}\big(A'_\infty\otimes_{\F\ppar{T_{K,\sigma_0}},\iota_0}D_{K,\sigma_0}(\rhobar)\big)\\
&\simeq &\F\ppar{T_{K,\sigma_0}^{p^{-\infty}}}\otimes_{A'_\infty}\big(A'_\infty \otimes_{A_\infty}D_{A_\infty}^{(0)}(\rhobar)\big)\\
&\simeq &\F\ppar{T_{K,\sigma_0}^{p^{-\infty}}}\otimes_{\iota,\F\ppar{T^{p^{-\infty}}}}\big(\F\ppar{T^{p^{-\infty}}}\otimes_{A_\infty}D_{A_\infty}^{(0)}(\rhobar)\big).
\end{eqnarray*}
By Proposition \ref{compare} we also have
\begin{multline*}
\F\ppar{T_{K,\sigma_0}^{p^{-\infty}}}\otimes_{\F\ppar{T_{K,\sigma_0}}}D_{K,\sigma_0}(\rhobar) \simeq \F\ppar{T_{K,\sigma_0}^{p^{-\infty}}}\otimes_{\F\ppar{T^{p^{-\infty}}}}\big(\F\ppar{T^{p^{-\infty}}}\otimes_{\F\ppar{T}}D_{\sigma_0}(\rhobar)\big)\\
\simeq  
\F\ppar{T_{K,\sigma_0}^{p^{-\infty}}}\otimes_{\iota, \F\ppar{T^{p^{-\infty}}}}\big(\F\ppar{T^{p^{-\infty}}}\otimes_{\iota^{-1}, \F\ppar{T^{p^{-\infty}}}}(\F\ppar{T^{p^{-\infty}}}\otimes_{\F\ppar{T}}D_{\sigma_0}(\rhobar))\big).
\end{multline*}
Since the action of $\Gal(K_\infty/K(\!\sqrt[p^\infty]{1}))\simeq \Gal(\F\ppar{T_{K,\sigma_0}^{p^{-\infty}}}/\F\ppar{T^{p^{-\infty}}})$ is trivial on both $\F\ppar{T^{p^{-\infty}}}\otimes_{A_\infty}D_{A_\infty}^{(0)}(\rhobar)$ and $\F\ppar{T^{p^{-\infty}}}\otimes_{\F\ppar{T}}D_{\sigma_0}(\rhobar)$, we deduce by Galois descent an isomorphism of $(\varphi_q,\Zp^\times)$-modules over $\F\ppar{T^{p^{-\infty}}}$
\[\F\ppar{T^{p^{-\infty}}}\otimes_{A_\infty}D_{A_\infty}^{(0)}(\rhobar)\!\simeq \F\ppar{T^{p^{-\infty}}}\otimes_{\iota^{-1}, \F\ppar{T^{p^{-\infty}}}}(\F\ppar{T^{p^{-\infty}}}\otimes_{\F\ppar{T}}D_{\sigma_0}(\rhobar)).\]
This easily gives the statement, using that all the above isomorphisms are functorial in $\rhobar$ (note that, if $\iota$ is given by $[1]\mapsto [\lambda]$ as above, then $d$ is the unique integer in $\{0,\dots, f-1\}$ congruent to $\val(\lambda)$ modulo $f$).
\end{proof}

\begin{rem1}\label{shift}
Using Theorem \ref{explicitK} below together with Lemma \ref{PS} and \cite[Prop.~3.5]{breuil-IL}, one can compute that $d=f-1$. We won't need this fact.
\end{rem1}
 
We can also consider the tensor product $\F\ppar{T}\otimes_AD_A^\otimes(\rhobar)$ for $\tr : A\twoheadrightarrow\F\ppar{T}$. It is obviously an \'etale $(\varphi,\Zp^\times)$-module.

\begin{cor1}\label{determine}
The $(\varphi,\Zp^\times)$-module $\F\ppar{T}\otimes_AD_A^\otimes(\rhobar)$ is the $(\varphi,\Zp^\times)$-module of the tensor induction $\ind_K^{\otimes\Qp}\!\rhobar$.
\end{cor1}
\begin{proof}
This easily follows from (\ref{daotimesbis}), Proposition \ref{lift}, Corollary \ref{descenti} and the ``tensor product version'' of \cite[Lemma 3.6]{breuil-IL} (which we leave to the reader).
\end{proof}

Proposition \ref{lift} also enables to prove the following full faithfulness statement.

\begin{cor1}\label{fully}
For $i\in \{0,\dots, f-1\}$ the functor $\rhobar \mapsto D_{A}^{(i)}(\rhobar)$ from continuous representations of $\Gal(\overline{K}/K)$ on finite-dimensional $\F$-vector spaces to \'etale $(\varphi_q,\cO_K^\times)$-modules over $A$ is exact and fully faithful.
\end{cor1}
\begin{proof}
By Corollary \ref{descenti} it is enough to prove the full faithfulness for $i=0$. We have morphisms:
{\scriptsize{\begin{equation*}\Hom_{\gK}(\rhobar, \rhobar')\longrightarrow \Hom_{(\varphi_q,\oK^\times)}(D_{A}^{(0)}(\rhobar), D_{A}^{(0)}(\rhobar'))\longrightarrow \Hom_{(\varphi_q,\Zp^\times)}(D_{\sigma_d}(\rhobar), D_{\sigma_d}(\rhobar'))
\end{equation*}}}
\!\!where we use Proposition \ref{lift} for the second. By the theory of $(\varphi_q,\Zp^\times)$-modules (see e.g.~\S\ref{LT}), we know that the composition of the two morphisms is bijective. Hence the first morphism is injective. It is enough to prove that the second morphism is also injective. Let $f:D_{A}^{(0)}(\rhobar) \rightarrow D_{A}^{(0)}(\rhobar')$ mapping to $0$, i.e.~$f(D_{A}^{(0)}(\rhobar))\subseteq {\mathfrak p}D_{A}^{(0)}(\rhobar')$, where ${\mathfrak p}\defeq \Ker(\tr : A\twoheadrightarrow\F\ppar{T})$ (a maximal ideal of the noetherian domain $A$). Using the fact that $D_{A}^{(0)}(\rhobar)$ is \'etale and that $f$ commutes with $\varphi_q$, we derive $f(D_{A}^{(0)}(\rhobar))\subseteq \varphi_q^n({\mathfrak p})D_{A}^{(0)}(\rhobar')$ for any $n\geq 0$. For those $n$ such that $x\mapsto x^{q^n}$ is $\F$-linear, the map $\varphi_q^n$ on $A$ is just $x\mapsto x^{q^n}$, hence $\varphi_q^n({\mathfrak p})\subseteq {\mathfrak p}^{q^n}$ for those $n$, and thus $f(D_{A}^{(0)}(\rhobar))\subseteq (\bigcap_{m\geq 0}{\mathfrak p}^m)D_{A}^{(0)}(\rhobar')=0$. This finishes the proof.
\end{proof}

\vspace{.1cm}

\begin{rem1}\ 
  \begin{enumerate}
  \item We do not expect the functor $\rhobar \mapsto D_{A}^{(i)}(\rhobar)$ to be essentially surjective (for any $i$). It is probably an interesting question to characterize its essential image.
  \item 
    It is {\it not true} that the functor $\rhobar \mapsto D_{A}^{\otimes}(\rhobar)$ is fully faithful, as in general the
    isomorphism class of the $(\varphi,\oK^\times)$-module $D_{A}^{\otimes}(\rhobar)$ does not determine the one of the Galois
    representation $\rhobar$. For instance, one can check by an explicit computation using Theorem \ref{explicitK} below that,
    if $f=2$ and $\rhobar\simeq (\ind\omega_{4}^h)\otimes {\rm unr}(\lambda)$ is irreducible, $D_A^\otimes(\rhobar)$ only
    sees $\lambda^4$, i.e.\ does not distinguish $\rhobar$ and $(\ind\omega_{4}^h)\otimes {\rm unr}(\lambda')$, where
    $\lambda'^2=-\lambda^2$. However, one can also check (again using Theorem \ref{explicitK}) that, at least when $\rhobar$
    is $2$-dimensional and semi-simple, $D_A^\otimes(\rhobar)$ determines $\rhobar$ if $\rhobar$ is split or if
    ${\det}(\rhobar)(p)=1$.
  \end{enumerate}
\end{rem1}

\subsection{An explicit computation in the semi-simple case}\label{explicitss}

When $\rhobar$ is semi-simple we show that the explicit \'etale $(\varphi_q,\cO_K^\times)$-module $D_{A,\sigma_0}(\rhobar)$ defined in \S\ref{oK} is isomorphic to the $(\varphi_q,\cO_K^\times)$-module $D_{A}^{(0)}(\rhobar)$ defined in \S\ref{arbitrary}.

It follows from (\ref{xiyi}) that
for all $a\in\cO_K^\times$ and $0\leq i\leq f-1$, we
have (using as usual $\sigma_0:\Fq\hookrightarrow \F$) $a(X_i)=\overline{a}^{p^i}X_i$ modulo terms of degree
$\geq2$. Therefore we have the following result:

\begin{lem1}\label{lemm:a_action}
  For $0\leq i\leq f-1$, we have
  $a(X_i)\in \overline{a}^{p^i}X_i(1+A^{\ronron})$.
\end{lem1}

We define
\[ f^X_{a,0}\defeq\frac{\overline{a}X_0}{a(X_0)}\in
  1+F_{-1}A=1+A^{\ronron}\subset 1+A_\infty^{\ronron}\]
  (note that by (\ref{xiyi}) $f^X_{a,0}$ in fact coincides with $ f_{a,\sigma_0}$ in (\ref{1-p}) up to a factor in $1+F_{-2}A$).
  
\begin{lem1}\label{lemm:existence_u}
  There exists $u\in\cO(U_{\underline{n}_0})^{(1+p\cO_K)^f\cap \Delta_1}$ such that
  \[ u^{q-1}=\frac{X_{f-1}^p}{X_0}\buildrel (\ref{phionainfini}) \over =\frac{\varphi(X_0)}{X_0}\in A\subseteq A_\infty = \cO(U_{\underline{n}_0})^{\Delta_1}\subset \cO(U_{\underline{n}_0})^{(1+p\cO_K)^f\cap \Delta_1}.\]
  Moreover we have
  \begin{gather*}
  \left\{
  \begin{array}{cl}
    \forall\  \underline{a}=(a_0,\dots,a_{f-1})\in\Delta_1& \quad
    \underline{a}(u)=\overline{a_0}u\\
    \forall\  a \in\cO_K^\times& \quad
    (a,1,\dots,1)(u)=\overline{a}\left(\frac{f^X_{a,0}}{\varphi(f^X_{a,0})}\right)^{\frac{1}{q-1}}u
    \end{array}
    \right.
  \end{gather*}
  noting that $\left(\frac{f^X_{a,0}}{\varphi(f^X_{a,0})}\right)^{\frac{1}{q-1}}$ is well-defined in $1+F_{-1}A\subset 1+A_\infty^{\ronron}$ since $\frac{f^X_{a,0}}{\varphi(f^X_{a,0})}\in 1+F_{-1}A$.
 \end{lem1}

\begin{proof}
  Let $\abs$ be a multiplicative norm on $A'_\infty=\cO(U_{\underline{n}_0})$ such that
  $\vabs{T_{K,i}}=\vabs{T_{K,0}^{p^i}}=p^{-p^i}$ for $0\leq i\leq f-1$ whose existence comes
  from Lemma \ref{lemm:norm}(iii)\&(i). Let $\abs_1$ be the associated norm on
  $\B^+(A'_\infty)$ defined in Remark \ref{rema:norm1}. As $\abs$ is
  multiplicative, the same proof as in \cite[Prop.~1.4.9]{FF} shows
  that $\abs_1$ is multiplicative.

  By definition of the map $m_{A'_\infty}$ in (\ref{mR}), we have the relation in $\B^+(A'_\infty)$
  \begin{equation}\label{prodm}
    \prod_{i=0}^{f-1}\left(\sum_{n\in\Z}[T_{K,i}^{q^{-n}}]p^n\right)=\sum_{n\in\Z}\sum_{i=0}^{f-1}[X_i^{p^{-nf-i}}]p^{nf+i}=F(X_0,\dots,X_{f-1}).
    \end{equation}
  For $c\in \R_{>0}$ let $\mathfrak{p}_c$ be the ideal of $\B^+(A'_\infty)$
  \[ \mathfrak{p}_c\defeq\set{x\in\B^+(A'_\infty),\ \vabs{x}_1<p^{-c}}\subset \B^+(A'_\infty)\]
  (note that it is an ideal as $\abs_1$ is multiplicative and with values in $[0,1]\subset \R_{\geq 0}$).
  Let $c=1+p+\cdots+\cdots+p^{f-1}$. As $\vabs{T_{K,i}^{q^{n}}}=p^{-p^iq^n}\leq p^{-q^n} < p^{-c}$
  for $n\geq1$, we have $\vabs*{\sum_{n\leq -1}[T_{K,i}^{q^{-n}}]p^n}_1\leq p^{-q} < p^{-c}$, see Remark
  \ref{rema:norm1}, hence we obtain from (\ref{prodm})
  \[\prod_{i=0}^{f-1}\left(\sum_{n\geq0}[T_{K,i}^{q^{-n}}]p^n\right)-F(X_0,\dots,X_{f-1})\in\mathfrak{p}_c\]
  and we deduce from Lemma \ref{lemm:technical_norm1} below applied to the element
  \[x\defeq \prod_{i=0}^{f-1}(\sum_{n\geq0}[T_{K,i}^{q^{-n}}]p^n)-F(X_0,\dots,X_{f-1})=\sum_{n\in\Z}[x_n]p^n \in \B^+(A'_\infty)\]
  (where $\sum_{n\geq 0}[x_n]p^n=\prod_{i=0}^{f-1}(\sum_{n\geq0}[T_{K,i}^{q^{-n}}]p^n)-\sum_{n\geq 0}\sum_{i=0}^{f-1}[X_i^{p^{-nf-i}}]p^{nf+i}$ and where $\sum_{n<0}[x_n]p^n=-\sum_{n<0}\sum_{i=0}^{f-1}[X_i^{p^{-nf-i}}]p^{nf+i}$) that we have
  \[\sum_{n\geq 0}[x_n]p^n=\prod_{i=0}^{f-1}\left(\sum_{n\geq0}[T_{K,i}^{q^{-n}}]p^n\right)-\sum_{n\geq0}\sum_{i=0}^{f-1}[X_i^{p^{-nf-i}}]p^{nf+i}\in\mathfrak{p}_c.\]
  Note that the left-hand side is now in $W((A'_{\infty})^\circ)$.
 As a consequence, we have
  \[\vabs{x_0}= \vabs{T_{K,0}\cdots T_{K,f-1}-X_0}<p^{-c}\]
  so that we can write in $(A'_{\infty})^\circ$
  \begin{equation}\label{w0}
  X_0=T_{K,0}\cdots T_{K,f-1}(1+w_0)
  \end{equation}
  with $\vabs{w_0}<p^{-c+(1+p+\cdots +p^{f-1})}=1$,
  i.e.~$w_0\in (A'_\infty)^{\ronron}$. Applying the automorphism $\varphi$ of $A'_\infty$ to (\ref{w0}) and since $\varphi$ respects $(A'_{\infty})^\circ$ and $(A'_\infty)^{\ronron}$ (as it is continuous) we obtain in $(A'_{\infty})^\circ$
  \[ X_{f-1}^p=T_{K,1}T_{K,2}\cdots T_{K,f-1}T_{K,0}^q(1+w_1)\]
  with $w_1\defeq \varphi(w_0)\in (A'_\infty)^{\ronron}$. We deduce the equality
  \[ X_{f-1}^pX_0^{-1}\in T_{K,0}^{q-1}(1+(A'_\infty)^{\ronron}).\]
 Using that $x\mapsto x^{q-1}$ is bijective on $1+(A'_\infty)^{\ronron}$, we see that there exists a unique $u\in T_{K,0}(1+(A'_\infty)^{\ronron})$ such that $u^{q-1}=X_{f-1}^pX_0^{-1}$.

  As $\Delta_1$ acts trivially on $A_\infty$, we have
  $\underline{a}(u)^{q-1}=u^{q-1}$ for all
  $\underline{a}\in\Delta_1$. Therefore there exists a character
  $\chi$ of $\Delta_1$ with values in $\F_q^\times\buildrel\sigma_0\over\hookrightarrow \F$ such that
  \[ \forall\  \underline{a}\in \Delta_1, \quad
    \underline{a}(u)=\chi(\underline{a})u.\]
  Writing $u=T_{K,0}(1+w)$ with $w\in(A'_\infty)^{\ronron}$, this gives
  $\overline{a_0}T_{K,0}f_{a_0}^{\LT}(T_{K,0})^{-1}(1+\underline{a}(w))=\chi(\underline{a})u$, where $f_a^{\LT}(T_{K,0})=\overline{a}T_{K,0}(a_{\LT}(T_{K,0}))^{-1}\in 1+(A'_\infty)^{\circ\circ}$. As $u\in T_{K,0}(1+(A'_\infty)^{\ronron})$ this implies
  \[ \chi(\underline{a})\overline{a_0}^{-1}\in
    (1+(A'_\infty)^{\ronron})\cap\F_q^\times=\set{1}\] which proves
  $\chi(\underline{a})=\overline{a_0}$.

  For the last relation, we have
  \begin{multline*}
    \left((a,1,\dots,1)(u)\right)^{q-1}=(a,1,\dots,1)(u^{q-1})=a(X_{f-1}^pX_0^{-1})=\frac{a(\varphi(X_0))}{a(X_0)}=\frac{\varphi(a(X_0))}{a(X_0)}\\
   =\frac{f^X_{a,0}}{\varphi(f^X_{a,0})}\frac{\varphi(X_0)}{X_0}=\frac{f^X_{a,0}}{\varphi(f^X_{a,0})}u^{q-1}= \left(\left(\frac{f^X_{a,0}}{\varphi(f^X_{a,0})}\right)^{\frac{1}{q-1}}u\right)^{q-1}\end{multline*}
   so that as above there is a character $\chi:\cO_K^\times\rightarrow \Fq^\times\subset \F^\times$ such that $(a,1,\dots,1)(u)=\chi(a)\left(\frac{f^X_{a,0}}{\varphi(f^X_{a,0})}\right)^{\frac{1}{q-1}}u$. But $u\in T_{K,0}(1+(A'_\infty)^{\ronron})$ implies $(a,1\dots,1)(u)\in\overline{a}T_{K,0}(1+(A'_\infty)^{\ronron})$ so that $\chi(a)=\overline a$ since $\left(\frac{f^X_{a,0}}{\varphi(f^X_{a,0})}\right)^{\frac{1}{q-1}}\in 1+A_\infty^{\ronron}$. This finishes the proof.
\end{proof}

\begin{lem1}\label{lemm:technical_norm1}
  Let $R$ be a perfectoid $\F$-algebra and let $(x_n)_{\in\Z}$ a
  family of elements of $R^{\ronron}$ such that the series
  $\sum_{n\in\Z}[x_n]p^n$ converges to an element $x$ in
  $\B^+(R)$. Assume that $\vabs{x}_1< c$ for some $c\in [0,1[$. Then we have
  \[ \vabs*{\sum_{n<0}[x_n]p^n}_1 < c.\]
\end{lem1}

\begin{proof}
  Recall that $\vabs{x}_1=\displaystyle{\lim_{\substack{\rho <1 \\
        \rho\rightarrow1}}\vabs{x}_\rho}$ (see the reference in Remark \ref{rema:norm1}). Therefore we can find
  $0<\rho<1$ such that $\vabs{x}_\rho<c$. This implies
  $\sup_{n\in\Z}\set{\vabs{x_n}\rho^n}<c$ and thus for $n\leq -1$,
  $\vabs{x_n}\rho^n<c$ which implies $\vabs{x_n}<c\rho <c$. The claim then follows from
  Remark \ref{rema:norm1} applied to $c\rho$.
\end{proof}

Let $v\defeq uT_{K,0}^{-1}$. We have
$v\in 1+\cO(U_{\underline{n}_0})^{\ronron}$ from the proof of Lemma \ref{lemm:existence_u}, so that, for each
$r\in\Z_{(p)}$ ($=\Z$ localized at the prime ideal $(p)$), the element
\[ v^r\defeq\sum_{n\geq0}\binom{r}{n}(v-1)^n\in
  1+\cO(U_{\underline{n}_0})^{\ronron}\]
exists. Writing $\underline a(v)=\underline a(u)a_0(T_{K,0})^{-1}$ and using the formula for $\underline a(u)$ in Lemma \ref{lemm:existence_u} and the fact that $f_{a_0}^{\LT}(T_{K,0})\in 1+\cO(U_{\underline{n}_0})^{\ronron}$, we have
\begin{equation}\label{actionu}
\forall\  \underline{a}\in \Delta_1 \ \forall\  r\in\Z_{(p)},
  \quad \underline{a}(v^r)=f_{a_0}^{\LT}(T_{K,0})^rv^r.
\end{equation}
We also have $\varphi^f(v)/v^q\in
1+\cO(U_{\underline{n}_0})^{\ronron}$ and
$(\varphi^f(v)/v^q)^{q-1}=\varphi^f(u^{q-1})/u^{q(q-1)}=1$ as
$u^{q-1}=\varphi(X_0)/X_0$. It follows that $\varphi^f(v)=v^q$ and
\begin{equation}\label{u^q}
  \varphi^f(u)=u^q.
\end{equation}

Now, let $\rhobar$ be an absolutely irreducible continuous representation of
$\Gal(\overline{K}/K)$ on a finite-dimensional $\F$-vector space and choose a basis $(e_0,\dots,e_{d-1})$ of the $\F\ppar{T_{K,\sigma_0}^{q-1}}$-module
$D_{K,\sigma_0}(\rhobar)^{[\F_q^\times]}$ as in (\ref{D0rhobar}). We consider the associated \'etale $(\varphi_q,\cO_K^\times)$-module $D_{A,\sigma_0}(\rhobar)=A\otimes_{\F\ppar{T_{K,\sigma_0}^{q-1}}}D_{K,\sigma_0}(\rhobar)^{[\F_q^\times]}$ defined in Lemma \ref{PSnatural},
where $A$ has the structure of $\F\ppar{T_{K,\sigma_0}^{q-1}}$-algebra given by (\ref{recette}).

\begin{thm1}\label{explicitK}
  Assume that $\rhobar$ is absolutely irreducible. The
  \'etale $(\varphi_q,\cO_K^\times)$-module $D_{A}^{(0)}(\rhobar)$ in Corollary \ref{descenti} is
  isomorphic to $D_{A,\sigma_0}(\rhobar)$.
\end{thm1}

\begin{proof}
  First, replacing the variable $Y_{\sigma_0}$ by the variable $X_0$
  in Lemma \ref{PSnatural} and using (\ref{xiyi}), it is easily
  checked that one obtains an isomorphic \'etale
  $(\varphi_q,\cO_K^\times)$-module. By Corollary \ref{descentoK}, it
  is enough to prove the statement of the theorem after extending
  scalars everywhere from $A$ to $A_\infty$. Recall we have
  $D_{K,\sigma_0}(\rhobar)\!=\!\F\ppar{T_{K,\sigma_0}}\otimes_{\F\ppar{T_{K,\sigma_0}^{q-1}}}D_{K,\sigma_0}(\rhobar)^{[\F_q^\times]}$
  with basis $(1\otimes e_i)_{0\leq i\leq d-1}$ as in Lemma \ref{PS},
  let $u\in T_{K,0}(1+(A'_\infty)^{\ronron})$ be as in Lemma
  \ref{lemm:existence_u} and let again $v = uT_{K,0}^{-1}$. Then using (\ref{dai}), (\ref{actionu}) and the action of $\cO_K^\times$ in (\ref{D0rhobar}) we obtain
  \[
    A_\infty \otimes_A
    D_{A}^{(0)}(\rhobar)=(A'_\infty\otimes_{\F\ppar{T_{K,\sigma_0}^{q-1}},\iota_0}D_{K,\sigma_0}(\rhobar)^{[\F_q^\times]})^{\Delta_1}=\bigoplus_{i=0}^{d-1}A_\infty
    v^{-\frac{hq^i(q-1)}{q^d-1}}(1\otimes
    e_i).\] Moreover it follows again from the last equality in Lemma \ref{lemm:existence_u} that we have in $A'_\infty\otimes_{\F\ppar{T_{K,\sigma_0}^{q-1}},\iota_0}D_{K,\sigma_0}(\rhobar)^{[\F_q^\times]}$ for $a\in \cO_K^\times$
  \begin{multline*}
    a(v^{-\frac{hq^i(q-1)}{q^d-1}}(1\otimes
    e_i))=\left(\left(\frac{f^X_{a,0}}{\varphi(f^X_{a,0})}\right)^{\!\!\frac{1}{q-1}}\!\frac{\overline{a}T_{K,0}}{a(T_{K,0})}v\right)^{\!\!\!\!-\frac{hq^i(q-1)}{q^d-1}}\!\!\!\!\!\!\!\! f_a^{\LT}(T_{K,0})^{\!\frac{hq^i(q-1)}{q^d-1}}\!(1\otimes
    e_i)\\
    =\left(\frac{\varphi(f^X_{a,0})}{f^X_{a,0}}\right)^{\frac{hq^i}{q^d-1}}v^{-\frac{hq^i(q-1)}{q^d-1}}(1\otimes
    e_i).
  \end{multline*}
We define an $A_\infty$-linear isomorphism $A_\infty\otimes_A D_{A,\sigma_0}(\rhobar)= A_\infty\otimes_{\F\ppar{T_{K,\sigma_0}^{q-1}}}D_{K,\sigma_0}(\rhobar)^{[\F_q^\times]}\buildrel\sim\over\longrightarrow A_\infty \otimes_A D_{A}^{(0)}(\rhobar)$ by
  $1\otimes e_i\mapsto
  v^{-\frac{hq^i(q-1)}{q^d-1}}\otimes e_i$ for $i\in \{0,\dots, d-1\}$. This
  isomorphism commutes with the actions of $\cO_K^\times$ on both
  sides by the above computation (together with Lemma \ref{PSnatural}). It also commutes with $\varphi_q$, namely we have in
  $A'_\infty\otimes_{\F\ppar{T_{K,\sigma_0}^{q-1}},\iota_0}D_{K,\sigma_0}(\rhobar)^{[\F_q^\times]}$ (using \eqref{u^q}):
  \[ \varphi_q\big(v^{-\frac{hq^i(q-1)}{q^d-1}}\otimes
    e_i\big)=v^{-\frac{hq^{i+1}(q-1)}{q^d-1}}\otimes e_{i+1}\ \ {\rm for}\ \ i< d-1\]
    and (using the formula for $u^{q-1}$ in Lemma \ref{lemm:existence_u})
  \begin{eqnarray*}
    \varphi_q\big(v^{-\frac{hq^{d-1}(q-1)}{q^d-1}}\otimes
    e_{d-1}\big)&= &v^{-h(q-1)}v^{-\frac{h(q-1)}{q^d-1}}\otimes \lambda^d
    T_{K}^{-h(q-1)}e_0\\
    & =& u^{-h(q-1)}\lambda^d\big(v^{-\frac{h(q-1)}{q^d-1}}\otimes
    e_0\big)\\
    & =&\lambda^d\left(\frac{\varphi(X_0)}{X_0}\right)^{-h}\big(v^{-\frac{h(q-1)}{q^d-1}}\otimes
    e_0\big).\qedhere
  \end{eqnarray*}
\end{proof}

\begin{rem1}\label{nonsemisimplebad}
Theorem \ref{explicitK} shows that, when $\rhobar$ is a direct sum of absolutely irreducible representations, one can obtain the \'etale $\varphi_q$-module $D_{A}^{(0)}(\rhobar)$ from the Lubin--Tate $(\varphi_q,\cO_K^\times)$-module $D_{K,\sigma_0}(\rhobar)=\F\ppar{T_{K,\sigma_0}}\otimes_{\F\ppar{T_{K,\sigma_0}^{q-1}}}D_{K,\sigma_0}(\rhobar)^{[\F_q^\times]}$ by the simple recipe (\ref{recette}). However, we do not expect this recipe to work in general when $\rhobar$ is not semi-simple.
\end{rem1}

Define $D_{A,\sigma}(\rhobar)$ as $D_{A,\sigma_0}(\rhobar)$ (see \S\ref{oK}) but using the embedding $\sigma$ instead of $\sigma_0$. From \S\ref{oK} one easily checks that there are canonical $A$-linear isomorphisms for $i\in \Z$
\begin{equation}\label{isoi}
\Id\otimes\varphi:A\otimes_{\varphi,A}D_{A,\sigma_{i}}(\rhobar)\buildrel\sim\over\longrightarrow D_{A,\sigma_{i-1}}(\rhobar)
\end{equation}
which commute with $\oK^\times$ and $\varphi_q$ on both sides. Comparing the isomorphism $\phi_i$ in Corollary \ref{descenti} with the isomorphism (\ref{isoi}) we see that we have for $i\in \{0,\dots, f-1\}$
\begin{equation}\label{comparei}
D_{A}^{(i)}(\rhobar)\simeq D_{A,\sigma_{f-i}}(\rhobar).
\end{equation}
Using (\ref{comparei}) and (\ref{daotimesbis}) we have therefore
\begin{equation}\label{daotimes}
D_A^\otimes(\rhobar)\simeq D_{A, \sigma_0}(\rhobar)\otimes_AD_{A, \sigma_1}(\rhobar)\otimes_A \cdots \otimes_A D_{A,\sigma_{f-1}}(\rhobar).
\end{equation}

When $\dim_{\F}\rhobar=1$, i.e.~for $\chi:\Gal(\overline K/K)\rightarrow \F^\times$ a continuous character, we will need in \S\ref{chapterGL2} the (very simple) description of $D_A^\otimes(\chi)$.

\begin{lem1}\label{twist1}
Viewing $\chi$ as a character of $K^\times$ via the local reciprocity map, we have {\upshape(}for $a\in \oK^\times${\upshape)}:
\begin{equation*}
\left\{\begin{array}{cll}
D_{A}^\otimes(\chi)&=&A F_\chi\\
\varphi(F_\chi)&=& \chi(p) F_\chi\\
a(F_\chi)&=&\chi(a)F_\chi.
\end{array}\right.
\end{equation*}
In particular, $D_A^\otimes(\rhobar\otimes \chi)$ equals $D_A^\otimes(\rhobar)$, but with the action of $\varphi$ multiplied by $\chi(p)$ and the action of $a\in \oK^\times$ multiplied by $\chi(a)$.
\end{lem1}
\begin{proof}
By (\ref{isoi}) and (\ref{daotimes}) replacing $\rhobar$ by $\chi$ we can describe $D_A^\otimes(\chi)$ as $AE_\chi$, where $E_\chi\defeq e_{\chi}\otimes \varphi(e_{\chi})\otimes\cdots\otimes \varphi^{f-1}(e_\chi)$ with $\varphi^{j}(e_{\chi})\in D_{A,{\sigma_{f-j}}}(\chi)$ (noting $e_\chi$ instead of $1\otimes e_\chi$). Write $\chi=\omega_f^{h_{\chi}}{\rm unr}(\lambda_{\chi})$ for $h_\chi\in \Z_{\geq 0}$ and $\lambda_\chi\in \F^\times$. Set $F_\chi\defeq Y_{\sigma_0}^{h_\chi}E_\chi$, then one computes:
\begin{multline*}
\varphi(F_\chi)=\varphi(Y_{\sigma_0})^{h_\chi} \varphi(E_\chi)=\varphi(Y_{\sigma_0})^{h_\chi} \varphi^f(e_\chi)\otimes \varphi(e_{\chi})\otimes\cdots\otimes \varphi^{f-1}(e_\chi)\\
=\lambda_\chi \varphi(Y_{\sigma_0})^{h_\chi}\Big(\frac{Y_{\sigma_0}}{\varphi(Y_{\sigma_0})}\Big)^{h_\chi}E_\chi=\lambda_\chi F_\chi=\chi(p)F_\chi
\end{multline*}
where the third equality follows from (\ref{dachi}). An analogous computation using $a(Y_{\sigma_0}^{h_\chi})=\sigma(\overline a)^{h_\chi}Y_{\sigma_0}^{h_\chi}f_{a,\sigma_0}^{-{h_\chi}}$ and $a(\varphi^j(e_\chi))=\big(\frac{\varphi^j(f_{a,\sigma_0})}{\varphi^{j+1}(f_{a,\sigma_0})}\big)^{\frac{h_\chi}{1-q}}\varphi^j(e_\chi)$ (see again (\ref{dachi})) gives $a(F_\chi)=\sigma_0(\overline a)^{h_\chi}F_\chi$. But $\sigma_0(\overline a)^{h_\chi}=\chi(a)$ (see (\ref{omegaf})). The rest of the statement follows from the discussion after (\ref{dachi}).
\end{proof}

\clearpage{}

\clearpage{}\section{\'Etale \texorpdfstring{$(\varphi,\oK^\times)$}{(phi,O\_K\^{}x)}-modules and modular representations of \texorpdfstring{$\GL_2$}{GL\_2}}\label{chapterGL2}

In this section we prove that the \'etale $(\varphi,\cO_{K}^\times)$-module $D_A(\pi)$ over $A$ associated in \cite[\S3]{BHHMS2} to certain automorphic admissible smooth representations $\pi$ of $\GL_2(K)$ over $\F$ is isomorphic to (a certain twist of) the \'etale $(\varphi,\cO_{K}^\times)$-module $D_A^\otimes(\rhobar)$ of \S\ref{chapterGalois}, where $\rhobar$ is the underlying $2$-dimensional representation of ${\rm Gal}(\overline K/K)$ over $\F$, which is assumed semi-simple and sufficiently generic. We conjecture that an analogous statement holds without these assumptions and for any automorphic admissible smooth representation of $\GL_2(K)$ over $\F$.

We let $I\defeq \smatr {\cO_K^\times}{\cO_K}{p\cO_K}{\cO_K^\times}$ be the Iwahori subgroup of $\GL_2(\cO_K)$, $K_1\defeq \smatr {1+p\cO_K}{p\cO_K}{p\cO_K}{1+p\cO_K}$ the first congruence subgroup, $I_1\defeq \smatr {1+p\cO_K}{\cO_K}{p\cO_K}{1+p\cO_K}$ the pro-$p$ radical of $I$ and $Z_1$ the center of $I_1$. We recall from \S\ref{oK} that $N_0= \smatr {1}{\mathcal{O}_K}{0}{1}\subseteq I_1$. If $C$ is a pro-$p$ group we denote by $\F\bbra{C}$ its Iwasawa algebra over $\F$ (a local ring), and $\m_{C}$ the maximal ideal of $\F\bbra{C}$. If $M$ is a filtered module in the sense of \cite[\S I.2]{LiOy} with $(F_nM)_{n\in \Z}$ its ascending filtration, we define $\gr(M)\defeq \bigoplus_{n\in \Z}F_nM/F_{n-1}M$. When $R=\F\bbra{C}$ and $M$ is an $R$-module, the filtration $F_nM=\m_R^{-n}M$ if $n\leq 0$ and $F_nM=M$ if $n\geq 0$ is called the $\m_R$-adic filtration on $M$.

\subsection{A local-global compatibility conjecture for \texorpdfstring{$(\varphi,\oK^\times)$}{(phi,O\_K\^{}x)}-modu\-les over \texorpdfstring{$A$}{A}}\label{conjecture}

We conjecture that any automorphic smooth representation of $\GL_2(K)$ over $\F$ gives rise to an \'etale $(\varphi,\cO_K^\times)$-module over $A$ which is (up to twist) a direct sum of copies of the module $D_A^\otimes$ in \S\ref{otimesarbitrary} of the corresponding local Galois representation at $p$. We state our main results.

First, we quickly review the construction of the $A$-module $D_A(\pi)$ associated to certain smooth representations $\pi$ of $\GL_2(K)$ over $\F$ in \cite[\S 3.1]{BHHMS2}.

Let $\pi$ be an admissible smooth representation of $\GL_2(K)$ over $\F$ with a central character and endow the $\F$-linear dual $\pi^\vee$ with the $\m_{I_1}$-adic filtration, or equivalently the $\m_{I_1/Z_1}$-adic filtration (which, in general, {\it strictly} contains the $\m_{N_0}$-adic filtration). We endow
\[(\pi^\vee)_{(Y_{\sigma_0}\cdots Y_{\sigma_{f-1}})}\defeq \F\bbra{N_0}_{(Y_{\sigma_0}\cdots Y_{\sigma_{f-1}})}\otimes_{\F\bbra{N_0}}\pi^\vee\]
with the tensor product filtration (where the localization $\F\bbra{N_0}_{(Y_{\sigma_0}\cdots Y_{\sigma_{f-1}})}$ is endowed with the 
filtration described by \eqref{eq:def:fil:Aq}, replacing $\F_q$ by $\F$)
and define $D_A(\pi)$ as the completion of $(\pi^\vee)_{(Y_{\sigma_0}\cdots Y_{\sigma_{f-1}})}$ for this filtration (\cite[\S I.3.4]{LiOy}). Then $D_A(\pi)$ is a complete filtered $A$-module and the action of $\cO_K^\times$ on $\pi^\vee$ extends by continuity to $D_A(\pi)$. Moreover the action $f\mapsto f\circ \smatr{p}{0}{0}{1}$ on $\pi^\vee$ gives rise to a continuous $A$-linear morphism (see \cite[\S 3.1.2]{BHHMS2})
\begin{equation}\label{betapsi}
\beta:D_A(\pi)\longrightarrow A\otimes_{\varphi,A}D_A(\pi),
\end{equation}
where $\varphi$ on $A$ is as in \S\ref{oK}. We let $\mathcal C$ be the abelian category of those $\pi$ such that $\gr(D_A(\pi))$ is a finitely generated $\gr(A)$-module. Then for any $\pi\in \mathcal C$, the $A$-module $D_A(\pi)$ is finite free (see \cite[Cor.~3.1.2.9]{BHHMS2} and Remark \ref{freeness}).

The following straightforward lemma will be used. For $\chi:K^\times\rightarrow \F^\times$ a smooth character, denote by $D_A(\chi)$ the rank $1$ \'etale $(\varphi,\oK^\times)$-module over $A$ defined by $Ae_\chi$ with $\varphi(e_\chi)\defeq \chi(p)e_\chi$ and $a(e_\chi)\defeq \chi(a)e_\chi$ for $a\in \oK^\times$.
{(Note that this is an ad hoc definition, as $D_A(\pi) = 0$ if $\pi = \chi \circ \det$.)}

\begin{lem1}\label{twist}
Let $\chi:K^\times\rightarrow \F^\times$ be a smooth character and $\pi$ in the category $\mathcal C$, then $D_A(\pi\otimes \chi)\cong D_A(\pi)\otimes_A D_A(\chi^{-1})$ with diagonal $\varphi$ and action of $\oK^\times$.
\end{lem1}
\begin{proof}
This directly follows from the definitions of $D_A(\pi)$ and of the actions of $\varphi$ and $\oK^\times$ on $D_A(\pi)$.
\end{proof}

For $\pi$ in $\mathcal C$, when $\beta$ is moreover an isomorphism, its inverse $\beta^{-1}=\Id\otimes \varphi$ makes $D_A(\pi)$ an \'etale $(\varphi,\cO_K^\times)$-module.

We now go to the global setting.

We fix a totally real number field $F$ that is unramified at $p$. We fix a quaternion algebra $D$ of center $F$ which is split at all places above $p$ and at not more than one infinite place. When $D$ is split at one infinite place we say that we are in the \emph{indefinite case}, and in the \emph{definite case} otherwise. For a compact open subgroup $U= \prod U_w\subset (D\otimes_F{\mathbb A}_F^{\infty})^\times$ we let $X_U$ be the associated smooth projective algebraic Shimura curve over $F$ (see \cite[\S 8.1]{BHHMS2} and the references therein for more details). 

Fix an absolutely irreducible continuous representation $\rbar:\gF\ra \GL_2(\F)$ and for a finite place $w$ of $F$ we write $\overline r_w\defeq \overline r\vert_{\gFw}$. We let $S_D$ be the set of finite places where $D$ ramifies, $S_{\rbar}$ the set of (finite) places where $\rbar$ is ramified and $S_p$ the set of (finite) places above $p$. Finally, we fix a place $v\in S_p$.
Let $\omega = \omega_1$ denote the mod $p$ cyclotomic character.

For any compact open subgroup $U^v= \prod_{w \ne v} U_w\subset (D\otimes_F{\mathbb A}_F^{\infty,v})^\times$ we consider the following admissible smooth representation $\pi$ of $\GL_2(F_v)$ over $\F$ with central character $(\omega {\det}(\overline r_v))^{-1}$:
\begin{equation}\label{piindef}
\pi \defeq\varinjlim_{U_v}\Hom_{\gF}\big(\rbar, H^1_{{\rm \acute et}}(X_{U^vU_v} \times_{F} \overline F, \F)\big)
\end{equation}
where the inductive limit runs over the compact open subgroups $U_v$ of $(D\otimes_FF_v)^\times\cong \GL_2(F_v)$.
In the definite case, we replace $\Hom_{\gF}(\rbar, H^1_{{\rm \acute et}}(X_{U} \times_{F} \overline F, \F))$ by the Hecke eigen\-space $S(U,\F)[{\mathfrak m}]\subseteq S(U,\F)\defeq \{f:D^\times\backslash (D\otimes_F{\mathbb A}_F^\infty)^\times/U \rightarrow \F\}$ associated to $\overline r$ (see \cite[\S 8.1]{BHHMS2}) and define analogously
\begin{equation}\label{pidef}
\pi \defeq\varinjlim_{U_v}S(U^vU_v,\F)[{\mathfrak m}].
\end{equation}
We also need the ``multiplicity $1$'' variants of the representations $\pi$. For that, we need to assume that $p\geq 5$, that $\rbar |_{G_{F(\!\sqrt[p]{1})}}$ is absolutely irreducible, that the image of $\rbar(G_{F(\!\sqrt[5]{1})})$ in ${\rm PGL}_2(\F)$ is not isomorphic to $A_5$, that $\rbar_w$ for $w\in S_p$ is generic in the sense of \cite[Def.\ 11.7]{BP} (which implies $S_p\subseteq S_{\rbar}$) and that $\rbar_w$ is non-scalar if $w\in S_D$. Under these assumptions, a so-called ``local factor'' is defined in \cite[\S3.3]{BD} (in the indefinite case and when $\rbar_w$ is reducible for all $w\in S_p$) and in \cite[\S 6.5]{EGS} (without these two conditions):
\begin{eqnarray}
\label{locpiindef} \!\!\!\!\!\!\pi\!&\!\!\defeq \!\!&\Hom_{U^v}\!\!\Big(\overline M^v,\Hom_{\gF}\big(\rbar, \varinjlim_{V} H^1_{{\rm \acute et}}(X_{V} \times_{F} \overline F, \F)\big)\Big)[{\mathfrak m}']\ \textrm{(indefinite case)}\\
\label{locpidef} \!\!\!\!\!\!\pi\!&\!\!\defeq \!\!&\Hom_{U^v}\!\!\Big(\overline M^v,\varinjlim_{V} S(V,\F)[{\mathfrak m}]\Big)[{\mathfrak m}']\ \textrm{(definite case)}
\end{eqnarray}
where the inductive limits run over the compact open subgroups $V$ of $(D\otimes_F{\mathbb A}_F^{\infty})^\times$, and where we refer to {\emph{loc.~cit.}\ for the definitions of the compact open subgroup $U^v\subseteq (D\otimes_F{\mathbb A}_F^{\infty,v})^\times$, of the (finite-dimensional) irreducible smooth representation $\overline M^v$ of $U^v$ over $\F$ and of the maximal ideal ${\mathfrak m}'$ in a certain Hecke algebra. 

\begin{conj1}\label{dAPi}
Let $\pi$ be as in (\ref{piindef}), (\ref{pidef}), (\ref{locpiindef}) or (\ref{locpidef}) and assume $\pi\ne 0$. Then $\pi$ is in the category $\mathcal C$, $\beta$ in (\ref{betapsi}) is a bijection and we have an isomorphism of \'etale $(\varphi,\oK^\times)$-modules $D_A(\pi)\simeq D_A^\otimes(\overline r_v(1))^{\oplus r}$ for some integer $r\geq 1$ which is equal to $1$ when $\pi$ is as in (\ref{locpiindef}) or (\ref{locpidef}).
\end{conj1}

In the sequel, we prove Conjecture \ref{dAPi} for $\pi$ as in (\ref{locpiindef}) or (\ref{locpidef}) when $\overline r_v$ is semi-simple and satisfies a strong genericity hypothesis (as defined below). We actually prove a purely local result for certain smooth representations $\pi$, that will ultimately include the representations in (\ref{locpiindef}) and (\ref{locpidef}).

Let first $\rhobar:\Gal(\overline K/K)\rightarrow \GL_2(\F)$ be a continuous representation satisfying the genericity assumption of \cite[Def.~11.7]{BP}. Let $\pi$ be a smooth representation of $\GL_2(K)$ over $\F$ satisfying the following two conditions:
\begin{enumerate}
\item there is an isomorphism of diagrams $(\pi^{I_1}\hookrightarrow \pi^{K_1})\cong D(\rhobar)^{\oplus r}$ for some $r\in \Z_{\geq 1}$, where $D(\rhobar)$ is a diagram associated to $\rhobar$ as in \cite{BP} or \cite[\S 3.2.1]{BHHMS2} with the constants $\nu_\sigma$ for $\sigma\in W(\rhobar)$ as in Remark \ref{rem:lambda-DL} below;
\item for any character $\chi:I\rightarrow \F^\times$ appearing in $\pi[{\mathfrak m}_{I_1}]$ there is an equality of multiplicities $[\pi[{\mathfrak m}^3_{I_1}]:\chi]=[\pi[{\mathfrak m}_{I_1}]:\chi]$.
\end{enumerate}
We moreover assume that $\rhobar$ is of the following form {\it up to twist}:
\begin{equation}\label{semisimpler}
\rhobar|_{I_K} \cong
\begin{cases}
\omega_{f}^{\sum_{j=0}^{f-1}(r_j+1)p^j}\oplus 1&\text{ if $\rhobar$ is reducible}
\\
\omega_{2f}^{\sum_{j=0}^{f-1} (r_{j}+1)p^j}\oplus \omega_{2f}^{\sum_{j=0}^{f-1}(r_j+1)p^{j+f}}
&\text{ if $\rhobar$ is irreducible}
\end{cases}
\end{equation}
where the integers $r_i$ satisfy the following (strong) genericity condition:
\begin{equation}\label{stronggen}
\begin{aligned}
\max\{12,2f-1\} &\le &\!\! r_j &\le & \! p-\max\{15,2f+2\} & \text{ if $j > 0$ or $\rhobar$ is reducible}\\
\max\{13,2f\} \!&\le &\!\! r_0 &\le & \! p-\max\{14,2f+1\} & \text{ if $\rhobar$ is irreducible.}
\end{aligned}
\end{equation}

The following is the main result of \S\ref{chapterGL2}.

\begin{thm1}[See \S\ref{sec:structure-d_ap}]\label{mainbis}
Assume that $\rhobar$ and $\pi$ are as above with moreover $(\pi^{I_1}\hookrightarrow \pi^{K_1})\!\cong D(\rhobar)$, i.e.~$r=1$. Then $\pi$ is in the category $\mathcal C$, $\beta$ in (\ref{betapsi}) is a bijection and we have an isomorphism of \'etale $(\varphi,\oK^\times)$-modules $D_A(\pi)\simeq D_A^\otimes(\rhobar^\vee(1))$, where $\rhobar^\vee(1)$ is the Cartier dual of $\overline \rho$.
\end{thm1}

It implies the following special cases of Conjecture \ref{dAPi}.

\begin{cor1}\label{maincor}
Let $\pi$ be as in (\ref{locpiindef}) or (\ref{locpidef}) and assume moreover that $\rbar_v$ satisfies (\ref{stronggen}), (\ref{semisimpler}), and that the framed deformation ring $R_{\rbar_w}$ of $\rbar_w$ over $W(\F)$ is formally smooth if $w\in (S_D \cup S_{\rbar})\setminus S_p$. Then Conjecture \ref{dAPi} is true for $\pi$.
\end{cor1}
\begin{proof}
By \cite[Thm.~5.36]{DoLe} (and the references therein) $\pi$ satisfies condition (i) above with $\rhobar=\rbar_v^\vee$ and $r=1$. By \cite[Thm.~8.3.14]{BHHMS1}, \cite[Thm.~1.5]{BHHMS1} and \cite[Rem.~8.4.5]{BHHMS1} $\pi$ satisfies condition (ii). Hence we can apply Theorem \ref{mainbis}.
\end{proof}

\begin{rem1}\ 
\begin{enumerate}
\item Under the assumptions of Theorem \ref{mainbis}, we already know that $\pi$ is in $\mathcal C$ (see \cite[Thm.~3.3.2.1]{BHHMS2}) and that $\beta$ is a bijection (see \cite[Rem.~3.3.5.4(ii)]{BHHMS2} noting that we do not need here the assumption (iii) in \cite[\S 3.3.5]{BHHMS2}). Hence we only need to prove $D_A(\pi)\simeq D_A^\otimes(\rhobar^\vee(1))$. In that direction, we already know the \'etale $(\varphi,\Zp^\times)$-module $\F\ppar{T}\otimes_AD_A(\pi)$. Indeed, it follows from \cite[Cor.~3.3.2.4]{BHHMS2}, \cite[Thm.~3.1.3.7]{BHHMS2}, Remark \ref{freeness} -- and some unravelling of the definition of the functor $V_{\GL_2}$ of \cite[\S2.1.1]{BHHMS2} using Lemma \ref{twist} and Lemma \ref{twist1} -- that $D_A(\pi)$ is free of rank $2^f$ and $\F\ppar{T}\otimes_AD_A(\pi)$ is isomorphic to the $(\varphi,\Zp^\times)$-module of the tensor induction $\ind_K^{\otimes\Qp}\!(\rhobar^\vee(1))$ (compare with Corollary \ref{determine}).\\
\item  {Under similar hypothesis but assuming that $\rhobar$ is reducible non-split with only one Serre weight, Conjecture \ref{dAPi} is proven in \cite{YitongWang3} (using the results of \cite{YitongWang2}).}
\end{enumerate}
\end{rem1}

The rest of this paper is devoted to the proof of the isomorphism $D_A(\pi)\simeq D_A^\otimes(\rhobar^\vee(1))$ in Theorem \ref{mainbis} and to the necessary material that needs to be introduced for that.

We fix $\rhobar$ and $\pi$ as in Theorem \ref{mainbis}. Twisting both $\rhobar$ and $\pi$ using Lemma \ref{twist1} and Lemma \ref{twist}, we can and do assume {\it from now on} $\rhobar\simeq (\ind\omega_{2f}^h)\otimes {\rm unr}(\lambda)$ or $\rhobar\simeq \begin{pmatrix}\omega_f^h{\rm unr}(\lambda_0)& 0\\ 0 & {\rm unr}(\lambda_1)\end{pmatrix}$ with $h=\sum_{j=0}^{f-1}p^j(r_j+1)$.

\subsection{Duality for \'etale \texorpdfstring{$(\varphi,\cO_K^\times)$}{(phi,O\_K\^{}x)}-modules over \texorpdfstring{$A$}{A}}
\label{sec:dual-vp-co_ks}

If $D$ is an \'etale $(\vp, \cO_K\s)$-module over $A$ we equip $\Hom_A(D,A)$ with the structure of an \'etale $(\vp,\cO_K\s)$-module over $A$.

Fix $D$ an \'etale $(\vp, \cO_K\s)$-module over $A$. 

We first equip $D$ with a left inverse $\psi : D \to D$ of $\vp$, as follows.
Fix a set of representatives $\{n\}$ of $N_0/N_0^p$ including 1.
Note that as $D$ is \'etale, every element $x$ of $D$ can be uniquely written as $x = \sum_{N_0/N_0^p} \delta_n \vp(x_n)$, where $\delta_n$ denotes the image of the element $[n]\in\F\bbra{N_0}$ in $A$.
Let $\psi : D \to D$ be defined by $\psi(x) \defeq x_1$. The following easy lemma is left to the reader.

\begin{lem1}\label{lm:psi}
  The map $\psi : D \to D$ is a left inverse of $\vp$ that is independent of any choices. We have
  $x = \sum_{N_0/N_0^p} \delta_n \vp(\psi(\delta_n^{-1} x))$ for any $x \in D$.
  Moreover, the actions of $\psi$ and $\cO_K\s$ commute.
\end{lem1}

To define $\varphi$ on $\Hom_A(D,A)$ recall that we have
\begin{equation}\label{eq:3f}
  \begin{aligned}
    \beta : D &\congto A \otimes_{\vp} D\\
    x &\mapsto \sum_n \delta_n \otimes_{\vp} \psi(\delta_n^{-1} x),
  \end{aligned}
\end{equation}
where the sum runs over representatives $\{n\}$ of $N_0/N_0^p$. Now if $M$, $N$ are $A$-modules with $M$ projective, we have for any $A$-algebra $B$
a canonical isomorphism $B \otimes_A \Hom_A(M,N) \cong \Hom_B(B
\otimes_A M,B \otimes_A N)$,
hence the $A$-linear dual of (\ref{eq:3f}) gives rise to
\begin{equation*}
  A \otimes_{\vp} \Hom_A(D,A) \congto \Hom_A(D,A), 
\end{equation*}
in other words we get a $\vp$-linear endomorphism of $\Hom_A(D,A)$ that we also call $\vp$ (an \'etale
Frobenius). Explicitly, this endomorphism is given by the formula
\begin{align}
  \Hom_A(D,A) &\to \Hom_A(D,A)\notag\\
  h &\mapsto \vp(h) = (x \mapsto \sum_{N_0/N_0^p} \delta_n \vp(h(\psi(\delta_n^{-1} x)))).\label{eq:5f}
\end{align}
By construction, it is independent of the choice of representatives.

Using Lemma~\ref{lm:psi} we can rewrite formula~\eqref{eq:5f} as follows:
\begin{equation}
  \vp(h) : \sum_n \delta_n \vp(x_n) \mapsto \sum_n \delta_n \vp(h(x_n)).\label{eq:6f}
\end{equation}
We also define the action of $a \in \cO_K\s$ by the formula $a(h) \defeq a \circ h \circ a^{-1}$.

\begin{lem1}\label{lm:dual-etale}
  With the definitions above, $\Hom_A(D,A)$ is an \'etale $(\vp,\cO_K\s)$-module. Moreover, the natural
  pairing $D \times \Hom_A(D,A) \to A$ is equivariant for the actions of $\vp$ and $\cO_K\s$.
\end{lem1}

Concretely, the $(\varphi,\cO_K^\times)$-module structure on
$\Hom_A(D,A)$ is uniquely characterized by the relations, for
$x\in D$, $y\in\Hom_A(D,A)$, $a\in\cO_K^\times$:
\begin{equation}
  \langle\varphi(x),\varphi(y)\rangle=\varphi(\langle x,y\rangle), \qquad \langle a(x),a(y)\rangle=a(\langle x,y\rangle),\label{eq:pairing-phi-gamma}
\end{equation}
where $\langle\cdot,\cdot\rangle$ is the natural pairing
$D\times\Hom_A(D,A)\rightarrow A$.

Fix now a smooth representation $\pi$ of $\GL_2(K)$ in the category $\cC$ and endow the finite
free $A$-module $D_A(\pi)$ with its filtration coming from the
$\fm_{I_1}$-adic filtration on $\pi^\vee$, cf.~\S\ref{conjecture}. If
$D$ is an \'etale $(\varphi,\cO_K^\times)$-module (endowed with its
natural topology of finite free $A$-module), recall that
$\Hom^{\cont}_\F(D,\F)$ is the vector space of continuous $\F$-linear morphisms $D\rightarrow \F$, or equivalently ($\F$ being endowed with the discrete topology) the $\F$-linear locally constant morphisms $D\rightarrow \F$. We give $\F$ the filtration such that $F_d\, \F = 0$ if and only if $d < 0$.

We write now $Y_i$ for $Y_{\sigma_i}$ (as in \cite[\S 3.1.1]{BHHMS2}, note that there will be no confusion with the variables $Y_i\in A_q$ in \S\ref{oK} which are not used here) and $\un Y^{(i_0,\dots,i_{f-1})}$ for $Y_0^{i_0}\cdots Y_{f-1}^{i_{f-1}}\in A$ (as in \cite[\S 3.2.2]{BHHMS2}).
We also sometimes use the shorthand $\un Y$ for $\un Y^{\un 1} = \prod_{j=0}^{f-1} Y_j$.

\begin{prop1}\label{prop:F-dual}
  There \ is \ an \ isomorphism \ of \ \ $\F\bbra{N_0}$-modules \ between \ $\Hom_\F^{\cont}(D_A(\pi),\F)$ and the set of sequences $(x_{k})_{k \ge 0}$ such that $x_k \in \pi$ and
  \begin{enumerate}
  \item $\un Y x_{k} = x_{k-1}$ for all $k \ge 1$;
  \item there exists $d \in \Z$ such that $x_{k} \in \pi[\m_{I_1}^{fk+d+1}]$ for all $k \ge 0$ {\upshape(}where $\pi[\m_{I_1}^j] \defeq 0$ for $j \le 0${\upshape)}.
  \end{enumerate}
  A continuous $\F$-linear map $h: D_A(\pi) \to \F$ corresponds to a sequence $(x_k)_{k \ge 0}$ as above if and only if
  \begin{equation*}
    h(\un Y^{-\un k}y) = \ang{x_{k},y}
  \end{equation*}
  for all $y\in \pi^\vee$, $k \ge 0$ and where we denote again by $y$ the image
  of $y$ in $D_A(\pi)$. Moreover, $h$ is filtered of degree $d$ if and only if $x_{k} \in \pi[\m_{I_1}^{fk+d+1}]$ for all $k \ge 0$.
\end{prop1}
\begin{proof}
  Let $S$ denote the multiplicative subset of $\F\bbra{N_0}$ generated by $Y_0\cdots Y_{f-1}$.
  Then from the definitions we have 
\begin{equation*}
    (\pi^\vee)_S \cong \vilim_{\substack{k \ge 0\\Y_0\cdots Y_{f-1}}} \pi^\vee \text{\ \ and\ \ } 
    F_{-d-1}(\pi^\vee)_S \cong \vilim_{\substack{k \ge 0\\Y_0\cdots Y_{f-1}}} \m_{I_1}^{fk+d+1}\pi^\vee,
  \end{equation*}
  so
  \begin{equation*}
    (\pi^\vee)_S/F_{-d-1}(\pi^\vee)_S \cong \vilim_{\substack{k \ge 0\\Y_0\cdots Y_{f-1}}} \pi^\vee/\m_{I_1}^{fk+d+1}\pi^\vee.
  \end{equation*}
  (Explicitly, the $k$-th map $\pi^\vee \to (\pi^\vee)_S$ is given by
  multiplication by $(Y_0\cdots Y_{f-1})^{-k}$.)
  Therefore, we have
  \begin{multline*}
    \Hom_\F^{\cont}(D_A(\pi),\F) = \Hom_\F^{\cont}((\pi^\vee)_S,\F) = \bigcup_{d \ge 0} \Hom_\F((\pi^\vee)_S/F_{-d-1}(\pi^\vee)_S,\F) \\ 
    = \bigcup_{d \ge 0} \Hom_\F((\pi^\vee)_S/F_{-d-1}(\pi^\vee)_S,\F)
    = \bigcup_{d \ge 0} \vplim_{\substack{k \ge 0\\Y_0\cdots Y_{f-1}}} \Hom_\F(\pi^\vee/\m_{I_1}^{fk+d+1}\pi^\vee,\F)\\
    = \bigcup_{d \ge 0} \vplim_{\substack{k \ge 0\\Y_0\cdots Y_{f-1}}} \pi[\m_{I_1}^{fk+d+1}].
  \end{multline*}
  The final claims follow by unravelling these identifications.
\end{proof}

We now make explicit the actions of $A$ and $\cO_K^\times$ on $\Hom_\F^{\cont}(D_A(\pi),\F)$, where the definitions of these actions in the following lemma are a posteriori motivated by Lemma~\ref{lm:mu-star} (namely, the map $\mu_*$ in~\eqref{eq:mu-star} becomes $A$ and $\cO_K^\times$-linear).

\begin{lem1}\label{lm:action}
  Suppose that $h : D_A(\pi) \to \F$ is continuous of degree $d$, i.e.\ sending $F_{-d-1}D_A(\pi)$ to 0.
  Let $h$ correspond to the sequence $(x_k)_{k \ge 0}$ as in Proposition~\ref{prop:F-dual}, so $\un Y x_{k+1} = x_k$ and
  $x_k \in \pi[\m_{I_1}^{kf+d+1}]$.
  \begin{enumerate}
  \item If $a \in A$, then $ah \defeq h \circ a$ corresponds to the sequence $(y_k)_{k \ge 0}$, where 
    \begin{equation}\label{eq:10f}
      y_k = \un Y^{\un\ell-\un k} a x_\ell
    \end{equation}
    for $\ell \gg_k 0$.
  \item If $a \in \cO_K\s$, then $a(h) \defeq N_{\F_q/\Fp}(\o a)^{-1} (h \circ \diag(a^{-1},1))$ corresponds to the sequence $(z_k)_{k \ge 0}$, where
    \begin{equation}\label{eq:11f}
      z_k = N_{\F_q/\Fp}(\o a)^{-1} \smatr{a}{}{}{1} \frac{\un Y^{\un\ell}}{a^{-1}(\un Y^{\un k})} x_\ell 
      = N_{\F_q/\Fp}(\o a)^{-1}\frac{a(\un Y^{\un \ell})}{\un Y^{\un k}} \smatr{a}{}{}{1} x_\ell
    \end{equation}
    for $\ell \gg_k 0$.
  \end{enumerate}
\end{lem1}

\begin{rem1}\label{rk:convention-action}
  To explain the notation in equations (\ref{eq:10f}), (\ref{eq:11f}) we note that for $x \in \pi[\m_{I_1}^e]$ ($e \ge 0$) we
  can extend the action of $\F\bbra{N_0}$ on $x$ to an action of the ring $\F\bbra{N_0}+F_{-e}A$ such that $F_{-e}A$ kills $x$
  (because $F_{-e}\F\bbra{N_0} = \F\bbra{N_0} \cap F_{-e}A$ kills $x$, by assumption). For (\ref{eq:10f}) we note that
  $\un Y^{-\un k}a \in A = \F\bbra{N_0}_S + F_{-d-1}A$ (where $S$ is generated by $\un Y$), so $\un Y^{\un \ell-\un k} a \in \F\bbra{N_0} + F_{-\ell f-d-1}A$ for $\ell \gg_k 0$ and
  $x_\ell \in \pi[\m_{I_1}^{\ell f+d+1}]$. Similarly for (\ref{eq:11f}) we note that $\frac{a(\un Y^{\un \ell})}{\un Y^{\un k}} \in \F\bbra{N_0} + F_{-\ell f-d-1}A$
  for $\ell \gg_k 0$ (and $\smatr{a}{}{}{1}$ normalizes $I_1$).
\end{rem1}

\begin{proof}
  For (i) we first note that $h(F_{-d-1}A \cdot \pi\dual) \subset h(F_{-d-1}D_A(\pi)) = 0$, so $h \circ a'|_{\pi\dual}$ only
  depends on $a'$ modulo $F_{-d-1}A$. Writing $\un Y^{-\un k}a \in \un Y^{-\un \ell} b + F_{-d-1}A$ as above with $b \in \F\bbra{N_0}$
  and $\ell \gg_k 0$, we compute for $k \ge 0$,
  \begin{equation}
    h \circ a \circ \un Y^{-\un k}|_{\pi\dual} = h \circ \un Y^{-\un \ell} \circ b|_{\pi\dual} = \ang{x_\ell,b(-)} = \ang{b x_\ell,-}
    = \ang{\un Y^{\un \ell-\un k}a x_\ell,-}\label{eq:12f}
  \end{equation}
  as functions $\pi\dual \to \F$, as desired (keeping in mind Remark~\ref{rk:convention-action}).

  For (ii), first note that $a(h) \circ \un Y^{-\un k} = N_{\F_q/\Fp}(\o a)^{-1} \big(h \circ a^{-1}(\un Y^{-\un k}) \circ \diag(a^{-1},1)\big)$.  By
  (\ref{eq:12f}) (applied with $k = 0$),
  $h \circ a^{-1}(\un Y^{-\un k})|_{\pi\dual} = \ang{\un Y^{\un \ell} \cdot a^{-1}(\un Y^{-\un k}) x_\ell,-}$ for $\ell \gg_k 0$
  and the result follows.
\end{proof}

\subsection{The continuous morphism \texorpdfstring{$\mu : A \to \F$}{mu : A -> F}}
\label{sec:finding-mu}

For $D$ an \'etale $(\vp, \cO_K\s)$-module over $A$ we relate $\Hom_A(D,A)$ to $\Hom^{\cont}_\F(D,\F)$ using a certain continuous morphism $\mu : A \to \F$.

Let us write $\F\bbra{N_0} = \F\bbra{T_0,\dots,T_{f-1}}$ with $T_j \defeq [\alpha_j]-1$, where $(\alpha_j)_{j\in \{0,\dots,f-1\}}$ is a fixed $\zp$-basis of $\cO_K$.
Recall that $A$ is endowed with a map $\psi:A\rightarrow A$ defined in \S \ref{sec:dual-vp-co_ks}, and which is a left inverse of $\vp: A\rightarrow A$.

\begin{prop1}\label{prop:mu}
  Up to scalar in $\F^{\times}$ there exists a unique $\mu \in
  \Hom^{\cont}_\F(A,\F)$ such that $\mu \circ \psi \in \F^\times \mu$,
  and we have
  $\mu \circ \psi = (-1)^{f-1}\mu$.
\end{prop1}

It will be convenient for the proof to avoid using the variables $Y_j$.
To obtain $A$ from $\F\bbra{N_0}$ it suffices to invert elements $Z_j$ ($0 \le j \le f-1$) such that $\gr(Z_j) = \gr(Y_j)$ in the graded ring and
then complete.
We will let $Z_j$ be the unique linear combination of the $T_{j'}$ such that $\gr(Z_j) = \gr(Y_j)$.
(Note that the $Z_j$ are not canonical but depend on the choice of $T_0$, \dots, $T_{f-1}$.)
There exists an element of $\GL_f(\F)$ that
relates the $Z_j$ and the $T_j$. Hence we get the same description of $A$ as in \cite[Rk.\ 3.1.1.3(iii)]{BHHMS2} with $Z_j$ instead of
$Y_j$, and also the valuation of an element of $A$ is still given by the minimal total degree as a series in $\un Z$. We note that
$\vp(T_j) = T_{j}^p$ and $\vp(Z_j) = Z_{j-1}^p$ (because $\vp(Z_j)$ is a homogeneous polynomial of degree $p$ in the $T_{j'}$ and hence in the $Z_{j'}$,
and since $\vp(Y_j) = Y_{j-1}^p$).

Before starting the proof of Proposition~\ref{prop:mu} we note that $\mu \circ \psi = c \mu$ (with $c \in \F^\times$) is equivalent to the two conditions
\begin{align}
  \mu &= c (\mu \circ \vp),\label{eq:mu1} \\
  \mu(\delta_n\vp(x)) &= 0\qquad\qquad\forall\ n \in N_0 \setminus N_0^p, \forall\ x \in A. \label{eq:mu2}
\end{align}
This follows immediately from the definition of $\psi : D \to D$ in \S\ref{sec:dual-vp-co_ks}.

\begin{proof}[Proof {\upshape(}Uniqueness{\upshape)}]
  Suppose that $\mu \circ \psi = c\mu$ for some $c \in \F^\times$.
  For the representatives $\{n\}$ of $N_0/N_0^p$ we take $n = \smatr{1}{\sum_j i_j\alpha_j}{}{1}$ ($0 \le i_j \le p-1$), so $\delta_n = \prod_j (1+T_j)^{i_j}$.
  By induction and~\eqref{eq:mu1}--\eqref{eq:mu2} we have
  for any $\un 0 \le \un i \le \un{p-1}$ that
  \begin{equation}\label{eq:1f'}
    \begin{aligned}
      \mu(\un T^{\un i} \vp(x)) &= (-1)^{\|\un i\|} \mu(\vp(x)) \\ &= (-1)^{\|\un i\|} c^{-1} \mu(x).
    \end{aligned}
  \end{equation}
  Take now $x \in F_{f-1} A$. Then by iterating~\eqref{eq:1f'} we have
  \[\mu(x) = c \mu(\un T^{\un{p-1}}\vp(x)) = \cdots = c^n \mu(\un T^{\un{p^n-1}}\vp^n(x)) = 0\] 
  for $n \gg 0$, since $\un T^{\un{p^n-1}}\vp^n(x) \to 0$ in $A$ as $n \to \infty$ if $x \in F_{f-1} A$ and $\mu$ is continuous. Hence
  \begin{equation}
    \mu(F_{f-1} A) = 0.\label{eq:2f}
  \end{equation}

  We claim that $\mu(\un Z^{\un i})$ for $\un i \in \Z^f$ is an explicit multiple of $\mu(\un Z^{-\un 1})$, only depending on $c$. 
To prove the claim, we may suppose that $\|\un i\| \le -f$ by~\eqref{eq:2f} and we will argue by descending induction on $\|\un i\|$.
  Write $\un i = \un r+p\un s$ with
  $\un 0 \le \un r \le \un{p-1}$ and $\un s \in \Z^f$. Hence $\mu(\un Z^{\un i}) = \mu(\un Z^{\un r} \un Z^{p\un s})$ and
  that can be expressed in terms of various $\mu(\un T^{\un{r'}} \un Z^{p\un s})$ with $\un{r'} \ge \un 0$ and
  $\|\un{r'}\| = \|\un r\|$. Fix now one such term and write $\un{r'} = \un{r''} + p\un{r'''}$ with
  $\un 0 \le \un{r''} \le \un{p-1}$ and $\un 0 \le \un{r'''}$. Then we can express $\mu(\un T^{\un{r'}} \un Z^{p\un s}) = \mu(\un T^{\un{r''}} \un T^{p\un{r'''}}\un Z^{p\un s})$ in
  terms of various $\mu(\un T^{\un{r''}}\un Z^{p\un t})$ with $\|\un t\| = \|\un s\|+\|\un{r'''}\|$.  By~\eqref{eq:1f'} we are
  reduced to $\pm \mu(\un Z^{p\un t}) = \pm c^{-1} \mu(\un Z^{\un{t'}})$, where $\un{t'}$ is a cyclic permutation of $\un t$
  and hence $\|\un{t'}\| = \|\un t\| = \|\un s\|+\|\un{r'''}\| = (\|\un i\|-\|\un{r''}\|)/p$.

  From $\|\un{r''}\| \le (p-1)f$ and $\|\un i\| \le -f$ it follows that $\|\un i\| \le \|\un{t'}\|$ and moreover that equality can only
  hold if $\un{r''} = \un{p-1}$ and $\|\un i\| = -f$, in which case $\un r = \un{r'} = \un{p-1}$ and $\un{r'''} = \un 0$ (as $\|\un r\| = \|\un{r''}\|+p\|\un{r'''}\| \le (p-1)f$).
  Thus $\|\un i\| < \|\un{t'}\|$ and we are done by induction, except possibly when $\|\un i\| = -f$ and $\un i \equiv -\un 1 \pmod p$.
  Applying the same argument to $\mu(\un Z^{\un{t'}})$, we are done in the exceptional case except if $\un{t'} \equiv -\un 1 \pmod p$, which implies $\un s = \un t \equiv -\un 1 \pmod p$
  and hence $\un i = \un{p-1}+p\un s \equiv -\un 1 \pmod{p^2}$. By iterating we are left with the case $\un i = -\un 1$, which completes the proof of the claim.

  Finally we show that $c$ is uniquely determined (assuming $\mu \ne 0$). Consider $\un i = -\un 1$ above. Then
  \begin{equation*}
    \mu(\un Z^{-\un1}) = \mu(\un Z^{\un{p-1}}\un Z^{-\un p}) = c' \mu(\un T^{\un{p-1}}\un Z^{-\un p})= c' c^{-1} \mu(\un Z^{-\un1}),
  \end{equation*}
  where $c'$ is the coefficient of $\un T^{\un{p-1}}$ in $\un Z^{\un{p-1}}$. Here, the second equality follows from the
  analysis in the preceding paragraph (the case $\|\un i\| = -f$) that all other intervening terms $\un T^{\un{r'}}\un Z^{-\un p}$ with $\un{r'} \ge \un 0$ and
  $\|\un{r'}\| = \|\un{p-1}\|$ lie in the kernel of $\mu$ (by~\eqref{eq:2f}). The third equality follows from~\eqref{eq:1f'} with $\un i = \un{p-1}$. 
  Hence $c = c'$ is uniquely determined.
\end{proof}

\begin{proof}[Proof {\upshape(}Existence{\upshape)}]
  We define
  \begin{equation}
    \mu(x) \defeq \ve_{-\un 1}(x \prod_j (1+T_j)^{-1})\label{eq:mu-explicit}
  \end{equation}
  for $x \in A$, where $\ve_{-\un 1}(y)$ is the coefficient of $\un Z^{-\un{1}}$ in $y$ for $y \in A$ (expanded in terms of
  the $\un Z^{\un i}$ as in \cite[Rk.\ 3.1.1.3(iii)]{BHHMS2}. Then $\mu \in \Hom^{\cont}_\F(A,\F)$, as $\mu(F_0 A) = \{0\}$.

  By~\eqref{eq:mu1}--\eqref{eq:mu2} it suffices to show that for $\un 1 \le \un i \le \un p$ we have
  \begin{equation}
    \ve_{-\un 1}( \prod_j (1+T_j)^{i_j-1} \vp(x)) = 0 \text{\ \ if $\un i \ne \un p$}\label{eq:eps1}
  \end{equation}
  and
  \begin{equation}
    \ve_{-\un 1}( \prod_j (1+T_j)^{p-1} \vp(x)) = (-1)^{f-1}\ve_{-\un 1}(x).\label{eq:eps2}
  \end{equation}
  (This time we take representatives $n = \smatr{1}{\sum_j i_j\alpha_j}{}{1}$ with $1 \le i_j \le p$.)

  Recalling that we can write 
  \begin{equation}\label{eq:4f}
    Z_j = \sum_i a_{ij} T_i \text{\ \ for some $(a_{ij}) \in \GL_f(\F)$},
  \end{equation}
  we deduce~\eqref{eq:eps1} and reduce~\eqref{eq:eps2} to showing that 
the coefficient of $\un Z^{\un{p-1}}$ in $\un T^{\un{p-1}}$ equals $(-1)^{f-1}$.
  From~\eqref{eq:4f}, by considering the action of $\vp$ and letting $a_i \defeq a_{i0}$, we obtain that
  \begin{equation*}
    Z_j = \sum_i a_i^{p^j} T_i \text{\ \ with $(a_i^{p^j}) \in \GL_f(\F)$}.
  \end{equation*}
  As $a_i^{p^f} = a_i$, the $a_i$ are in the image of $\F_q$ in $\F$ and in fact they form an $\Fp$-basis of $\F_q$. 
  (If not, then $\sum_i \lambda_i a_i = 0$ for some $\lambda_i \in \Fp$ that are not all zero. This implies that
  $\sum_i \lambda_i a_i^{p^j} = 0$ for all $0 \le j \le f-1$, contradicting that $(a_i^{p^j}) \in \GL_f(\F)$.)
  
  Let us now work with formal variables $\un x \defeq (x_i)_{0 \le i \le f-1}$ and $b_i$ ($0 \le i \le f-1$).

  \begin{lem1}\label{lem:coeff}
    The coefficient of $\un x^{\un{p-1}} (= \prod_j x_j^{p-1})$ in $\prod_j (\sum_i b_i^{p^j} x_i)^{p-1}$ equals
    \[\prod_{\un c \in (\F_p^f - \{\un 0\})/\F_p^\times} (\sum_i c_i b_i)^{p-1}
    = (-1)^{(p^f-1)/(p-1)} \prod_{\un c \in \F_p^f - \{\un 0\}} (\sum_i c_i b_i).\]
  \end{lem1}

  (Note that the first product does not depend on the choice of representatives, and for the equality note that $\prod_{x \in \Fp\s} x = -1$.)

  This lemma implies what we want: as the $a_i$ form an $\Fp$-basis of $\F_q$, the lemma (applied with $x_i = T_i$, $b_i = a_i$) shows
  that the coefficient of $\un T^{\un{p-1}}$ in $\un Z^{\un{p-1}}$ equals $-(-1)^{(p^f-1)/(p-1)} = (-1)^{f-1}$, as $\prod_{x
    \in k\s} x = -1$.

  To prove Lemma~\ref{lem:coeff}, we use the following.

  \begin{sublem1}\label{sublem:coeff}
    Suppose $h \in \F[x_0,\dots,x_{f-1}]$. Then the coefficient of $\un x^{\un{p-1}}$ in $h$ is invariant under any linear change of variables over $\Fp$, i.e.\
    is equal to the coefficient of $\un y^{\un{p-1}}$ in $h$ if $\un x$ and $\un y$ are related by an element $\gamma$ of $\GL_f(\Fp)$.
  \end{sublem1}
  
  (This is presumably well known. For the proof we may assume that $h$ is a monomial and that $\gamma$ is an elementary transformation, in which case
  it follows from the facts that $\Fp^\times$ is of order $p-1$ and that $\binom{r}{p-1} = 0$ for $p \le r \le 2p-2$.)

  Let $C$ denote the coefficient of $\un x^{\un{p-1}} (= \prod_j x_j^{p-1})$ in $\prod_j (\sum_i b_i^{p^j} x_i)^{p-1}$.
  Then $C \in \F[b_0,\dots,b_{f-1}]$ is a homogeneous polynomial of degree $p^f-1$, which is clearly divisible by $\un b^{\un{p-1}}$.
  For any linear change of variables $x_i = \sum_i \lambda_{ij} y_j$ with $\lambda_{ij} \in \Fp$,
  Sublemma~\ref{sublem:coeff} then implies that $\prod_j (\sum_i b_i \lambda_{ij})^{p-1}$ divides $C$.
  In particular, $(\sum_i c_i b_i)^{p-1}$ divides $C$ for each $c \in (\F_p^f - \{\un 0\})/\Fp^\times$.
  But the product of such polynomials is already of degree $p^f-1$ and they are pairwise relatively prime, hence we are done by
  remarking that the coefficient of $\prod_i b_i^{p^i(p-1)}$ is the same on both sides. 
\end{proof}

\begin{rem1}\label{rk:mu-a}
  Fix $\mu \ne 0$ as in Proposition~\ref{prop:mu}.
  By uniqueness we must have $\mu \circ a^{-1} \in \F \mu$ for any $a \in \cO_K\s$.
  But it is easy to compute the scalar: by applying the explicit formula~\eqref{eq:mu-explicit} to the element $\prod_j (1+T_j)^{p-1}\un Z^{-\un{1}}$
  we obtain 
  \[\mu \circ a^{-1} = N_{\F_q/\Fp}(\o a) \mu \quad\forall\ a \in \cO_K\s.\]
\end{rem1}

Suppose $\mu \in \Hom^{\cont}_\F(A,\F)$ is nonzero such that $\mu \circ \psi = (-1)^{f-1}\mu$ and $D$ is an \'etale $(\vp,
\cO_K\s)$-module over $A$. Then composition with $\mu$ induces an $A$-linear map
\begin{equation}\label{eq:mu-star}
  \mu_* : \Hom_A(D,A) \to \Hom^{\cont}_\F(D,\F).
\end{equation}
Recall from Lemma~\ref{lm:dual-etale} that $\Hom_A(D,A)$ is naturally an \'etale $(\vp,\cO_K\s)$-module. The following
lemma will allow us to calculate this structure on the level of $\Hom^{\cont}_\F(D,\F)$.

\begin{lem1}\label{lm:mu-star}\ 
\begin{enumerate}
  \item The map $\mu_*$ in~\eqref{eq:mu-star} is injective.
  \item We have $\mu_*(\vp(h)) = (-1)^{f-1} \mu_*(h) \circ \psi$.
  \item We have $\mu_*(a(h)) = N_{\F_q/\Fp}(\o a)^{-1} \mu_*(h) \circ a^{-1}$ for $a \in \cO_K\s$.
  \end{enumerate}
\end{lem1}
\begin{proof}
  Part (iii) follows immediately from Remark~\ref{rk:mu-a}.

  For (i) we can reduce to $D = A$ by using that $D$ is finite projective. Observe then that the kernel of $\mu_*$
  is an $\cO_K\s$-stable ideal of $A$ by (iii); but by \cite[Cor.\ 3.1.1.7]{BHHMS2} it is zero, as it cannot be all of $A$.
  (Alternatively part (i) also follows from the explicit
  formula for $\mu$ above.)

  Part (ii) follows from the explicit formula (\ref{eq:6f}) for $\varphi$ on $\Hom_A(D,A)$ as well as the two conditions 
  at the beginning of the proof of Proposition~\ref{prop:mu}.
\end{proof}

We make part (ii) more explicit. Suppose that $h \in \Hom_\F^{\cont}(D_A(\pi),\F)$ corresponds to a sequence $(x_k)_{k \ge 0}$
as in Proposition~\ref{prop:F-dual}. Then $(-1)^{f-1} h \circ \psi$ corresponds to a sequence $(x'_k)_{k \ge 0}$ determined by
the relation
\begin{equation}\label{eq:7f}
  x'_{pk} = (-1)^{f-1} \smatr{p}{}{}{1} x_k,
\end{equation}
since $\psi \circ \un Y^{-\un{pk}} = \un Y^{-\un k} \circ \psi$ on $D_A(\pi)$.

\begin{lem1}\label{lm:descent-criterion}
  Suppose that $D$ is a finite projective $A$-module.
  Then the image of $\mu_* : \Hom_A(D,A) \into \Hom^{\cont}_\F(D,\F)$ consists precisely of all continuous $\F$-linear maps $h : D \to \F$ such that
  for all $M \in \Z$ and all $x \in D$ the set $X'_M \defeq \{\un i \in \Z^f : h(\un Z^{\un i} x) \ne 0, \|\un i \| = M\}$ is finite.

  Equivalently, the image of $\mu_* : \Hom_A(D,A) \into \Hom^{\cont}_\F(D,\F)$ consists precisely of all continuous $\F$-linear maps $h : D \to \F$ such that
  for all $M \in \Z$ and all $x \in D$ the set $X_M \defeq \{\un i \in \Z^f : h(\un Y^{\un i} x) \ne 0, \|\un i \| = M\}$ is finite.
\end{lem1}
\begin{proof}
  For the first part it is easy to reduce to the case where $D = A$, using the compatibility of $\mu_*$ with direct sums
  $D = D_1 \oplus D_2$.  If $h = \mu_*(a)$ for some $a \in A$ and $x \in A$, then we write
  $ax \prod_j (1+T_j)^{-1} = \sum_{\un i} \lambda_{\un i} \un Z^{\un i}$ for $\lambda_{\un i} \in \F$.  Then
  $h(\un Z^{\un i} x) = \lambda_{-\un i-\un 1}$ (by the explicit formula for $\mu_*$ in \S\ref{sec:finding-mu}), so
  $h(\un Z^{\un i} x) \ne 0$ can only happen for finitely many $\un i$ of any fixed degree $\|\un i \| = M$.
  Conversely, if $h : A \to \F$ is continuous such that for all $M \in \Z$ the set
  $\{\un i : h(\un Z^{\un i}) \ne 0, \|\un i \| = M\}$ is finite, then by continuity of $h$ and the finiteness assumption it
  follows that $a \defeq (\prod_j (1+T_j))\sum_{\un i} h(\un Z^{\un i})\un Z^{-\un i-\un 1} \in A$, and by the explicit formula
  for $\mu_*$ we have $\mu_*(a) = h$.

  To justify the second part, recall that $Y_j$, $Z_j \in \F\bbra{N_0}$ with $\gr(Y_j) = \gr(Z_j)$, so
  $Z_j = Y_j \sum_{d=0}^\infty F_{d,j}$, where $Y_j F_{d,j}$ is a homogeneous polynomial in $Y_0,\dots,Y_{f-1}$ of degree $d+1$
  and $F_{0,j} = 1$.
  Define the subring
  \begin{equation*}
    A_0 \defeq \Big\{ \sum_{d=0}^{\infty} \frac{F_d}{\un Y^{\un{d}}} : \text{$F_d$ a homog.\ poly.\ in $Y_0,\dots,Y_{f-1}$ of degree $d(f+1)$} \Big\}
  \end{equation*}
  of $A$ with maximal ideal $\fm_0$ defined by the condition $F_0 = 0$. The above observation then implies that 
  $Z_j \in Y_j (1+\fm_0)$ for any $j$, hence 
  \begin{equation}\label{eq:15f}
    \un Z^{\un i} \in \un Y^{\un i} (1+\fm_0) \quad\forall\ \un i \in \Z^f.
  \end{equation}
  Also note that
  \begin{equation}\label{eq:16f}
    A_0 = \bigg\{ \sum_{\un k \in \Z^f;\ k_j \ge -\|\un k\|\; \forall j} \lambda_{\un k} \un Y^{\un k} : \lambda_{\un k} \in \F \bigg\}
  \end{equation}
  and that the condition $k_j \ge -\|\un k\|$ for all $j$ implies $k_j \le f \|\un k\|$ for all $j$ (and $\|\un k\| \ge 0$),
  so that there are only finitely many terms of any fixed degree.

  Fix now $x \in D$ and suppose that the set $X_N = \{\un i \in \Z^f : h(\un Y^{\un i} x) \ne 0, \|\un i\| = N\}$ is finite for any $N \in \Z$.
  By continuity of $h$ we know that $h(\un Y^{\un i} x) = 0$ for all $\|\un i\| \ge e$ (some $e \in \Z$).
  Fix any $M \in \Z$ and suppose that $h(\un Z^{\un i} x) \ne 0$ and $\|\un i\| = M$. 
  By equations \eqref{eq:15f}--\eqref{eq:16f} we get that $h(\un Y^{\un i + \un k}x) \ne 0$ for some $\un k \in \Z^f$ such that $k_j \ge -\|\un k\|$ for all $j$.
  In particular, $\|\un i\| \le \|\un i\| + \|\un k\| < e$, so
  \begin{equation*}
    X'_M \subset \bigcup_{\stackrel{k_j \ge -\|\un k\|\; \forall j}{0 \le \|\un k\| < e-M}} \big(X_{M + \|\un k\|} - \un k\big),
  \end{equation*}
  a finite union of finite sets. The converse direction follows by reversing the roles of $Y_j$ and $Z_j$.
\end{proof}

\subsection{Some combinatorial lemmas and computations}
\label{ss:degree}

We \ give \ several \ technical \ but \ important \ lemmas \ (some \ generalizing \ results \ in \cite[\S 3.2]{BHHMS2}) involving the combined action of $\underline{Y}^{\underline{k}}$ (for some $\underline{k}\in \Z_{\geq 0}^f$) and $\smatr{p}{0}{0}{1}$ in a representation $\pi$ as at the end of \S\ref{conjecture}.

We recall some notation and results from \cite{BHHMS2}. 
Let $H\defeq \begin{pmatrix}\F_q^\times&0\\0&\F_q^\times\end{pmatrix} \le \GL_2(\F_q)$.
As in \cite{BP} we write $(s_0,s_1,\dots, s_{f-1})\otimes\eta$ for the Serre weight
\[
\mathrm{Sym}^{s_0}\F^2\otimes_{\F}(\mathrm{Sym}^{s_1}\F^2)^{\rm Fr}\otimes\cdots\otimes_{\F}(\mathrm{Sym}^{s_{f-1}}\F^2)^{{\rm Fr}^{f-1}}\otimes_{\F}\eta\circ\det,\] 
where the $s_i$ are integers between $0$ and $p-1$, $\eta$ is a character $\F_q^{\times}\rightarrow \F^{\times}$ and $\GL_2(\Fq)$ acts on $(\mathrm{Sym}^{s_i}\F^2)^{{\rm Fr}^i}$ via $\sigma_i:\F_q\hookrightarrow \F$. 
We fix $\rhobar$ as at the end of \S\ref{conjecture}. We identify $W(\brho)$ with the subsets of $\{0,\dots,f-1\}$ as in \cite[\S2]{breuil-IL} and let $J_{\sigma}$ be the subset associated to $\sigma$.  More precisely,  $W(\brho)$ is exactly the set of Serre weights of the form 
\[(\lambda_0(r_0),\dots,\lambda_{f-1}(r_{f-1}))\otimes{\det}^{e(\lambda)(r_0,\dots,r_{f-1})}\]
where $\lambda\in \mathcal{ID}(x_0,\dots,x_{f-1})$ (resp.\ $\lambda\in \mathcal{RD}(x_0,\dots,x_{f-1})$) if $\brho$ is irreducible (resp.\ $\brho$ is reducible), see \cite[\S11]{BP} or \cite[\S2]{breuil-IL}.  If $\sigma\in W(\brho)$ corresponds to $\lambda$, then  
we have $\lambda_j(x_j)\in\{p-2-x_j,p-3-x_j\}$ if and only if $j\in J_{\sigma}$ when $j>0$ or $\brho$ is reducible, $\lambda_0(x_0)\in \{p-2-x_0,p-1-x_0\}$ if and only if $0\in J_{\sigma}$ when $\brho$ is irreducible. }

Let $\sigma\in W(\brho)$. Let $\delta(\sigma)\defeq\delta_{\rm irr}(\sigma)$ if $\brho$ is irreducible and $\delta(\sigma)\defeq\delta_{\rm red}(\sigma)$ if $\brho$ is reducible the Serre weights defined in \cite[\S5]{breuil-IL}. 
 Then $\delta(\sigma)\in W(\brho)$ and we have the following    
 explicit description of $J_{\delta(\sigma)}$ (see \cite[\S5]{breuil-IL}): 
 \[\begin{array}{rlll}  j\in J_{\delta(\sigma)}, j<f-1~ (\mathrm{resp}.~  f-1\in J_{\delta(\sigma)}) &\!\!\!\Longleftrightarrow\!\!\!& j+1\in J_{\sigma}~(\mathrm{resp}.~0\notin J_{\sigma})& \!\!\mathrm{if}\ \delta=\delta_{\rm irr}\\
 j\in J_{\delta(\sigma)}&\!\!\!\Longleftrightarrow \!\!\!& j+1\in J_{\sigma}& \!\!\mathrm{if}\ \delta=\delta_{\rm red}.
 \end{array} \]
We fix a nonzero vector $v_{\sigma}\in\sigma^{N_0}$, and let
$\chi_{\sigma}:H\ra \F^{\times}$ be the
$H$-eigencharacter
of $v_{\sigma}$. Let $\chi_{\sigma}^s:H\ra \F^{\times}$ denote the conjugate of $\chi_{\sigma}$ by $\smatr{0}110$. As in \cite[\S2]{BP} we identify the irreducible constituents of $\Ind_I^{\GL_2(\cO_K)}(\chi_{\sigma}^s)$ with the subsets of $\{0,\dots,f-1\}$ (for example $\emptyset$ corresponds to the socle $\sigma$ of $\Ind_I^{\GL_2(\cO_K)}(\chi_{\sigma}^s)$). We know that $\delta(\sigma)$ occurs in $\Ind_{I}^{\GL_2(\cO_K)}(\chi_{\sigma}^s)$ and we denote by $J^{\rm max}(\sigma)\subseteq \{0,\dots,f-1\}$ the associated subset. Precisely, using \cite[Lemma 2.7]{BP} one checks that 
\begin{equation*}J^{\rm max}(\sigma)=(J_{\sigma}\cup J_{\delta(\sigma)})\setminus (J_{\sigma}\cap J_{\delta(\sigma)}).
\end{equation*}
By \cite[Lemma 3.2.3.2]{BHHMS2}, we have $|J^{\rm max}(\sigma)|=|J^{\rm max}(\delta(\sigma))|$. As a consequence, the quantity 
\[m\defeq|J^{\rm max}(\sigma)|\in\{0,\dots,f-1\}\]
depends only on the orbit of $\sigma$.  By the proof of \cite[Lemma 19.5]{BP}, the vector $\smatr{0}1p0 v_{\sigma}$ generates a $\GL_2(\cO_K)$-subrepresentation of $\pi$ isomorphic to the unique quotient of $\Ind_I^{\GL_2(\cO_K)}(\chi_{\sigma}^s)$ with irreducible socle parametrized by $J^{\rm max}(\sigma)$, which in particular yields an embedding of $\delta(\sigma)$ in $\mathrm{soc}_{\GL_2(\cO_K)}(\pi)$.

Write \begin{equation*}\sigma=(s_0,\dots,s_{f-1})\otimes\eta, \ \ \delta(\sigma)=(s_0',\dots,s_{f-1}')\otimes\eta'.\end{equation*}
Define $\underline c \in \Z^f$ by $c_{j}\defeq s_{j}'$ if $j\in J^{\rm max}(\sigma)$, and $c_{j}\defeq p-1$ otherwise.

The following lemma explicitly determines $s_{j-1}'$ and $c_{j-1}$ in terms of  $s_{j}$.  We remark that if $f=1$ and $\brho$ is irreducible, some formulas need to be modified, e.g.\ Lemma \ref{lem:s-s'}(i). But the main result (Theorem \ref{mainbis}) is known in this case, so it is harmless to exclude it.

\begin{lem1}\label{lem:s-s'}\ 
\begin{enumerate}
  \item Assume that $\brho$ is irreducible and $f\geq 2$. Then{\begin{center}
        \begin{tabular}{ |c||c|c|c|c|c|c|}
          \hline  
          $s_0$& $s_{f-1}'$ &     & $c_{f-1}$     \\
          \hline \hline
          $r_0$ & $p-2-r_{f-1}$&$f-1\in J^{\rm max}(\sigma)$ & $p-2-r_{f-1}$  \\
          \hline
          $r_0-1$&$p-3-r_{f-1}$&$f-1 \notin J^{\rm max}(\sigma)$  & $p-1$\\
          \hline
          $p-2-r_0$&$r_{f-1}$&$f-1 \notin J^{\rm max}(\sigma)$  &$p-1$\\
          \hline
          $p-1-r_0$&$r_{f-1}+1$&$f-1 \in J^{\rm max}(\sigma)$  &$r_{f-1}+1$ \\
          \hline
          \end{tabular}
        \end{center}}
      \noindent while if $1\leq j\leq f-1$ we have {\begin{center}
          \begin{tabular}{ |c||c|c|c|c|c|c|}
            \hline  
            $s_j$&  $s_{j-1}',~ j=1$& $s_{j-1}', ~j>1$ &    &$c_{j-1}$     \\
            \hline \hline
            $r_j$  & $r_0-1$&$r_{j-1}$&$j-1\notin J^{\rm max}(\sigma)$ &$p-1$  \\
            \hline
            $r_j+1$ &$r_0$&$r_{j-1}+1$&$j-1\in J^{\rm max}(\sigma)$  & $s_{j-1}'$\\
            \hline
            $p-2-r_j$ &$p-1-r_0$&$p-2-r_{j-1}$&$j-1\in J^{\rm max}(\sigma)$  &$s_{j-1}'$\\
            \hline
            $p-3-r_j$ &$p-2-r_0$&$p-3-r_{j-1}$&$j-1\notin J^{\rm max}(\sigma)$  &$p-1$\\
            \hline
            \end{tabular}
          \end{center}}
      \item 
        Assume that $\brho$ is {\upshape(}split{\upshape)} reducible. Then for any $0\leq j\leq f-1$ we have {\begin{center}
            \begin{tabular}{ |c||c|c|c|c|c|c|}
              \hline  
              $s_j$& $s_{j-1}'$ &    &$c_{j-1}$      \\
              \hline \hline
              $r_j$  &$r_{j-1}$&  $j-1\notin J^{\rm max}(\sigma)$   & $p-1$\\
              \hline
              $r_j+1$ &$r_{j-1}+1$ &$j-1\in J^{\rm max}(\sigma)$ &$r_{j-1}+1$ \\
              \hline
              $p-2-r_j$  &$p-2-r_{j-1}$&$j-1\in J^{\rm max}(\sigma)$ & $p-2-r_{j-1}$\\
              \hline
              $p-3-r_j$  &$p-3-r_{j-1}$&$j-1\notin J^{\rm max}(\sigma)$&  $p-1$\\
              \hline
              \end{tabular}
            \end{center}}
        \end{enumerate}
\end{lem1}
 
\begin{proof}
This is an easy exercise using the relation between $J_{\sigma}$ and $J_{\delta(\sigma)}$.  
Note also that $j \not\in J_\sigma$ if and only if $s_{j+1} \in \{r_{j+1},p-2-r_{j+1}\}$.
\end{proof}

\begin{rem1}\label{rk:s-s'}
  Strictly speaking, we should state Lemma~\ref{lem:s-s'} in terms of $\lambda,\lambda'$, which \ are \ the \ elements \ in \ $\mathcal{ID}(x_0,\dots,x_{f-1})$ \ or \ $\mathcal{RD}(x_0,\dots,x_{f-1})$ \ (depending \ on \ whether $\brho$ is irreducible or reducible) corresponding to $\sigma,\delta(\sigma)$ respectively.
  Lemma \ref{lem:s-s'} determines $\lambda_{j-1}'(x_{j-1})$ in terms of $\lambda_{j}(x_{j})$, not $s'_{j-1}$ in terms of $s_j$ (because ambiguities arise when $r_0 = \frac{p-1}2$ in the first table, and when $r_j = \frac{p-3}2$ in the second and third tables.) The same comment applies to Lemma~\ref{lem:a-d'-sigma}.
\end{rem1}

\begin{lem1}\label{lem:BP}
The vector $\underline{Y}^{\underline c}\smatr{p}{}{}{1}(v_{\sigma})$ 
spans $\delta(\sigma)^{N_0}$ as an $\F$-vector space.
Hence there is a unique scalar $\mu_\sigma \in \F\s$ such that
\begin{equation}\label{eq:x-delta}
v_{\delta(\sigma)}=\mu_\sigma\cdot \underline{Y}^{\underline c}\smatr{p}{}{}{1}(v_{\sigma})\end{equation}
\end{lem1}
\begin{proof}
This is \cite[Prop.~3.2.3.1(i)]{BHHMS2}.
\end{proof}

By \cite[Lemma 3.2.2.6(ii)]{BHHMS2}, if $\underline{0}\leq \underline{i}\leq \underline{s}$, there is a unique $H$-eigenvector $\underline{Y}^{-\underline{i}}v_{\sigma}\in \sigma$ that is sent by $\underline{Y}^{\underline{i}}$ to $v_{\sigma}$. The  following result is a generalization  of \cite[Lemma 3.2.3.5]{BHHMS2}.

\begin{lem1}\label{lem:BP-florian}
Assume \ $m>0$. \ Let \ $\underline{k},\underline{i}\in\Z_{\geq 0}^f$ \ such \ that \ $\|\underline{i}\|\leq f-1$ \ and $\underline{Y}^{\underline{k}}\smatr{p}{}{}{1}(\underline{Y}^{-\underline{i}}v_{\sigma})\neq 0$.
\begin{enumerate}
\item We have
\begin{equation*}\|\underline{k}\|\leq p\|\underline{i}\|+\|\underline c\|.
\end{equation*}
\item If $\|\underline{k}\|\geq p\|\underline{i}\|-(f-1) +\|\underline c\|$, then \[\mu_\sigma\cdot \underline{Y}^{\underline{k}}\smatr{p}{}{}{1}(\underline{Y}^{-\underline{i}}v_{\sigma})=\underline{Y}^{-\underline{\ell}}v_{\delta(\sigma)} \in \delta(\sigma) \]
for some $\underline{\ell}\geq \underline{0}$ with $\|\underline{\ell}\|\leq f-1$. 
More precisely, $\ell_j= i_{j+1}p + c_j-k_j$ for all $j$. 
\end{enumerate}
\end{lem1}
\begin{proof}
 Before starting the proof, we first remark that Lemma 3.2.3.3 and Lemma 3.2.3.4 of  \cite{BHHMS2} remain true if we replace the assumption $\|\underline{i}\|\leq m-1$ by $\|\underline{i}\|\leq f-1$ in the  statements. Indeed, for Lemma 3.2.3.3,  this new assumption $\|\underline{i}\|\leq f-1$ implies  $i_j\leq f-1$, and so \[2i_j+1\leq 2f-1\leq s_j\]  for all $j$ ($s_j$ is denoted by $t_j$ in \emph{loc.~cit.}) by the genericity assumption.   Hence, \cite[Prop.~6.2.2]{BHHMS1} still applies and the rest of the proof of Lemma 3.2.3.3 works without change. The proof of Lemma 3.2.3.4 of \cite{BHHMS2} also works through, because one checks that besides the citation to Lemma 3.2.3.3 the condition $\|\underline{i}\|\leq m-1$ is only used to   deduce   $\|\underline{i}\|\leq f-1$.

Now we prove the lemma, following the proof of \cite[Lemma 3.2.3.5]{BHHMS2}. We first prove by induction on $\|\underline{i}\|\leq f-1$ the following fact: if 
\[\|\underline{k}\|\geq p\|\underline{i}\|-(f-1)+\|\underline c\|\defeq B\] 
and $\underline{Y}^{\underline{k}}\smatr{p}{}{}{1}(\underline{Y}^{-\underline{i}}v_{\sigma})\neq0$, then
$\underline{Y}^{\underline{k}}\smatr{p}{}{}{1}(\underline{Y}^{-\underline{i}}v_{\sigma})=\underline{Y}^{\underline{k}'}\smatr{p}{}{}{1}(v_{\sigma})$ for some $\underline{k}'\in \Z_{\geq 0}^f$ such that $k_j'=k_j-i_{j+1}p$ for all $j$.   This is trivial if $\underline{i}=\underline{0}$, so we can assume $\underline{i}\neq \underline{0}$. Moreover, as in \emph{loc.~cit.}, by induction we are reduced to the case $k_j<p $ for all $j$.  We make this assumption and derive below a contradiction (so this case cannot happen). 

Define a set $J$ as in \emph{loc.~cit.}, i.e.\
\begin{equation}\label{eq:J}
J\defeq\{j\in J^{\rm max}(\sigma), i_{j+1}=0\}.
\end{equation}
As in \emph{loc.~cit.} we have 
\[\|\underline{k}\|\leq (p-1)(f-|J|)+\sum_{j\in J}(s_j'+2i_j)+|J\setminus (J^{\rm max}(\sigma)+1)|\defeq A\]
and to get a  contradiction it is enough to show $A< B$, which is equivalent to 
\begin{equation}\label{eq:A<B}
mp +|J\setminus (J^{\rm max}(\sigma)+1)|< (p-2)\|\underline{i}\|  +(p-1)|J|+C+D, \end{equation}
where \[C\defeq m-(f-1),\ \  D\defeq 2\sum_{j\notin J}i_j+\sum_{j\in J^{\rm max}(\sigma)\setminus J}s_j'. \]
We have the following two cases.
\begin{itemize}
\item If $|J^{\rm max}(\sigma)\setminus J|>0$, then as in \emph{loc.~cit.} $m\leq \|\underline{i}\|+|J|$, hence \eqref{eq:A<B} is implied by 
\[mp+|J\setminus (J^{\rm max}(\sigma)+1)|<(p-2)\|\underline{i}\|+(p-2)(m-\|\underline{i}\|)+|J|+C+D,\]
or equivalently
\[m+(f-1)+|J\setminus (J^{\rm max}(\sigma)+1)|<|J|+D.\]
This is slightly stronger than (140) of \cite{BHHMS2}, but one checks that the argument in \emph{loc.~cit.} still  allows to conclude. 
\item If $J^{\rm max}(\sigma)=J$, then as in \emph{loc.~cit.} we have $|J\setminus(J^{\rm max}(\sigma)+1)|\leq f-m$ and $|J|=m$, and \eqref{eq:A<B} is implied by 
\[mp+(f-m)<(p-2)\|\underline{i}\|+(p-1)m+C+D\]
or equivalently
\[2f-1<(p-2)\|\underline{i}\|+m+D.\]
As $\|\underline{i}\| > 0$ and $D \ge 0$, the last inequality holds
  by our genericity condition (i.e.\ $p>4f$).
\end{itemize} 
This proves the desired fact. The rest of the proof is the same as the proof of  \cite[Lemma 3.2.3.5]{BHHMS2} and we omit the details.
(Several times $f-m = (f-1)-(m-1)$ has to be added or subtracted from expressions in the last three paragraphs of the proof in \emph{loc.~cit.} to account for the weaker lower bound in Lemma~\ref{lem:BP-florian}(ii).)
\end{proof}

\begin{rem1}\label{rem:BP-florian}
Taking $\underline{i}=\underline{0}$ in Lemma \ref{lem:BP-florian}, we get the following. If $\underline{Y}^{\underline{k}}\smatr{p}{}{}{1}(v_\sigma)\neq0$ for some $\underline{k}\in\Z_{\geq 0}^f$ and if \[\|\underline{k}\|\geq \|\underline{c}\| -(f-1),\] then
$\mu_\sigma\cdot \underline{Y}^{\underline{k}}\smatr{p}{}{}{1}(v_{\sigma})=\underline{Y}^{-\underline{\ell}}v_{\delta(\sigma)} \in \delta(\sigma)$ for some $\|\underline{\ell}\|\leq f-1$. 
More precisely, $\underline{\ell}=\underline{c}-\underline{k}$.
\end{rem1}

We will need the following analogue of Lemma \ref{lem:BP-florian}.  

\begin{lem1}\label{lem:BP-florian-v}
Assume $m>0$. Let $\underline{i}\in\Z_{\geq 0}^f$ such that $\|\underline{i}\|\leq f-1$. Let $\underline{k}\in \Z_{\geq0}^{f}$ and assume that there exists  $0\leq j_0\leq f-1$ such that 
\begin{enumerate}
\item[\upshape(a)] $k_{j_0}\leq p(i_{j_0+1}-1)$ (hence $i_{j_0+1}\geq 1$);
\item[\upshape(b)]  
$\|\underline{k}\|> p\|\underline{i}\| + \|\underline{c}\| - c_{j_0}$.
\end{enumerate}
Then $\underline{Y}^{\underline{k}}\smatr{p}{}{}{1}(\underline{Y}^{-\underline{i}}v_{\sigma})=0$. 
\end{lem1}
 \begin{proof}
Assume for a contradiction that $\underline{Y}^{\underline{k}}\smatr{p}{}{}{1}(\underline{Y}^{-\underline{i}}v_{\sigma})\neq0$.
As in the proof  of \cite[Lemma 3.2.3.5]{BHHMS2}, by induction we are reduced to the case  $k_{j}<p$ for all $j$; we make this assumption from now on.
Note $\|\underline{c}\| = \sum_{j\in J^{\rm max}(\sigma)}s_j'+\sum_{j\notin J^{\rm max}(\sigma)}(p-1)$.

Let $J$ be the set defined by \eqref{eq:J}.   
 Then by (a) we have $j_0\notin J$. As explained in the proof of Lemma \ref{lem:BP-florian},   \cite[Lemma 3.2.3.4]{BHHMS2} still applies, and we get (see the fourth paragraph of the proof of Lemma 3.2.3.5 of \emph{loc.~cit.})
\[\sum_{j\neq j_0}k_j\leq (f-1-|J|)(p-1)+\sum_{j\in J}(s_j'+2i_j)+|J\setminus (J^{\rm max}(\sigma)+1)|\defeq A.\]
On the other hand, letting $\gamma \defeq 1$ if $i_{j_0+1} > 1$ and $\gamma \defeq 0$ if $i_{j_0+1} = 1$ we see that $k_{j_0} \le (p-1)\gamma$ (using (a) when $\gamma = 0$),
which together with condition (b) implies
\begin{equation*}
\sum_{j\neq j_0}k_j>p\|\underline{i}\|-(p-1)\gamma+\sum_{j\in J^{\rm max}(\sigma)}s_j'+\sum_{j\notin J^{\rm max}(\sigma)}(p-1)-c_{j_0}\defeq B.
\end{equation*}
To get a contradiction it is enough to show $A\leq B$.

A computation shows that $A\leq B$   is equivalent to
\begin{equation}\label{eq:A<Bagain}
  \begin{aligned}
    mp+|J\setminus (J^{\rm max}(\sigma)+1)|&\leq (p-2)\|\underline{i}\| + (p-1)|J|+(p-1)(1-\gamma)+C+D\\
    &= (p-2)(\|\underline{i}\|+|J|+1-\gamma)+|J|+1-\gamma+C+D,
  \end{aligned}
\end{equation}
where \[C\defeq m-c_{j_0},\ \ D\defeq 2\sum_{j\notin J}i_j+\sum_{j\in J^{\rm max}(\sigma)\setminus J}s_j'.\] 
If $j\in J^{\rm max}(\sigma)\setminus J$, then $i_{j+1}>0$, so we obtain
\begin{align*}
|J^{\rm max}(\sigma)\setminus J| &\leq \sum_{J^{\rm max}(\sigma)\setminus (J \cup \{j_0\})} i_{j+1} + 1 = \bigg(\sum_{J^{\rm max}(\sigma)\setminus (J \cup \{j_0\})} i_{j+1}\bigg) + i_{j_0+1} + (1-i_{j_0+1}) \\ &\le \|\underline{i}\|+(1-i_{j_0+1}).
\end{align*}
As $|J^{\rm max}(\sigma)\setminus J|=m-|J|$ and $i_{j_0+1} \ge \gamma+1$, this means 
\[m\leq \|\underline{i}\|+|J|+(1-i_{j_0+1})\leq \|\underline{i}\|+|J|-\gamma.\]
Thus, to show \eqref{eq:A<Bagain} it is enough to show 
\[mp +|J\setminus(J^{\rm max}(\sigma)+1)|\leq (p-2)(m+1)+|J|+1-\gamma+C+D\]
or equivalently
\[2m+|J\setminus(J^{\rm max}(\sigma)+1)|\leq |J|+(p-1-\gamma+C)+D.\]
 If $|J^{\rm max}(\sigma)\setminus J|>0$, then it is true by  \cite[Eq.~(140)]{BHHMS2} (and using  $p-2+m-c_{j_0}\geq 0$ as $m\geq 1$). If $J^{\rm max}(\sigma)=J$, then again as in \emph{loc.~cit.}, we have $|J\setminus (J^{\rm max}(\sigma)+1)|\leq f-m$ and $|J|=m$, and  \eqref{eq:A<Bagain} is implied by
 \[mp+f-m\leq (p-2)\|\underline{i}\|+(p-1)(m+1-\gamma)+ (m-c_{j_0})+D,\]
 equivalently, 
 \begin{align*}
f&\leq (p-2)\|\underline{i}\|+(p-1)(1-\gamma)+(m-c_{j_0})+D\\
&= (p-2)(\|i\|-\gamma) + (p-1-c_{j_0}) +(m-\gamma) +D.
 \end{align*}
This is true by our genericity condition: indeed, as $\|\underline{i}\|\geq \gamma+1$, $m \ge 1$, $c_{j_0} \le p-1$, and $D\geq 0$, the above inequality is implied by 
$f \le p-2 \le p-1-\gamma$.
\end{proof}

Now, fix $\sigma\in W(\brho)$ and define $\sigma_i\in W(\brho)$ inductively by $\sigma_0\defeq \sigma$ and $\sigma_i\defeq \delta(\sigma_{i-1})$ for $i\ge 1$. Let $d\geq 1$ be the smallest integer such that $\sigma_{d}\cong \sigma_0$. For convenience, if $i\geq0$ we set $\sigma_i\defeq \sigma_{i'}$, where $i'\in\{0,\dots,d-1\}$ is the unique integer such that $i\equiv i'\ (\mathrm{mod}\ d)$. Write 
\begin{equation*}\sigma_i=(s_0^{(i)},\dots,s_{f-1}^{(i)})\otimes\eta_i.
\end{equation*}
To make the notation consistent, we also write $s_j=s_j^{(0)}$.

For convenience, we introduce the following notation. For $i\geq1$, define $\underline{c_i^\sigma}\in\Z_{\geq 0}^{f}$ by 
\begin{equation}
\label{eq:c_sigma}
c_{i,j}^\sigma\defeq\begin{cases}
s_j^{(i)}&\text{if $j\in J^{\rm max}(\sigma_{i-1})$,}\\
p-1&\text{otherwise}
\end{cases}
\end{equation}
(in particular $\underline{0}\leq \underline{c_{i}^\sigma}\leq \underline{p-1}$).
Define a shift function $\delta : \Z^{f} \to \Z^f$ by setting  \[\delta(\underline{i})_j\defeq i_{j+1}, \ \ \underline{i}=(i_j)\in\Z^f.\]
Note that $\delta$ does not change $\|\cdot\|$ and that 
$\underline{Y}^{p\delta(\underline i)} \smatr{p}{}{}{1} = \smatr{p}{}{}{1} \underline{Y}^{\underline i}.$ 
We inductively define $\underline{a_n^\sigma}\in\Z_{\geq 0}^{f}$  for $n\geq 0$ as follows: $\underline{a_0^\sigma}\defeq \underline{0}$ and for $n\geq 1$, 
\begin{equation}\label{eq:abc}
\underline{a_n^\sigma}\defeq p\delta(\underline{a_{n-1}^\sigma})+\underline{c_{n}^\sigma}.
\end{equation}
For example,  $a_{1,j}^\sigma=c_{1,j}^\sigma=s_j^{(1)}$ if $j\in J^{\rm max}(\sigma)$ and $a_{1,j}^\sigma=c_{1,j}^\sigma=p-1$ if $j\notin J^{\rm max}(\sigma)$.

The following result determines $\un{a_{d'}^{\sigma}}$ in terms of the $s_j$ (where $d' \defeq df$).
For $0\leq j\leq f-1$, recall that $h_j = r_j+1$ and define 
\begin{equation}\label{eq:h[i]}
h^{[j]}\defeq h_j+ph_{j+1}+\cdots +p^{f-1-j}h_{f-1}
\end{equation}
(thus $h^{[0]}=h$).   

\begin{lem1}\label{lem:a-d'-sigma}\ 
\begin{enumerate}
  \item Assume that $\brho$ is irreducible and $f\geq 2$. Then {\begin{center}
        \begin{tabular}{ |c||c|c|c|c|c|c|}
          \hline  
          $s_0$& $r_0$ & $r_0-1$& $p-2-r_0$&$p-1-r_0$          \\
          \hline  
          $\frac{a_{d',0}^{\sigma}}{1-p^{d'}}$ & $-1+\frac{h}{1+q}$&$-1$&$-1$&$-\frac{h}{1+q}$   \\
          \hline
          \end{tabular}
        \end{center}}
      \noindent while if $1\leq j\leq f-1$ we have {\begin{center}
          \begin{tabular}{ |c||c|c|c|c|c|c|}
            \hline  
            $s_j$& $r_j$ & $r_j+1$& $p-2-r_j$&$p-3-r_j$          \\
            \hline  
            $\frac{a_{d',j}^{\sigma}}{1-p^{d'}}$ & $-1$&$h^{[j]}-\frac{hp^{f-j}}{1+q}$&$-1-h^{[j]}+\frac{hp^{f-j}}{1+q}$&$-1$   \\
            \hline
            \end{tabular}
          \end{center}}
      \item Assume that $\brho$ is {\upshape(}split{\upshape)} reducible. Then for any $0\leq j\leq f-1$:
      {\begin{center}
            \begin{tabular}{ |c||c|c|c|c|c|c|}
              \hline  
              $s_j$& $r_j$ & $r_j+1$& $p-2-r_j$&$p-3-r_j$          \\
              \hline  
              $\frac{a_{d',j}^{\sigma}}{1-p^{d'}}$ & $-1$&$h^{[j]}+\frac{hp^{f-j}}{1-q}$&$ -1-h^{[j]} - \frac{hp^{f-j}}{1-q}$&$-1$   \\
              \hline
              \end{tabular}
            \end{center}}
        \end{enumerate}
\end{lem1}
\begin{proof}
(i) Note that we always have $2|d$ (as $d \nmid f$ but $d|(2f)$) and so $(2f)|d'$. Thus it suffices to prove the formulas for $\frac{a_{2f,j}^{\sigma}}{1-p^{2f}}$; we choose to work with $2f$ because $d|(2f)$ by \cite[Lem.~5.2]{breuil-IL}. Using Lemma \ref{lem:s-s'}, we can inductively determine $c_{n,j}^{\sigma}$ for $1\leq n\leq 2f$, and then compute $a_{2f,j}^{\sigma}$ using the formula $a_{2f,j}^{\sigma}=\sum_{k=0}^{2f-1}p^kc_{2f-k,j+k}^{\sigma}$, where $c^{\sigma}_{n,j}$ is understood to be $c^{\sigma}_{n, j~(\mathrm{mod}\ f)}$ if $j\geq f$.  

We do this in the case $j=0$ and $s_0=r_0$, and leave the other cases to the reader. In this case, we obtain using Lemma \ref{lem:s-s'} that
\[c^{\sigma}_{1,f-1}=p-2-r_{f-1}, ~  \dots,  ~ c^{\sigma}_{f-1,1}=p-2-r_1,~ c^{\sigma}_{f,0}=p-1-r_0,\]
\[c^{\sigma}_{f+1,f-1}=r_{f-1}+1, ~\dots, ~ c^{\sigma}_{2f-1,1}=r_1+1,~ c^{\sigma}_{2f,0}=r_0,\]
and so
\[\begin{array}{rll}a_{2f,0}^{\sigma}&=&r_0+p(r_1+1)+\cdots +p^{f-1}(r_{f-1}+1)\\
&&+p^f(p-1-r_0)+p^{f+1}(p-2-r_{1})+\cdots+p^{2f-1}(p-2-r_{f-1})\\
&=&(h-1)+p^f(p^f-h)\\
&=&(1-p^{2f})(-1+\frac{h}{1+p^f}),\end{array} \]
proving the result. 

(ii) In this case it suffices to prove the formulas for $\frac{a_{f,j}^{\sigma}}{1-p^f}$. The computation is similar to (i) and is easier, and we leave it to the reader.
\end{proof}

For $i \ge 0$ let
\[v_i \defeq v_{\delta^i(\sigma)} \in \delta^i(\sigma)^{N_0}\setminus \{0\}\]
and
$\mu_i \defeq \mu_{\delta^i(\sigma)} \in \F\s$, as defined in Lemma \ref{lem:BP}. Then by~\eqref{eq:x-delta} we have
\begin{equation}
v_i = \mu_{i-1}\cdot \underline{Y}^{\underline{c_i^\sigma}}\smatr{p}{}{}{1}(v_{i-1}) \quad\forall\ i \ge 1.\label{eq:3y}
\end{equation}
Let
\begin{equation}\label{eq:lambda-sigma}
\lambda_\sigma \defeq (-1)^{d(f-1)}\big(\prod_{0\leq i'\leq d-1}\prod_{j\in J^{\rm max}(\sigma_{i'})}(p-1-s_j^{(i'+1)})!\big)^{-1}\nu_\sigma,
\end{equation}
where $\nu_\sigma \in \F\s$ is defined as before \cite[Prop.~3.2.4.2]{BHHMS2}, i.e.\ the eigenvalue of the operator $S^d$ defined in \cite[\S 4]{breuil-IL} acting on $\sigma^{I_1}$. Note that $\nu_\sigma$ depends only on the orbit of $\sigma$, and hence the same is true for $\lambda_\sigma$.

\begin{lem1}\label{lem:n+d:lambda}
We have
\begin{equation*}
  \prod_{i=0}^{d-1} \mu_i = \lambda_\sigma^{-1}.
\end{equation*}
\end{lem1}
\begin{proof}
This follows from \cite[Lemma 3.2.2.5]{BHHMS2} and the definition of $\nu_\sigma$. 
\end{proof}

\begin{rem1}\label{rem:lambda-DL}
When $\pi$ moreover comes from cohomology, i.e.~is as in (\ref{piindef}) or (\ref{pidef}), it is conjectured in \cite[\S6]{breuil-IL} and proved in \cite[Thm.~5.36]{DoLe} that 
\begin{itemize}
\item if $\brho$ is irreducible, then $\nu_\sigma=(-1)^{\frac{dh}{2f}(1+\sum_{j=0}^{f-1}r_j)}(-\det(\brho)(p))^{\frac{d}{2}}$;
\item if $\brho$ is reducible, then $\nu_\sigma=(-1)^{\frac{dh}{f}\sum_{j=0}^{f-1}r_j}\lambda_0^{|\overline{J_{\sigma}}|\frac{d}{f}}\lambda_1^{|J_{\sigma}|\frac{d}{f}}$, where $J_{\sigma}\subset\{0,1,\dots,f-1\}$ is the set corresponding to $\sigma$ and $\overline{J_{\sigma}}$ denotes its complement.
\end{itemize}
Here, $h$ is the number attached to $\sigma$ in \cite[Lemma 6.2]{breuil-IL} (it is not the integer $h$ of \S\ref{conjecture}). By the proof of \cite[Lemma 6.2]{breuil-IL}, we deduce 
\begin{equation}\label{eq:lambda-DL}
\lambda_\sigma=\left\{\begin{array}{lll}(-1)^{d(f-1)}(-\det(\brho)(p))^{\frac{d}{2}}&\text{if $\brho$ irreducible,}\\
(-1)^{d(f-1)}\lambda_0^{|\overline{J_{\sigma}}|\frac{d}{f}}\lambda_1^{|J_{\sigma}|\frac{d}{f}}&\text{if $\brho$ reducible}.\end{array}\right.
\end{equation}
\end{rem1}
 
The following result follows by induction from~\eqref{eq:3y}, as well as Lemma~\ref{lem:n+d:lambda}.
\begin{lem1}\label{lem:a(n)}
For all $n\geq 0$, we have $\big(\prod\limits_{i=0}^{n-1} \mu_i\big)\cdot \underline{Y}^{\underline{a_n^\sigma}}\smatr{p}{}{}{1}^n(v_0) = v_{n}$. 
In particular, for all $n \ge 0$, $\underline{Y}^{\underline{a_{nd}^\sigma}}\smatr{p}{}{}{1}^{nd}(v_{\sigma}) = \lambda_\sigma^n v_\sigma$.
\end{lem1}

\begin{prop1}\label{prop:deg-phin}
 Let  $\underline{k}\in \Z_{\geq0}^{f}$ and $n \ge 0$.  If $\|\underline{k}\|\geq \|\underline{a_n^\sigma}\|-(f-1)$ and  $\underline{Y}^{\underline{k}}\smatr{p}{}{}{1}^n(v_0)\neq0$, then $\underline{k}=\underline{a_n^\sigma}-\underline{\ell}$ for some $\underline{\ell}\geq\underline{0}$ satisfying $\|\underline{\ell}\|\leq f-1$ and 
\[\big(\prod\limits_{i=0}^{n-1} \mu_i\big)\cdot \underline{Y}^{\underline{k}}\smatr{p}{}{}{1}^n(v_0) = \underline{Y}^{-\underline{\ell}}v_{n} \in \sigma_n.\]
\end{prop1}
 
\begin{proof} 
If $n=0$, we necessarily have $\underline k = \underline{a_0^\sigma} = \underline \ell = \underline 0$ and there is nothing to prove. Assume $n\geq 1$ and that the statement holds for $n-1$. 

Let $\underline{k}\in \Z_{\geq0}^{f}$ with $\|\underline{k}\|\geq \|\underline{a_n^\sigma}\|-(f-1)$.  Write $\underline{k}=p\delta(\underline{k}')+\underline{k}''$, with $\underline{k}' \ge \underline 0$
and $\underline{0}\leq \underline{k}''\leq \underline{p-1}$.
Recalling that $\|\delta(\cdot)\|=\|\cdot\|$, the assumption implies the following inequalities 
\[p\|\underline{k}'\|+(p-1)f\geq  \|\underline{k}\|> \|\underline{a}_n^\sigma\|-f \geq p\|\underline{a_{n-1}^\sigma}\|-f, \]
from which we deduce $\|\underline{k}'\|>\|\underline{a_{n-1}^\sigma}\|-f$, equivalently \[\|\underline{k}'\|\geq \|\underline{a_{n-1}^\sigma}\|-(f-1). \]

We clearly have 
\begin{equation}
  \underline{Y}^{\underline{k}}\smatr{p}{}{}{1}^n(v_0)=\underline{Y}^{\underline{k}''}\smatr{p}{}{}{1}\big(\underline{Y}^{\underline{k}'}\smatr{p}{}{}{1}^{n-1}(v_0)\big),\label{eq:2y}
\end{equation}
so in particular $\underline{Y}^{\underline{k}'}\smatr{p}{}{}{1}^{n-1}(v_0)\neq 0$. As $\|\underline{k}'\|\geq \|\underline{a_{n-1}^\sigma}\|-(f-1)$, by the inductive hypothesis there exists $\underline{\ell}'\geq\underline{0}$ with $\|\underline{\ell}'\|\leq f-1$ such that 
\begin{equation}\label{eq:k'ind}
\underline{k}'=\underline{a_{n-1}^\sigma}-\underline{\ell}'
\text{\ \ and\ \ } \big(\prod\limits_{i=0}^{n-2} \mu_i\big)\cdot \underline{Y}^{\underline{k}'}\smatr{p}{}{}{1}^{n-1}(v_0) = \underline{Y}^{-\underline{\ell}'}v_{n-1}\in \sigma_{n-1}.
\end{equation}

We first assume $m>0$ and claim that $\underline{\ell}'=\underline{0}$. Indeed, the relation $\|\underline{k}\|\geq \|\underline{a_n^\sigma}\|-(f-1)$ together with \eqref{eq:k'ind} gives 
\[\|\underline{k}''\|\geq p\|\underline{\ell}'\|-(f-1)+\|\underline{c_n^\sigma}\|.\]
Lemma \ref{lem:BP-florian}(ii) applied with $\sigma = \sigma_{n-1}$ (and genericity) shows that $k_j''\geq \ell'_{j+1}p$ for all $j$. However, by definition $0\leq k_j''\leq p-1$, so we must have   $\ell'_{j+1}=0$ for all $j$. This proves the claim.

By the claim and by equations \eqref{eq:2y}--\eqref{eq:k'ind} we have $\underline{k}'=\underline{a_{n-1}^\sigma}$ and
\[\big(\prod\limits_{i=0}^{n-2} \mu_i\big)\cdot \underline{Y}^{\underline{k}'}\smatr{p}{}{}{1}^{n-1}(v_0) = v_{n-1}, \text{\ \ so\ \ } \underline{Y}^{\underline{k}''}\smatr{p}{}{}{1}(v_{n-1}) \ne 0.\]
By the previous paragraph we have moreover that $\|\underline{k}''\|\geq \|\underline{c_n^\sigma}\|-(f-1)$.
Remark \ref{rem:BP-florian} applied with $\sigma = \sigma_{n-1}$ gives $\mu_{n-1}\cdot\underline{Y}^{\underline{k}''}\smatr{p}{}{}{1}(v_{n-1}) = \underline{Y}^{-\underline{\ell}}v_{n} \in \sigma_n$
for some $\underline{\ell}\geq\underline{0}$ satisfying $\|\underline{\ell}\|\leq f-1$ and $\underline{\ell} = \underline{c_n^\sigma}-\underline{k}''$. As $\underline{k}'=\underline{a_{n-1}^\sigma}$ we deduce $\underline{\ell} = \underline{a_n^\sigma}-\underline{k}$ and the result follows.

Now we assume $m=0$, equivalently $\sigma\cong\delta(\sigma)$. It is easy to see that this case only happens when $\brho$ is (split) reducible and either $J_{\sigma}=\emptyset$ or $J_{\sigma}=\{0,\dots,f-1\}$. In this case we have $\underline{a_n^{\sigma}}=\underline{p^n-1}$ for any $n\geq 0$, and Lemma \ref{lem:BP} implies that $\underline{Y}^{\underline{a_n^{\sigma}}}\smatr{p}{}{}1^n(v_0)\neq0$. Using \eqref{eq:2y} and the fact $Y_jv_0=0$ for all $j$, an induction shows that if  $k_{j}\geq p^{n}$ for some $0\leq j\leq f-1$, then \[\underline{Y}^{\underline{k}} \smatr{p}{}{}{1}^{n}(v_0)=\underline{Y}^{\underline{k}-p^{n}\underline{k}'}\smatr{p}{}{}1^{n}\big(\underline{Y}^{\delta^{-n}(\underline{k}')}v_0\big)=0,\]
where $k'\in\Z_{\geq0}^f$ is defined as: $k'_j=1$ and $k'_{j'}=0$ for $j'\neq j$. We deduce that $\underline{Y}^{\underline{k}}\smatr{p}{}{}1^{n}(v_0)\neq0 $ if and only if $\underline{k}\leq \underline{a_n^{\sigma}}$,  which implies the first assertion. The second assertion can be proved as above, noting that Remark \ref{rem:BP-florian} remains true when $m=0$.
\end{proof}

\begin{cor1}\label{cor:deg-phin}
  Let  $\underline{k}\in \Z_{\geq0}^{f}$ and $n \ge 0$. 
  \begin{enumerate}
  \item If $\|\underline{k}\|>\|\underline{a_n^\sigma}\|$, then $\underline{Y}^{\underline{k}}\smatr{p}{}{}{1}^n(v_0)=0$.
  \item If $\|\underline{k}\|=\|\underline{a_n^\sigma}\|$ and if $\underline{Y}^{\underline{k}}\smatr{p}{}{}{1}^n(v_0)\neq0$, then
    $\underline{k}=\underline{a_n^\sigma}$.
  \end{enumerate}
\end{cor1}
\begin{proof}
It is a direct consequence of Proposition \ref{prop:deg-phin}.
\end{proof}

\subsection{The degree function on an admissible smooth representation of \texorpdfstring{$\GL_2(K)$}{GL\_2(K)}}
\label{ss:gr-pi}

We define and study a ``degree function'' on representations $\pi$ as at the end of \S\ref{conjecture}.

Let $\rhobar$ and $\pi$ be as in \emph{loc.~cit.}  {We define $\gr(\pi)\defeq \bigoplus_{n\geq 0}\pi[\m_{I_1}^{n+1}]/\pi[\m_{I_1}^{n}]$.} For $v\in \pi$, we define \[\deg(v)\defeq \min\{n \ge -1:\ v\in \pi[\m_{I_1}^{n+1}]\} \in \Z_{\ge -1}.\]
 {We let $\gr(v)$ be the image of $v$ in $\pi[\m_{I_1}^{\deg(v)+1}]/\pi[\m_{I_1}^{\deg(v)}]$ if $v\ne 0$ and $\gr(v)=0$ if $v=0$.}

Fix $\sigma\in W(\brho)$ and let $v_{\sigma}\in \sigma^{N_0}\setminus\{0\}$. Define $\underline{a_n^\sigma}\in \Z_{\geq 0}^f$ as in \S\ref{ss:degree}.

\begin{prop1}\label{prop:degree-of-phi-n-v}
For all $n \ge 0$ we have
\[\deg\big(\smatr{p}{}{}{1}^{n}(v_{\sigma})\big)=\|\underline{a_n^\sigma}\|.\]
\end{prop1}
\begin{proof}
Put $u_n\defeq\smatr{p}{}{}{1}^n(v_{\sigma})$ for simplicity.  First, by the proof of \cite[Cor.~5.3.5]{BHHMS1}, we know that as a $\gr(\F\bbra{I_1/Z_1})$-module $\gr(\pi)$ is annihilated by the ideal $J$ defined by $J\defeq (y_jz_j,z_jy_j; 0\leq j\leq f-1)$, so that $\gr(\pi)$ becomes a graded module over $R\defeq \gr(\F\bbra{I_1/Z_1})/J$ which is commutative, isomorphic to $\F[y_j,z_j]/(y_jz_j ; 0 \le j\le f-1)$, with $y_j$, $z_j$ of degree $-1$. On the other hand, as $v_{\sigma}\in \sigma^{N_0}=\sigma^{I_1}$, it is direct to check that $u_n$ is fixed by $\smatr{1}0{p\cO_K}1$ for all $n\geq 0$, hence 
\[\sum_{\lambda\in\F_q}\lambda^{-p^j}\smatr{1}0{p[\lambda]}1u_n=\big(\sum_{\lambda\in\F_q}\lambda^{-p^j}\big)u_n=0.\]
Namely,  
 $u_n$ is annihilated by $\sum_{\lambda\in\F_q}\lambda^{-p^j}\smatr{1}0{p[\lambda]}1\in\F\bbra{I_1/Z_1}$ (a lifting of $z_j$), 
  hence $z_j\gr(u_n)=0$ and consequently we observe that any element in $\langle R\cdot \gr(u_n)\rangle$ is annihilated by $z_j$. 

Next we note the following fact: if $v\in \pi$ with $\deg(v) > 0$ and if $\gr(v)$ is annihilated by all $z_j$, then there exists some $i\in\{0,\dots,f-1\}$ such that $y_i\gr(v)\neq0$. (If not,  then $R_{-1}$, the degree $-1$ part of $R$, annihilates $\gr(v)$. Suppose $v \in \pi[\m_{I_1}^{n+1}] \setminus \pi[\m_{I_1}^{n}]$ for some $n \ge 1$. Since $R_{-1} = \m_{I_1}/\m_{I_1}^{2}$,  we deduce  $\m_{I_1}v \subseteq \pi[\m_{I_1}^{n-1}]$, i.e.\ $v \in \pi[\m_{I_1}^{n}]$, contradiction.)
As a consequence, $Y_iv\neq0$ and 
\[\deg(Y_iv)=\deg(v)-1;\]
moreover we have $\gr(Y_iv)=y_i\gr(v)\in \langle R\cdot\gr(v)\rangle$. Applying this fact to $u_n$ (and to $Y_iu_n$, etc.) and using the observation of the last paragraph, we find that there exists $\underline{a_n'^\sigma}\in\Z_{\geq 0}^f$ such that $\underline{Y}^{\underline{a_n'^\sigma}}u_n$ is of degree 0, i.e.\ $\underline{Y}^{\underline{a_n'^\sigma}}u_n\in \pi^{I_1}\setminus\{0\}$ and 
\[\deg(u_n)=\|\underline{a_n'^\sigma}\|.\]
On the one hand,  we have $\|\underline{a_n'^\sigma}\|\leq \|\underline{a_n^\sigma}\|$ by Corollary \ref{cor:deg-phin}(i) (as $\underline{Y}^{\underline{a_n'^\sigma}}u_n\neq 0$). On the other hand, we  have $\deg(\smatr{p}{}{}{1}^n(v_{\sigma}))\geq \|\underline{a_n^\sigma}\|$ by Lemma \ref{lem:a(n)}, so the result follows.   
\end{proof}

If $V$ is any admissible  smooth  representation of $\GL_2(K)$ over $\F$, we define $\deg(v)$ for $v\in V$ as above. On the other hand, by restricting $V$ to $N_0$, we can also define
\[\deg'(v)\defeq\min\{n\geq -1: v\in V[\frak{m}_{N_0}^{n+1}]\}.\]
This is well-defined as $V$ is smooth. It is clear that $\deg(v)\geq \deg'(v)$. 

We note the following  consequence of the proof of Proposition \ref{prop:degree-of-phi-n-v} (it will not be used in this paper). 
  
\begin{cor1} 
Let $V$ be in the category $\mathcal C$ of \S\ref{conjecture} and assume that $\gr(V)$ is annihilated by the ideal $J$ defined in the proof of Proposition \ref{prop:degree-of-phi-n-v}. If $v\in V$ is an element fixed by $\smatr10{p\cO_K}1$, then 
 there exists $\underline{k}\in\Z_{\geq 0}^f$ with $|\!|\underline{k}|\!|=\deg(v)$ such that $0\neq \underline{Y}^{\underline{k}}v\in V^{I_1}$. Moreover, we have 
 $\deg(v)=\deg'(v)$.
\end{cor1}

\subsection{A crucial finiteness result}
\label{ss:finiteness}

We prove an important finiteness result (Proposition \ref{prop:cond-ii}) which will be crucially used in \S\ref{sec:basis-d_api} to construct elements of $\Hom_A(D_A(\pi),A)$.

Fix $\sigma\in W(\brho)$ and define $\sigma_i\in W(\brho)$, $v_i\in\sigma_i$ and $d\in\Z_{\geq1}$ as in \S\ref{ss:degree} (before Lemma \ref{lem:n+d:lambda}).
We have elements $\un{c_n^\sigma},\underline{a_n^\sigma}\in \Z_{\geq 0}^f$ defined for $n \ge 1$ (resp.\ $n \ge 0$) in \eqref{eq:c_sigma} (resp.\ \eqref{eq:abc}). 
By induction we have 
\[\underline{a_n^\sigma}=\sum_{i=0}^{n-1}p^{i}\delta^i(\underline{c_{n-i}^\sigma})\]
and as   $\underline{c_{n}^\sigma}$ is periodic with period $d$, we deduce 
\begin{equation}\label{eq:a-ndp}
\underline{a_{nd'}^\sigma}=p^{d'}\underline{a_{(n-1)d'}^\sigma}+\underline{a_{d'}^\sigma},
\end{equation}
where we recall that $d' = df$ (so $\delta^{d'}$ is the identity). 

We consider the following elements for $\underline{i}\in\Z^{f}$:
\begin{equation*}x_{\sigma,\underline{i}}\defeq \lambda_\sigma^n\underline{Y}^{\underline{a_{nd}^\sigma}-\underline{i}}\smatr{p}{}{}{1}^{nd}(v_\sigma),
\end{equation*}
where $\lambda_{\sigma}$ is defined in \eqref{eq:lambda-sigma} and $n\geq 0$ is chosen large enough so that $\underline{a_{nd}^\sigma}-\underline{i}\geq \underline 0$. By Lemma~\ref{lem:a(n)} the definition is independent of the choice of $n$.

The following finiteness result is the main result of this section.

\begin{prop1}\label{prop:cond-ii}
For any $M\in\Z$ the set $\{\underline{i}\in\Z^f:x_{\sigma,\underline{i}}\neq0, \|\underline{i}\|=M\}$ is finite.
\end{prop1}
For Lemmas \ref{lem:d-case} and \ref{lem:n'nd-case} below, we  assume $m=|J^{\rm max}(\sigma)|>0$.
\begin{lem1}\label{lem:d-case}
Let $\underline{k}\in \Z_{\geq0}^{f}$ and $n\geq 1$. Assume that for some $0\leq j_0\leq f-1$, 
\begin{enumerate}
\item[\upshape(a)] $k_{j_0} \le a_{n,j_0}^\sigma-p-c_{n,j_0}^\sigma$, 
\item[\upshape(b)] $\|\underline{k}\|>\|\underline{a_n^\sigma}\|-c_{n,j_0}^\sigma$. 
\end{enumerate}
Then $\underline{Y}^{\underline{k}}\smatr{p}{}{}{1}^{n}(v_0)=0$. 
\end{lem1}
\begin{proof}
  Write $\underline{k}=p\delta(\underline{k}')+\underline{k}''$ with $\underline{k}',\underline{k}''\geq \underline{0}$ and $\underline{k}''\leq \underline{p-1}$.  
Condition (b) implies 
\[p\|\underline{k}'\|+(p-1)f\geq \|\underline{k}\|> p\|\underline{a_{n-1}^\sigma}\|+\sum_{j\neq j_0}c_{n,j}^\sigma\]
and consequently
\[p\|\underline{k}'\|+pf>p\|\underline{a_{n-1}^\sigma}\|.\] Thus, we have $\|\underline{k}'\|>\|\underline{a_{n-1}^\sigma}\|-f$. 

Assume $\underline{Y}^{\underline{k}}\smatr{p}{}{}{1}^n(v_0)\neq0$ for a contradiction. Then by the proof of Proposition \ref{prop:deg-phin} we also have $\underline{Y}^{\underline{k}'}\smatr{p}{}{}{1}^{n-1}(v_0)\neq0$. Moreover, by  Proposition \ref{prop:deg-phin},  there exists $\underline{i}\geq \underline{0}$ with $\|\underline{i}\|\leq f-1$ such that  $\underline{k}'=\underline{a_{n-1}^\sigma}-\underline{i}$ and
\[\big(\prod\limits_{i=0}^{n-2} \mu_i\big)\cdot \underline{Y}^{\underline{{k}}'}\smatr{p}{}{}{1}^{n-1}(v_0)= \underline{Y}^{-\underline{i}}v_{n-1}\in \sigma_{n-1}.\]
Thus, condition (a)  translates to
\[p(a_{n-1,j_0+1}^\sigma-i_{j_0+1})+k_{j_0}''\le a_{n,j_0}^\sigma-p-c_{n,j_0}^\sigma\] 
from which we deduce
$k_{j_0}''\leq p(i_{j_0+1}-1)$ using \eqref{eq:abc}, 
and  we get a contradiction by Lemma \ref{lem:BP-florian-v} applied to $\underline{k}''$. Indeed, $\underline{Y}^{\underline{{k}}''}\smatr{p}{}{}{1}(\underline{Y}^{-\underline{i}}v_{n-1}) \ne 0$ and the equality $\underline{k}'=\underline{a_{n-1}^\sigma}-\underline{i}$ together with condition (b) imply 
\[\|\underline{k}''\|>p\|\underline{i}\|+\|c_{n}^\sigma\|-c_{n,j_0}^\sigma\]
which verifies the corresponding condition (b) of Lemma \ref{lem:BP-florian-v} (with $\sigma = \sigma_{n-1}$). 
\end{proof}

Recall that $d' = df$, that $\delta^{d'}$ is the identity, and that $\underline{c_n^\sigma}$ is periodic with period $d$.

\begin{lem1}\label{lem:n'nd-case}
Let $\underline{k}\in \Z_{\geq0}^{f}$ and $n'>n\geq 0$. Assume that for some $0\leq j_0\leq f-1$,
\begin{enumerate}
\item[\upshape(a)] $k_{j_0} \le a_{n'd',j_0}^\sigma-a_{nd',j_0}^\sigma-p^{nd'}(p+c_{d,j_0}^\sigma)$ and 
\item[\upshape(b)] $\|\underline{k}\|>  \|\underline{a_{n'd'}^\sigma}\|-\|\underline{a_{nd'}^\sigma}\|-p^{nd'}c_{d,j_0}^\sigma+f(p^{nd'}-1)$.
\end{enumerate}Then $\underline{Y}^{\underline{k}}\smatr{p}{}{}{1}^{n'd'}(v_0)=0$.
\end{lem1}
\begin{proof}
Applying Lemma \ref{lem:d-case} with $n\defeq (n'-n)d'$, we see that if $\underline{k}'\in \Z_{\geq0}^{f}$ such that $k'_{j_0} \le a_{(n'-n)d',j_0}^\sigma-p-c_{d,j_0}^\sigma$ (recall that $c_{d,j_0}^\sigma$ is periodic) and if $\|\underline{k}'\|>\|\underline{a_{(n'-n)d'}^\sigma}\|-c_{d,j_0}^\sigma $, then $\underline{Y}^{\underline{k}'}\smatr{p}{}{}{1}^{(n'-n)d'}(v_0)=0$. 

 Write $\underline{k}=p^{nd'}\underline{k}'+\underline{k}''$ with $\underline{k}',\underline{k}''\geq\underline{0}$ and $\underline{k}''\leq \underline{p^{nd'}-1}$. Note that $\underline{a_{n'd'}^\sigma}-\underline{a_{nd'}^\sigma}=p^{nd'}\underline{a_{(n'-n)d'}^\sigma}$ by \eqref{eq:a-ndp}. Firstly, by condition (a)  we have
 \[p^{nd'}  k'_{j_0}\leq k_{j_0} \le p^{nd'} a_{(n'-n)d',j_0}^\sigma-p^{nd'}(p+c_{d,j_0}^\sigma)\]
 and so $k_{j_0}' \le a_{(n'-n)d',j_0}^\sigma-p-c_{d,j_0}^\sigma$. Secondly,  as $f(p^{nd'}-1)-\|\underline{k}''\|\ge 0$,  condition (b) implies that \[p^{nd'}\|\underline{k}'\|> p^{nd'}\|\underline{a_{(n'-n)d'}^\sigma}\|-p^{nd'}c_{d,j_0}^\sigma\]
so that 
\[\|\underline{k}'\|> \|\underline{a_{(n'-n)d'}^\sigma}\|- c_{d,j_0}^\sigma.\]
We then conclude that $\underline{Y}^{\underline{k}'}\smatr{p}{}{}{1}^{(n'-n)d'}(v_0)=0$ as explained above, hence
\[ \underline{Y}^{\underline{k}}\smatr{p}{}{}{1}^{n'd'}(v_0)=\underline{Y}^{\underline{k}''}\smatr{p}{}{}{1}^{nd'}\underline{Y}^{\underline{k}'}\smatr{p}{}{}{1}^{(n'-n)d'}(v_0) = 0.\qedhere \]
\end{proof}

\begin{proof}[Proof of Proposition \ref{prop:cond-ii}]
If $m=0$, then the end of the proof of Proposition \ref{prop:deg-phin}  implies $x_{\sigma,\underline{i}}=0$ if $i_j<0$ for some $0\leq j\leq f-1$, from which the result easily follows.

 Assume $m>0$ from now on, so that Lemma \ref{lem:n'nd-case} applies. Fix any $M \in \Z$. We will show that the set $\{\underline{i} \in \Z^f: x_{\sigma,\underline{i}}\neq 0, \|\underline{i}\|=M\} $ is finite. 
Choose $n$  large enough such that for all $0\leq j\leq f-1$:  
\begin{equation}
\|\underline{a_{nd'}^\sigma}\|+p^{nd'}c_{d,j}^\sigma-f(p^{nd'}-1)>M;\label{eq:1y}
\end{equation}
this is always possible because the left-hand side tends to infinity when 
$n\ra\infty$ (recall that $c_{d,j}^\sigma>f$, by genericity). 

Now pick any $\underline{i} \in \Z^f$ such that $\|\underline{i}\| = M$. Choose $n'>n$ large enough such that $\underline{a_{n'd'}^\sigma}\geq \underline{i}$,
hence $x_{\sigma,\underline{i}} \in \F^\times \underline{Y}^{\underline{a_{n'd'}^\sigma}-\underline{i}} \smatr{p}{}{}{1}^{n'd'}(v_0)$. By \eqref{eq:1y} and as $\|\underline{i}\| = M$, we get for all $0\leq j\leq f-1$:
\[\|\underline{a_{n'd'}^\sigma}\|-\|\underline{i}\| > \|\underline{a_{n'd'}^\sigma}\|-(\|\underline{a_{nd'}^\sigma}\|+p^{nd'}c_{d,j}^\sigma-f(p^{nd'}-1)).\]
There are two cases: 
\begin{itemize}
\item If $i_{j_0} \ge a_{nd',j_0}^\sigma+p^{nd'}(p+c_{d,j_0}^\sigma)$ for some $j_0$, then $x_{\sigma,\underline{i}}=0$  by Lemma \ref{lem:n'nd-case} (applied to $\underline{k}\defeq \underline{a_{n'd'}^\sigma}-\underline{i}$).
\item Otherwise, we must have $i_j < a_{nd',j}^\sigma+p^{nd'}(p+c_{d,j}^\sigma)$ for all $j$, and such a set (together with the restriction $\|\underline{i}\|=M$) is automatically finite. Note that the quantities $a_{nd',j}^\sigma+p^{nd'}(p+c_{d,j}^\sigma)$ depend only on our fixed $M$, as $n$ does. \qedhere
\end{itemize}
\end{proof}

\subsection{An explicit basis of \texorpdfstring{$\Hom_A(D_A(\pi),A)$}{Hom\_A(D\_A(pi),A)}}
\label{sec:basis-d_api}

We exhibit an $A$-basis of $\Hom_A(D_A(\pi),A)$ and explicitly describe its image in the vector space $\Hom_\F^{\cont}(D_A(\pi),\F)$ via the embedding \eqref{eq:mu-star}.

Recall $\pi$ and $\rhobar$ are as in Theorem \ref{mainbis} with $\rhobar$ as at the end of \S\ref{conjecture}, in particular $\pi^{I_1}$ is multiplicity-free for the action of $I$. For any $\sigma \in W(\rhob)$ and our fixed choice of $v_\sigma \in \sigma^{N_0}\setminus \{0\}$ we define:
\begin{equation}\label{eq:13f}
  x_{\sigma,k} \defeq \lambda_\sigma^n \un Y^{\un{a_{nd}^\sigma}-\un{k}} \smatr{p}{}{}{1}^{nd} v_\sigma
\end{equation}
for $k \ge 0$ and any $n \gg_k 0$.
This is well-defined by Lemma~\ref{lem:a(n)}.

Recall from Proposition~\ref{prop:degree-of-phi-n-v} that 
\begin{equation*}
\text{$\smatr{p}{}{}{1}^{nd} v_\sigma \in \pi[\m_{I_1}^{\|\un{a_{nd}^\sigma}\|+1}]$,\ \ \ so $x_{\sigma,k}\in \pi[\m_{I_1}^{kf+1}]$},
\end{equation*}
hence \ by \ Proposition~\ref{prop:F-dual} \ the \ sequence \ $(x_{\sigma,k})_{k \ge 0}$ \ defines \ an \ element \ $x_\sigma$ \ of $\Hom_\F^{\cont}(D_A(\pi),\F)$ of degree 0.

\begin{thm1}\label{thm:descent}
  The elements $\{x_\sigma: \sigma \in W(\rhob)\}$ are contained in the image of the injection
  \begin{equation*}
    \mu_* : \Hom_A(D_A(\pi),A) \into \Hom_\F^{\cont}(D_A(\pi),\F)
  \end{equation*}
  and form an $A$-basis of $\Hom_A(D_A(\pi),A)$.
\end{thm1}

We first need a lemma. Note that $\pi^{I_1}$ is multiplicity-free  for the action of $I$, so there
exist unique $I$-eigenvectors $v_\sigma^* \in (\pi^{I_1})\dual = \gr_0(\pi\dual)$ such that $\ang{v_\sigma,v_{\sigma'}^*} = \delta_{\sigma,\sigma'}$
(for $\sigma$, $\sigma' \in W(\rhob)$). We already know that $D_A(\pi)$ is free by Remark \ref{freeness}. The following result only applies to our current $\pi$ but is more precise.

\begin{lem1}\label{lm:freeness}
  Suppose that $\pi$ is as above. Then $\gr(D_A(\pi))$ is a free $\gr(A)$-module with basis $(v_\sigma^*)_{\sigma \in W(\rhob)}$.
  In particular, $D_A(\pi)$ is a filtered free $A$-module of rank $2^f$.
\end{lem1}

\begin{proof}
  Recall from \cite[\S3.1]{BHHMS2} that $\gr(D_A(\pi))$ is obtained from $\gr(\pi\dual)$ by localizing at $\prod_j y_j$.
  By localizing the surjection in \cite[Thm.\ 3.3.2.1]{BHHMS2} at $\prod_j y_j$ and using \cite[Lemma 3.3.1.3(i)]{BHHMS2} we obtain a surjection
  $\bigoplus_{\sigma \in W(\rhob)} \gr(A) \onto \gr(D_A(\pi))$ of $\gr(A)$-modules, sending the standard basis element indexed by $\sigma$ on the left to 
  $v_\sigma^*$. But $\rank_{\gr(A)}(\gr(D_A(\pi))) = \rank_A(D_A(\pi)) = 2^f$ by \cite[Lemma 3.1.4.1]{BHHMS2} and \cite[Cor.~3.3.2.4]{BHHMS2}, hence the surjection $\bigoplus_{\sigma \in W(\rhob)} \gr(A) \onto \gr(D_A(\pi))$ is an isomorphism. By \cite[Thm.\ I.4.2.4(5)]{LiOy}
  we can lift it to an isomorphism $\bigoplus_{\sigma \in W(\rhob)} A \congto D_A(\pi)$ of filtered $A$-modules.
\end{proof}

\begin{proof}[Proof of Theorem~\ref{thm:descent}]
  Fix any $\sigma \in W(\rhob)$ and consider the continuous $\F$-linear map $h_\sigma(\defeq x_\sigma) : D_A(\pi) \to \F$ of degree 0 corresponding to the sequence $(x_{\sigma,k})_{k \ge 0}$.
  We endow $D_A(\pi)$ with its natural good filtration (coming from the $\fm_{I_1}$-adic filtration on $\pi^\vee$, cf.\ \cite[\S3.1.2]{BHHMS2}).
  Let again $S$ denote the multiplicative subset of $\F\bbra{N_0}$ generated by $Y_0\cdots Y_{f-1}$.
  To descend $h_\sigma$ to $\Hom_A(D_A(\pi),A)$ we now check the second criterion in Lemma~\ref{lm:descent-criterion}. Thus fix any $x \in D_A(\pi)$ and $M \in \Z$.
  By continuity there exists $e \in \Z$ such that $h_\sigma(F_e D_A(\pi)) = 0$. As $(\pi\dual)_S$ is dense in $D_A(\pi)$ we can find $\un\ell \in \Z^f$ and $x^* \in \pi\dual$
  such that $x - \un Y^{\un \ell} x^* \in F_{e+\|\un i\|} D_A(\pi)$. Then $h_\sigma((x - \un Y^{\un \ell} x^*) \un Y^{\un i}) = 0$ for all $\un i \in \Z^f$ such
  that $\|\un i\| = M$, so we may assume that $x = \un Y^{\un \ell} x^* \in (\pi\dual)_S$.

  As in \S\ref{ss:finiteness} we define $x_{\sigma,\un i} \defeq \lambda_\sigma^n\underline{Y}^{\underline{a_{nd}^\sigma}-\underline{i}}\smatr{p}{}{}{1}^{nd}(v_\sigma)$ for $\un i \in \Z^f$,
  where $n \gg_{\un i} 0$. (In particular, $x_{\sigma,(k,\dots,k)} = x_{\sigma,k}$ for $k \ge 0$ and $\un Y^{\un j} x_{\sigma,\un i} = x_{\sigma,\un i-\un j}$ for any $\un j \ge \un 0$.)
  Explicitly,
  \begin{equation*}
    h_\sigma \circ \un Y^{-\un k} = \ang{x_{\sigma,k},-} \quad\text{on $\pi\dual$}
  \end{equation*}
  for all $k \ge 0$, from which it follows from the properties of $(x_{\sigma,\un i})_{\un i}$ that
  \begin{equation}\label{eq:17f}
    h_\sigma \circ \un Y^{-\un i} = \ang{x_{\sigma,\un i},-} \quad\text{on $\pi\dual$}
  \end{equation}
  for all $\un i \in \Z^f$.
  This implies that
  \begin{equation*}
    h_\sigma(\un Y^{\un i} x) = h_\sigma(\un Y^{\un i+\un \ell} x^*) = x^*(x_{\sigma,-(\un i+\un \ell)})
  \end{equation*}
  which can be nonzero for only finitely many $\un i$ by Proposition~\ref{prop:cond-ii}. Thus $h_\sigma$ indeed descends to an element
  $H_\sigma$ of $\Hom_A(D_A(\pi),A)$.

  For the final claim, first note that
  \[\gr\big(\Hom_A(D_A(\pi),A)\big) \cong \Hom_{\gr(A)}\big(\gr(D_A(\pi)),\gr(A)\big)\]
  by \cite[Lemma I.6.9]{LiOy} and Lemma~\ref{lm:freeness}.
  By \cite[Cor.\ I.4.2.5(2)]{LiOy} it then suffices to show that the $\gr(H_\sigma)$ ($\sigma \in W(\rhob)$) form a basis of
  $\Hom_{\gr(A)}(\gr(D_A(\pi)),\gr(A))$. By Lemma~\ref{lm:freeness}, the $\gr(A)$-module $\gr(D_A(\pi))$ has basis
  $v_\sigma^*$ ($\sigma \in W(\rhob)$), so it will be enough to establish
  $\ang{\gr (H_\sigma),v_{\sigma'}^*} = \delta_{\sigma,\sigma'} \un y^{-\un 1}$ for all $\sigma$, $\sigma' \in W(\rhob)$.

  By the explicit formula from the proof of Lemma~\ref{lm:descent-criterion} we know that
  \begin{equation*}
    H_\sigma(x) = (\prod_j (1+T_j))\sum_{\un i} h_\sigma(\un Z^{\un i}x)\un Z^{-\un i-\un 1} \quad\forall\ x \in D_A(\pi).
  \end{equation*}
Consider the equality $\mu \circ H_\sigma = h_\sigma$. 
  Note that $H_\sigma$ is a filtered map of degree $f$, since 
  $h_\sigma$ is of degree 0, $\un Z^{\un i} \in F_{-\|\un i\|} A$, and $\prod_j (1+T_j) \in F_0 A$. Similarly,
  $\mu$ is a filtered map of degree $-f$. Therefore
  \begin{equation}\label{eq:18f}
    \gr(\mu) \circ \gr(H_\sigma) = \gr(h_\sigma).
  \end{equation}

  Recall that $\gr(A) = \F[y_0^{\pm 1},\dots,y_{f-1}^{\pm 1}]$. Let $\o\ve_{\un i} : \gr(A) \to \F$ be the map sending
  $\sum_{\un j \in \Z^f} \lambda_{\un j} \un y^{\un j}$ to $\lambda_{\un i}$; it is $\F$-linear and of degree $\|\un i\|$.
  By definition, $\gr(\mu) : \gr_f A \to \F$ sends $\gr(\prod_j (1+T_j) \sum_{\|\un i\| \ge -f} \lambda_{\un i} \un Z^{\un i})$ to $\lambda_{-\un 1}$.
  As $\gr(Y_j) = \gr(Z_j)$, it follows that 
  \begin{equation*}\gr(\mu) = \o\ve_{-\un 1}.
  \end{equation*}
  On the other hand, relation~\eqref{eq:17f} implies that
  \begin{equation}\label{eq:20f}
    \gr(h_\sigma) \circ \un y^{-\un i} = \ang{\gr(x_{\sigma,\un i}),-} \quad\text{on $\gr(\pi\dual)$}
  \end{equation}
  for all $\un i \in \Z^f$. (They are graded maps of degree $\|\un i\|$; we filter $\pi$ as in \S\ref{ss:gr-pi}.)

  Using equations~\eqref{eq:18f}--\eqref{eq:20f} we compute that
  \begin{align*}
    \o\ve_{\un i-\un 1} \circ \gr(H_\sigma) &= \o\ve_{-\un 1} \circ \un y^{-\un i} \circ \gr(H_\sigma) = \gr(\mu) \circ \gr(H_\sigma) \circ \un y^{-\un i}\\
    &= \gr(h_\sigma) \circ \un y^{-\un i} = \ang{\gr(x_{\sigma,\un i}),-} \quad\text{on $\gr(\pi\dual)$.}
  \end{align*}
  As $\o\ve_{\un i-\un 1} \circ \gr(H_\sigma)$ is a map of degree $\|\un i\|$, if $(\o\ve_{\un i-\un 1} \circ \gr(H_\sigma))(v_{\sigma'}^*) \ne 0$, then
  $\|\un i\| = 0$. By the definition of $x_{\sigma,\un i}$ and by Corollary~\ref{cor:deg-phin} we know that $x_{\sigma,\un i} = 0$
  if $\|\un i\| = 0$ and $\un i \ne \un 0$. Therefore,
  \begin{equation*}
    \gr(H_\sigma) = \ang{\gr(x_{\sigma,\un 0}),-} \un y^{-\un 1} = \ang{\gr(v_{\sigma}),-} \un y^{-\un 1} \quad\text{on $\gr(\pi\dual)$,}
  \end{equation*}
  as desired.
\end{proof}

\subsection{The \texorpdfstring{$(\varphi,\cO_K\s)$}{(phi,O\_K\^{}x)}-action on \texorpdfstring{$\Hom_A(D_A(\pi),A)$}{Hom\_A(D\_A(pi),A)}}
\label{sec:o_ks-action}

We determine the $\varphi$- and $\okt$-actions on the elements $x_\sigma \in \Hom_A(D_A(\pi),A)$ ($\sigma \in W(\rhob)$), as defined in~\S\ref{sec:basis-d_api}.

We first determine the $\varphi$-action on $x_{\sigma}$. Let $v_{\sigma}\in \sigma^{N_0}\setminus \{0\}$ be as in~\S\ref{sec:basis-d_api} which defines $x_{\sigma}=(x_{\sigma,k})_{k\geq 0}$ via \eqref{eq:13f}. 
By Lemma \ref{lem:BP} there exists a constant $\mu_{\sigma}\in \F^{\times}$ such that 
\begin{equation}\label{eq:const-mu}
v_{\delta(\sigma)}=\mu_{\sigma}\cdot \un{Y}^{\un{c_1^{\sigma}}}\smatr{p}{}{}{1}(v_{\sigma}),
\end{equation}
where $\un{c_{1}^\sigma}$ is defined in \eqref{eq:c_sigma}.

\begin{prop1}\label{prop:phi-x}
 For any $\sigma \in W(\rhob)$ we have 
 \[\varphi(x_{\sigma})=(-1)^{f-1}\mu_{\sigma}^{-1}\un{Y}^{-\un{c_1^{\sigma}}}x_{\delta(\sigma)}.\]
\end{prop1}

\begin{proof}
Equivalently, we need to check that  
\begin{equation}\label{eq:relation-check}
x_{\delta(\sigma),k}=(-1)^{f-1}\mu_{\sigma} \un{Y}^{\un{c_1^{\sigma}}} (\varphi(x_{\sigma}))_k
\end{equation}
for any $k\geq0$. 

We have 
\begin{equation}\label{eq:and-c1}
\un{a_{nd+1}^\sigma}=\un{c_1^\sigma}+p\delta(\un{a_{nd}^\sigma})=\un{a_{nd}^{\delta(\sigma)}}+p^{nd}\delta^{nd}(\un{c_1^\sigma})
\end{equation}
by \eqref{eq:abc}. By definition \eqref{eq:13f} and using \eqref{eq:const-mu}, we have
\[\begin{array}{rll}
x_{\delta(\sigma),k}&=&\lambda_{\delta(\sigma)}^n \un{Y}^{\un{a_{nd}^{\delta(\sigma)}}-\un{k}}\smatr{p}{}{}{1}^{nd}(v_{\delta(\sigma)})\\
&=&\lambda_{\delta(\sigma)}^n\un{Y}^{\un{a_{nd}^{\delta(\sigma)}}-\un{k}}\cdot \mu_{\sigma}\un{Y}^{p^{nd}\delta^{nd}(\un{c_1^\sigma})}\smatr{p}{}{}{1}^{nd+1}(v_{\sigma})\\
&=&\mu_{\sigma}\lambda_{\delta(\sigma)}^{n}\un{Y}^{\un{a_{nd+1}^\sigma}-\un{k}} \smatr{p}{}{}{1}^{nd+1}(v_{\sigma}),
\end{array}\]
where we applied \eqref{eq:and-c1}. On the other hand, by Lemma~\ref{lm:mu-star}(ii) and~\eqref{eq:7f} the action of $\varphi$ on $x_{\sigma}$ can be computed on
sequences as follows: for any $k \ge 0$,
\begin{equation*}(\varphi(x_{\sigma}))_{k}=(-1)^{f-1}\un{Y}^{\un{p\ell}-\un{k}}\smatr{p}{}{}{1}(x_{\sigma,\ell}),
\end{equation*}
where $\ell$ is chosen arbitrarily so that $p\ell \ge k$. Thus  
we have  (for $\ell$ large enough)
\[\begin{array}{rll} 
\un{Y}^{\un{c_1^\sigma}}(\varphi(x_{\sigma}))_k&=&(-1)^{f-1} \un{Y}^{\un{c_1^\sigma}}\un{Y}^{\un{p\ell}-\un{k}}\smatr{p}{}{}{1}(x_{\sigma,\ell})\\
&=&(-1)^{f-1} \un{Y}^{\un{c_1^\sigma}}\un{Y}^{\un{p\ell}-\un{k}}\cdot \lambda_{\sigma}^n\un{Y}^{p\delta(\un{a_{nd}^\sigma})-\un{p\ell}}\smatr{p}{}{}{1}^{nd+1}(v_{\sigma})\\
&=&(-1)^{f-1} \lambda_{\sigma}^n \un{Y}^{\un{a_{nd+1}^\sigma}-\un{k}}\smatr{p}{}{}{1}^{nd+1}(v_{\sigma}),
 \end{array}\]
where we used again \eqref{eq:13f} and~\eqref{eq:and-c1}.
As $\lambda_{\sigma}=\lambda_{\delta(\sigma)}$ by the discussion after \eqref{eq:lambda-sigma}, relation \eqref{eq:relation-check} is verified.  
\end{proof}

We now determine the $\okt$-action on  $x_\sigma$. By Lemma~\ref{lm:mu-star}(iii) 
we can compute this action on the image of $x_\sigma$ in $\Hom_\F^{\cont}(D_A(\pi),\F)$ (i.e.\ before descending).

For $a \in \cO_K\s$ and $0 \le i \le f-1$ we put
\[f_{a,i} \defeq \frac{\o a^{p^i} Y_i}{a(Y_i)} \in 1+F_{1-p}A,\]
where we follow the convention in \S\ref{sec:dual-vp-co_ks} of just writing an index $i$ instead of an index $\sigma_i$ (in particular $f_{a,0}=f_{a,\sigma_0}$ in \eqref{1-p}). Note that $\vp(f_{a,i}) = f_{a,i-1}^p$. We also let $\chi_\sigma : \F_q\s \to \F\s$ denote the eigencharacter
of $\diag(-,1)$ on $\sigma^{I_1}$.

\begin{prop1}\label{prop:gamma-action}
  For any $\sigma \in W(\rhob)$ and $a \in \cO_K\s$ we have
  \begin{equation*}
    a(x_\sigma) = N_{\F_q/\Fp}(\o a)^{-1} \chi_\sigma(\o a) \bigg(\prod_{i=0}^{f-1} f_{a,i}^{-a_{d',i}^\sigma/(1-p^{d'})}\bigg) x_\sigma
  \end{equation*}
  in $\Hom_\F^{\cont}(D_A(\pi),\F)$, where $d' = df$.
\end{prop1}
\begin{proof}
  First note that we may apply any element of $\F\bbra{N_0} + F_{-k f-1}A$ to \eqref{eq:13f} (with $F_{-k f-1}A$ killing both
  sides) by applying our convention in Remark~\ref{rk:convention-action} to both $x_{\sigma,k}\in \pi[\m_{I_1}^{kf+1}]$ and
  $\smatr{p}{}{}{1}^{nd'} v_\sigma \in \pi[\m_{I_1}^{\|\un{a_{nd'}^\sigma}\|+1}]$.

To simplify the notation we set $M \defeq \prod_{i=0}^{f-1} f_{a,i}^{a_{d',i}^\sigma/(1-p^{d'})}$.  Let us now consider $N_{\F_q/\Fp}(\o a) \chi_\sigma(\o a)^{-1} M a(x_\sigma)$.
  Combining both parts of Lemma~\ref{lm:action} and the previous paragraph, we obtain that its $k$-th component is given by the following formulas (where $\ell \gg_k 0$
  and $n \gg_\ell 0$):
  \begin{align*}
    &\chi_\sigma(\o a)^{-1} M \frac{a(\un Y^{\un\ell})}{\un Y^{\un k}} \smatr{a}{}{}{1} x_{\sigma,\ell} \\
    &= \chi_\sigma(\o a)^{-1} M \frac{a(\un Y^{\un\ell})}{\un Y^{\un k}} \smatr{a}{}{}{1} \lambda_\sigma^n \un Y^{\un{a_{nd'}^\sigma}-\un \ell} \smatr{p}{}{}{1}^{nd'} v_\sigma \\
    &= \chi_\sigma(\o a)^{-1} M \lambda_\sigma^n \frac{a(\un Y^{\un{a_{nd'}^\sigma}})}{\un Y^{\un k}} \smatr{p}{}{}{1}^{nd'} \smatr{a}{}{}{1} v_\sigma \\
    &= M \lambda_\sigma^n \frac{a(\un Y^{\un{a_{nd'}^\sigma}})}{\un Y^{\un k}} \smatr{p}{}{}{1}^{nd'} v_\sigma.
  \end{align*}
Recalling that $a(Y_i) = \o a^{p^i} Y_i f_{a,i}^{-1}$ and $\un{a_{nd'}^\sigma} = \un{a_{d'}^\sigma} \frac{p^{nd'}-1}{p^{d'}-1}$ the formula simplifies to
\begin{align*}
    &=  M \lambda_\sigma^n \bigg(\prod_{i=0}^{f-1} \o a^{p^i a_{nd',i}^\sigma}\bigg) \un Y^{\un{a_{nd'}^\sigma}-\un k} \bigg(\prod_{i=0}^{f-1} f_{a,i}^{-a_{nd',i}^\sigma}\bigg) \smatr{p}{}{}{1}^{nd'} v_\sigma\\
    &= M^{p^{nd'}} \bigg(\prod_{i=0}^{f-1} \o a^{p^i a_{nd',i}^\sigma}\bigg) \lambda_\sigma^n \un Y^{\un{a_{nd'}^\sigma}-\un k} \smatr{p}{}{}{1}^{nd'} v_\sigma.
  \end{align*}
  Now $M^{p^{nd'}}$ only matters modulo $F_{-kf-1}A$. But as $f_{a,i} \in 1+F_{-(p-1)}A$ we have 
  $ M^{p^{nd'}} \in 1+F_{-p^{nd'}(p-1)}A$, so for $n$ sufficiently large we can omit this factor.
  In summary, the $k$-th component of $N_{\F_q/\Fp}(\o a)\chi_\sigma(\o a)^{-1} M a(x_\sigma)$ is
  given by $$\big(\prod_{i=0}^{f-1} \o a^{p^i a_{nd',i}^\sigma}\big) \lambda_\sigma^n \un Y^{\un{a_{nd'}^\sigma}-\un k} \smatr{p}{}{}{1}^{nd'} v_\sigma =
  \big(\prod_{i=0}^{f-1} \o a^{p^i a_{nd',i}^\sigma}\big) x_{\sigma,k}.$$

  Finally notice that $\sum p^i a_{nd',i}^\sigma$ = $(1+p^{d'}+\cdots+p^{(n-1)d'}) \sum p^i a_{d',i}^\sigma \equiv n \sum p^i a_{d',i}^\sigma \pmod {q-1}$, as $f\mid d'$.
  Since $n$ (sufficiently large) was arbitrary above, we deduce that $\sum p^i a_{d',i}^\sigma \equiv 0 \pmod{q-1}$, and the result follows.
\end{proof}
 
Let $\{x_{\sigma}^{*} : \sigma\in W(\brho)\}$ denote the $A$-basis of $D_A(\pi)$ that is dual to $\{x_\sigma : \sigma \in W(\brho)\}$.
By~\eqref{eq:pairing-phi-gamma} we deduce the $\varphi$- and $\cO_K^{\times}$-actions on the elements $x_{\sigma}^{*}$ from Propositions \ref{prop:phi-x} and \ref{prop:gamma-action}.
\begin{cor1}\label{cor:phi-gamma-action}
Fix $\sigma \in W(\rhob)$. 
We have 
\[\varphi(x_{\sigma}^*)=(-1)^{f-1}\mu_{\sigma}\un{Y}^{\un{c_1^{\sigma}}} x_{\delta(\sigma)}^*,\]
and for $a\in\cO_K^{\times}$,
\[a(x_\sigma^*) = N_{\F_q/\Fp}(\o a) \chi_\sigma(\o a)^{-1}  \bigg(\prod_{i=0}^{f-1} f_{a,i}^{a_{d',i}^\sigma/(1-p^{d'})}\bigg) x_\sigma^*.\]
\end{cor1}

\subsection{The main theorem on \texorpdfstring{$D_A(\pi)$}{D\_A(pi)}}
\label{sec:structure-d_ap}

We prove $D_A(\pi)\simeq D_A^\otimes(\rhobar^\vee(1))$ which finishes the proof of Theorem \ref{mainbis}.

Recall that $\rhobar$ is as at the end of \S\ref{conjecture}.

\begin{thm1}\label{main}
There is an isomorphism of \'etale $(\varphi,\oK^\times)$-modules \[D_A(\pi)\simeq D_A^\otimes(\rhobar^\vee(1)).\]
\end{thm1}
\begin{proof}
We write $D_{A, \sigma_0}(\rhobar)=Ae_0\oplus Ae_1$ with $(e_0,e_1)$ as in Lemma \ref{PSnatural} (for $d=2$ and noting $e_i$ instead of $1\otimes e_i$) when $\rhobar$ is absolutely irreducible and where
\begin{equation}\label{actionsplit}
\left\{\begin{array}{cll}
\varphi_q(e_0)&=& \displaystyle{\lambda_0\Big(\frac{Y_{0}}{\varphi(Y_{0})}\Big)^h e_0}\\
\varphi_q(e_1)&=&\lambda_1 e_1\\
a(e_0)&=&\displaystyle{\Big(\frac{f_{a,0}}{\varphi(f_{a,0})}\Big)^{\frac{h}{1-q}} e_0}\\
a(e_1)&=&e_1.
\end{array}\right.
\end{equation}
when $\rhobar$ is (split) reducible. Let $I\defeq \{0,1\}^{f}$ and denote by $\underline i=(i_j)_j$ an element of $I$. By (\ref{daotimes}) and since $\varphi^{f-1-j}(e_{i_j})\in D_{A,{\sigma_{j+1}}}(\rhobar)$ (see (\ref{isoi})) we have $D_A^\otimes(\rhobar)=\bigoplus_{\underline i\in I}AE_{\underline i}$, where
\[E_{\underline i}\defeq \bigotimes_{j=0}^{f-1}\varphi^{f-1-j}(e_{i_j}).\]
We will define an explicit $A$-linear isomorphism from $D_{A}^{\otimes}(\brho^\vee(1))$ to $D_{A}(\pi)$ and check that 
it is a morphism of $(\varphi,\cO_K^{\times})$-modules. Twisting $\rhobar$ and $\pi$ by the same unramified \ character \ and \ using \ Lemma \ref{twist1} \ and \ Lemma \ref{twist}, \ we \ can \ assume ${\det}(\rhobar)(p)=1$, i.e.\ ${\det}(\rhobar)=\omega_f^{\sum_{i=0}^{f-1}p^i(r_i+1)}$. Then
\[D_{A}^{\otimes}(\brho^\vee(1))\simeq D_{A}^{\otimes}(\brho\otimes {\det}(\rhobar)^{-1}\omega)=D_{A}^{\otimes}\big(\brho\otimes \omega_f^{-\sum_{i=0}^{f-1}p^ir_i}\big),\]
and using Lemma \ref{twist1} and Lemma \ref{twist} again, it is equivalent to define an isomorphism {$\vartheta$} of $(\varphi,\cO_K^{\times})$-modules 
\begin{equation}\label{eq:morphism-twist}
D_{A}^{\otimes}(\brho)\longrightarrow D_{A}\big(\pi\otimes \omega_f^{-\sum_{i=0}^{f-1}p^ir_i}\big)\cong D_A(\pi)\otimes_A D_A\big(\omega_f^{\sum_{i=0}^{f-1}p^ir_i}\big).
\end{equation}

We \ know \ by \ Theorem~\ref{thm:descent} \ that \ $\{x_\sigma: \sigma\in W(\brho)\}$ \ form \ an \ $A$-basis \ of $\Hom_{A}(D_A(\pi),A)$.  
Let $\{x_{\sigma}^{*} : \sigma\in W(\brho)\}$ denote the $A$-basis of $D_A(\pi)$ that is dual to $\{x_\sigma : \sigma \in W(\brho)\}$, as in~\S\ref{sec:o_ks-action}. 
For convenience, below we write $x_J^{*}$ instead of $x_{\sigma}^{*} \otimes 1$ in~\eqref{eq:morphism-twist}, where $J = J_\sigma$.   

Write $\sigma=(s_0,\dots,s_{f-1})\otimes\eta$ and $\delta(\sigma)=(s_0',\dots,s_{f-1}')\otimes\eta'$. Below when we write, for example, $s_i=r_i+1$, we actually mean that  $\lambda_i(x_i)=x_i+1$, where $\lambda\in \mathcal{ID}(x_0,\dots,x_{f-1})$ or $\mathcal{RD}(x_0,\dots,x_{f-1})$ is the element corresponding to $\sigma$; see Remark \ref{rk:s-s'}. 

(i) Assume first $\brho$ is absolutely irreducible.  
For $J\subset \{0,1,\dots,f-1\}$, with corresponding Serre weight $\sigma\in W(\rhobar)$, define 
\begin{equation*}\vartheta: E_{\underline{i_J}} \longmapsto \alpha_J \underline{Y}^{\underline{b_J}-\underline{1}}x_J^{*},
\end{equation*}
where $\alpha_J\in\F^{\times}$ are suitable constants, $\underline{i_J}\defeq \mathbf{1}_J$  (i.e.\ $i_{J,j}=1$ if $j\in J$ and $i_{J,j}=0$ if $j\notin J$), and
\begin{itemize}
\item $b_{J,i}\defeq 0$ if either $i=0$ and $s_0\in \{r_0,r_0-1\}$, or $i>0$ and $s_i\in\{r_i,p-3-r_i\}$;
\item $b_{J,0}\defeq -h^{[0]}+1$ if $s_0=p-1-r_0$;
\item $b_{J,i}\defeq h^{[i]}+1$ if $i>0$ and $s_i=r_i+1$;
\item $b_{J,i}\defeq -h^{[i]}$ if $s_i=p-2-r_i$;
\end{itemize}
where $h^{[i]}$ was defined in \eqref{eq:h[i]}.
Below 
we check that for well-chosen $\alpha_J$, $\vartheta$ commutes with  $\varphi$, i.e.\ $\vartheta\big(\varphi(E_{\un{i_J}})\big)=\varphi(\alpha_J \un{Y}^{\un{b_J}-\un{1}}x_J^*)$.
Writing $J'=J_{\delta(\sigma)}$, Corollary~\ref{cor:phi-gamma-action} implies 
\begin{equation}\label{eq:phi-xJdual}
\varphi(x_J^*) = (-1)^{f-1}\mu_{J}\underline{Y}^{\underline{c_{J'}}}x_{J'}^{*},
\end{equation}
where $\mu_{J}\defeq \mu_{\sigma_J}$, and $\underline{c_{J'}}$ is defined as in \S\ref{ss:degree} with respect to the pair $(\sigma,\delta(\sigma))$. Also, using Lemma \ref{PSnatural} it is easy to check that  \begin{equation*}\varphi(E_{\underline{i_J}})=\left\{\begin{array}{rll} 
 E_{\underline{i_{J'}}}& \text{if $i_{J,0}=0,$}\\
-(\frac{Y_{\sigma_0}}{Y_{\sigma_{f-1}}^{p}})^hE_{\underline{i_{J'}}}&\text{if $i_{J,0}=1$.}\end{array}\right.
\end{equation*}
Thus, we are reduced to checking:
\begin{itemize}
\item if $i_{J,0}=0$ then
\begin{equation}\label{eq:bJ-irr1}
\left\{\begin{array}{ll}\alpha_J\cdot\mu_{J}=(-1)^{f-1}\alpha_{J'},\\
p\delta(\underline{b_J})+\underline{1-p}+\underline{c_{J'}}=\underline{b_{J'}};\end{array}\right.
\end{equation}
\item if $i_{J,0}=1$ then  
\begin{equation}\label{eq:bJ-irr2}
\left\{\begin{array}{ll} \alpha_J\cdot\mu_{J}=(-1)^{f}\alpha_{J'},\\p\delta(\underline{b_J})+\underline{1-p}+\underline{c_{J'}}=\underline{b_{J'}}+(h,0,\dots,-ph).\end{array}\right.
\end{equation}
\end{itemize}
 
First assume $0\notin J$, i.e.\ $s_0\in\{r_0,r_0-1\}$; note that this implies $s_1\in \{r_1,p-2-r_1\}$ by the property of $W(\brho)$. 
We need to check 
\begin{equation*}pb_{J,i}+1-p+c_{J',i-1}=b_{J',i-1}
\end{equation*}
for any $0\leq i\leq f-1$. 
It is a direct check using Lemma \ref{lem:s-s'}. We do it for  $i=0,1$ and leave the other cases as an exercise. Recall that $c_{J',i-1}=s_{i-1}'$ if $i-1\in J^{\rm max}(\sigma_{J})$ and $c_{J',i-1}=p-1$ otherwise. 
 \begin{itemize}
 \item If $i=0$ and $s_0=r_0$, then $b_{J,0}=0$ by definition and $c_{J',f-1}=s_{f-1}'=p-2-r_{f-1}$ by Lemma \ref{lem:s-s'}, so we obtain
 \[p\cdot 0+(1-p)+(p-2-r_{f-1})=-h^{[f-1]},\]
 which is equal to $b_{J',f-1}$. 
 \item If $i=0$ and $s_0=r_0-1$, then $b_{J,0}=0$ by definition and $c_{J',f-1}=p-1$ by Lemma \ref{lem:s-s'}, so  we obtain
 \[p\cdot 0+(1-p)+(p-1)=0,\]
 which is equal to $b_{J',f-1}$ (as $s_{f-1}'=p-3-r_{f-1}$). 
 \item If $i=1$ and $s_1=r_1$, then $b_{J,1}=0$ by definition and $c_{J',0}=p-1$ by Lemma \ref{lem:s-s'}, so we obtain
 \[p\cdot 0+(1-p)+(p-1)=0,\]
 which is equal to $b_{J',0}$ (as $s_{0}'=r_0-1$).
\item If $i=1$ and $s_1=p-2-r_1$, then $b_{J,1}=-h^{[1]}$ by definition and $c_{J',0}=s_{0}'=p-1-r_0$, so we obtain
\[p(-h^{[1]})+(1-p)+(p-1-r_0)=-h^{[0]}+1,\]
which is equal to $b_{J',0}$. 
 \end{itemize}
 
Assume $0\in J$, i.e.\ $s_0\in \{p-2-r_0,p-1-r_0\}$; note that this implies $s_1\in \{r_1+1,p-3-r_1\}$. We check \eqref{eq:bJ-irr2} for $i=0$ and leave the other cases as an exercise. 
\begin{itemize}
\item If $s_0=p-2-r_0$, then $b_{J,0}=-h^{[0]}$ by definition and $c_{J',f-1}=p-1$ by Lemma \ref{lem:s-s'}, so we obtain
\[p(-h^{[0]})+(1-p)+(p-1)=-ph,\]
which equals to $b_{J',f-1}-ph$ (as $b_{J',f-1}=0$, since $s'_{f-1} = r_{f-1}$).
\item If $s_0=p-1-r_0$, then $b_{J,0}=-h^{[0]}+1$ by definition and $c_{J',f-1}=s_{f-1}'=r_{f-1}+1$, so we obtain
\[p(-h^{[0]}+1)+(1-p)+(r_{f-1}+1)=(r_{f-1}+1)+1-ph,\]
which  is equal to $b_{J',f-1}-ph$ (as $b_{J',f-1}=h^{[f-1]}+1$).
\end{itemize}

Now we show that  the constants $\alpha_{J}$ can be  compatibly chosen so that $\vartheta$ is $\varphi$-equivariant.  Using \eqref{eq:bJ-irr1} and \eqref{eq:bJ-irr2}  it suffices to check, for any  $J$ whose orbit has length $d$, that 
\[(-1)^{d(f-1)}\prod_{j=0}^{d-1}\mu_{\delta^j(J)}=\prod_{j=0}^{d-1}(-1)^{i_{\delta^j(J),0}}.\]
As the left-hand side is equal to $(-1)^{-\frac{d}{2}}=(-1)^{\frac{d}{2}}$ by Lemma \ref{lem:n+d:lambda} and \eqref{eq:lambda-DL} (and ${\det}(\rhobar)(p)=1$), it suffices to show that
\begin{equation}\label{eq:number-d}
\sharp \big\{0\leq j\leq d-1, 0\in \delta^j(J)\big\}=\frac{d}{2}.
\end{equation}
By the proof of \cite[Lemma 5.2]{breuil-IL}, letting $J'=J\cup \{f+j, j\in\overline{J}\}$ (where $\overline J$ is the complement of $J$), 
then $d$ is also the smallest positive integer such that $J' = J'-d$ as subsets of $\Z/2f\Z$, and in particular $d$ divides $2f$. Since $|J'|=f$ and $J'\cap \{0,1,\dots,f-1\}=J$, it is easy to see that
\[\sharp\big\{0\leq j\leq 2f-1,~ 0\in \delta^{j}(J)\big\}=f\]
from which we deduce \eqref{eq:number-d}.

We now check that $\vartheta$ is $\cO_K^{\times}$-equivariant. By Lemma~\ref{PSnatural} we know that
\begin{equation*}
  a(E_{\un{i_J}}) = \prod_{i=0}^{f-1} \varphi^{f-1-i}(f_{a,0}^{h(1-\varphi)q^{i_{J,i}}/(1-q^2)}) E_{\un{i_J}} 
\end{equation*}
and by Corollary~\ref{cor:phi-gamma-action} we have
\begin{equation*}
  a(x_J^*) = N_{\F_q/\Fp}(\o a) \chi_\sigma(\o a)^{-1} \overline a^{\sum_{i=0}^{f-1}p^ir_i}\bigg(\prod_{i=0}^{f-1} f_{a,i}^{a_{d',i}^\sigma/(1-p^{d'})}\bigg) x_J^*,
\end{equation*}
where $d' = df$ and recall the twist $D_A(\omega_f^{\sum_{i=0}^{f-1}p^ir_i})$ in \eqref{eq:morphism-twist}. Thus it suffices to show that
\begin{multline*}
  a(\un Y^{\un{b_J}-\un 1}) N_{\F_q/\Fp}(\o a) \chi_\sigma(\o a)^{-1} \overline a^{\sum_{i=0}^{f-1}p^ir_i}\bigg(\prod_{i=0}^{f-1} f_{a,i}^{a_{d',i}^\sigma/(1-p^{d'})}\bigg)\\
   =  \prod_{i=0}^{f-1}f_{a,i+1}^{p^{f-1-i}hq^{i_{J,i}}/(1-q^2)}f_{a,i}^{-p^{f-i}hq^{i_{J,i}}/(1-q^2)} \un Y^{\un{b_J}-\un 1},
\end{multline*}
which is implied by the following claims (where we use that $2f\mid d'$):
\begin{enumerate}
\item[\upshape(a)] $\chi_\sigma(\o a) = \o a^{\sumfi p^i b_{J,i}+\sum_{i=0}^{f-1}p^ir_i}$,
\item[\upshape(b)] 
  \begin{equation*}
  \displaystyle \frac{a_{d',i}^\sigma}{1-p^{d'}}+1-b_{J,i} = 
  \begin{cases}
    \displaystyle\frac{hp^{f-i}}{1-q^2} [q^{\mathbf{1}_J(i-1)}-q^{\mathbf{1}_J(i)}] & \text{if $1 \le i \le f-1$,}\\
    \displaystyle\frac h{1-q^2} [q^{\mathbf{1}_J(f-1)}-q\cdot q^{\mathbf{1}_J(0)}] & \text{if $i = 0$}.
  \end{cases}
\end{equation*}
\end{enumerate}

To verify the first claim, note from \cite[\S2]{breuil-IL} that 
\begin{equation*}
\chi_\sigma(\o a) = \o a^{\frac 12\big(\sumfi p^i (r_i+s_i)+(q-1)\mathbf 1_J(f-1)\big)}.
\end{equation*}
It then suffices to show that
\begin{equation*}
  \frac 12\big(\sumfi p^i (s_i-r_i) + (q-1)\mathbf 1_J(f-1)\big)  \equiv \sumfi p^i b_{J,i}   \pmod{q-1}.
\end{equation*}
First assume $f-1\notin J$ (so that $\mathbf{1}_J(f-1)=0$), equivalently $s_0\in\{r_0,p-2-r_0\}$. Then $(s_0,\dots,s_{f-1})$ consists of subsequences  of the form $p-2-r_j,p-3-r_{j+1},\dots,p-3-r_{j'-1},r_{j'}+1$ for some $0\leq j<j'\leq f-1$ (and $r_i$ for $i\notin\{j,\dots,j'\}$). Since $b_{J,i}=0$ if $s_i=r_i$, we are reduced to prove that for $0\leq j< j'\leq f-1$, 
\begin{equation}\label{eq:claim-a}
\frac{1}{2}\big((p-2-2r_{j})+\sum_{j<i<j'}p^i(p-3-2r_i)+p^{j'}\big)\equiv\sum_{j\leq i\leq j'}p^ib_{J,i} \pmod{q-1}.
\end{equation}
It is direct to check that the left-hand side of \eqref{eq:claim-a} is equal to 
\[p^j(p-1-r_j)+\sum_{j<i<j'}p^{i}(p-2-r_i)=p^j(p-h_j)+\sum_{j<i<j'}p^i(p-1-h_i).\]
On the other hand, by the definition of $b_{J,i}$ the right-hand side of \eqref{eq:claim-a} is equal to  
\[\begin{array}{rll}p^j(-h^{[j]})+p^{j'}(h^{[j']}+1)&=&p^{j'}-p^jh_j-\cdots p^{j'-1}h_{j'-1}\\
&=&p^j(p-h_j)+\sum_{j<i<j'}p^i(p-1-h_i),\end{array}\]
hence \eqref{eq:claim-a} is verified in this case (we actually have an equality). Now assume $f-1\in J$ (so that $\mathbf{1}_J(f-1)=1$), equivalently $s_0\in \{r_0-1,p-1-r_0\}$. 
\begin{itemize}
\item[$\bullet$] If $s_0=r_0-1$, then $(s_0,\dots,s_{f-1}) $ contains a subsequence of the form $p-2-r_j, p-3-r_{j+1},\dots, r_0-1$  for some $0< j\leq f-1$ (note that the case $j=f-1$ is allowable), and one computes
\[\begin{array}{rll}&\frac{1}{2}\big(p^j(p-2-2r_j)+\sum_{j<i\leq f-1}p^i(p-3-2r_{i})+(-1)+q-1\big)\\
=&p^j(p-1-r_j)+\sum_{j<i\leq f-1}p^i(p-2-r_i)+(-1)\\
\equiv&p^j(-h^{[j]}) \pmod{q-1}.
\end{array}\] 
\item[$\bullet$] If $s_0=p-1-r_0$, then $(s_0,\dots,s_{f-1}) $ contains a subsequence of the form $p-2-r_{j}, p-3-r_{j+1},\cdots,p-1-r_0,p-3-r_1,\dots,p-3-r_{j'-1},r_{j'}+1 $ for some   $0< j'<j\leq f-1$, and  one checks   the following congruence relation mod $q-1$:
\[\frac{1}{2}\big((p-2-2r_{j})+\sum_{j<i<j',i\neq0}p^i(p-3-2r_i)+(p-1-2r_0)+p^{j'}+(q-1)\big)\equiv \sum_{j\leq i\leq j'}p^ib_{J,i},\]    where $\sum_{j<i<j'}$ means $\sum_{j<i\leq f-1}+\sum_{0\leq i< j'}$ and similarly for $\sum_{j\leq i\leq j'}$.  
\end{itemize}
Together with \eqref{eq:claim-a}, claim (a) is  verified in this case.

We check claim (b). Using Lemma \ref{lem:a-d'-sigma} and the definition of $\un{b_J}$ one checks that  
{\begin{center}
\begin{tabular}{ |c||c|c|c|c|c|c|}
\hline  
$s_0$& $r_0$ & $r_0-1$& $p-2-r_0$&$p-1-r_0$          \\
\hline  
 $\frac{a_{d',0}^{\sigma}}{1-p^{d'}}+1-b_{J,0}$ & $\frac{h}{1+q}$&$0$&$h$&$\frac{hq}{1+q}$   \\
\hline
 \end{tabular}
 \end{center}}
\noindent while if $1\leq i\leq f-1$ we have
  {\begin{center}
\begin{tabular}{ |c||c|c|c|c|c|c|}
\hline  
$s_i$& $r_i$ & $r_i+1$& $p-2-r_i$&$p-3-r_i$          \\
\hline  
 $\frac{a_{d',i}^{\sigma}}{1-p^{d'}}+1-b_{J,i}$ & $0$& $-\frac{hp^{f-i}}{1+q}$&$\frac{hp^{f-i}}{1+q}$& $0$   \\
\hline
 \end{tabular}
\end{center}}
\noindent Then (b) can easily be checked case by case.

(ii) Assume $\brho$ is (split) reducible.
For $J\subset \{0,1,\dots,f-1\}$, with corresponding Serre weight $\sigma\in W(\rhobar)$, define
\[\vartheta:\ E_{\un{i_J}}\longmapsto \alpha_{J}Y^{\un{b_{J}}-\un{1}}x_J^{*}, \]
where $\alpha_{J}\in \F^{\times}$ are suitable constants, $\un{i_J}\defeq \mathbf{1}_{J^c}$ (i.e.\ $i_{J,j}=1$ if $j\notin J$ and $i_{J,j}=0$ if $j\in J$), and
\begin{itemize}
\item $b_{J,i}=0$ if $s_i=r_i$;
\item $b_{J,i}=-h^{[i]}$ if $s_i=p-2-r_i$;
\item $b_{J,i}=h^{[i]}+1$ if $s_i=r_i+1$ and $i>0$ (resp.\ $b_{J,0}=1$ if $i=0$);
\item $b_{J,i}=0$ if $s_i=p-3-r_i$ and $i>0$ (resp.\ $b_{J,0}=-h^{[0]}$ if $i=0$).
\end{itemize}

Write $J'=J_{\delta(\sigma)}$. Then \eqref{eq:phi-xJdual} remains true,  and it is easy to check that \begin{equation*}\varphi(E_{\underline{i_J}})=\left\{\begin{array}{rll} \lambda_0(\frac{Y_{\sigma_0}}{Y_{\sigma_{f-1}}^{p}})^hE_{\underline{i_{J'}}}& \text{if $i_{J,0}=0$,}\\
\lambda_1E_{\underline{i_{J'}}}&\text{if $i_{J,0}=1$}.\end{array}\right.
\end{equation*}
Thus, to check that $\vartheta$ is $\varphi$-equivariant it is equivalent to check
\begin{itemize}
\item if $i_{J,0}=0$ then
\begin{equation}\label{eq:bJ-red1}
\left\{\begin{array}{ll}
\alpha_J\cdot\mu_{J}=(-1)^{f-1}\lambda_0\alpha_{J'},\\
p\delta(\underline{b_J})+\underline{1-p}+\underline{c_{J'}}=\underline{b_{J'}}+(h,0,\dots,-ph);\end{array}\right.\end{equation}
\item if $i_{J,0}=1$ then
\begin{equation}\label{eq:bJ-red2}
\left\{\begin{array}{ll}
\alpha_J\cdot\mu_{J}=(-1)^{f-1}\lambda_1\alpha_{J'},\\
p\delta(\underline{b_J})+\underline{1-p}+\underline{c_{J'}}=\underline{b_{J'}}. \end{array}\right.
\end{equation}
\end{itemize}

We leave it as an exercise to check the second equation of \eqref{eq:bJ-red1}, resp.\ \eqref{eq:bJ-red2}, using Lemma \ref{lem:s-s'}. Thus,  to show that  the constants $\alpha_{J}$ can be compatibly chosen so that  $\vartheta$ is $\varphi$-equivariant, it suffices to check, for any $J$ whose orbit has length $d$, 
\[(-1)^{d(f-1)}\prod_{j=0}^{d-1}\mu_{\delta^j(J)}=\lambda_0^{|J|\frac{d}{f}}\lambda_0^{-|\overline{J}|\frac{d}{f}},\]
where we have used ${\det}(\rhobar)(p)=1$ and the fact that
\[\sharp\big\{0\leq j\leq d-1, 0\in \delta^j(J)\big\}=|J|\frac{d}{f},\ \ \sharp\big\{0\leq j\leq d-1, 0\notin \delta^j(J)\big\}=|\overline{J}|\frac{d}{f}.\]
We conclude by Lemma \ref{lem:n+d:lambda} and \eqref{eq:lambda-DL}.

We now check that $\vartheta$ is $\cO_K^{\times}$-equivariant. Using (\ref{actionsplit}) we know that
\begin{equation*}
  a(E_{i_J}) = \prod_{i : i_{J,i} = 0} \varphi^{f-1-i}(f_{a,0}^{h(1-\varphi)/(1-q)}) E_{i_J}
\end{equation*}
and by Corollary~\ref{cor:phi-gamma-action} we have
\begin{equation*}
  a(x_J^*) = N_{\F_q/\Fp}(\o a) \chi_\sigma(\o a)^{-1} \overline a^{\sum_{i=0}^{f-1}p^ir_i}\bigg(\prod_{i=0}^{f-1} f_{a,i}^{a_{d',i}^\sigma/(1-p^{d'})}\bigg) x_J^*.
\end{equation*}
Thus it suffices to show that
\begin{multline*}
  a(\un Y^{\un{b_J}-\un 1}) N_{\F_q/\Fp}(\o a) \chi_\sigma(\o a)^{-1}\overline a^{\sum_{i=0}^{f-1}p^ir_i} \bigg(\prod_{i=0}^{f-1} f_{a,i}^{a_{d',i}^\sigma/(1-p^{d'})}\bigg)\\
   = \prod_{\stackrel{0 \le i \le f-1}{i_{J,i} = 0}} (f_{a,i+1}^{hp^{f-i-1}/(1-q)}f_{a,i}^{-hp^{f-i}/(1-q)}) \un Y^{\un{b_J}-\un 1},
\end{multline*}
which is implied by the following claims (where we use that $f\mid d'$):
\begin{enumerate}
\item[\upshape(a)] $\chi_\sigma(\o a) = \o a^{\sumfi p^i b_{J,i}+\sum_{i=0}^{f-1}p^ir_i}$,
\item[\upshape(b)] 
  \begin{equation*}
  \displaystyle \frac{a_{d',i}^\sigma}{1-p^{d'}}+1-b_{J,i} = 
  \begin{cases}
    \displaystyle\frac{hp^{f-i}}{1-q} [\mathbf{1}_J(i-1)-\mathbf{1}_J(i)] & \text{if $1 \le i \le f-1$,}\\
    \displaystyle\frac h{1-q} [\mathbf{1}_J(f-1)-q\mathbf{1}_J(0)] & \text{if $i = 0$}.
  \end{cases}
\end{equation*}
\end{enumerate}
Both claims are checked as in the irreducible case (we omit the details).
\end{proof}

\clearpage{}

\newpage

\bibliography{Biblio}

\bibliographystyle{amsalpha}

\end{document}